\noindent


\edef\resetatcatcode{\catcode`\noexpand\@\the\catcode`\@\relax}
\ifx\miniltx\undefined\else\endinput\fi
\let\miniltx\box

\def\makeatletter{\catcode`\@11\relax}
\def\makeatother{\catcode`\@12\relax}
\makeatletter

\def\@makeother#1{\catcode`#1=12\relax}

\def\@ifnextchar#1#2#3{%
  \let\reserved@d=#1%
  \def\reserved@a{#2}\def\reserved@b{#3}%
  \futurelet\@let@token\@ifnch}
\def\@ifnch{%
  \ifx\@let@token\@sptoken
    \let\reserved@c\@xifnch
  \else
    \ifx\@let@token\reserved@d
      \let\reserved@c\reserved@a
    \else
      \let\reserved@c\reserved@b
    \fi
  \fi
  \reserved@c}
\begingroup
\def\:{\global\let\@sptoken= } \:  
\def\:{\@xifnch} \expandafter\gdef\: {\futurelet\@let@token\@ifnch}
\endgroup

\def\@ifstar#1{\@ifnextchar *{\@firstoftwo{#1}}}
\long\def\@dblarg#1{\@ifnextchar[{#1}{\@xdblarg{#1}}}
\long\def\@xdblarg#1#2{#1[{#2}]{#2}}

\long\def \@gobble #1{}
\long\def \@gobbletwo #1#2{}
\long\def \@gobblefour #1#2#3#4{}
\long\def\@firstofone#1{#1}
\long\def\@firstoftwo#1#2{#1}
\long\def\@secondoftwo#1#2{#2}

\def\NeedsTeXFormat#1{\@ifnextchar[\@needsf@rmat\relax}
\def\@needsf@rmat[#1]{}
\def\ProvidesPackage#1{\@ifnextchar[%
    {\@pr@videpackage{#1}}{\@pr@videpackage#1[]}}
\def\@pr@videpackage#1[#2]{\wlog{#1: #2}}
\let\ProvidesFile\ProvidesPackage

\let\DeclareOption\@gobbletwo
\def\ProcessOptions{\@ifstar\relax\relax}

\def\RequirePackage{%
  \@fileswithoptions\@pkgextension}
\def\@fileswithoptions#1{%
  \@ifnextchar[
    {\@fileswith@ptions#1}%
    {\@fileswith@ptions#1[]}}
\def\@fileswith@ptions#1[#2]#3{%
  \@ifnextchar[
  {\@fileswith@pti@ns#1[#2]#3}%
  {\@fileswith@pti@ns#1[#2]#3[]}}

\def\@fileswith@pti@ns#1[#2]#3[#4]{%
    \def\reserved@b##1,{%
      \ifx\@nil##1\relax\else
        \ifx\relax##1\relax\else
         \noexpand\@onefilewithoptions##1[#2][#4]\noexpand\@pkgextension
        \fi
        \expandafter\reserved@b
      \fi}%
      \edef\reserved@a{\zap@space#3 \@empty}%
      \edef\reserved@a{\expandafter\reserved@b\reserved@a,\@nil,}%
  \reserved@a}

\def\zap@space#1 #2{%
  #1%
  \ifx#2\@empty\else\expandafter\zap@space\fi
  #2}

\let\@empty\empty
\def\@pkgextension{sty}

\def\@onefilewithoptions#1[#2][#3]#4{%
  \input #1.#4 }

\def\typein{%
  \let\@typein\relax
  \@testopt\@xtypein\@typein}
\def\@xtypein[#1]#2{%
  \message{#2}%
  \advance\endlinechar\@M
  \read\@inputcheck to#1%
  \advance\endlinechar-\@M
  \@typein}
\def\@namedef#1{\expandafter\def\csname #1\endcsname}
\def\@nameuse#1{\csname #1\endcsname}
\def\@cons#1#2{\begingroup\let\@elt\relax\xdef#1{#1\@elt #2}\endgroup}
\def\@car#1#2\@nil{#1}
\def\@cdr#1#2\@nil{#2}
\def\@carcube#1#2#3#4\@nil{#1#2#3}
\def\@preamblecmds{}

\def\@star@or@long#1{%
  \@ifstar
   {\let\l@ngrel@x\relax#1}%
   {\let\l@ngrel@x\long#1}}

\let\l@ngrel@x\relax
\def\newcommand{\@star@or@long\new@command}
\def\new@command#1{%
  \@testopt{\@newcommand#1}0}
\def\@newcommand#1[#2]{%
  \@ifnextchar [{\@xargdef#1[#2]}%
                {\@argdef#1[#2]}}
\long\def\@argdef#1[#2]#3{%
   \@ifdefinable #1{\@yargdef#1\@ne{#2}{#3}}}
\long\def\@xargdef#1[#2][#3]#4{%
  \@ifdefinable#1{%
     \expandafter\def\expandafter#1\expandafter{%
          \expandafter
          \@protected@testopt
          \expandafter
          #1%
          \csname\string#1\expandafter\endcsname
          {#3}}%
       \expandafter\@yargdef
          \csname\string#1\endcsname
           \tw@
           {#2}%
           {#4}}}
\def\@testopt#1#2{%
  \@ifnextchar[{#1}{#1[#2]}}
\def\@protected@testopt#1{
  \ifx\protect\@typeset@protect
    \expandafter\@testopt
  \else
    \@x@protect#1%
  \fi}
\long\def\@yargdef#1#2#3{%
  \@tempcnta#3\relax
  \advance \@tempcnta \@ne
  \let\@hash@\relax
  \edef\reserved@a{\ifx#2\tw@ [\@hash@1]\fi}%
  \@tempcntb #2%
  \@whilenum\@tempcntb <\@tempcnta
     \do{%
         \edef\reserved@a{\reserved@a\@hash@\the\@tempcntb}%
         \advance\@tempcntb \@ne}%
  \let\@hash@##%
  \l@ngrel@x\expandafter\def\expandafter#1\reserved@a}
\long\def\@reargdef#1[#2]#3{%
  \@yargdef#1\@ne{#2}{#3}}
\def\renewcommand{\@star@or@long\renew@command}
\def\renew@command#1{%
  {\escapechar\m@ne\xdef\@gtempa{{\string#1}}}%
  \expandafter\@ifundefined\@gtempa
     {\@latex@error{\string#1 undefined}\@ehc}%
     {}%
  \let\@ifdefinable\@rc@ifdefinable
  \new@command#1}
\long\def\@ifdefinable #1#2{%
      \edef\reserved@a{\expandafter\@gobble\string #1}%
     \@ifundefined\reserved@a
         {\edef\reserved@b{\expandafter\@carcube \reserved@a xxx\@nil}%
          \ifx \reserved@b\@qend \@notdefinable\else
            \ifx \reserved@a\@qrelax \@notdefinable\else
              #2%
            \fi
          \fi}%
         \@notdefinable}
\let\@@ifdefinable\@ifdefinable
\long\def\@rc@ifdefinable#1#2{%
  \let\@ifdefinable\@@ifdefinable
  #2}
\def\newenvironment{\@star@or@long\new@environment}
\def\new@environment#1{%
  \@testopt{\@newenva#1}0}
\def\@newenva#1[#2]{%
   \@ifnextchar [{\@newenvb#1[#2]}{\@newenv{#1}{[#2]}}}
\def\@newenvb#1[#2][#3]{\@newenv{#1}{[#2][#3]}}
\def\renewenvironment{\@star@or@long\renew@environment}
\def\renew@environment#1{%
  \@ifundefined{#1}%
     {\@latex@error{Environment #1 undefined}\@ehc
     }{}%
  \expandafter\let\csname#1\endcsname\relax
  \expandafter\let\csname end#1\endcsname\relax
  \new@environment{#1}}
\long\def\@newenv#1#2#3#4{%
  \@ifundefined{#1}%
    {\expandafter\let\csname#1\expandafter\endcsname
                         \csname end#1\endcsname}%
    \relax
  \expandafter\new@command
     \csname #1\endcsname#2{#3}%
     \l@ngrel@x\expandafter\def\csname end#1\endcsname{#4}}

\def\providecommand{\@star@or@long\provide@command}
\def\provide@command#1{%
  {\escapechar\m@ne\xdef\@gtempa{{\string#1}}}%
  \expandafter\@ifundefined\@gtempa
    {\def\reserved@a{\new@command#1}}%
    {\def\reserved@a{\renew@command\reserved@a}}%
   \reserved@a}%

\def\@ifundefined#1{%
  \expandafter\ifx\csname#1\endcsname\relax
    \expandafter\@firstoftwo
  \else
    \expandafter\@secondoftwo
  \fi}

\chardef\@xxxii=32
\mathchardef\@Mi=10001
\mathchardef\@Mii=10002
\mathchardef\@Miii=10003
\mathchardef\@Miv=10004

\newcount\@tempcnta
\newcount\@tempcntb
\newif\if@tempswa\@tempswatrue
\newdimen\@tempdima
\newdimen\@tempdimb
\newdimen\@tempdimc
\newbox\@tempboxa
\newskip\@tempskipa
\newskip\@tempskipb
\newtoks\@temptokena

\long\def\@whilenum#1\do #2{\ifnum #1\relax #2\relax\@iwhilenum{#1\relax
     #2\relax}\fi}
\long\def\@iwhilenum#1{\ifnum #1\expandafter\@iwhilenum
         \else\expandafter\@gobble\fi{#1}}
\long\def\@whiledim#1\do #2{\ifdim #1\relax#2\@iwhiledim{#1\relax#2}\fi}
\long\def\@iwhiledim#1{\ifdim #1\expandafter\@iwhiledim
        \else\expandafter\@gobble\fi{#1}}
\long\def\@whilesw#1\fi#2{#1#2\@iwhilesw{#1#2}\fi\fi}
\long\def\@iwhilesw#1\fi{#1\expandafter\@iwhilesw
         \else\@gobbletwo\fi{#1}\fi}
\def\@nnil{\@nil}
\def\@empty{}
\def\@fornoop#1\@@#2#3{}
\long\def\@for#1:=#2\do#3{%
  \expandafter\def\expandafter\@fortmp\expandafter{#2}%
  \ifx\@fortmp\@empty \else
    \expandafter\@forloop#2,\@nil,\@nil\@@#1{#3}\fi}
\long\def\@forloop#1,#2,#3\@@#4#5{\def#4{#1}\ifx #4\@nnil \else
       #5\def#4{#2}\ifx #4\@nnil \else#5\@iforloop #3\@@#4{#5}\fi\fi}
\long\def\@iforloop#1,#2\@@#3#4{\def#3{#1}\ifx #3\@nnil
       \expandafter\@fornoop \else
      #4\relax\expandafter\@iforloop\fi#2\@@#3{#4}}
\def\@tfor#1:={\@tf@r#1 }
\long\def\@tf@r#1#2\do#3{\def\@fortmp{#2}\ifx\@fortmp\space\else
    \@tforloop#2\@nil\@nil\@@#1{#3}\fi}
\long\def\@tforloop#1#2\@@#3#4{\def#3{#1}\ifx #3\@nnil
       \expandafter\@fornoop \else
      #4\relax\expandafter\@tforloop\fi#2\@@#3{#4}}
\long\def\@break@tfor#1\@@#2#3{\fi\fi}
\def\@removeelement#1#2#3{%
  \def\reserved@a##1,#1,##2\reserved@a{##1,##2\reserved@b}%
  \def\reserved@b##1,\reserved@b##2\reserved@b{%
    \ifx,##1\@empty\else##1\fi}%
  \edef#3{%
    \expandafter\reserved@b\reserved@a,#2,\reserved@b,#1,\reserved@a}}

\let\ExecuteOptions\@gobble

\def\PackageError#1#2#3{%
  \errhelp{#3}\errmessage{#1: #2}}
\def\@latex@error#1#2{%
  \errhelp{#2}\errmessage{#1}}

\bgroup\uccode`\!`\%\uppercase{\egroup
\def\@percentchar{!}}

\ifx\@@input\@undefined
 \let\@@input\input
\fi

\def\input{\@ifnextchar\bgroup\@iinput\@@input}
\def\@iinput#1{\input{#1} }

\ifx\filename@parse\@undefined
  \def\reserved@a{./}\ifx\@currdir\reserved@a
    \wlog{^^JDefining UNIX/DOS style filename parser.^^J}
    \def\filename@parse#1{%
      \let\filename@area\@empty
      \expandafter\filename@path#1/\\}
    \def\filename@path#1/#2\\{%
      \ifx\\#2\\%
         \def\reserved@a{\filename@simple#1.\\}%
      \else
         \edef\filename@area{\filename@area#1/}%
         \def\reserved@a{\filename@path#2\\}%
      \fi
      \reserved@a}
  \else\def\reserved@a{[]}\ifx\@currdir\reserved@a
    \wlog{^^JDefining VMS style filename parser.^^J}
    \def\filename@parse#1{%
      \let\filename@area\@empty
      \expandafter\filename@path#1]\\}
    \def\filename@path#1]#2\\{%
      \ifx\\#2\\%
         \def\reserved@a{\filename@simple#1.\\}%
      \else
         \edef\filename@area{\filename@area#1]}%
         \def\reserved@a{\filename@path#2\\}%
      \fi
      \reserved@a}
  \else\def\reserved@a{:}\ifx\@currdir\reserved@a
    \wlog{^^JDefining Mac style filename parser.^^J}
    \def\filename@parse#1{%
      \let\filename@area\@empty
      \expandafter\filename@path#1:\\}
    \def\filename@path#1:#2\\{%
      \ifx\\#2\\%
         \def\reserved@a{\filename@simple#1.\\}%
      \else
         \edef\filename@area{\filename@area#1:}%
         \def\reserved@a{\filename@path#2\\}%
      \fi
      \reserved@a}
  \else
    \wlog{^^JDefining generic filename parser.^^J}
    \def\filename@parse#1{%
      \let\filename@area\@empty
      \expandafter\filename@simple#1.\\}
  \fi\fi\fi
  \def\filename@simple#1.#2\\{%
    \ifx\\#2\\%
       \let\filename@ext\relax
    \else
       \edef\filename@ext{\filename@dot#2\\}%
    \fi
    \edef\filename@base{#1}}
  \def\filename@dot#1.\\{#1}
\else
  \wlog{^^J^^J%
    \noexpand\filename@parse was defined in texsys.cfg:^^J%
    \expandafter\strip@prefix\meaning\filename@parse.^^J%
    }
\fi

\long\def \IfFileExists#1#2#3{%
  \openin\@inputcheck#1 %
  \ifeof\@inputcheck
    \ifx\input@path\@undefined
      \def\reserved@a{#3}%
    \else
      \def\reserved@a{\@iffileonpath{#1}{#2}{#3}}%
    \fi
  \else
    \closein\@inputcheck
    \edef\@filef@und{#1 }%
    \def\reserved@a{#2}%
  \fi
  \reserved@a}
\long\def\@iffileonpath#1{%
  \let\reserved@a\@secondoftwo
  \expandafter\@tfor\expandafter\reserved@b\expandafter
             :\expandafter=\input@path\do{%
    \openin\@inputcheck\reserved@b#1 %
    \ifeof\@inputcheck\else
      \edef\@filef@und{\reserved@b#1 }%
      \let\reserved@a\@firstoftwo%
      \closein\@inputcheck
      \@break@tfor
    \fi}%
  \reserved@a}
\long\def \InputIfFileExists#1#2{%
  \IfFileExists{#1}%
    {#2\@addtofilelist{#1}\@@input \@filef@und}}

\chardef\@inputcheck0

\let\@addtofilelist \@gobble

\def\@defaultunits{\afterassignment\remove@to@nnil}
\def\remove@to@nnil#1\@nnil{}

\newdimen\leftmarginv
\newdimen\leftmarginvi

\newdimen\@ovxx
\newdimen\@ovyy
\newdimen\@ovdx
\newdimen\@ovdy
\newdimen\@ovro
\newdimen\@ovri
\newdimen\@xdim
\newdimen\@ydim
\newdimen\@linelen
\newdimen\@dashdim

\long\def\mbox#1{\leavevmode\hbox{#1}}

\let\@onlypreamble\@gobble

\def\AtBeginDocument#1{#1}
\let\protect\relax

\newdimen\fboxsep
\newdimen\fboxrule

\fboxsep = 3pt
\fboxrule = .4pt

\def\@height{height} \def\@depth{depth} \def\@width{width}
\def\@minus{minus}
\def\@plus{plus}
\def\hb@xt@{\hbox to}

\long\def\@begin@tempboxa#1#2{%
   \begingroup
     \setbox\@tempboxa#1{\color@begingroup#2\color@endgroup}%
     \def\width{\wd\@tempboxa}%
     \def\height{\ht\@tempboxa}%
     \def\depth{\dp\@tempboxa}%
     \let\totalheight\@ovri
     \totalheight\height
     \advance\totalheight\depth}
\let\@end@tempboxa\endgroup

\let\set@color\relax
\let\color@begingroup\relax
\let\color@endgroup\relax
\let\color@setgroup\relax

\let\color@hbox\relax
\let\color@vbox\relax
\let\color@endbox\relax


\begingroup
  \catcode`P=12
  \catcode`T=12
  \lowercase{
    \def\x{\def\rem@pt##1.##2PT{##1\ifnum##2>\z@.##2\fi}}}
  \expandafter\endgroup\x
\def\strip@pt{\expandafter\rem@pt\the}


\def\@input#1{%
  \IfFileExists{#1}{\@@input\@filef@und}{\message{No file #1.}}}

\def\@warning{\immediate\write16}
\ifx\eplain\undefined
  \let\next\relax
\else
  \expandafter\let\expandafter\next\csname endinput\endcsname
\fi
\next
\ifx\ProvidesPackage\undefined
\def\next#1#2[#3]{\wlog{#2 #3}}
\expandafter\next\fi
\ProvidesPackage{ifpdf}
[2016/04/04 v3.0 Provides the ifpdf switch]
\expandafter\ifx\csname ifpdf\endcsname\relax
 \csname newif\expandafter\endcsname\csname ifpdf\endcsname
\else
 \ifx\pdftrue\undefined
  \ifx\PackageError\undefined
  \begingroup\def\PackageError#1#2#3{\endgroup\errmessage{#2}}
  \fi
  \PackageError{ifpdf}{incompatible ifpdf definition}{}
  \expandafter\expandafter\expandafter
 \fi
\fi
\let\ifpdf\iffalse
\ifx\directlua\undefined
\begingroup\expandafter\expandafter\expandafter\endgroup
\expandafter\ifx\csname pdfoutput\endcsname\relax
\else
  \ifnum\pdfoutput>0 %
    \pdftrue
  \fi
\fi
\else
\directlua{%
if (tex.outputmode or tex.pdfoutput or 0) > 0 then
  tex.print('\string\\pdftrue')
end
}
\fi
\def\makeactive#1{\catcode`#1 = \active \ignorespaces}%
\chardef\letter = 11
\chardef\other = 12
\def\makeatletter{%
  \edef\resetatcatcode{\catcode`\noexpand\@\the\catcode`\@\relax}%
  \catcode`\@11\relax
}%
\def\makeatother{%
  \edef\resetatcatcode{\catcode`\noexpand\@\the\catcode`\@\relax}%
  \catcode`\@12\relax
}%
\edef\leftdisplays{\the\catcode`@}%
\catcode`@ = \letter
\let\@eplainoldatcode = \leftdisplays
\toksdef\toks@ii = 2
\def\uncatcodespecials{%
   \def\do##1{\catcode`##1 = \other}%
   \dospecials
}%
{%
   \makeactive\^^M %
   \long\gdef\letreturn#1{\let^^M = #1}%
}%
\let\@eattoken = \relax  
\def\eattoken{\let\@eattoken = }%
\def\gobble#1{}%
\def\gobbletwo#1#2{}%
\def\gobblethree#1#2#3{}%
\def\@emptymarkA{\@emptymarkB}
\def\ifempty#1{\@@ifempty #1\@emptymarkA\@emptymarkB}%
\def\@@ifempty#1#2\@emptymarkB{\ifx #1\@emptymarkA}%
\def\@gobbleminus#1{\ifx-#1\else#1\fi}%
\def\ifinteger#1{\ifcat_\ifnum9<1\@gobbleminus#1 _\else A\fi}%
\def\isinteger{TT\fi\ifinteger}%
\def\@gobblemeaning#1:->{}%
\def\sanitize{\expandafter\@gobblemeaning\meaning}%
\def\ifundefined#1{\expandafter\ifx\csname#1\endcsname\relax}%
\def\csn#1{\csname#1\endcsname}%
\def\ece#1#2{\expandafter#1\csname#2\endcsname}%
\def\expandonce{\expandafter\noexpand}%
\let\@plainwlog = \wlog
\let\wlog = \gobble
\newlinechar = `^^J
\def\loggingall{\tracingcommands\tw@\tracingstats\tw@
   \tracingpages\@ne\tracingoutput\@ne\tracinglostchars\@ne
   \tracingmacros\tw@\tracingparagraphs\@ne\tracingrestores\@ne
   \showboxbreadth\maxdimen\showboxdepth\maxdimen
}%
\def\tracingoff{\tracingonline\z@\tracingcommands\z@\tracingstats\z@
  \tracingpages\z@\tracingoutput\z@\tracinglostchars\z@
  \tracingmacros\z@\tracingparagraphs\z@\tracingrestores\z@
  \showboxbreadth5 \showboxdepth3
}%
\begingroup
  \catcode`\{ = 12 \catcode`\} = 12
  \catcode`\[ = 1 \catcode`\] = 2
  \gdef\lbracechar[{]%
  \gdef\rbracechar[}]%
  \catcode`\% = \other
  \gdef\percentchar[
\def\vpenalty{\ifhmode\par\fi \penalty}%
\def\hpenalty{\ifvmode\leavevmode\fi \penalty}%
\def\iterate{%
  \let\loop@next\relax
  \body
  \let\loop@next\iterate
  \fi
  \loop@next
}%
\def\edefappend#1#2{%
  \toks@ = \expandafter{#1}%
  \edef#1{\the\toks@ #2}%
}%
\def\allowhyphens{\nobreak\hskip\z@skip}%
\long\def\hookprepend{\@hookassign{\the\toks@ii \the\toks@}}%
\long\def\hookappend{\@hookassign{\the\toks@ \the\toks@ii}}%
\let\hookaction = \hookappend 
\long\def\@hookassign#1#2#3{%
  \expandafter\ifx\csname @#2hook\endcsname \relax
    \toks@ = {}%
  \else
    \expandafter\let\expandafter\temp \csname @#2hook\endcsname
    \toks@ = \expandafter{\temp}%
  \fi
  \toks2 = {#3}
  \ece\edef{@#2hook}{#1}%
}%
\long\def\hookactiononce#1#2{%
  \edefappend#2{\global\let\noexpand#2\relax}
  \hookaction{#1}#2%
}%
\def\hookrun#1{%
  \expandafter\ifx\csname @#1hook\endcsname \relax \else
    \def\temp{\csname @#1hook\endcsname}%
    \expandafter\temp
  \fi
}%
\def\setpropertyglobal#1#2#3{\ece\xdef{#1@p#2}{#3}}%
\def\getproperty#1#2{%
  \expandafter\ifx\csname#1@p#2\endcsname\relax
  \else \csname#1@p#2\endcsname
  \fi
}%
\ifx\@undefinedmessage\@undefined
  \def\@undefinedmessage
    {No .aux file; I won't warn you about undefined labels.}%
\fi
\edef\cite{\the\catcode`@}%
\catcode`@ = 11
\let\@oldatcatcode = \cite
\chardef\@letter = 11
\chardef\@other = 12
\def\@innerdef#1#2{\edef#1{\expandafter\noexpand\csname #2\endcsname}}%
\@innerdef\@innernewcount{newcount}%
\@innerdef\@innernewdimen{newdimen}%
\@innerdef\@innernewif{newif}%
\@innerdef\@innernewwrite{newwrite}%
\def\@gobble#1{}%
\ifx\inputlineno\@undefined
   \let\@linenumber = \empty 
\else
   \def\@linenumber{\the\inputlineno:\space}%
\fi
\long\def\@futurenonspacelet#1{\def\cs{#1}%
   \afterassignment\@stepone\let\@nexttoken=
}%
\begingroup 
\def\\{\global\let\@stoken= }%
\\ 
\endgroup
\def\@stepone{\expandafter\futurelet\cs\@steptwo}%
\def\@steptwo{\expandafter\ifx\cs\@stoken\let\@@next=\@stepthree
   \else\let\@@next=\@nexttoken\fi \@@next}%
\def\@stepthree{\afterassignment\@stepone\let\@@next= }%
\def\@getoptionalarg#1{%
   \let\@optionalusercs = #1%
   \let\@optionalnext = \relax
   \@futurenonspacelet\@optionalnext\@bracketcheck
}%
\def\@bracketcheck{%
   \ifx [\@optionalnext
      \expandafter\@@getoptionalarg 
   \else
      \let\@optionalarg = \empty    
      \expandafter\@optionalusercs
   \fi
}%
\def\@@getoptionalarg[#1]{%
   \def\@optionalarg{#1}%
   \let\@optdummy=\relax 
   \@futurenonspacelet\@optdummy\@optionalusercs
}%
\def\@nnil{\@nil}%
\def\@fornoop#1\@@#2#3{}%
\def\@for#1:=#2\do#3{%
   \edef\@fortmp{#2}%
   \ifx\@fortmp\empty \else
      \expandafter\@forloop#2,\@nil,\@nil\@@#1{#3}%
   \fi
}%
\def\@forloop#1,#2,#3\@@#4#5{\def#4{#1}\ifx #4\@nnil \else
       #5\def#4{#2}\ifx #4\@nnil \else#5\@iforloop #3\@@#4{#5}\fi\fi
}%
\def\@iforloop#1,#2\@@#3#4{\def#3{#1}\ifx #3\@nnil
       \let\@nextwhile=\@fornoop \else
      #4\relax\let\@nextwhile=\@iforloop\fi\@nextwhile#2\@@#3{#4}%
}%
\@innernewif\if@fileexists
\def\@testfileexistence{\@getoptionalarg\@finishtestfileexistence}%
\def\@finishtestfileexistence#1{%
   \begingroup
      \def\extension{#1}%
      \immediate\openin0 =
         \ifx\@optionalarg\empty\jobname\else\@optionalarg\fi
         \ifx\extension\empty \else .#1\fi
         \space
      \ifeof 0
         \global\@fileexistsfalse
      \else
         \global\@fileexiststrue
      \fi
      \immediate\closein0
   \endgroup
}%
\toks0 = {%
\def\bibliographystyle#1{%
   \@readauxfile
   \@writeaux{\string\bibstyle{#1}}%
}%
\let\bibstyle = \@gobble
\let\bblfilebasename = \jobname
\def\bibliography#1{%
   \@readauxfile
   \@writeaux{\string\bibdata{#1}}%
   \@testfileexistence[\bblfilebasename]{bbl}%
   \if@fileexists
      \nobreak
      \@readbblfile
   \fi
}%
\let\bibdata = \@gobble
\def\nocite#1{%
   \@readauxfile
   \@writeaux{\string\citation{#1}}%
}%
\@innernewif\if@notfirstcitation
\def\cite{\@getoptionalarg\@cite}%
\def\@cite#1{%
   \let\@citenotetext = \@optionalarg
   \printcitestart
   \nocite{#1}%
   \@notfirstcitationfalse
   \@for \@citation :=#1\do
   {%
      \expandafter\@onecitation\@citation\@@
   }%
   \ifx\empty\@citenotetext\else
      \printcitenote{\@citenotetext}%
   \fi
   \printcitefinish
}%
\def\@onecitation#1\@@{%
   \if@notfirstcitation
      \printbetweencitations
   \fi
   \expandafter \ifx \csname\@citelabel{#1}\endcsname \relax
      \if@citewarning
         \message{\@linenumber Undefined citation `#1'.}%
      \fi
      \expandafter\gdef\csname\@citelabel{#1}\endcsname{%
         {\tt
            \escapechar = -1
            \nobreak\hskip0pt
            \expandafter\string\csname#1\endcsname
            \nobreak\hskip0pt
         }%
      }%
   \fi
   \printcitepreitem{#1}%
   \csname\@citelabel{#1}\endcsname
   \printcitepostitem
   \@notfirstcitationtrue
}%
\def\@citelabel#1{b@#1}%
\def\@citedef#1#2{{\expandafter\gdef\csname\@citelabel{#1}\endcsname{#2}}}%
\def\@readbblfile{%
   \ifx\@itemnum\@undefined
      \@innernewcount\@itemnum
   \fi
   \begingroup
      \ifx\begin\@undefined
         \def\begin##1##2{%
            \setbox0 = \hbox{\biblabelcontents{##2}}%
            \biblabelwidth = \wd0
         }%
         \let\end = \@gobble 
      \fi
      \@itemnum = 0
      \def\bibitem{\@getoptionalarg\@bibitem}%
      \def\@bibitem{%
         \ifx\@optionalarg\empty
            \expandafter\@numberedbibitem
         \else
            \expandafter\@alphabibitem
         \fi
      }%
      \def\@alphabibitem##1{%
         \expandafter \xdef\csname\@citelabel{##1}\endcsname {\@optionalarg}%
         \ifx\biblabelprecontents\@undefined
            \let\biblabelprecontents = \relax
         \fi
         \ifx\biblabelpostcontents\@undefined
            \let\biblabelpostcontents = \hss
         \fi
         \@finishbibitem{##1}%
      }%
      \def\@numberedbibitem##1{%
         \advance\@itemnum by 1
         \expandafter \xdef\csname\@citelabel{##1}\endcsname{\number\@itemnum}%
         \ifx\biblabelprecontents\@undefined
            \let\biblabelprecontents = \hss
         \fi
         \ifx\biblabelpostcontents\@undefined
            \let\biblabelpostcontents = \relax
         \fi
         \@finishbibitem{##1}%
      }%
      \def\@finishbibitem##1{%
         \bblitemhook{##1}%
         \biblabelprint{\csname\@citelabel{##1}\endcsname}%
         \@writeaux{\string\@citedef{##1}{\csname\@citelabel{##1}\endcsname}}%
         \ignorespaces
      }%
      \ifx\undefined\em \let\em=\bblem \fi
      \ifx\undefined\emph \let\emph=\bblemph \fi
      \ifx\undefined\mbox \let\mbox=\bblmbox \fi
      \ifx\undefined\newblock \let\newblock=\bblnewblock \fi
      \ifx\undefined\sc \let\sc=\bblsc \fi
      \ifx\undefined\textbf \let\textbf=\bbltextbf \fi
      \frenchspacing
      \clubpenalty = 4000 \widowpenalty = 4000
      \tolerance = 10000 \hfuzz = .5pt
      \everypar = {\hangindent = \biblabelwidth
                      \advance\hangindent by \biblabelextraspace}%
      \bblrm
      \parskip = 1.5ex plus .5ex minus .5ex
      \biblabelextraspace = .5em
      \bblhook
      \input \bblfilebasename.bbl
   \endgroup
}%
\@innernewdimen\biblabelwidth
\@innernewdimen\biblabelextraspace
\def\biblabelprint#1{%
   \noindent
   \hbox to \biblabelwidth{%
      \biblabelprecontents
      \biblabelcontents{#1}%
      \biblabelpostcontents
   }%
   \kern\biblabelextraspace
}%
\def\biblabelcontents#1{{\bblrm [#1]}}%
\def\bblrm{\rm}%
\def\bblem{\it}%
\def\bblemph#1{{\bblem #1\/}}
\def\bbltextbf#1{{\bf #1}}
\def\bblmbox{\leavevmode\hbox}
\def\bblsc{\ifx\@scfont\@undefined
              \font\@scfont = cmcsc10
           \fi
           \@scfont
}%
\def\bblnewblock{\hskip .11em plus .33em minus .07em }%
\let\bblhook = \empty
\let\bblitemhook = \@gobble
\def\printcitestart{[}
\def\printcitefinish{]}
\def\printbetweencitations{, }
\let\printcitepreitem\@gobble 
\let\printcitepostitem\empty
\def\printcitenote#1{, #1}
\let\citation = \@gobble
\@innernewcount\@btxnumparams
\ifx\newcommand\undefined
\long\def\newcommand#1{%
   \def\@btxcommandname{#1}%
   \@getoptionalarg\@btxcontinuenewcommand
}%
\fi
\ifx\renewcommand\undefined
  \let\renewcommand = \newcommand
\fi
\ifx\providecommand\undefined
\long\def\providecommand#1{%
   \def\@btxcommandname{#1}%
   \expandafter\ifx\@btxcommandname \@undefined
     \let\cs=\@continuenewcommand  
   \else
     \let\cs=\@gobble              
   \fi
   \@getoptionalarg\cs
}%
\fi
\def\@btxcontinuenewcommand{%
   \@btxnumparams = \ifx\@optionalarg\empty 0\else\@optionalarg \fi \relax
   \@btxnewcommand
}%
\def\@btxnewcommand#1{%
   \def\@btxstartdef{\expandafter\def\@btxcommandname}%
   \ifnum\@btxnumparams=0
      \let\@btxparamdef = \empty
   \else
      \ifnum\@btxnumparams>9
         \errmessage{\the\@btxnumparams\space is too many parameters}%
      \else
         \ifnum\@btxnumparams<0
            \errmessage{\the\@btxnumparams\space is too few parameters}%
         \else
            \edef\@btxparamdef{%
               \ifcase\@btxnumparams
                  \empty  No arguments.
               \or ####1%
               \or ####1####2%
               \or ####1####2####3%
               \or ####1####2####3####4%
               \or ####1####2####3####4####5%
               \or ####1####2####3####4####5####6%
               \or ####1####2####3####4####5####6####7%
               \or ####1####2####3####4####5####6####7####8%
               \or ####1####2####3####4####5####6####7####8####9%
               \fi
            }%
         \fi
      \fi
   \fi
   \expandafter\@btxstartdef\@btxparamdef{#1}%
}%
}%
\ifx\nobibtex\@undefined \the\toks0 \fi
\def\@readauxfile{%
   \if@auxfiledone \else 
      \global\@auxfiledonetrue
      \@testfileexistence{aux}%
      \if@fileexists
         \begingroup
            \endlinechar = -1
            \catcode`@ = 11
            \input \jobname.aux
         \endgroup
      \else
         \message{\@undefinedmessage}%
         \global\@citewarningfalse
      \fi
      \immediate\openout\@auxfile = \jobname.aux
   \fi
}%
\newif\if@auxfiledone
\ifx\noauxfile\@undefined \else \@auxfiledonetrue\fi
\@innernewwrite\@auxfile
\def\@writeaux#1{\ifx\noauxfile\@undefined \write\@auxfile{#1}\fi}%
\ifx\@undefinedmessage\@undefined
   \def\@undefinedmessage{No .aux file; I won't give you warnings about
                          undefined citations.}%
\fi
\@innernewif\if@citewarning
\ifx\noauxfile\@undefined \@citewarningtrue\fi
\catcode`@ = \@oldatcatcode
\let\auxfile = \@auxfile
\let\for = \@for
\let\futurenonspacelet = \@futurenonspacelet
\def\iffileexists{\if@fileexists}%
\let\innerdef = \@innerdef
\let\innernewcount = \@innernewcount
\let\innernewdimen = \@innernewdimen
\let\innernewif = \@innernewif
\let\innernewwrite = \@innernewwrite
\let\linenumber = \@linenumber
\let\readauxfile = \@readauxfile
\let\spacesub = \@spacesub
\let\testfileexistence = \@testfileexistence
\let\writeaux = \@writeaux
\def\innerinnerdef#1{\expandafter\innerdef\csname inner#1\endcsname{#1}}%
\innerinnerdef{newbox}%
\innerinnerdef{newfam}%
\innerinnerdef{newhelp}%
\innerinnerdef{newinsert}%
\innerinnerdef{newlanguage}%
\innerinnerdef{newmuskip}%
\innerinnerdef{newread}%
\innerinnerdef{newskip}%
\innerinnerdef{newtoks}%
\def\immediatewriteaux#1{%
  \ifx\noauxfile\@undefined
    \immediate\write\@auxfile{#1}%
  \fi
}%
\def\bblitemhook#1{\gdef\@hlbblitemlabel{#1}}%
\def\biblabelprint#1{%
   \noindent
   \hbox to \biblabelwidth{%
      \hldest@impl{bib}{\@hlbblitemlabel}%
      \biblabelprecontents
      \biblabelcontents{#1}%
      \biblabelpostcontents
   }%
   \kern\biblabelextraspace
}%
\def\eplainprintcitepreitem#1{\hlstart@impl{cite}{#1}}%
\def\eplainprintcitepostitem{\hlend@impl{cite}}%
\def\printcitepreitem#1{%
  \testfileexistence[\bblfilebasename]{bbl}%
  \iffileexists
    \global\let\printcitepreitem\eplainprintcitepreitem
    \global\let\printcitepostitem\eplainprintcitepostitem
  \else
    \global\let\printcitepreitem\gobble
    \global\let\printcitepostitem\relax
  \fi
  \printcitepreitem{#1}%
}%
\def\@Nnil{\@Nil}%
\def\@Fornoop#1\@@#2#3{}%
\def\For#1:=#2\do#3{%
   \edef\@Fortmp{#2}%
   \ifx\@Fortmp\empty \else
      \expandafter\@Forloop#2,\@Nil,\@Nil\@@#1{#3}%
   \fi
}%
\def\@Forloop#1,#2,#3\@@#4#5{\@Fordef#1\@@#4\ifx #4\@Nnil \else
       #5\@Fordef#2\@@#4\ifx #4\@Nnil \else#5\@iForloop #3\@@#4{#5}\fi\fi
}%
\def\@iForloop#1,#2\@@#3#4{\@Fordef#1\@@#3\ifx #3\@Nnil
       \let\@Nextwhile=\@Fornoop \else
      #4\relax\let\@Nextwhile=\@iForloop\fi\@Nextwhile#2\@@#3{#4}%
}%
\def\@Forspc{ }%
\def\@Fordef{\futurelet\@Fortmp\@@Fordef}
\def\@@Fordef{%
  \expandafter\ifx\@Forspc\@Fortmp 
    \expandafter\@Fortrim
  \else
    \expandafter\@@@Fordef
  \fi
}%
\expandafter\def\expandafter\@Fortrim\@Forspc#1\@@{\@Fordef#1\@@}%
\def\@@@Fordef#1\@@#2{\def#2{#1}}%
\def\tmpfileextension{.tmp}%
\let\tmpfilebasename = \jobname
\ifx\eTeXversion\undefined
  \innernewwrite\eplain@tmpfile
  \def\scantokens#1{%
    \toks@={#1}%
    \immediate\openout\eplain@tmpfile=\tmpfilebasename\tmpfileextension
    \immediate\write\eplain@tmpfile{\the\toks@}%
    \immediate\closeout\eplain@tmpfile
    \input \tmpfilebasename\tmpfileextension\relax
  }%
\fi
\begingroup
   \makeactive\^^M \makeactive\ 
\gdef\obeywhitespace{%
\makeactive\^^M\def^^M{\par\futurelet\next\@finishobeyedreturn}%
\makeactive\ \let =\ %
\aftergroup\@removebox%
\futurelet\next\@finishobeywhitespace%
}%
\gdef\@finishobeywhitespace{{%
\ifx\next %
\aftergroup\@obeywhitespaceloop%
\else\ifx\next^^M%
\aftergroup\gobble%
\fi\fi}}%
\gdef\@finishobeyedreturn{%
\ifx\next^^M\vskip\blanklineskipamount\fi%
\indent%
}%
\endgroup
\def\@obeywhitespaceloop#1{\futurelet\next\@finishobeywhitespace}%
\def\@removebox{%
  \ifhmode
    \setbox0 = \lastbox
    \ifdim\wd0=\parindent
      \setbox2 = \hbox{\unhcopy0}
      \ifdim\wd2=0pt
        \ignorespaces
      \else
        \box0 
      \fi
    \else
       \box0 
    \fi
  \fi
}%
\newskip\blanklineskipamount
\blanklineskipamount = 0pt
\def\frac#1/#2{\leavevmode
   \kern.1em \raise .5ex \hbox{\the\scriptfont0 #1}%
   \kern-.1em $/$%
   \kern-.15em \lower .25ex \hbox{\the\scriptfont0 #2}%
}%
\newdimen\hruledefaultheight  \hruledefaultheight = 0.4pt
\newdimen\hruledefaultdepth   \hruledefaultdepth = 0.0pt
\newdimen\vruledefaultwidth   \vruledefaultwidth = 0.4pt
\def\ehrule{\hrule height\hruledefaultheight depth\hruledefaultdepth}%
\def\evrule{\vrule width\vruledefaultwidth}%
\ifx\sc\undefined
    \def\sc{%
      \expandafter\ifx\the\scriptfont\fam\nullfont
        \font\temp = cmr7 \temp
      \else
        \the\scriptfont\fam
      \fi
      \def\uppercasesc{\char\uccode`}%
    }%
\fi
\ifx\uppercasesc\undefined
  \let\uppercasesc = \relax
\fi
\def\TeX{T\kern-.1667em\lower.5ex\hbox{E}\kern-.125emX\spacefactor1000 }%
\ifx\AmS\undefined
    \def\AmS{{\the\textfont2 A}\kern-.1667em\lower.5ex\hbox
        {\the\textfont2 M}\kern-.125em{\the\textfont2 S}}
\fi
\ifx\AMS\undefined \let\AMS=\AmS \fi
\ifx\AmSLaTeX\undefined
    \def\AmSLaTeX{\AmS-\LaTeX}
\fi
\ifx\AMSLaTeX\undefined \let\AMSLaTeX=\AmSLaTeX \fi
\ifx\AmSTeX\undefined
    \def\AmSTeX{$\cal A$\kern-.1667em\lower.5ex\hbox{$\cal M$}%
            \kern-.125em$\cal S$-\TeX}%
\fi
\ifx\AMSTEX\undefined \let\AMSTEX=\AmSTeX \fi
\ifx\AMSTeX\undefined \let\AMSTeX=\AmSTeX \fi
\ifx\BibTeX\undefined
    \def\BibTeX{B{\sc \uppercasesc i\kern-.025em \uppercasesc b}\kern-.08em
                \TeX}%
\fi
\ifx\BIBTeX\undefined \let\BIBTeX=\BibTeX \fi
\ifx\BIBTEX\undefined \let\BIBTEX=\BibTeX \fi
\ifx\LAMSTeX\undefined
    \def\LAMSTeX{L\raise.42ex\hbox{\kern-.3em\the\scriptfont2 A}%
                 \kern-.2em\lower.376ex\hbox{\the\textfont2 M}%
                 \kern-.125em {\the\textfont2 S}-\TeX}%
\fi
\ifx\LamSTeX\undefined \let\LamSTeX=\LAMSTeX \fi
\ifx\LAmSTeX\undefined \let\LAmSTeX=\LAMSTeX \fi
\ifx\LaTeX\undefined
    \def\LaTeX{L\kern-.36em\raise.3ex\hbox{\sc \uppercasesc a}\kern-.15em\TeX}%
\fi
\ifx\LATEX\undefined \let\LATEX=\LaTeX \fi
\ifx\LaTeXe\undefined
    \def\LaTeXe{\LaTeX{}\kern.05em2$_{\textstyle\varepsilon}$}
\fi
\ifx\MF\undefined
    \ifx\manfnt\undefined
            \font\manfnt=logo10
    \fi
    \ifx\manfntsl\undefined
            \font\manfntsl=logosl10
    \fi
    \def\MF{{\ifdim\fontdimen1\font>0pt \let\manfnt = \manfntsl \fi
      {\manfnt META}\-{\manfnt FONT}}\spacefactor1000 }%
\fi
\ifx\METAFONT\undefined \let\METAFONT=\MF \fi
\ifx\SLITEX\undefined
    \def\SLITEX{S\kern-.065em L\kern-.18em\raise.32ex\hbox{i}\kern-.03em\TeX}%
\fi
\ifx\SLiTeX\undefined \let\SLiTeX=\SLITEX \fi
\ifx\SliTeX\undefined \let\SliTeX=\SLITEX \fi
\ifx\SLITeX\undefined \let\SLITeX=\SLITEX \fi
\edef \path {\the \catcode `\@}%
\catcode `\@ = 11
\let \oldc@tcode = \path
\catcode `\@ = 11
\newcount \c@tcode
\newcount \c@unter
\newif \ifspecialpathdelimiters
\let \pathafterhook = \relax
\begingroup
\catcode `\ = 10
\gdef \passivesp@ce { }
\catcode `\ = 13\relax%
\gdef\activesp@ce{ }%
\endgroup
\def \discretionaries 
    {\begingroup
        \c@tcodes = 13
        \discr@tionaries
    }
\def \discr@tionaries #1
    {\def \discr@ti@naries ##1#1
         {\endgroup
          \def \discr@ti@n@ries ####1
              {\if   \noexpand ####1\noexpand #1%
                     \let \n@xt = \relax
               \else
                     \catcode `####1 = 13
                     \def ####1{\discretionary
                                  {\char `####1}{}{\char `####1}}%
                     \let \n@xt = \discr@ti@n@ries
               \fi
               \n@xt
              }%
          \def \discr@ti@n@ri@s {\discr@ti@n@ries ##1#1}%
         }%
     \discr@ti@naries
    }
\def \path
    {\ifspecialpathdelimiters
        \begingroup
        \c@tcodes = 12
        \def \endp@th {\endgroup \endgroup \pathafterhook}%
     \else
        \def \endp@th {\endgroup \pathafterhook}%
     \fi
     \p@th
    }
\def \p@th #1
    {\begingroup
        \tt
        \c@tcode = \catcode `#1
        \discr@ti@n@ri@s
        \catcode `\ = \active
        \expandafter \edef \activesp@ce {\passivesp@ce \hbox {}}%
        \catcode `#1 = \c@tcode
        \def \p@@th ##1#1
            {\leavevmode \hbox {}##1%
             \endp@th
            }%
     \p@@th
    }
\def \c@tcodes {\afterassignment \c@tc@des \c@tcode}
\def \c@tc@des
    {\c@unter = 0
     \loop
            \ifnum \catcode \c@unter = \c@tcode
            \else
                \catcode \c@unter = \c@tcode
            \fi
     \ifnum \c@unter < 255
            \advance \c@unter by 1
     \repeat
     \catcode `\ = 10
    }
\catcode `\@ = \oldc@tcode
\discretionaries |~!@$
\ifx\eTeX\undefined
  \def\eTeX{\hbox{\mathsurround=0pt $\varepsilon$-\kern-.125em\TeX}}%
\fi
\ifx\ExTeX\undefined
  \def\ExTeX{\hbox{\mathsurround=0pt
    $\textstyle\varepsilon_{\kern-0.15em\cal{X}}$\kern-.2em\TeX}}%
\fi
\def\eplain@Xe@reflect#1{%
  \ifx\reflectbox\undefined
    \errmessage{A graphics package must be loaded for \string\XeTeX}%
  \else
    \ifdim \fontdimen1\font>0pt
      \raise 1.75ex \hbox{\kern.1em\rotatebox{180}{#1}}\kern-.1em
    \else
      \reflectbox{#1}%
    \fi
  \fi
}%
\def\eplain@Xe#1{\leavevmode
  \smash{\hbox{X%
    \setbox0=\hbox{\TeX}\setbox2=\hbox{E}%
    \lower\dp0\hbox{\raise\dp2\hbox{\kern-.125em\eplain@Xe@reflect{E}}}%
    \kern-.1667em #1}}}%
\ifx\XeTeX\undefined
  \def\XeTeX{\eplain@Xe\TeX}%
\fi
\ifx\XeLaTeX\undefined
  \def\XeLaTeX{\eplain@Xe{\thinspace\LaTeX}}%
\fi
\def\blackbox{\vrule height .8ex width .6ex depth -.2ex \relax}
\def\makeblankbox#1#2{%
  \ifvoid0
    \errhelp = \@makeblankboxhelp
    \errmessage{Box 0 is void}%
  \fi
  \hbox{\lower\dp0
    \vbox{\hidehrule{#1}{#2}%
      \kern -#1
      \hbox to \wd0{\hidevrule{#1}{#2}%
        \raise\ht0\vbox to #1{}
        \lower\dp0\vtop to #1{}
        \hfil\hidevrule{#2}{#1}%
      }%
      \kern-#1\hidehrule{#2}{#1}%
    }%
  }%
}%
\newhelp\@makeblankboxhelp{Assigning to the dimensions of a void^^J%
  box has no effect.  Do `\string\setbox0=\string\null' before you^^J%
  define its dimensions.}%
\def\hidehrule#1#2{\kern-#1\hrule height#1 depth#2 \kern-#2}%
\def\hidevrule#1#2{%
  \kern-#1%
  \dimen@=#1\advance\dimen@ by #2%
  \vrule width\dimen@
  \kern-#2%
}%
\newdimen\boxitspace \boxitspace = 3pt
\long\def\boxit#1{%
  \vbox{%
    \ehrule
    \hbox{%
      \evrule
      \kern\boxitspace
      \vbox{\kern\boxitspace \parindent = 0pt #1\kern\boxitspace}%
      \kern\boxitspace
      \evrule
    }%
    \ehrule
  }%
}%
\def\numbername#1{\ifcase#1%
   zero%
   \or one%
   \or two%
   \or three%
   \or four%
   \or five%
   \or six%
   \or seven%
   \or eight%
   \or nine%
   \or ten%
   \or #1%
   \fi
}%
\let\@plainnewif = \newif
\let\@plainnewdimen = \newdimen
\let\newif = \innernewif
\let\newdimen = \innernewdimen
\edef\@eplainoldandcode{\the\catcode`& }%
\catcode`& = 11
\toks0 = {%
\edef\thinlines{\the\catcode`@ }%
\catcode`@ = 11
\let\@oldatcatcode = \thinlines
\def\smash@@{\relax 
  \ifmmode\def\next{\mathpalette\mathsm@sh}\else\let\next\makesm@sh
  \fi\next}
\def\makesm@sh#1{\setbox\z@\hbox{#1}\finsm@sh}
\def\mathsm@sh#1#2{\setbox\z@\hbox{$\m@th#1{#2}$}\finsm@sh}
\def\finsm@sh{\ht\z@\z@ \dp\z@\z@ \box\z@}
\edef\@oldandcatcode{\the\catcode`& }%
\catcode`& = 11
\def\&whilenoop#1{}%
\def\&whiledim#1\do #2{\ifdim #1\relax#2\&iwhiledim{#1\relax#2}\fi}%
\def\&iwhiledim#1{\ifdim #1\let\&nextwhile=\&iwhiledim
        \else\let\&nextwhile=\&whilenoop\fi\&nextwhile{#1}}%
\newif\if&negarg
\newdimen\&wholewidth
\newdimen\&halfwidth
\font\tenln=line10
\def\thinlines{\let\&linefnt\tenln \let\&circlefnt\tencirc
  \&wholewidth\fontdimen8\tenln \&halfwidth .5\&wholewidth}%
\def\thicklines{\let\&linefnt\tenlnw \let\&circlefnt\tencircw
  \&wholewidth\fontdimen8\tenlnw \&halfwidth .5\&wholewidth}%
\def\drawline(#1,#2)#3{\&xarg #1\relax \&yarg #2\relax \&linelen=#3\relax
  \ifnum\&xarg =0 \&vline \else \ifnum\&yarg =0 \&hline \else \&sline\fi\fi}%
\def\&sline{\leavevmode
  \ifnum\&xarg< 0 \&negargtrue \&xarg -\&xarg \&yyarg -\&yarg
  \else \&negargfalse \&yyarg \&yarg \fi
  \ifnum \&yyarg >0 \&tempcnta\&yyarg \else \&tempcnta -\&yyarg \fi
  \ifnum\&tempcnta>6 \&badlinearg \&yyarg0 \fi
  \ifnum\&xarg>6 \&badlinearg \&xarg1 \fi
  \setbox\&linechar\hbox{\&linefnt\&getlinechar(\&xarg,\&yyarg)}%
  \ifnum \&yyarg >0 \let\&upordown\raise \&clnht\z@
  \else\let\&upordown\lower \&clnht \ht\&linechar\fi
  \&clnwd=\wd\&linechar
  \&whiledim \&clnwd <\&linelen \do {%
    \&upordown\&clnht\copy\&linechar
    \advance\&clnht \ht\&linechar
    \advance\&clnwd \wd\&linechar
  }%
  \advance\&clnht -\ht\&linechar
  \advance\&clnwd -\wd\&linechar
  \&tempdima\&linelen\advance\&tempdima -\&clnwd
  \&tempdimb\&tempdima\advance\&tempdimb -\wd\&linechar
  \hskip\&tempdimb \multiply\&tempdima \@m
  \&tempcnta \&tempdima \&tempdima \wd\&linechar \divide\&tempcnta \&tempdima
  \&tempdima \ht\&linechar \multiply\&tempdima \&tempcnta
  \divide\&tempdima \@m
  \advance\&clnht \&tempdima
  \ifdim \&linelen <\wd\&linechar \hskip \wd\&linechar
  \else\&upordown\&clnht\copy\&linechar\fi}%
\def\&hline{\vrule height \&halfwidth depth \&halfwidth width \&linelen}%
\def\&getlinechar(#1,#2){\&tempcnta#1\relax\multiply\&tempcnta 8
  \advance\&tempcnta -9 \ifnum #2>0 \advance\&tempcnta #2\relax\else
  \advance\&tempcnta -#2\relax\advance\&tempcnta 64 \fi
  \char\&tempcnta}%
\def\drawvector(#1,#2)#3{\&xarg #1\relax \&yarg #2\relax
  \&tempcnta \ifnum\&xarg<0 -\&xarg\else\&xarg\fi
  \ifnum\&tempcnta<5\relax \&linelen=#3\relax
    \ifnum\&xarg =0 \&vvector \else \ifnum\&yarg =0 \&hvector
    \else \&svector\fi\fi\else\&badlinearg\fi}%
\def\&hvector{\ifnum\&xarg<0 \rlap{\&linefnt\&getlarrow(1,0)}\fi \&hline
  \ifnum\&xarg>0 \llap{\&linefnt\&getrarrow(1,0)}\fi}%
\def\&vvector{\ifnum \&yarg <0 \&downvector \else \&upvector \fi}%
\def\&svector{\&sline
  \&tempcnta\&yarg \ifnum\&tempcnta <0 \&tempcnta=-\&tempcnta\fi
  \ifnum\&tempcnta <5
    \if&negarg\ifnum\&yarg>0                   
      \llap{\lower\ht\&linechar\hbox to\&linelen{\&linefnt
        \&getlarrow(\&xarg,\&yyarg)\hss}}\else 
      \llap{\hbox to\&linelen{\&linefnt\&getlarrow(\&xarg,\&yyarg)\hss}}\fi
    \else\ifnum\&yarg>0                        
      \&tempdima\&linelen \multiply\&tempdima\&yarg
      \divide\&tempdima\&xarg \advance\&tempdima-\ht\&linechar
      \raise\&tempdima\llap{\&linefnt\&getrarrow(\&xarg,\&yyarg)}\else
      \&tempdima\&linelen \multiply\&tempdima-\&yarg 
      \divide\&tempdima\&xarg
      \lower\&tempdima\llap{\&linefnt\&getrarrow(\&xarg,\&yyarg)}\fi\fi
  \else\&badlinearg\fi}%
\def\&getlarrow(#1,#2){\ifnum #2 =\z@ \&tempcnta='33\else
\&tempcnta=#1\relax\multiply\&tempcnta \sixt@@n \advance\&tempcnta
-9 \&tempcntb=#2\relax\multiply\&tempcntb \tw@
\ifnum \&tempcntb >0 \advance\&tempcnta \&tempcntb\relax
\else\advance\&tempcnta -\&tempcntb\advance\&tempcnta 64
\fi\fi\char\&tempcnta}%
\def\&getrarrow(#1,#2){\&tempcntb=#2\relax
\ifnum\&tempcntb < 0 \&tempcntb=-\&tempcntb\relax\fi
\ifcase \&tempcntb\relax \&tempcnta='55 \or
\ifnum #1<3 \&tempcnta=#1\relax\multiply\&tempcnta
24 \advance\&tempcnta -6 \else \ifnum #1=3 \&tempcnta=49
\else\&tempcnta=58 \fi\fi\or
\ifnum #1<3 \&tempcnta=#1\relax\multiply\&tempcnta
24 \advance\&tempcnta -3 \else \&tempcnta=51\fi\or
\&tempcnta=#1\relax\multiply\&tempcnta
\sixt@@n \advance\&tempcnta -\tw@ \else
\&tempcnta=#1\relax\multiply\&tempcnta
\sixt@@n \advance\&tempcnta 7 \fi\ifnum #2<0 \advance\&tempcnta 64 \fi
\char\&tempcnta}%
\def\&vline{\ifnum \&yarg <0 \&downline \else \&upline\fi}%
\def\&upline{\hbox to \z@{\hskip -\&halfwidth \vrule width \&wholewidth
   height \&linelen depth \z@\hss}}%
\def\&downline{\hbox to \z@{\hskip -\&halfwidth \vrule width \&wholewidth
   height \z@ depth \&linelen \hss}}%
\def\&upvector{\&upline\setbox\&tempboxa\hbox{\&linefnt\char'66}\raise
     \&linelen \hbox to\z@{\lower \ht\&tempboxa\box\&tempboxa\hss}}%
\def\&downvector{\&downline\lower \&linelen
      \hbox to \z@{\&linefnt\char'77\hss}}%
\def\&badlinearg{\errmessage{Bad \string\arrow\space argument.}}%
\thinlines
\countdef\&xarg     0
\countdef\&yarg     2
\countdef\&yyarg    4
\countdef\&tempcnta 6
\countdef\&tempcntb 8
\dimendef\&linelen  0
\dimendef\&clnwd    2
\dimendef\&clnht    4
\dimendef\&tempdima 6
\dimendef\&tempdimb 8
\chardef\@arrbox    0
\chardef\&linechar  2
\chardef\&tempboxa  2           
\let\lft^%
\let\rt_
\newif\if@pslope 
\def\@findslope(#1,#2){\ifnum#1>0
  \ifnum#2>0 \@pslopetrue \else\@pslopefalse\fi \else
  \ifnum#2>0 \@pslopefalse \else\@pslopetrue\fi\fi}%
\def\generalsmap(#1,#2){\getm@rphposn(#1,#2)\plnmorph\futurelet\next\addm@rph}%
\def\sline(#1,#2){\setbox\@arrbox=\hbox{\drawline(#1,#2){\sarrowlength}}%
  \@findslope(#1,#2)\d@@blearrfalse\generalsmap(#1,#2)}%
\def\arrow(#1,#2){\setbox\@arrbox=\hbox{\drawvector(#1,#2){\sarrowlength}}%
  \@findslope(#1,#2)\d@@blearrfalse\generalsmap(#1,#2)}%
\newif\ifd@@blearr
\def\bisline(#1,#2){\@findslope(#1,#2)%
  \if@pslope \let\@upordown\raise \else \let\@upordown\lower\fi
  \getch@nnel(#1,#2)\setbox\@arrbox=\hbox{\@upordown\@vchannel
    \rlap{\drawline(#1,#2){\sarrowlength}}%
      \hskip\@hchannel\hbox{\drawline(#1,#2){\sarrowlength}}}%
  \d@@blearrtrue\generalsmap(#1,#2)}%
\def\biarrow(#1,#2){\@findslope(#1,#2)%
  \if@pslope \let\@upordown\raise \else \let\@upordown\lower\fi
  \getch@nnel(#1,#2)\setbox\@arrbox=\hbox{\@upordown\@vchannel
    \rlap{\drawvector(#1,#2){\sarrowlength}}%
      \hskip\@hchannel\hbox{\drawvector(#1,#2){\sarrowlength}}}%
  \d@@blearrtrue\generalsmap(#1,#2)}%
\def\adjarrow(#1,#2){\@findslope(#1,#2)%
  \if@pslope \let\@upordown\raise \else \let\@upordown\lower\fi
  \getch@nnel(#1,#2)\setbox\@arrbox=\hbox{\@upordown\@vchannel
    \rlap{\drawvector(#1,#2){\sarrowlength}}%
      \hskip\@hchannel\hbox{\drawvector(-#1,-#2){\sarrowlength}}}%
  \d@@blearrtrue\generalsmap(#1,#2)}%
\newif\ifrtm@rph
\def\@shiftmorph#1{\hbox{\setbox0=\hbox{$\scriptstyle#1$}%
  \setbox1=\hbox{\hskip\@hm@rphshift\raise\@vm@rphshift\copy0}%
  \wd1=\wd0 \ht1=\ht0 \dp1=\dp0 \box1}}%
\def\@hm@rphshift{\ifrtm@rph
  \ifdim\hmorphposnrt=\z@\hmorphposn\else\hmorphposnrt\fi \else
  \ifdim\hmorphposnlft=\z@\hmorphposn\else\hmorphposnlft\fi \fi}%
\def\@vm@rphshift{\ifrtm@rph
  \ifdim\vmorphposnrt=\z@\vmorphposn\else\vmorphposnrt\fi \else
  \ifdim\vmorphposnlft=\z@\vmorphposn\else\vmorphposnlft\fi \fi}%
\def\addm@rph{\ifx\next\lft\let\temp=\lftmorph\else
  \ifx\next\rt\let\temp=\rtmorph\else\let\temp\relax\fi\fi \temp}%
\def\plnmorph{\dimen1\wd\@arrbox \ifdim\dimen1<\z@ \dimen1-\dimen1\fi
  \vcenter{\box\@arrbox}}%
\def\lftmorph\lft#1{\rtm@rphfalse \setbox0=\@shiftmorph{#1}%
  \if@pslope \let\@upordown\raise \else \let\@upordown\lower\fi
  \llap{\@upordown\@vmorphdflt\hbox to\dimen1{\hss 
    \llap{\box0}\hss}\hskip\@hmorphdflt}\futurelet\next\addm@rph}%
\def\rtmorph\rt#1{\rtm@rphtrue \setbox0=\@shiftmorph{#1}%
  \if@pslope \let\@upordown\lower \else \let\@upordown\raise\fi
  \llap{\@upordown\@vmorphdflt\hbox to\dimen1{\hss
    \rlap{\box0}\hss}\hskip-\@hmorphdflt}\futurelet\next\addm@rph}%
\def\getm@rphposn(#1,#2){\ifd@@blearr \dimen@\morphdist \advance\dimen@ by
  .5\channelwidth \@getshift(#1,#2){\@hmorphdflt}{\@vmorphdflt}{\dimen@}\else
  \@getshift(#1,#2){\@hmorphdflt}{\@vmorphdflt}{\morphdist}\fi}%
\def\getch@nnel(#1,#2){\ifdim\hchannel=\z@ \ifdim\vchannel=\z@
    \@getshift(#1,#2){\@hchannel}{\@vchannel}{\channelwidth}%
    \else \@hchannel\hchannel \@vchannel\vchannel \fi
  \else \@hchannel\hchannel \@vchannel\vchannel \fi}%
\def\@getshift(#1,#2)#3#4#5{\dimen@ #5\relax
  \&xarg #1\relax \&yarg #2\relax
  \ifnum\&xarg<0 \&xarg -\&xarg \fi
  \ifnum\&yarg<0 \&yarg -\&yarg \fi
  \ifnum\&xarg<\&yarg \&negargtrue \&yyarg\&xarg \&xarg\&yarg \&yarg\&yyarg\fi
  \ifcase\&xarg \or  
    \ifcase\&yarg    
      \dimen@i \z@ \dimen@ii \dimen@ \or 
      \dimen@i .7071\dimen@ \dimen@ii .7071\dimen@ \fi \or
    \ifcase\&yarg    
      \or 
      \dimen@i .4472\dimen@ \dimen@ii .8944\dimen@ \fi \or
    \ifcase\&yarg    
      \or 
      \dimen@i .3162\dimen@ \dimen@ii .9486\dimen@ \or
      \dimen@i .5547\dimen@ \dimen@ii .8321\dimen@ \fi \or
    \ifcase\&yarg    
      \or 
      \dimen@i .2425\dimen@ \dimen@ii .9701\dimen@ \or\or
      \dimen@i .6\dimen@ \dimen@ii .8\dimen@ \fi \or
    \ifcase\&yarg    
      \or 
      \dimen@i .1961\dimen@ \dimen@ii .9801\dimen@ \or
      \dimen@i .3714\dimen@ \dimen@ii .9284\dimen@ \or
      \dimen@i .5144\dimen@ \dimen@ii .8575\dimen@ \or
      \dimen@i .6247\dimen@ \dimen@ii .7801\dimen@ \fi \or
    \ifcase\&yarg    
      \or 
      \dimen@i .1645\dimen@ \dimen@ii .9864\dimen@ \or\or\or\or
      \dimen@i .6402\dimen@ \dimen@ii .7682\dimen@ \fi \fi
  \if&negarg \&tempdima\dimen@i \dimen@i\dimen@ii \dimen@ii\&tempdima\fi
  #3\dimen@i\relax #4\dimen@ii\relax }%
\catcode`\&=4  
}%
\catcode`& = 4
\toks2 = {%
\catcode`\&=4  
\def\generalhmap{\futurelet\next\@generalhmap}%
\def\@generalhmap{\ifx\next^ \let\temp\generalhm@rph\else
  \ifx\next_ \let\temp\generalhm@rph\else \let\temp\m@kehmap\fi\fi \temp}%
\def\generalhm@rph#1#2{\ifx#1^
    \toks@=\expandafter{\the\toks@#1{\rtm@rphtrue\@shiftmorph{#2}}}\else
    \toks@=\expandafter{\the\toks@#1{\rtm@rphfalse\@shiftmorph{#2}}}\fi
  \generalhmap}%
\def\m@kehmap{\mathrel{\smash@@{\the\toks@}}}%
\def\mapright{\toks@={\mathop{\vcenter{\smash@@{\drawrightarrow}}}\limits}%
  \generalhmap}%
\def\mapleft{\toks@={\mathop{\vcenter{\smash@@{\drawleftarrow}}}\limits}%
  \generalhmap}%
\def\bimapright{\toks@={\mathop{\vcenter{\smash@@{\drawbirightarrow}}}\limits}%
  \generalhmap}%
\def\bimapleft{\toks@={\mathop{\vcenter{\smash@@{\drawbileftarrow}}}\limits}%
  \generalhmap}%
\def\adjmapright{\toks@={\mathop{\vcenter{\smash@@{\drawadjrightarrow}}}\limits}%
  \generalhmap}%
\def\adjmapleft{\toks@={\mathop{\vcenter{\smash@@{\drawadjleftarrow}}}\limits}%
  \generalhmap}%
\def\hline{\toks@={\mathop{\vcenter{\smash@@{\drawhline}}}\limits}%
  \generalhmap}%
\def\bihline{\toks@={\mathop{\vcenter{\smash@@{\drawbihline}}}\limits}%
  \generalhmap}%
\def\drawrightarrow{\hbox{\drawvector(1,0){\harrowlength}}}%
\def\drawleftarrow{\hbox{\drawvector(-1,0){\harrowlength}}}%
\def\drawbirightarrow{\hbox{\raise.5\channelwidth
  \hbox{\drawvector(1,0){\harrowlength}}\lower.5\channelwidth
  \llap{\drawvector(1,0){\harrowlength}}}}%
\def\drawbileftarrow{\hbox{\raise.5\channelwidth
  \hbox{\drawvector(-1,0){\harrowlength}}\lower.5\channelwidth
  \llap{\drawvector(-1,0){\harrowlength}}}}%
\def\drawadjrightarrow{\hbox{\raise.5\channelwidth
  \hbox{\drawvector(-1,0){\harrowlength}}\lower.5\channelwidth
  \llap{\drawvector(1,0){\harrowlength}}}}%
\def\drawadjleftarrow{\hbox{\raise.5\channelwidth
  \hbox{\drawvector(1,0){\harrowlength}}\lower.5\channelwidth
  \llap{\drawvector(-1,0){\harrowlength}}}}%
\def\drawhline{\hbox{\drawline(1,0){\harrowlength}}}%
\def\drawbihline{\hbox{\raise.5\channelwidth
  \hbox{\drawline(1,0){\harrowlength}}\lower.5\channelwidth
  \llap{\drawline(1,0){\harrowlength}}}}%
\def\generalvmap{\futurelet\next\@generalvmap}%
\def\@generalvmap{\ifx\next\lft \let\temp\generalvm@rph\else
  \ifx\next\rt \let\temp\generalvm@rph\else \let\temp\m@kevmap\fi\fi \temp}%
\toksdef\toks@@=1
\def\generalvm@rph#1#2{\ifx#1\rt 
    \toks@=\expandafter{\the\toks@
      \rlap{$\vcenter{\rtm@rphtrue\@shiftmorph{#2}}$}}\else 
    \toks@@={\llap{$\vcenter{\rtm@rphfalse\@shiftmorph{#2}}$}}%
    \toks@=\expandafter\expandafter\expandafter{\expandafter\the\expandafter
      \toks@@ \the\toks@}\fi \generalvmap}%
\def\m@kevmap{\the\toks@}%
\def\mapdown{\toks@={\vcenter{\drawdownarrow}}\generalvmap}%
\def\mapup{\toks@={\vcenter{\drawuparrow}}\generalvmap}%
\def\bimapdown{\toks@={\vcenter{\drawbidownarrow}}\generalvmap}%
\def\bimapup{\toks@={\vcenter{\drawbiuparrow}}\generalvmap}%
\def\adjmapdown{\toks@={\vcenter{\drawadjdownarrow}}\generalvmap}%
\def\adjmapup{\toks@={\vcenter{\drawadjuparrow}}\generalvmap}%
\def\vline{\toks@={\vcenter{\drawvline}}\generalvmap}%
\def\bivline{\toks@={\vcenter{\drawbivline}}\generalvmap}%
\def\drawdownarrow{\hbox to5pt{\hss\drawvector(0,-1){\varrowlength}\hss}}%
\def\drawuparrow{\hbox to5pt{\hss\drawvector(0,1){\varrowlength}\hss}}%
\def\drawbidownarrow{\hbox to5pt{\hss\hbox{\drawvector(0,-1){\varrowlength}}%
  \hskip\channelwidth\hbox{\drawvector(0,-1){\varrowlength}}\hss}}%
\def\drawbiuparrow{\hbox to5pt{\hss\hbox{\drawvector(0,1){\varrowlength}}%
  \hskip\channelwidth\hbox{\drawvector(0,1){\varrowlength}}\hss}}%
\def\drawadjdownarrow{\hbox to5pt{\hss\hbox{\drawvector(0,-1){\varrowlength}}%
  \hskip\channelwidth\lower\varrowlength
  \hbox{\drawvector(0,1){\varrowlength}}\hss}}%
\def\drawadjuparrow{\hbox to5pt{\hss\hbox{\drawvector(0,1){\varrowlength}}%
  \hskip\channelwidth\raise\varrowlength
  \hbox{\drawvector(0,-1){\varrowlength}}\hss}}%
\def\drawvline{\hbox to5pt{\hss\drawline(0,1){\varrowlength}\hss}}%
\def\drawbivline{\hbox to5pt{\hss\hbox{\drawline(0,1){\varrowlength}}%
  \hskip\channelwidth\hbox{\drawline(0,1){\varrowlength}}\hss}}%
\def\commdiag#1{\null\,
  \vcenter{\commdiagbaselines
  \m@th\ialign{\hfil$##$\hfil&&\hfil$\mkern4mu ##$\hfil\crcr
      \mathstrut\crcr\noalign{\kern-\baselineskip}
      #1\crcr\mathstrut\crcr\noalign{\kern-\baselineskip}}}\,}%
\def\commdiagbaselines{\baselineskip15pt \lineskip3pt \lineskiplimit3pt }%
\def\gridcommdiag#1{\null\,
  \vcenter{\offinterlineskip
  \m@th\ialign{&\vbox to\vgrid{\vss
    \hbox to\hgrid{\hss\smash@@{$##$}\hss}}\crcr
      \mathstrut\crcr\noalign{\kern-\vgrid}
      #1\crcr\mathstrut\crcr\noalign{\kern-.5\vgrid}}}\,}%
\newdimen\harrowlength \harrowlength=60pt
\newdimen\varrowlength \varrowlength=.618\harrowlength
\newdimen\sarrowlength \sarrowlength=\harrowlength
\newdimen\hmorphposn \hmorphposn=\z@
\newdimen\vmorphposn \vmorphposn=\z@
\newdimen\morphdist  \morphdist=4pt
\dimendef\@hmorphdflt 0       
\dimendef\@vmorphdflt 2       
\newdimen\hmorphposnrt  \hmorphposnrt=\z@
\newdimen\hmorphposnlft \hmorphposnlft=\z@
\newdimen\vmorphposnrt  \vmorphposnrt=\z@
\newdimen\vmorphposnlft \vmorphposnlft=\z@

\newdimen\hgrid \hgrid=15pt
\newdimen\vgrid \vgrid=15pt
\newdimen\hchannel  \hchannel=0pt
\newdimen\vchannel  \vchannel=0pt
\newdimen\channelwidth \channelwidth=3pt
\dimendef\@hchannel 0         
\dimendef\@vchannel 2         
\catcode`& = \@oldandcatcode
\catcode`@ = \@oldatcatcode
}%
\let\newif = \@plainnewif
\let\newdimen = \@plainnewdimen
\ifx\noarrow\@undefined \the\toks0 \the\toks2 \fi
\catcode`& = \@eplainoldandcode
\def\environment#1{%
   \ifx\@groupname\@undefined\else
      \errhelp = \@unnamedendgrouphelp
      \errmessage{`\@groupname' was not closed by \string\endenvironment}%
   \fi
   \edef\@groupname{#1}%
   \begingroup
      \let\@groupname = \@undefined
}%
\def\endenvironment#1{%
   \endgroup
   \edef\@thearg{#1}%
   \ifx\@groupname\@thearg
   \else
      \ifx\@groupname\@undefined
         \errhelp = \@isolatedendenvironmenthelp
         \errmessage{Isolated \string\endenvironment\space for `#1'}%
      \else
         \errhelp = \@mismatchedenvironmenthelp
         \errmessage{Environment `#1' ended, but `\@groupname' started}%
         \endgroup 
      \fi
   \fi
   \let\@groupname = \@undefined
}%
\newhelp\@unnamedendgrouphelp{Most likely, you just forgot an^^J%
   \string\endenvironment.  Maybe you should try inserting another^^J%
   \string\endgroup to recover.}%
\newhelp\@isolatedendenvironmenthelp{You ended an environment X, but^^J%
   no \string\environment{X} to start it is anywhere in sight.^^J%
   You might also be at an \string\endenvironment\space that would match^^J%
   a \string\begingroup, i.e., you forgot an \string\endgroup.}%
\newhelp\@mismatchedenvironmenthelp{You started an environment named X, but^^J%
   you ended one named Y.  Maybe you made a typo in one^^J%
   or the other of the names?}%
\newif\ifenvironment
\def\checkenv{\ifenvironment \errhelp = \@interwovenenvhelp
   \errmessage{Interwoven environments}%
   \egroup \fi
}%
\newhelp\@interwovenenvhelp{Perhaps you forgot to end the previous^^J%
   environment? I'm finishing off the current group,^^J%
   hoping that will fix it.}%
\newtoks\previouseverydisplay
\let\@leftleftfill\relax 
\newdimen\leftdisplayindent \leftdisplayindent=\parindent
\newif\if@leftdisplays
\def\leftdisplays{%
  \if@leftdisplays\else
    \previouseverydisplay = \everydisplay
    \everydisplay = {\the\previouseverydisplay \leftdisplaysetup}%
    \let\@save@maybedisableeqno = \@maybedisableeqno
    \let\@saveeqno = \eqno
    \let\@saveleqno = \leqno
    \let\@saveeqalignno = \eqalignno
    \let\@saveleqalignno = \leqalignno
    \let\@maybedisableeqno = \relax
    \def\eqno{\hfill\textstyle\enspace}%
    \def\leqno{%
      \hfill
      \hbox to0pt\bgroup
        \kern-\displaywidth
        \kern-\leftdisplayindent    
        $\aftergroup\@leftleqnoend  
    }%
    \@redefinealignmentdisplays
    \@leftdisplaystrue
  \fi
}%
\newbox\@lignbox
\newdimen\disprevdepth
\def\centereddisplays{%
  \if@leftdisplays
    \everydisplay = \previouseverydisplay
    \let\@maybedisableeqno = \@save@maybedisableeqno
    \let\eqno = \@saveeqno
    \let\leqno = \@saveleqno
    \let\eqalignno = \@saveeqalignno
    \let\leqalignno = \@saveleqalignno
    \@leftdisplaysfalse
  \fi
}%
\def\leftdisplaysetup{%
   \dimen@ = \leftdisplayindent
   \advance\dimen@ by \leftskip
   \advance\displayindent by \dimen@
   \advance\displaywidth by -\dimen@
   \halign\bgroup##\cr \noalign\bgroup
      \disprevdepth = \prevdepth
      \setbox\z@ = \hbox to\displaywidth\bgroup
      $\displaystyle
      \aftergroup\@lefteqend 
}
\def\@lefteqend{
   \hfil\egroup
   \@putdisplay}
\def\@leftleqnoend{\hss \egroup $}
\def\@putdisplay{%
   \ifvoid\@lignbox 
     \moveright\displayindent\box\z@
   \else 
     \prevdepth = \dp\@lignbox 
     \unvbox\@lignbox
   \fi
   \egroup\egroup 
   $
}
\def\@redefinealignmentdisplays{%
  \def\displaylines##1{
    \global\setbox\@lignbox\vbox{%
      \prevdepth = \disprevdepth
      \displ@y
      \tabskip\displayindent
      \halign{\hbox to\displaywidth
        {$\@lign\displaystyle####\hfil$\hfil}\crcr
              ##1\crcr}}}%
  \def\eqalignno##1{%
    \def\eqno{&}%
    \def\leqno{&}%
    \global\setbox\@lignbox\vbox{%
      \prevdepth = \disprevdepth
      \displ@y
      \advance\displaywidth by \displayindent
      \tabskip\displayindent
      \halign to\displaywidth{%
         \hfil $\@lign\displaystyle{####}$\@leftleftfill\tabskip\z@skip
        &$\@lign\displaystyle{{}####}$\hfil\tabskip\centering
        &\llap{$\@lign####$}\tabskip\z@skip\crcr
        ##1\crcr}}}%
  \def\leqalignno##1{%
    \def\eqno{&}%
    \def\leqno{&}%
    \global\setbox\@lignbox\vbox{%
      \prevdepth = \disprevdepth
      \displ@y
      \advance\displaywidth by \displayindent
      \tabskip\displayindent
      \halign to\displaywidth{%
         \hfil $\@lign\displaystyle{####}$\@leftleftfill\tabskip\z@skip
        &$\@lign\displaystyle{{}####}$\hfil\tabskip\centering
        &\kern-\displaywidth
         \rlap{\kern\displayindent \kern-\leftdisplayindent$\@lign####$}%
         \tabskip\displaywidth\crcr
        ##1\crcr}}}%
}%
\let\@primitivenoalign = \noalign
\newtoks\@everynoalign
\def\@lefteqalignonoalign#1{%
  \@primitivenoalign{%
    \advance\leftskip by -\parindent
    \advance\leftskip by -\leftdisplayindent
    \parskip = 0pt
    \parindent = 0pt
    \the\@everynoalign
    #1%
  }%
}%
\def\monthname{%
   \ifcase\month
      \or Jan\or Feb\or Mar\or Apr\or May\or Jun%
      \or Jul\or Aug\or Sep\or Oct\or Nov\or Dec%
   \fi
}%
\def\fullmonthname{%
   \ifcase\month
      \or January\or February\or March\or April\or May\or June%
      \or July\or August\or September\or October\or November\or December%
   \fi
}%
\def\timestring{\begingroup
   \count0 = \time
   \divide\count0 by 60
   \count2 = \count0   
   \count4 = \time
   \multiply\count0 by 60
   \advance\count4 by -\count0   
   \ifnum\count4<10
      \toks1 = {0}%
   \else
      \toks1 = {}%
   \fi
   \ifnum\count2<12
      \toks0 = {a.m.}%
   \else
      \toks0 = {p.m.}%
      \advance\count2 by -12
   \fi
   \ifnum\count2=0
      \count2 = 12
   \fi
   \number\count2:\the\toks1 \number\count4 \thinspace \the\toks0
\endgroup}%
\def\today{\the\day\ \fullmonthname\ \the\year}%
\newskip\abovelistskipamount      \abovelistskipamount = .5\baselineskip
  \newcount\abovelistpenalty      \abovelistpenalty    = 10000
  \def\abovelistskip{\vpenalty\abovelistpenalty \vskip\abovelistskipamount}%
\newskip\interitemskipamount      \interitemskipamount = 0pt
  \newcount\belowlistpenalty      \belowlistpenalty    = -50
\newskip\belowlistskipamount      \belowlistskipamount = .5\baselineskip
  \newcount\interitempenalty      \interitempenalty    = 0
  \def\interitemskip{\vpenalty\interitempenalty \vskip\interitemskipamount}%
\newdimen\listleftindent    \listleftindent = 0pt
\newdimen\listrightindent   \listrightindent = 0pt
\let\listmarkerspace = \enspace
\newtoks\everylist
\newdimen\@listindent
\def\beginlist{%
  \abovelistskip
  \@listindent = \parindent
  \advance\@listindent by \listleftindent
  \advance\leftskip by \@listindent
  \advance\rightskip by \listrightindent
  \itemnumber = 1
  \the\everylist
}%
\def\li{\@getoptionalarg\@finli}%
\def\@finli{%
  \let\@lioptarg\@optionalarg
  \ifx\@lioptarg\empty \else
    \begingroup
      \@@hldestoff
      \expandafter\writeitemxref\expandafter{\@lioptarg}%
    \endgroup
  \fi
  \ifnum\itemnumber=1 \else \interitemskip \fi
  \begingroup
    \ifx\@lioptarg\empty \else
      \toks@ = \expandafter{\@lioptarg}%
      \let\li@nohldest@marker\marker
      \edef\marker{\noexpand\hldest@impl{li}{\the\toks@}\noexpand\li@nohldest@marker}%
    \fi
    \printitem
  \endgroup
  \advance\itemnumber by 1
  \advance\itemletter by 1
  \advance\itemromannumeral by 1
  \ignorespaces
}%
\def\writeitemxref#1{\definexref{#1}\marker{item}}%
\def\printitem{%
  \par
  \nobreak
  \vskip-\parskip
  \noindent
  \printmarker\marker
}%
\def\printmarker#1{\llap{\marker \enspace}}%
\newcount\numberedlistdepth
\newcount\itemnumber
\newcount\itemletter
\newcount\itemromannumeral
\def\numberedmarker{%
  \ifcase\numberedlistdepth
      (impossible)%
  \or \printitemnumber
  \or \printitemletter
  \or \printitemromannumeral
  \else *%
  \fi
}%
\def\printitemnumber{\number\itemnumber}%
\def\printitemletter{\char\the\itemletter}%
\def\printitemromannumeral{\romannumeral\itemromannumeral}%
\def\numberedprintmarker#1{\llap{#1) \listmarkerspace}}%
\def\numberedlist{\environment{@numbered-list}%
  \advance\numberedlistdepth by 1
  \itemletter = `a
  \itemromannumeral = 1
  \beginlist
  \let\marker = \numberedmarker
  \let\printmarker = \numberedprintmarker
}%

\newcount\unorderedlistdepth
\def\unorderedmarker{%
  \ifcase\unorderedlistdepth
      (impossible)%
  \or \blackbox
  \or ---%
  \else *%
  \fi
}%
\def\unorderedprintmarker#1{\llap{#1\listmarkerspace}}%
\def\unorderedlist{\environment{@unordered-list}%
  \advance\unorderedlistdepth by 1
  \beginlist
  \let\marker = \unorderedmarker
  \let\printmarker = \unorderedprintmarker
}%
\def\listing#1{%
   \par \begingroup
   \@setuplisting
   \setuplistinghook
   \input #1
   \endgroup
}%
\let\setuplistinghook = \relax
\def\linenumberedlisting{%
  \ifx\lineno\undefined \innernewcount\lineno \fi
  \lineno = 0
  \everypar = {\advance\lineno by 1 \printlistinglineno}%
}%
\def\printlistinglineno{\llap{[\the\lineno]\quad}}%
\def\nolastlinelisting{\aftergroup\@removeboxes}%
\def\@removeboxes{%
  \setbox0 = \lastbox
  \ifvoid0
    \ignorespaces 
  \else
    \expandafter\@removeboxes
  \fi
}%
{%
  \makeactive\^^L
  \let^^L = \relax
  \gdef\@setuplisting{%
     \uncatcodespecials
     \obeywhitespace
     \makeactive\`
     \makeactive\^^I
     \makeactive\^^L
     \def^^L{\vfill\break}%
     \parskip = 0pt
     \listingfont
  }%
}%
\def\listingfont{\tt}%
{%
   \makeactive\`
   \gdef`{\relax\lq}
}%
{%
   \makeactive\^^I
   \gdef^^I{\hskip8\fontdimen2}%
}%
\def\verbatimescapechar#1{%
  \gdef\@makeverbatimescapechar{%
    \@makeverbatimdoubleescape #1%
    \catcode`#1 = 0
  }%
}%
\def\@makeverbatimdoubleescape#1{%
  \catcode`#1 = \other
  \begingroup
    \lccode`\* = `#1%
    \lowercase{\endgroup \ece\def*{*}}%
}%
\verbatimescapechar\|  
\def\verbatim{\begingroup
  \uncatcodespecials
  \makeactive\` 
  \@makeverbatimescapechar
  \tt\obeywhitespace}

\def\definecontentsfile#1{%
  \ece\innernewwrite{#1file}%
  \ece\innernewif{if@#1fileopened}%
  \ece\let{#1filebasename} = \jobname
  \ece\def{open#1file}{\opencontentsfile{#1}}%
  \ece\def{write#1entry}{\writecontentsentry{#1}}%
  \ece\def{writenumbered#1entry}{\writenumberedcontentsentry{#1}}%
  \ece\def{writenumbered#1line}{\writenumberedcontentsline{#1}}%
  \ece\innernewif{ifrewrite#1file} \csname rewrite#1filetrue\endcsname
  \ece\def{read#1file}{\readcontentsfile{#1}}%
}%
\definecontentsfile{toc}%
\def\opencontentsfile#1{%
  \csname if@#1fileopened\endcsname \else
     \ece{\immediate\openout}{#1file} = \csname #1filebasename\endcsname.#1
     \ece\global{@#1fileopenedtrue}%
  \fi
}%
\def\writecontentsentry#1#2#3{\writenumberedcontentsentry{#1}{#2}{#3}{}}%
\def\writenumberedcontentsentry#1#2#3#4{%
  \csname ifrewrite#1file\endcsname
    \writenumberedcontents@cmdname{#1}{#2}%
    \def\temp{#3}
    \toks2 = \expandafter{#4}%
    \edef\cs{\the\toks2}%
    \edef\@wr{%
      \write\csname #1file\endcsname{%
        \the\toks0 
        {\sanitize\temp}
        \ifx\empty\cs\else {\sanitize\cs}\fi 
        {\noexpand\folio}
      }%
    }%
    \@wr
  \fi
  \ignorespaces
}%
\def\writenumberedcontentsline#1#2#3#4{%
  \csname ifrewrite#1file\endcsname
    \writenumberedcontents@cmdname{#1}{#2}%
    \def\temp{#4}
    \toks2 = \expandafter{#3}%
    \edef\cs{\the\toks2}%
    \edef\@wr{%
      \write\csname #1file\endcsname{%
        \the\toks0 
        \ifx\empty\cs\else {\sanitize\cs}\fi 
        {\sanitize\temp}
        {\noexpand\folio}
      }%
    }%
    \@wr
  \fi
  \ignorespaces
}%
\def\writenumberedcontents@cmdname#1#2{%
  \csname open#1file\endcsname
  \edef\temp{#2}
  \expandafter\if\expandafter\isinteger\expandafter{\temp}%
    \toks0 = {\expandafter\noexpand \csname #1entry\endcsname}%
    \edef\temp{\the\toks0{\temp}}%
    \toks0 = \expandafter{\temp}%
  \else
    \toks0 = {\expandafter\noexpand \csname #1#2entry\endcsname}%
  \fi
}%
\def\readcontentsfile#1{%
   \edef\temp{%
     \noexpand\testfileexistence[\csname #1filebasename\endcsname]{#1}%
   }\temp
   \if@fileexists
      \input \csname #1filebasename\endcsname.#1\relax
   \fi
}%
\let\ifxrefwarning = \iftrue
\def\xrefwarningtrue{\@citewarningtrue \let\ifxrefwarning = \iftrue}%
\def\xrefwarningfalse{\@citewarningfalse \let\ifxrefwarning = \iffalse}%
\begingroup
  \catcode`\_ = 8
  \gdef\xrlabel#1{#1_x}%
\endgroup
\def\xrdef#1{%
  \begingroup
    \hldest@impl{xrdef}{#1}%
    \begingroup
      \@@hldestoff
      \definexref{#1}{\noexpand\folio}{page}%
    \endgroup
  \endgroup
  \ignorespaces
}%
\def\definexref#1#2#3{%
  \hldest@impl{definexref}{#1}%
  \edef\temp{#1}%
  \readauxfile
  \edef\@wr{\noexpand\writeaux{\string\@definelabel{\temp}{#2}{#3}}}%
  \@wr
  \ignorespaces
}%
\def\@definelabel#1{
  \begingroup 
    \expandafter\ifx\csname\xrlabel{#1}\endcsname \relax
      \expandafter\@definelabel@nocheck
    \else
      \expandafter\@definelabel@warn
    \fi
    {#1}%
}%
\def\@definelabel@nocheck#1#2#3{%
    \expandafter\gdef\csname\xrlabel{#1}\endcsname{#2}%
    \setpropertyglobal{\xrlabel{#1}}{class}{#3}%
  \endgroup 
}%
\def\@definelabel@warn#1#2#3{%
  \message{^^J\linenumber Label `#1' multiply defined,
           value `#2' of class `#3' overwriting value
           `\csname\xrlabel{#1}\endcsname' of class
           `\getproperty{\xrlabel{#1}}{class}'.}%
  \@definelabel@nocheck{#1}{#2}{#3}%
}%
\def\reftie{\penalty\@M \ }
\let\refspace\
\def\xrefn{\@getoptionalarg\@finxrefn}%
\def\@finxrefn#1{%
  \hlstart@impl{ref}{#1}%
  \ifx\@optionalarg\empty \else
    \let\@xrefnoptarg\@optionalarg
    \readauxfile
    {\@@hloff\@xrefnoptarg}\reftie
  \fi
  \plain@xrefn{#1}%
  \hlend@impl{ref}%
}%
\def\plain@xrefn#1{%
  \readauxfile
  \expandafter \ifx\csname\xrlabel{#1}\endcsname\relax
    \if@citewarning
       \message{\linenumber Undefined label `#1'.}%
    \fi
    \expandafter\def\csname\xrlabel{#1}\endcsname{%
      `{\tt
        \escapechar = -1
        \expandafter\string\csname#1\endcsname
      }'%
    }%
  \fi
  \csname\xrlabel{#1}\endcsname 
}%
\let\refn = \xrefn
\def\xrefpageword{p.\thinspace}%
\def\xref{\@getoptionalarg\@finxref}%
\def\@finxref#1{%
  \hlstart@impl{xref}{#1}%
  \ifx\@optionalarg\empty \else
    {\@@hloff\@optionalarg}\refspace
  \fi
  \xrefpageword\plain@xrefn{#1}%
  \hlend@impl{xref}%
}%
\def\@maybewarnref{%
  \ifundefined{amsppt.sty}%
  \else
    \message{Warning: amsppt.sty and Eplain both define \string\ref. See
             the Eplain manual.}%
    \let\amsref = \ref
  \fi
  \let\ref = \eplainref
  \ref
}
\let\ref = \@maybewarnref
\def\eplainref{\@getoptionalarg\@fineplainref}%
\def\@fineplainref{\@generalref{1}{}}%
\def\refs{\let\@optionalarg\empty \@generalref{0}s}%
\def\@generalref#1#2#3{%
  \let\@generalrefoptarg\@optionalarg
  \readauxfile
  \ifcase#1 \else \hlstart@impl{ref}{#3}\fi
  \edef\@generalref@class{\getproperty{\xrlabel{#3}}{class}}%
  \expandafter\ifx\csname \@generalref@class word\endcsname\relax
    \ifx\@generalrefoptarg\empty \else {\@@hloff\@generalrefoptarg\reftie}\fi
  \else
    \begingroup
      \@@hloff
      \ifx\@generalrefoptarg\empty \else \@generalrefoptarg \refspace \fi
      \csname \@generalref@class word\endcsname
      #2\reftie
    \endgroup
  \fi
  \ifcase#1 \hlstart@impl{ref}{#3}\fi
  \plain@xrefn{#3}%
  \hlend@impl{ref}%
}%
\newcount\eqnumber
\newcount\subeqnumber
\def\eqdefn{\@getoptionalarg\@fineqdefn}%
\def\@fineqdefn#1{%
  \ifx\@optionalarg\empty
    \global\advance\eqnumber by 1
    \def\temp{\eqconstruct{\number\eqnumber}}%
  \else
    \def\temp{\@optionalarg}%
  \fi
  \global\subeqnumber = 0
  \gdef\@currenteqlabel{#1}%
  \toks0 = \expandafter{\@currenteqlabel}%
  \begingroup
    \def\eqrefn{\noexpand\plain@xrefn}%
    \def\xrefn{\noexpand\plain@xrefn}%
    \edef\temp{\noexpand\@eqdefn{\the\toks0}{\temp}}%
    \temp
  \endgroup
}%
\def\eqsubdefn#1{%
  \global\advance\subeqnumber by 1
  \toks0 = {#1}%
  \toks2 = \expandafter{\@currenteqlabel}%
  \begingroup
    \def\eqrefn{\noexpand\plain@xrefn}%
    \def\xrefn{\noexpand\plain@xrefn}%
    \def\eqsubreftext{\noexpand\eqsubreftext}%
    \edef\temp{%
      \noexpand\@eqdefn
        {\the\toks0}%
        {\eqsubreftext{\eqrefn{\the\toks2}}{\the\subeqnumber}}%
    }%
    \temp
  \endgroup
}%
\newcount\phantomeqnumber
\def\phantomeqlabel{PHEQ\the\phantomeqnumber}%
\def\@eqdefn#1#2{%
  \ifempty{#1}%
    \global\advance\phantomeqnumber by 1
    \edef\hl@eqlabel{\phantomeqlabel}%
    \readauxfile
  \else
    \def\hl@eqlabel{#1}%
    {\@@hldestoff \definexref{#1}{#2}{eq}}%
  \fi
  \hldest@impl{eq}{\hl@eqlabel}%
  \begingroup 
    \@definelabel@nocheck{#1}{#2}{eq}%
}%
\def\eqdef{\@getoptionalarg\@fineqdef}%
\def\@fineqdef{%
  \toks0 = \expandafter{\@optionalarg}%
  \edef\temp{\noexpand\@eqdef{\noexpand\eqdefn[\the\toks0]}}%
  \temp
}%
\def\eqsubdef{\@eqdef\eqsubdefn}%
\def\@eqdef#1#2{%
  \@maybedisableeqno
  \eqnum #1{#2}
        \let\@optionalarg\empty 
        {\@@hloff\@fineqref{#2}}
  \@mayberestoreeqno
  \ignorespaces
}%
\let\@mayberestoreeqno = \relax
\def\@maybedisableeqno{%
  \ifinner
    \global\let\eqno = \relax
    \global\let\leqno = \relax
    \global\let\@mayberestoreeqno = \@restoreeqno
  \fi
}%
\let\@primitiveeqno = \eqno
\let\@primitiveleqno = \leqno
\def\@restoreeqno{%
  \global\let\eqno = \@primitiveeqno
  \global\let\leqno = \@primitiveleqno
  \global\let\@mayberestoreeqno = \empty
}%
\def\righteqnumbers{%
  \def\eqnum{\eqno}%
  \def\eqalignnum{\eqalignno}%
}%
\righteqnumbers
\def\eqrefn{\@getoptionalarg\@fineqrefn}%
\def\@fineqrefn#1{%
  \eqref@start{#1}%
  \plain@xrefn{#1}%
  \hlend@impl{eq}%
}%
\def\eqref{\@getoptionalarg\@fineqref}%
\def\@fineqref#1{%
  \eqref@start{#1}%
  \eqprint{\plain@xrefn{#1}}%
  \hlend@impl{eq}%
}%
\def\eqref@start#1{%
  \let\@eqrefoptarg\@optionalarg
  \ifempty{#1}%
    \hlstart@impl{eq}{\phantomeqlabel}%
  \else
    \hlstart@impl{eq}{#1}%
  \fi
  \ifx\@eqrefoptarg\empty \else
    {\@@hloff\@eqrefoptarg\reftie}%
  \fi
}%
\let\eqconstruct = \identity
\def\eqprint#1{(#1)}%
\def\eqsubreftext#1#2{#1.#2}%
\let\extraidxcmdsuffixes = \empty
\def\defineindex#1{%
  \def\@idxprefix{#1}%
  \expandafter\innernewif\csname if\@idxprefix dx\endcsname
  \csname \@idxprefix dxtrue\endcsname
  \for\@idxcmd:=,marked,submarked,name%
                \extraidxcmdsuffixes\do
  {%
    \@defineindexcmd\@idxcmd
  }%
  \ece\innernewwrite{@#1indexfile}%
  \ece\innernewif{if@#1indexfileopened}%
}%
\newif\ifsilentindexentry
\def\@defineindexcmd#1{%
  \@defineoneindexcmd{s}{#1}\silentindexentrytrue
  \@defineoneindexcmd{}{#1}\silentindexentryfalse
}%
\def\@defineoneindexcmd#1#2#3{%
  \toks@ = {#3}%
  \edef\temp{%
    \def
      \expandonce\csname#1\@idxprefix dx#2\endcsname 
      {\def\noexpand\@idxprefix{\@idxprefix}
       \expandonce\csname @@#1idx#2\endcsname
      }%
    \def
      \expandonce\csname @@#1idx#2\endcsname{
        \the\toks@
        \noexpand\@idxgetrange\expandonce\csname @#1idx#2\endcsname
      }%
  }%
  \temp
}%
\let\indexfilebasename = \jobname
\def\@idxwrite#1#2{%
  \csname if\@idxprefix dx\endcsname
    \@openidxfile
    \def\temp{#1}%
    \edef\@wr{%
      \expandafter\write\csname @\@idxprefix indexfile\endcsname{%
        \string\indexentry
        {\sanitize\temp}%
        {\noexpand#2}%
      }%
    }%
    \@wr
  \else
    \write-1{}%
  \fi
  \ifindexproofing
    \def\temp{#1}%
    \edef\temp{%
      \insert\@indexproof{\noexpand\indexproofterm{\sanitize\temp}}%
    }%
    \temp
    \ifhmode\allowhyphens\fi
  \fi
  \hookrun{afterindexterm}%
  \ifsilentindexentry \expandafter\ignorespaces\fi
}%
\def\@openidxfile{%
  \csname if@\@idxprefix indexfileopened\endcsname \else
    \expandafter\immediate\openout\csname @\@idxprefix indexfile\endcsname =
      \indexfilebasename.\@idxprefix dx
    \expandafter\global\csname @\@idxprefix indexfileopenedtrue\endcsname
  \fi
}%
\newif\ifindexproofing
\newinsert\@indexproof
\dimen\@indexproof = \maxdimen                  
\count\@indexproof = 0  \skip\@indexproof = 0pt 
\font\indexprooffont = cmtt8
\def\indexproofterm#1{\hbox{\strut \indexprooffont #1}}%
\let\@plainmakeheadline = \makeheadline
\def\makeheadline{%
  \expandafter\ifx\csname\idxpageanchor{\folio}\endcsname\relax \else
    {\@@hldeston \hldest@impl{idx}{\hlidxpagelabel{\folio}}}%
  \fi
  \indexproofunbox
  \@plainmakeheadline
}%
\def\indexsetmargins{%
  \ifx\undefined\outsidemargin
    \dimen@ = 1truein
    \advance\dimen@ by \hoffset
    \edef\outsidemargin{\the\dimen@}%
    \let\insidemargin = \outsidemargin
  \fi
}%
\def\indexproofunbox{%
  \ifvoid\@indexproof\else
    \indexsetmargins
    \rlap{%
      \kern\hsize
      \ifodd\pageno \kern\outsidemargin \else \kern\insidemargin \fi
      \vbox to 0pt{\unvbox\@indexproof\vss}%
    }\nointerlineskip
  \fi
}%
\def\idxrangebeginword{begin}%
\def\idxbeginrangemark{(}
\def\idxrangeendword{end}%
\def\idxendrangemark{)}%
\def\idxseecmdword{see}%
\def\idxseealsocmdword{seealso}%
\newif\if@idxsee
\newif\if@hlidxencap
\let\@idxseenterm = \relax
\def\idxpagemarkupcmdword{pagemarkup}%
\let\@idxpagemarkup = \relax
\def\@idxgetrange#1{%
  \let\@idxrangestr = \empty
  \let\@afteridxgetrange = #1%
  \begingroup
    \catcode\idxargopen=1
    \@getoptionalarg\@finidxgetopt
}%
\def\@finidxgetopt{%
    \global\let\@idxgetrange@arg\@optionalarg
  \endgroup
  \@hlidxencaptrue
  \for\@idxarg:=\@idxgetrange@arg\do{%
    \expandafter\@idxcheckpagemarkup\@idxarg=,%
    \ifx\@idxarg\idxrangebeginword
      \def\@idxrangestr{\idxencapoperator\idxbeginrangemark}%
    \else
      \ifx\@idxarg\idxrangeendword
        \def\@idxrangestr{\idxencapoperator\idxendrangemark}%
        \@hlidxencapfalse
      \else
        \ifx\@idxarg\idxseecmdword
          \def\@idxpagemarkup{indexsee}%
          \@idxseetrue
          \@hlidxencapfalse
        \else
          \ifx\@idxarg\idxseealsocmdword
            \def\@idxpagemarkup{indexseealso}%
            \@idxseetrue
            \@hlidxencapfalse
          \else
             \ifx\@idxpagemarkup\relax
               \errmessage{Unrecognized index option `\@idxarg'}%
             \fi
          \fi
        \fi
      \fi
    \fi
  }%
  \ifnum\hldest@place@idx < 0 \else
    \if@hlidxencap
      \ifx\@idxpagemarkup\relax
        \let\@idxpagemarkup\empty
      \fi
      \ifcase\hldest@place@idx \relax
        \edef\@idxpagemarkup{hlidxpage{\@idxpagemarkup}}%
        \definepageanchor{\noexpand\folio}%
      \else
        \global\advance\hlidxlabelnumber by 1
        \edef\@idxpagemarkup{hlidx{\hlidxlabel}{\@idxpagemarkup}}%
        \hldest@impl{idx}{\hlidxlabel}%
      \fi
    \fi
  \fi
  \@afteridxgetrange
}%
\def\@idxcheckpagemarkup#1=#2,{%
  \def\temp{#1}%
  \ifx\temp\idxpagemarkupcmdword
    \if ,#2, 
      \errmessage{Missing markup command to `pagemarkup'}%
    \else
      \def\temp##1={##1}%
      \edef\@idxpagemarkup{\temp\string#2}%
    \fi
  \fi
}%
\def\hldest@place@idx{-1}%
\begingroup
  \catcode`\_ = 8
  \gdef\idxpageanchor#1{#1_p}%
\endgroup
\def\definepageanchor#1{%
  \readauxfile
  \edef\@wr{\noexpand\writeaux{\string\@definepageanchor{#1}}}%
  \@wr
  \ignorespaces
}%
\def\@definepageanchor#1{%
  \expandafter\gdef\csname\idxpageanchor{#1}\endcsname{}%
}%
\newcount\hlidxlabelnumber
\def\hlidxlabel{IDX\the\hlidxlabelnumber}%
\def\hlidxpagelabel#1{IDXPG#1}%
\def\hlidx#1#2#3{%
  \ifempty{#2}%
    \let\@idxpageencap\relax
  \else
    \expandafter\let\expandafter\@idxpageencap\csname #2\endcsname
  \fi
  \hlstart@impl{idx}{#1}%
  \@idxpageencap{#3}%
  \hlend@impl{idx}%
}%
\def\hlidxpage#1#2{%
  \@hlidxgetpages{#2}%
  \ifempty{#1}%
    \let\@idxpageencap\relax
  \else
    \expandafter\let\expandafter\@idxpageencap\csname #1\endcsname
  \fi
  \hlstart@impl{idx}{\hlidxpagelabel{\@idxpageiref}}%
  \expandafter\@idxpageencap\expandafter{\@idxpagei}%
  \hlend@impl{idx}%
  \ifx\@idxpageii\empty \else
    \@idxpagesep
    \hlstart@impl{idx}{\hlidxpagelabel{\@idxpageiiref}}%
    \expandafter\@idxpageencap\expandafter{\@idxpageii}%
    \hlend@impl{idx}%
  \fi
}%
\def\@hlidxgetpages#1{%
  \idxparselist{#1}%
  \ifx\idxpagei\empty
    \idxparserange{#1}%
    \ifx\idxpagei\empty
      \def\@idxpageiref{#1}
    \else
      \let\@idxpageiref\idxpagei 
    \fi
    \def\@idxpagei{#1}%
    \let\@idxpageii\empty
  \else
    \let\@idxpagei\idxpagei
    \let\@idxpageii\idxpageii
    \let\@idxpageiref\idxpagei 
    \let\@idxpageiiref\idxpageii 
    \let\@idxpagesep\idxpagelistdelimiter
  \fi
}%
\def\setidxpagelistdelimiter#1{%
  \gdef\idxpagelistdelimiter{#1}%
  \gdef\@removepagelistdelimiter##1#1{##1}%
  \gdef\@idxparselist##1#1##2\@@{%
    \ifempty{##2}%
      \let\idxpagei\empty
    \else
      \def\idxpagei{##1}%
      \edef\idxpageii{\@removepagelistdelimiter##2}%
    \fi
  }%
  \gdef\idxparselist##1{\@idxparselist##1#1\@@}%
}%
\def\setidxpagerangedelimiter#1{%
  \gdef\idxpagerangedelimiter{#1}%
  \gdef\@removepagerangedelimiter##1#1{##1}%
  \gdef\@idxparserange##1#1##2\@@{%
    \ifempty{##2}%
      \let\idxpagei\empty
    \else
      \def\idxpagei{##1}%
      \edef\idxpageii{\@removepagerangedelimiter##2}%
    \fi
  }%
  \gdef\idxparserange##1{\@idxparserange##1#1\@@}%
}%
\setidxpagelistdelimiter{, }%
\setidxpagerangedelimiter{--}%
\def\idxsubentryseparator{!}%
\def\idxencapoperator{|}%
\def\idxmaxpagenum{99999}%
\newtoks\@idxmaintoks
\newtoks\@idxsubtoks
\def\@idxtokscollect{%
  \edef\temp{\the\@idxsubtoks}%
  \edef\@indexentry{%
    \the\@idxmaintoks
    \ifx\temp\empty\else \idxsubentryseparator\the\@idxsubtoks \fi
    \@idxrangestr
  }%
  \if@idxsee
    \@idxseefalse 
    \edef\temp{\noexpand\idx@getverbatimarg
      {\noexpand\@finidxtokscollect{\idxmaxpagenum}}}%
  \else
    \def\temp{\@finfinidxtokscollect\folio}%
  \fi
  \temp
}%
\def\@finidxtokscollect#1#2{%
  \def\@idxseenterm{#2}%
  \@finfinidxtokscollect{#1}%
}%
\def\@finfinidxtokscollect#1{%
  \ifx\@idxpagemarkup\relax \else
    \toks@ = \expandafter{\@indexentry}%
    \edef\@indexentry{%
      \the\toks@
      \ifx\@idxrangestr\empty \idxencapoperator \fi
      \@idxpagemarkup
    }%
    \let\@idxpagemarkup = \relax
  \fi
  \ifx\@idxseenterm\relax \else
    \toks@ = \expandafter{\@indexentry}%
    \edef\@indexentry{\the\toks@{\sanitize\@idxseenterm}}%
    \let\@idxseenterm = \relax
  \fi
  \expandafter\@idxwrite\expandafter{\@indexentry}{#1}%
}%
\def\@idxcollect#1#2{%
  \@idxmaintoks = {#1}%
  \@idxsubtoks = {#2}%
  \@idxtokscollect
}%
\def\idxargopen{`\{}%
\def\idxargclose{`\}}%
\def\idx@getverbatimarg#1{%
  \def\idx@getverbatimarg@cmd{#1}
  \begingroup
    \uncatcodespecials
    \catcode\idxargopen=1
    \catcode\idxargclose=2
    \catcode`\ =10   
    \catcode`\^^I=10 
    \catcode`\^^M=5  
    \idx@fingetverbatimarg
}%
\def\idx@fingetverbatimarg#1{\endgroup\idx@getverbatimarg@cmd{#1}}%
\def\idx@getverboptarg#1{%
  \def\idx@optionaltemp{#1}
  \let\idx@optionalnext = \relax
  \begingroup
    \if@idxsee \catcode\idxargopen=1 \fi
    \@futurenonspacelet\idx@optionalnext\idx@bracketcheck
}%
\def\idx@bracketcheck{%
  \ifx [\idx@optionalnext
    \endgroup 
    \expandafter\idx@finbracketcheck
  \else
    \endgroup 
    \let\@optionalarg = \empty
    \expandafter\idx@optionaltemp
  \fi
}%
\def\idx@finbracketcheck{%
  \begingroup
    \uncatcodespecials
    \catcode`\ =10   
    \catcode`\^^I=10 
    \catcode`\^^M=5  
    \idx@@getoptionalarg
}%
\def\idx@@getoptionalarg[#1]{%
  \endgroup
  \def\@optionalarg{#1}%
  \idx@optionaltemp
}%
{\catcode`\&=0
\gdef\idx@scanterm#1{%
  \edef\idx@scanterm@nl@catcode{\the\catcode`\^^M\relax}%
  \catcode`\^^M=5
  \scantokens{#1&endinput}%
  \catcode`\^^M=\idx@scanterm@nl@catcode
}}%
\def\@idx{\idx@getverbatimarg\@finidx}%
\def\@finidx#1{%
  \idx@scanterm{#1}
  \@idxcollect{#1}{}%
}%
\def\@sidx{\idx@getverbatimarg\@finsidx}%
\def\@finsidx#1{\@idxmaintoks = {#1}\idx@getverboptarg\@finfinsidx}%
\def\@finfinsidx{%
  \@idxsubtoks = \expandafter{\@optionalarg}%
  \@idxtokscollect
}%
\def\idxsortkeysep{@}
\def\@idxconstructmarked#1#2#3{%
  \toks@ = {#2}
  \toks2 = {#3}
  \edef\temp{\the\toks2 \idxsortkeysep \the\toks@{\the\toks2}}%
  #1 = \expandafter{\temp}%
}%
\def\@idxmarked#1{\idx@getverbatimarg{\@finidxmarked{#1}}}%
\def\@finidxmarked#1#2{%
  \idx@scanterm{#1{#2}}
  \@idxconstructmarked\@idxmaintoks{#1}{#2}%
  \@idxsubtoks = {}%
  \@idxtokscollect
}%
\def\@sidxmarked#1{\idx@getverbatimarg{\@finsidxmarked{#1}}}%
\def\@finsidxmarked#1#2{%
  \@idxconstructmarked\toks@{#1}{#2}%
  \edef\temp{{\the\toks@}}%
  \expandafter\@finsidx\temp
}%
\def\@idxsubmarked{%
  \let\sidxsubmarked@print\idxsubmarked@print
  \idx@getverbatimarg\@finsidxsubmarked
}%
\def\idxsubmarked@print{\expandafter\idx@scanterm\expandafter}%
\def\@sidxsubmarked{%
  \let\sidxsubmarked@print\gobble
  \idx@getverbatimarg\@finsidxsubmarked
}%
\def\@finsidxsubmarked#1{\@idxmaintoks = {#1}\@@finsidxsubmarked}
\def\@@finsidxsubmarked#1{\idx@getverbatimarg{\@@@finsidxsubmarked{#1}}}
\def\@@@finsidxsubmarked#1#2{
  \sidxsubmarked@print 
    {\the\@idxmaintoks\space #1{#2}}
  \@idxconstructmarked\@idxsubtoks{#1}{#2}%
  \@idxtokscollect
}%
\def\idxnameseparator{, }
\def\@idxcollectname#1#2{%
  \def\temp{#1}%
  \ifx\temp\empty
    \toks@ = {}%
  \else
    \toks@ = \expandafter{\idxnameseparator #1}%
  \fi
  \toks2 = {#2}%
  \edef\temp{\the\toks2 \the\toks@}%
}%
\def\@idxname{\idx@getverbatimarg\@finidxname}%
\def\@finidxname#1{\idx@getverbatimarg{\@finfinidxname{#1}}}%
\def\@finfinidxname#1#2{%
  \idx@scanterm{#1 #2}
  \@idxcollectname{#1}{#2}%
  \expandafter\@idxcollect\expandafter{\temp}{}%
}%
\def\@sidxname{\idx@getverbatimarg\@finsidxname}%
\def\@finsidxname#1{\idx@getverbatimarg{\@finfinsidxname{#1}}}%
\def\@finfinsidxname#1#2{%
  \@idxcollectname{#1}{#2}%
  \expandafter\@finsidx\expandafter{\temp}%
}%
\let\indexfonts = \relax
\def\readindexfile#1{%
  \edef\@idxprefix{#1}%
  \testfileexistence[\indexfilebasename]{\@idxprefix nd}%
  \iffileexists \begingroup
    \ifx\begin\undefined
      \def\begin##1{\@beginindex}%
      \let\end = \@gobble
    \fi
    \input \indexfilebasename.\@idxprefix nd
    \singlecolumn
  \endgroup
  \else
    \message{No index file \indexfilebasename.\@idxprefix nd.}%
  \fi
}%
\def\@beginindex{%
  \let\item = \@indexitem
  \let\subitem = \@indexsubitem
  \let\subsubitem = \@indexsubsubitem
  \indexfonts
  \doublecolumns
  \parindent = 0pt
  \hookrun{beginindex}%
}%

\newskip\aboveindexitemskipamount  \aboveindexitemskipamount = 0pt plus2pt
\def\aboveindexitemskip{\vskip\aboveindexitemskipamount}%
\def\@indexitem{\begingroup
  \@indexitemsetup
  \leftskip = 0pt
  \aboveindexitemskip
  \penalty-100 
  \def\par{\endgraf\endgroup\nobreak}%
}%
\def\@indexsubitem{%
  \@indexitemsetup
  \leftskip = 1em
}%
\def\@indexsubsubitem{%
  \@indexitemsetup
  \leftskip = 2em
}%
\def\@indexitemsetup{%
  \par
  \hangindent = 1em
  \raggedright
  \hyphenpenalty = 10000
  \hookrun{indexitem}%
}%
\def\seevariant{\it}%
\def\indexseeword{see}%
\def\indexsee{\idx@getverbatimarg\@finindexsee}%
\def\@finindexsee#1#2{{\seevariant \indexseeword\/ }\idx@scanterm{#1}}%
\def\indexseealsowords{see also}%
\def\indexseealso{\idx@getverbatimarg\@finindexseealso}%
\def\@finindexseealso#1#2{{\seevariant \indexseealsowords\/ }\idx@scanterm{#1}}%
\defineindex{i}%
\begingroup
  \catcode `\^^M = \active %
  \gdef\flushleft{%
    \def\@endjustifycmd{\@endflushleft}%
    \def\@eoljustifyaction{\null\hfil\break}%
    \let\@firstlinejustifyaction = \relax
    \@startjustify %
  }%
  \gdef\flushright{%
    \def\@endjustifycmd{\@endflushright}%
    \def\@eoljustifyaction{\break\null\hfil}%
    \def\@firstlinejustifyaction{\hfil\null}%
    \@startjustify %
  }%
  \gdef\center{%
    \def\@endjustifycmd{\@endcenter}%
    \def\@eoljustifyaction{\hfil\break\null\hfil}%
    \def\@firstlinejustifyaction{\hfil\null}%
    \@startjustify %
  }%
  \gdef\@startjustify{%
    \parskip = 0pt
    \catcode`\^^M = \active %
    \def^^M{\futurelet\next\@finjustifyreturn}%
    \def\@eateol##1^^M{%
      \def\temp{##1}%
      \@firstlinejustifyaction %
      \ifx\temp\empty\else \temp^^M\fi %
    }%
    \expandafter\aftergroup\@endjustifycmd %
    \checkenv \environmenttrue %
    \par\noindent %
    \@eateol %
  }%
  \gdef\@finjustifyreturn{%
    \@eoljustifyaction %
    \ifx\next^^M%
      \def\par{\endgraf\vskip\blanklineskipamount \global\let\par = \endgraf}%
      \@endjustifycmd %
      \noindent %
      \@firstlinejustifyaction %
    \fi %
  }%
\endgroup
\def\@endflushleft{\unpenalty{\parfillskip = 0pt plus1fil\par}\ignorespaces}%
\def\@endflushright{
   \unskip \setbox0=\lastbox \unpenalty
   {\parfillskip = 0pt \par}\ignorespaces
}%
\def\@endcenter{
   \unskip \setbox0=\lastbox \unpenalty
   {\parfillskip = 0pt plus1fil \par}\ignorespaces
}%
\ifx\@undefined\raggedleft
\innernewskip\raggedleftskip \raggedleftskip=0pt plus2em
\def\raggedleft{\leftskip\raggedleftskip \parindent=0pt
  \spaceskip.3333em \xspaceskip.5em \parfillskip0pt \relax}
\fi 
\newcount\abovecolumnspenalty   \abovecolumnspenalty = 10000
\newcount\@linestogo         
\newcount\@linestogoincolumn 
\newcount\@columndepth       
\newdimen\@columnwidth       
\newtoks\crtok  \crtok = {\cr}%
\newcount\currentcolumn
\def\makecolumns#1/#2: {\par \begingroup
   \@columndepth = #1
   \advance\@columndepth by -1
   \divide \@columndepth by #2
   \advance\@columndepth by 1
   \@linestogoincolumn = \@columndepth
   \@linestogo = #1
   \currentcolumn = 1
   \def\@endcolumnactions{%
      \ifnum \@linestogo<2
         \the\crtok \egroup \endgroup \par 
      \else
         \global\advance\@linestogo by -1
         \ifnum\@linestogoincolumn<2
            \global\advance\currentcolumn by 1
            \global\@linestogoincolumn = \@columndepth
            \the\crtok
         \else
            &\global\advance\@linestogoincolumn by -1
         \fi
      \fi
   }%
   \makeactive\^^M
   \letreturn \@endcolumnactions
   \@columnwidth = \hsize
     \advance\@columnwidth by -\parindent
     \divide\@columnwidth by #2
   \penalty\abovecolumnspenalty
   \noindent 
   \valign\bgroup
     &\hbox to \@columnwidth{\strut \hsize = \@columnwidth ##\hfil}\cr
}%
\newcount\footnotenumber
\newcount\hlfootlabelnumber
\newdimen\footnotemarkseparation \footnotemarkseparation = .5em
\newskip\interfootnoteskip \interfootnoteskip = 0pt
\newtoks\everyfootnote
\newdimen\footnoterulewidth \footnoterulewidth = 2in
\newdimen\footnoteruleheight \footnoteruleheight = 0.4pt
\newdimen\belowfootnoterulespace \belowfootnoterulespace = 2.6pt
\let\@plainfootnote = \footnote
\let\@plainvfootnote = \vfootnote
\def\hlfootlabel{FOOT\the\hlfootlabelnumber}%
\def\hlfootbacklabel{FOOTB\the\hlfootlabelnumber}%
\def\@eplainfootnote#1{\let\@sf\empty 
  \ifhmode\edef\@sf{\spacefactor\the\spacefactor}\/\fi
  \global\advance\hlfootlabelnumber by 1
  \hlstart@impl{foot}{\hlfootlabel}%
  \hldest@impl{footback}{\hlfootbacklabel}%
  #1%
  \hlend@impl{foot}%
  \@sf\vfootnote{#1}%
}%
\let\footnote\@eplainfootnote
\def\vfootnote#1{\insert\footins\bgroup
  \interlinepenalty\interfootnotelinepenalty
  \splittopskip\ht\strutbox 
  \advance\splittopskip by \interfootnoteskip
  \splitmaxdepth\dp\strutbox
  \floatingpenalty\@MM
  \leftskip\z@skip \rightskip\z@skip \spaceskip\z@skip \xspaceskip\z@skip
  \everypar = {}%
  \parskip = 0pt 
  \ifnum\@numcolumns > 1 \hsize = \@normalhsize \fi
  \the\everyfootnote
  \vskip\interfootnoteskip
  \indent\llap{%
    \hlstart@impl{footback}{\hlfootbacklabel}%
    \hldest@impl{foot}{\hlfootlabel}%
    #1%
    \hlend@impl{footback}%
    \kern\footnotemarkseparation
  }%
  \footstrut\futurelet\next\fo@t
}%
\def\footnoterule{\dimen@ = \footnoteruleheight
  \advance\dimen@ by \belowfootnoterulespace
  \kern-\dimen@
  \hrule width\footnoterulewidth height\footnoteruleheight depth0pt
  \kern\belowfootnoterulespace
  \vskip-\interfootnoteskip
}%
\def\numberedfootnote{%
  \global\advance\footnotenumber by 1
  \@eplainfootnote{$^{\number\footnotenumber}$}%
}%
\newdimen\paperheight
\ifnum\mag=1000
  \paperheight = 11in 
\else
  \paperheight = 11truein 
\fi
\def\topmargin{\afterassignment\@finishtopmargin \dimen@}%
\def\@finishtopmargin{%
  \dimen2 = \voffset		
  \voffset = \dimen@ \advance\voffset by -1truein
  \advance\dimen2 by -\voffset	
  \advance\vsize by \dimen2	
}%
\def\advancetopmargin{%
  \dimen@ = 0pt \afterassignment\@finishadvancetopmargin \advance\dimen@
}%
\def\@finishadvancetopmargin{%
  \advance\voffset by \dimen@
  \advance\vsize by -\dimen@
}%
\def\bottommargin{\afterassignment\@finishbottommargin \dimen@}%
\def\@finishbottommargin{%
  \@computebottommargin		
  \advance\dimen2 by -\dimen@	
  \advance\vsize by \dimen2	
}%
\def\advancebottommargin{%
  \dimen@ = 0pt \afterassignment\@finishadvancebottommargin \advance\dimen@
}%
\def\@finishadvancebottommargin{%
  \advance\vsize by -\dimen@
}%
\def\@computebottommargin{%
  \dimen2 = \paperheight	
  \advance\dimen2 by -\vsize	
  \advance\dimen2 by -\voffset	
  \advance\dimen2 by -1truein	
}%
\newdimen\paperwidth
\ifnum\mag=1000
  \paperwidth = 8.5in 
\else
  \paperwidth = 8.5truein 
\fi
\def\leftmargin{\afterassignment\@finishleftmargin \dimen@}%
\def\@finishleftmargin{%
  \dimen2 = \hoffset		
  \hoffset = \dimen@ \advance\hoffset by -1truein
  \advance\dimen2 by -\hoffset	
  \advance\hsize by \dimen2	
}%
\def\advanceleftmargin{%
  \dimen@ = 0pt \afterassignment\@finishadvanceleftmargin \advance\dimen@
}%
\def\@finishadvanceleftmargin{%
  \advance\hoffset by \dimen@
  \advance\hsize by -\dimen@
}%
\def\rightmargin{\afterassignment\@finishrightmargin \dimen@}%
\def\@finishrightmargin{%
  \@computerightmargin		
  \advance\dimen2 by -\dimen@	
  \advance\hsize by \dimen2	
}%
\def\advancerightmargin{%
  \dimen@ = 0pt \afterassignment\@finishadvancerightmargin \advance\dimen@
}%
\def\@finishadvancerightmargin{%
  \advance\hsize by -\dimen@
}%
\def\@computerightmargin{%
  \dimen2 = \paperwidth		
  \advance\dimen2 by -\hsize	
  \advance\dimen2 by -\hoffset	
  \advance\dimen2 by -1truein	
}%
\let\@plainm@g = \m@g
\def\m@g{\@plainm@g \paperwidth = 8.5 true in \paperheight = 11 true in}%
\newskip\abovecolumnskip \abovecolumnskip = \bigskipamount
\newskip\belowcolumnskip \belowcolumnskip = \bigskipamount
\newdimen\gutter \gutter = 2pc
\newbox\@partialpage
\newdimen\@normalhsize
\newdimen\@normalvsize  
\newtoks\previousoutput
\def\quadcolumns{\@columns4}%
\def\triplecolumns{\@columns3}%
\def\doublecolumns{\@columns2}%
\def\begincolumns#1{\ifcase#1\relax \or \singlecolumn \or \@columns2 \or
                            \@columns3 \or \@columns4 \else \relax \fi}%
\let\@ndcolumns = \relax
\chardef\@numcolumns = 1
\mathchardef\@ejectpartialpenalty = 10141
\chardef\@col@minlines = 3
\chardef\@col@extralines = 3
\newdimen\@col@extraheight
\def\@columns#1{%
  \@ndcolumns
  \global\let\@ndcolumns = \@endcolumns
  \global\chardef\@numcolumns = #1
  \global\previousoutput = \expandafter{\the\output}%
  \global\output = {%
    \ifnum\outputpenalty = -\@ejectpartialpenalty
      \dimen@ = \vsize
      \advance\dimen@ by \@col@minlines\baselineskip
      \global\setbox\@partialpage =
        \vbox  \ifdim \pagetotal > \vsize  to \dimen@  \fi  {%
	  \unvbox255 \unskip
	}%
    \else
      \the\previousoutput
    \fi
  }%
  \vskip \abovecolumnskip
  \vskip \@col@minlines\baselineskip
  \penalty -\@ejectpartialpenalty
  \global\output = {\@columnoutput}%
  \global\@normalhsize = \hsize
  \global\@normalvsize = \vsize
  \count@ = \@numcolumns
  \advance\count@ by -1
  \global\advance\hsize by -\count@\gutter
  \global\divide\hsize by \@numcolumns
  \advance\vsize by -\ht\@partialpage
  \advance\vsize by -\ht\footins
  \ifvoid\footins\else \advance\vsize by -\skip\footins \fi
  \multiply\count\footins by \@numcolumns
  \advance\vsize by -\ht\topins
  \ifvoid\topins\else \advance\vsize by -\skip\topins \fi
  \multiply\count\topins by \@numcolumns
  \global\vsize = \@numcolumns\vsize
  \@col@extraheight=\@col@extralines\baselineskip
  \multiply\@col@extraheight by \@numcolumns
  \global\advance\vsize by \@col@extraheight
}%
\def\gutterbox{\vbox to \dimen0{\vfil\hbox{\hfil}\vfil}}%
\def\@columnsplit{%
  \splittopskip = \topskip
  \splitmaxdepth = \baselineskip
  \begingroup
    \vbadness = 10000
    \global\setbox1 = \vsplit255 to \dimen@  \global\wd1 = \hsize
    \global\setbox3 = \vsplit255 to \dimen@  \global\wd3 = \hsize
    \ifnum\@numcolumns > 2
      \global\setbox5 = \vsplit255 to \dimen@ \global\wd5 = \hsize
    \fi
    \ifnum\@numcolumns > 3
      \global\setbox7 = \vsplit255 to \dimen@ \global\wd7 = \hsize
    \fi
  \endgroup
  \setbox0 = \box255
  \global\setbox255 = \vbox{%
    \unvbox\@partialpage
    \ifcase\@numcolumns \relax\or\relax
      \or \hbox to \@normalhsize{\box1\hfil\gutterbox\hfil\box3}%
      \or \hbox to \@normalhsize{\box1\hfil\gutterbox\hfil\box3%
                                      \hfil\gutterbox\hfil\box5}%
      \or \hbox to \@normalhsize{\box1\hfil\gutterbox\hfil\box3%
                                      \hfil\gutterbox\hfil\box5%
                                      \hfil\gutterbox\hfil\box7}%
    \fi
  }%
  \setbox\@partialpage = \box0
}%
\def\@columnoutput{%
  \dimen@ = \ht255
    \advance\dimen@ by -\@col@extraheight
    \divide\dimen@ by \@numcolumns
  \@columnsplit
  \@recoverclubpenalty
  \hsize = \@normalhsize 
  \vsize = \@normalvsize
  \the\previousoutput
  \unvbox\@partialpage
  \penalty\outputpenalty
  \global\vsize = \@numcolumns\@normalvsize
  \global\advance\vsize by \@col@extraheight
}%
\def\singlecolumn{%
  \@ndcolumns
  \chardef\@numcolumns = 1
  \vskip\belowcolumnskip
  \nointerlineskip
}%
\newbox\@singlecolumnbox
\newdimen\column@pagegoal
\newdimen\column@vsize
\def\@endcolumns{%
  \global\let\@ndcolumns = \relax
  \par 
  \column@pagegoal = \pagegoal
  \advance\column@pagegoal by-\@col@extraheight
  \ifdim \pagetotal > \column@pagegoal
    \column@vsize = \column@pagegoal
  \else
    \column@vsize = \pagetotal
  \fi
  \global\output = {\global\setbox1 = \box255}%
  \pagegoal = \pagetotal
  \break                     
  \setbox2 = \box1           
  \global\output = \expandafter{\the\previousoutput}%
  \setbox\@singlecolumnbox = \box\@partialpage
  \@balancecolumns
}%
\def\@balancecolumns{%
  \global\setbox255 = \copy2  
  \dimen@ = \column@vsize
    \divide\dimen@ by \@numcolumns
  \@columnsplit
  \ifvoid\@partialpage
    \global\vsize = \@normalvsize
    \global\hsize = \@normalhsize
    \dump@balanced@columns
    \let\next\relax
  \else
    \advance \column@vsize by \@numcolumns pt
    \test@spill@columns
    \ifspill@columns
      \begingroup
        \vsize = \@normalvsize
        \hsize = \@normalhsize
        \dump@balanced@columns
        \break
        \@recoverclubpenalty
      \endgroup
      \unvbox\@partialpage
      \let\next\@endcolumns
    \else
      \setbox0=\box\@partialpage 
      \let\next\@balancecolumns
    \fi
  \fi
  \next
}%
\def\dump@balanced@columns{%
  \ifvoid\topins\else\topinsert\unvbox\topins\endinsert\fi
  \unvbox\@singlecolumnbox
  \nointerlineskip
  \box255
}%
\newif\ifspill@columns
\def\test@spill@columns{%
  \spill@columnsfalse
  \ifdim \column@vsize > \column@pagegoal
    \ifvoid\footins
      \ifvoid\topins
        \spill@columnstrue
      \fi
    \fi
  \fi
}%
\def\@saveclubpenalty{
  \edef\@recoverclubpenalty{%
     \global\clubpenalty=\the\clubpenalty\relax%
     \global\let\noexpand\@recoverclubpenalty\relax
  }
}%
\let\@recoverclubpenalty\relax
\newdimen\temp@dimen
\def\columnfill{%
  \par
  \dimen@ = \pagetotal  
  \temp@dimen = \vsize  
  \advance\temp@dimen by -\@col@extraheight
  \divide\temp@dimen by \@numcolumns  
  \loop
    \ifdim \dimen@ > \temp@dimen  
      \advance \dimen@ by -\temp@dimen
      \advance \dimen@ by \topskip 
  \repeat
  \advance \temp@dimen by -\dimen@
  \advance \temp@dimen by -\prevdepth
  \@saveclubpenalty
  \clubpenalty=10000\relax
  \hrule height\temp@dimen width0pt depth0pt\relax  
  \nointerlineskip
  \par
  \nointerlineskip
  \allowbreak \vfil 
  \relax
}%
\def\@hldest{%
  \def\hl@prefix{hldest}%
  \let\after@hl@getparam\hldest@aftergetparam
  \begingroup
    \hl@getparam
}%
\def\hldest@aftergetparam{%
  \ifvmode
    \hldest@driver
  \else
    \allowhyphens
    \smash{\ifx\hldest@opt@raise\empty \else \raise\hldest@opt@raise\fi
             \hbox{\hldest@driver}}%
    \allowhyphens
  \fi
  \endgroup
}%
\def\@hlstart{%
  \leavevmode
  \def\hl@prefix{hl}%
  \let\after@hl@getparam\hlstart@aftergetparam
  \begingroup
    \hl@getparam
}%
\def\hlstart@aftergetparam{%
  \ifx\color\undefined \else
    \ifx\hl@opt@color\empty \else
      \ifx\hl@opt@colormodel\empty
        \edef\temp{\noexpand\color{\hl@opt@color}}%
      \else
        \edef\temp{\noexpand\color[\hl@opt@colormodel]{\hl@opt@color}}%
      \fi
      \temp
    \fi
  \fi
  \hl@driver
}%
\def\hl@getparam#1#2{
  \edef\@hltype{#1}%
  \ifx\@hltype\empty
    \expandafter\let\expandafter\@hltype
      \csname \hl@prefix @type\endcsname
  \fi
  \expandafter\ifx\csname \hl@prefix @typeh@\@hltype\endcsname \relax
    \errmessage{Invalid hyperlink type `\@hltype'}%
  \fi
  \For\hl@arg:=#2\do{%
    \ifx\hl@arg\empty \else
      \expandafter\hl@set@opt\hl@arg=,%
    \fi
  }%
  \bgroup
    \uncatcodespecials
    \catcode`\{=1 \catcode`\}=2
    \@hl@getparam
}%
\def\@hl@getparam#1{%
  \egroup
  \edef\@hllabel{#1}%
  \after@hl@getparam
  \ignorespaces
}%
\def\hl@set@opt#1=#2,{%
  \expandafter\ifx\csname \hl@prefix @opt@#1\endcsname \relax
    \errmessage{Invalid hyperlink option `#1'}%
  \fi
  \if ,#2, 
    \errmessage{Missing value for option `#1'}%
  \else
    \def\temp##1={##1}%
    \expandafter\edef\csname \hl@prefix @opt@#1\endcsname{\temp#2}%
  \fi
}%
\def\hldest@impl#1{%
  \expandafter\ifcase\csname hldest@on@#1\endcsname
    \relax\expandafter\gobble
  \else
    \toks@=\expandafter{\csname hldest@type@#1\endcsname}%
    \toks@ii=\expandafter{\csname hldest@opts@#1\endcsname}%
    \edef\temp{\noexpand\hldest{\the\toks@}{\the\toks@ii}}%
    \expandafter\temp
  \fi
}%
\def\hlstart@impl#1{%
  \expandafter\ifcase\csname hl@on@#1\endcsname
    \leavevmode\expandafter\gobble
  \else
    \toks@=\expandafter{\csname hl@type@#1\endcsname}%
    \toks@ii=\expandafter{\csname hl@opts@#1\endcsname}%
    \edef\temp{\noexpand\hlstart{\the\toks@}{\the\toks@ii}}%
    \expandafter\temp
  \fi
}%
\def\hlend@impl#1{%
  \expandafter\ifcase\csname hl@on@#1\endcsname
  \else
    \hlend
  \fi
}%
\def\hl@asterisk@word{*}%
\def\hl@opts@word{opts}%
\newif\if@params@override
\def\hldest@groups{definexref,xrdef,li,eq,bib,foot,footback,idx}%
\def\hl@groups{ref,xref,eq,cite,foot,footback,idx,url,hrefint,hrefext}%
\def\hldesttype{%
  \def\hl@prefix{hldest}%
  \def\hl@param{type}%
  \let\hl@all@groups\hldest@groups
  \futurelet\hl@excl\hl@param@read@excl
}%
\def\hldestopts{%
  \def\hl@prefix{hldest}%
  \def\hl@param{opts}%
  \let\hl@all@groups\hldest@groups
  \futurelet\hl@excl\hl@param@read@excl
}%
\def\hltype{%
  \def\hl@prefix{hl}%
  \def\hl@param{type}%
  \let\hl@all@groups\hl@groups
  \futurelet\hl@excl\hl@param@read@excl
}%
\def\hlopts{%
  \def\hl@prefix{hl}%
  \def\hl@param{opts}%
  \let\hl@all@groups\hl@groups
  \futurelet\hl@excl\hl@param@read@excl
}%
\def\hl@param@read@excl{%
  \ifx!\hl@excl
    \let\next\hl@param@read@opt@arg
    \@params@overridetrue
  \else
    \def\next{\hl@param@read@opt@arg{!}}%
    \@params@overridefalse
  \fi
  \next
}%
\def\hl@param@read@opt@arg#1{
  \@getoptionalarg\hl@setparam
}%
\def\@hl@setparam#1{%
  \ifx\@optionalarg\empty
    \hl@setparam@default{#1}
  \else
    \let\hl@do@all@groups\gobble
    \For\hl@group:=\@optionalarg\do{%
      \ifx\hl@group\hl@asterisk@word
        \def\hl@do@all@groups{\let\@optionalarg\hl@all@groups \hl@setparam}%
      \else
        \hl@setparam@group{#1}%
      \fi
    }%
    \hl@do@all@groups{#1}%
  \fi
}%
\def\hl@setparam@group#1{%
  \ifx\hl@group\empty
    \hl@setparam@default{#1}%
  \else
    \expandafter\ifx\csname\hl@prefix @\hl@param @\hl@group\endcsname\relax
      \errmessage{Hyperlink group `\hl@prefix:\hl@param:\hl@group' is not defined}%
    \fi
    \ifx\hl@param\hl@opts@word
      \if@params@override
        \expandafter\let\csname\hl@prefix @\hl@param @\hl@group\endcsname\empty
      \fi
      \hl@update@opts@with@list{#1}
    \else
      \ece\def{\hl@prefix @\hl@param @\hl@group}{#1}%
    \fi
  \fi
}%
\def\hl@setparam@default#1{%
  \ifx\hl@param\hl@opts@word
    \For\hl@opt:=#1\do{%
      \ifx\hl@opt\empty \else
        \expandafter\hl@set@opt\hl@opt=,%
      \fi
    }%
  \else
    \expandafter\ifx\csname\hl@prefix @\hl@param\endcsname\relax
      \message{Default hyperlink parameter `\hl@prefix:\hl@param' is not defined}%
    \fi
    \ece\def{\hl@prefix @\hl@param}{#1}%
  \fi
}%
\def\hl@update@opts@with@list#1{%
  \global\expandafter\let\expandafter\hl@update@new@list
    \csname \hl@prefix @opts@\hl@group\endcsname
  \begingroup
    \For\hl@opt:=#1\do{%
      \hl@update@opts@with@opt
    }%
  \endgroup
  \ece\let{\hl@prefix @opts@\hl@group}\hl@update@new@list
}%
\def\hl@update@opts@with@opt{%
  \global\let\hl@update@old@list\hl@update@new@list
  \global\let\hl@update@new@list\empty
  \global\let\hl@update@new@opt\hl@opt
  \expandafter\hl@parse@opt@key\hl@opt=,%
  \let\hl@update@new@key\hl@update@key
  \global\let\hl@update@comma\empty
  \begingroup
    \for\hl@opt:=\hl@update@old@list\do{%
      \ifx\hl@opt\empty \else 
        \expandafter\hl@parse@opt@key\hl@opt=,%
        \toks@=\expandafter{\hl@update@new@list}%
        \ifx\hl@update@key\hl@update@new@key
          \ifx\hl@update@new@opt\empty \else 
            \toks@ii=\expandafter{\hl@update@new@opt}%
            \xdef\hl@update@new@list{\the\toks@\hl@update@comma\the\toks@ii}%
            \global\let\hl@update@new@opt\empty
            \global\def\hl@update@comma{,}%
          \fi
        \else
          \toks@ii=\expandafter{\hl@opt}%
          \xdef\hl@update@new@list{\the\toks@\hl@update@comma\the\toks@ii}%
          \global\def\hl@update@comma{,}%
        \fi
      \fi
    }%
  \endgroup
  \ifx\hl@update@new@opt\empty \else
    \toks@=\expandafter{\hl@update@new@list}%
    \toks@ii=\expandafter{\hl@update@new@opt}%
    \xdef\hl@update@new@list{\the\toks@\hl@update@comma\the\toks@ii}%
  \fi
}%
\def\hl@parse@opt@key#1=#2,{\def\hl@update@key{#1}}%
\def\hldest@opt@raise{\normalbaselineskip}%
\def\hl@opt@colormodel{cmyk}%
\def\hl@opt@color{0.28,1,1,0.35}%
\def\hldest@on@definexref{0}%
\def\hldest@on@xrdef{0}%
\def\hldest@on@li{0}%
\def\hldest@on@eq{0}
\def\hldest@on@bib{0}
\def\hldest@on@foot{0}
\def\hldest@on@footback{0}
\def\hldest@on@idx{0}
\let\hldest@type@definexref\empty
\let\hldest@type@xrdef\empty
\let\hldest@type@li\empty
\let\hldest@type@eq\empty 
\let\hldest@type@bib\empty 
\let\hldest@type@foot\empty 
\let\hldest@type@footback\empty 
\let\hldest@type@idx\empty 
\let\hldest@opts@definexref\empty
\let\hldest@opts@xrdef\empty
\let\hldest@opts@li\empty
\def\hldest@opts@eq{raise=1.7\normalbaselineskip}
\let\hldest@opts@bib\empty 
\let\hldest@opts@foot\empty 
\let\hldest@opts@footback\empty 
\let\hldest@opts@idx\empty 
\def\hl@on@ref{0}
\def\hl@on@xref{0}%
\def\hl@on@eq{0}
\def\hl@on@cite{0}
\def\hl@on@foot{0}
\def\hl@on@footback{0}
\def\hl@on@idx{0}%
\def\hl@on@url{0}
\def\hl@on@hrefint{0}
\def\hl@on@hrefext{0}
\let\hl@type@ref\empty 
\let\hl@type@xref\empty
\let\hl@type@eq\empty 
\let\hl@type@cite\empty 
\let\hl@type@foot\empty 
\let\hl@type@footback\empty 
\let\hl@type@idx\empty
\let\hl@type@url\empty 
\let\hl@type@hrefint\empty 
\let\hl@type@hrefext\empty 
\let\hl@opts@ref\empty 
\let\hl@opts@xref\empty
\let\hl@opts@eq\empty 
\let\hl@opts@cite\empty 
\let\hl@opts@foot\empty 
\let\hl@opts@footback\empty 
\let\hl@opts@idx\empty
\let\hl@opts@url\empty 
\let\hl@opts@hrefint\empty 
\let\hl@opts@hrefext\empty 
\def\@hlon{\@hlonoff@value@stub{hl}\@@hlon1 }%
\def\@hloff{\@hlonoff@value@stub{hl}\@@hloff0 }%
\def\@hldeston{\@hlonoff@value@stub{hldest}\@@hldeston1 }%
\def\@hldestoff{\@hlonoff@value@stub{hldest}\@@hldestoff0 }%
\def\@hlonoff@value@stub#1#2#3{%
  \def\hl@prefix{#1}%
  \let\hl@on@empty#2%
  \def\hl@value{#3}%
  \expandafter\let\expandafter\hl@all@groups
    \csname \hl@prefix @groups\endcsname
  \@getoptionalarg\@finhlswitch
}%
\def\@finhlswitch{%
  \ifx\@optionalarg\empty
    \hl@on@empty
  \fi
  \let\hl@do@all@groups\relax
  \For\hl@group:=\@optionalarg\do{%
    \ifx\hl@group\hl@asterisk@word
      \let\@optionalarg\hl@all@groups
      \let\hl@do@all@groups\@finhlswitch
    \else
      \ifx\hl@group\empty
        \hl@on@empty
      \else
        \expandafter\ifx\csname\hl@prefix @on@\hl@group\endcsname \relax
          \errmessage{Hyperlink group `\hl@prefix:on:\hl@group'
                      is not defined}%
        \fi
        \ece\edef{\hl@prefix @on@\hl@group}{\hl@value}%
      \fi
    \fi
  }%
  \hl@do@all@groups
}%
\def\@@hlon{%
  \let\hlstart\@hlstart
  \let\hlend\@hlend
}%
\def\@@hloff{%
  \def\hlstart##1##2##3{\leavevmode\ignorespaces}%
  \let\hlend\relax
}%
\def\@@hldeston{%
  \let\hldest\@hldest
}%
\def\@@hldestoff{%
  \def\hldest##1##2##3{\ignorespaces}%
}%
\def\hl@idxexact@word{idxexact}%
\def\hl@idxpage@word{idxpage}%
\def\hl@idxnone@word{idxnone}%
\def\hl@raw@word{raw}%
\def\enablehyperlinks{\@getoptionalarg\@finenablehyperlinks}%
\def\@finenablehyperlinks{%
  \let\hl@selecteddriver\empty
  \def\hldest@place@idx{0}%
  \for\hl@arg:=\@optionalarg\do{%
    \ifx\hl@arg\hl@idxexact@word
      \def\hldest@place@idx{1}%
    \else
      \ifx\hl@arg\hl@idxnone@word
        \def\hldest@place@idx{-1}%
      \else
        \ifx\hl@arg\hl@idxpage@word
          \def\hldest@place@idx{0}%
        \else
          \let\hl@selecteddriver\hl@arg
        \fi
      \fi
    \fi
  }%
  \ifx\hl@selecteddriver\empty
    \ifpdf
      \def\hl@selecteddriver{pdftex}%
      \message{^^JEplain: using `pdftex' hyperlink driver.}%
    \else
      \def\hl@selecteddriver{hypertex}%
      \message{^^JEplain: using `hypertex' hyperlink driver.}%
    \fi
  \else
    \expandafter\ifx\csname hldriver@\hl@selecteddriver\endcsname \relax
      \errmessage{No hyperlink driver `\hl@selecteddriver' available}%
    \fi
  \fi
  \let\hl@setparam\@hl@setparam
  \csname hldriver@\hl@selecteddriver\endcsname
  \def\@finenablehyperlinks{\errmessage{Hyperlink driver `\hl@selecteddriver'
                                        already selected}}%
  \let\hldriver@nolinks\undefined
  \let\hldriver@hypertex\undefined
  \let\hldriver@pdftex \undefined
  \let\hldriver@dvipdfm\undefined
  \let\hloff\@hloff
  \let\hlon\@hlon
  \let\hldestoff\@hldestoff
  \let\hldeston\@hldeston
  \hlon[*,]\hloff[foot,footback]%
  \hldeston[*,]\hldestoff[foot,footback]%
}%
\def\hldriver@nolinks{%
  \def\@hldest##1##2##3{%
    \edef\temp{\write-1{hldest: ##3}}%
    \ifvmode
      \temp
    \else
      \allowhyphens
      \expandafter\smash\expandafter{\temp}%
      \allowhyphens
    \fi
    \ignorespaces
  }%
  \def\@hlstart##1##2##3{%
    \leavevmode
    \begingroup 
    \edef\temp{\write-1{hlstart: ##3}}%
    \temp
    \ignorespaces
  }%
  \def\@hlend{%
    \edef\temp{\write-1{hlend}}%
    \temp
    \endgroup 
  }%
  \let\hl@setparam\gobble
}%
{\catcode`\#=\other
\gdef\hlhash{#}}%
\def\hldriver@hypertex{%
  \def\hldest@type{xyz}%
  \let\hldest@opt@cmd \empty
  \def\hldest@driver{%
    \ifx\@hltype\hl@raw@word
      \csname \hldest@opt@cmd \endcsname
    \else
    \fi
  }%
  \let\hldest@typeh@raw \empty
  \let\hldest@typeh@xyz \empty
  \def\hl@type{name}%
  \ifx\hl@type@url\empty
    \def\hl@type@url{url}%
  \fi
  \ifx\hl@type@hrefext\empty
    \def\hl@type@hrefext{url}%
  \fi
  \let\hl@opt@cmd  \empty
  \let\hl@opt@ext  \empty
  \let\hl@opt@file \empty
  \def\hl@driver{%
    \ifx\@hltype\hl@raw@word
      \csname \hl@opt@cmd \endcsname
    \else
      \def\hlstart@preamble{html:<a href="}%
      \csname hl@typeh@\@hltype\endcsname
    \fi
  }%
  \let\hl@typeh@raw \empty
  \def\hl@typeh@name{\special{\hlstart@preamble \hlhash\@hllabel">}}%
  \def\hl@typeh@filename{%
    \special{%
      \hlstart@preamble
        file:\hl@opt@file\hl@opt@ext
        \ifempty\@hllabel \else \hlhash\@hllabel\fi
      ">%
    }%
  }%
  \def\hl@typeh@url{%
    \special{%
      \hlstart@preamble
        \@hllabel
      ">%
    }%
  }%
  \def\@hlend{\endgroup}
}%
\def\hldriver@pdftex{%
\ifpdf 
  \def\hldest@type{xyz}%
  \let\hldest@opt@width  \empty
  \let\hldest@opt@height \empty
  \let\hldest@opt@depth  \empty
  \let\hldest@opt@zoom   \empty
  \let\hldest@opt@cmd    \empty
  \def\hldest@driver{%
    \ifx\@hltype\hl@raw@word
      \csname \hldest@opt@cmd \endcsname
    \else
      \pdfdest name{\@hllabel}\@hltype
        \csname hldest@typeh@\@hltype\endcsname
    \fi
  }%
  \let\hldest@typeh@raw   \empty
  \let\hldest@typeh@fit   \empty
  \let\hldest@typeh@fith  \empty
  \let\hldest@typeh@fitv  \empty
  \let\hldest@typeh@fitb  \empty
  \let\hldest@typeh@fitbh \empty
  \let\hldest@typeh@fitbv \empty
  \def\hldest@typeh@fitr{%
    \ifx\hldest@opt@width  \empty \else width  \hldest@opt@width  \fi
    \ifx\hldest@opt@height \empty \else height \hldest@opt@height \fi
    \ifx\hldest@opt@depth  \empty \else depth  \hldest@opt@depth  \fi
  }%
  \def\hldest@typeh@xyz{%
    \ifx\hldest@opt@zoom\empty \else zoom \hldest@opt@zoom \fi
  }%
  \def\hl@type{name}%
  \ifx\hl@type@url\empty
    \def\hl@type@url{url}%
  \fi
  \ifx\hl@type@hrefext\empty
    \def\hl@type@hrefext{url}%
  \fi
  \let\hl@opt@width   \empty
  \let\hl@opt@height  \empty
  \let\hl@opt@depth   \empty
  \def\hl@opt@bstyle  {S}%
  \def\hl@opt@bwidth  {1}%
  \let\hl@opt@bcolor  \empty
  \let\hl@opt@hlight  \empty
  \let\hl@opt@bdash   \empty
  \let\hl@opt@pagefit \empty
  \let\hl@opt@cmd     \empty
  \let\hl@opt@file    \empty
  \let\hl@opt@newwin  \empty
  \def\hl@driver{%
    \ifx\@hltype\hl@raw@word
      \csname \hl@opt@cmd \endcsname
    \else
      \let\hl@BSspec\relax 
      \ifx\hl@opt@bstyle \empty
        \ifx\hl@opt@bwidth \empty
          \ifx\hl@opt@bdash \empty
            \let\hl@BSspec\empty 
          \fi
        \fi
      \fi
      \def\hlstart@preamble{%
        \pdfstartlink
          \ifx\hl@opt@width  \empty \else width  \hl@opt@width  \fi
          \ifx\hl@opt@height \empty \else height \hl@opt@height \fi
          \ifx\hl@opt@depth  \empty \else depth  \hl@opt@depth \fi
          attr{%
            \ifx\hl@opt@bcolor\empty\else /C[\hl@opt@bcolor]\fi
            \ifx\hl@opt@hlight\empty\else /H/\hl@opt@hlight\fi
            \ifx\hl@BSspec\relax
              /BS<<%
                /Type/Border%
                \ifx\hl@opt@bstyle\empty\else /S/\hl@opt@bstyle\fi
                \ifx\hl@opt@bwidth\empty\else /W \hl@opt@bwidth\fi
                \ifx\hl@opt@bdash\empty \else /D[\hl@opt@bdash]\fi
              >>%
            \fi
          }%
      }%
      \csname hl@typeh@\@hltype\endcsname
    \fi
  }%
  \let\hl@typeh@raw\empty
  \def\hl@typeh@name{\hlstart@preamble goto name{\@hllabel}}%
  \def\hl@typeh@num{\hlstart@preamble  goto num \@hllabel}%
  \def\hl@typeh@page{%
    \count@=\@hllabel
    \advance\count@ by-1
    \hlstart@preamble
    user{%
      /Subtype/Link%
      /Dest%
        [\the\count@
          \ifx\hl@opt@pagefit\empty/Fit\else\hl@opt@pagefit\fi]%
    }%
  }%
  \def\hl@typeh@filename{\hl@file{(\@hllabel)}}%
  \def\hl@typeh@filepage{%
    \count@=\@hllabel
    \advance\count@ by-1
    \hl@file{%
      [\the\count@ \ifx\hl@opt@pagefit\empty/Fit\else\hl@opt@pagefit\fi]%
    }%
  }%
  \def\hl@file##1{%
    \hlstart@preamble
    user{%
      /Subtype/Link%
      /A<<%
        /Type/Action%
        /S/GoToR%
        /D##1%
        /F(\hl@opt@file)%
        \ifx\hl@opt@newwin\empty \else
          /NewWindow \ifcase\hl@opt@newwin false\else true\fi
        \fi
      >>%
    }%
  }%
  \def\hl@typeh@url{%
    \hlstart@preamble
    user{%
      /Subtype/Link%
      /A<<%
        /Type/Action%
        /S/URI%
        /URI(\@hllabel)%
      >>%
    }%
  }%
  \def\@hlend{\pdfendlink\endgroup}
\else 
  \message{Eplain warning: `pdftex' hyperlink driver: PDF output is^^J
           \space not enabled, falling back on `nolinks' driver.}%
  \hldriver@nolinks
\fi
}%
\def\hldriver@dvipdfm{%
  \def\hldest@type{xyz}%
  \let\hldest@opt@left   \empty
  \let\hldest@opt@top    \empty
  \let\hldest@opt@right  \empty
  \let\hldest@opt@bottom \empty
  \let\hldest@opt@zoom   \empty
  \let\hldest@opt@cmd    \empty
  \def\hldest@driver{%
    \ifx\@hltype\hl@raw@word
      \csname \hldest@opt@cmd \endcsname
    \else
      \def\hldest@preamble{%
        pdf: dest (\@hllabel) [@thispage
      }%
      \csname hldest@typeh@\@hltype\endcsname
    \fi
  }%
  \let\hldest@typeh@raw\empty
  \def\hldest@typeh@fit{%
    \special{\hldest@preamble /Fit]}%
  }%
  \def\hldest@typeh@fith{%
    \special{\hldest@preamble /FitH
      \ifx\hldest@opt@top\empty @ypos \else \hldest@opt@top \fi]}%
  }%
  \def\hldest@typeh@fitv{%
    \special{\hldest@preamble /FitV
      \ifx\hldest@opt@left\empty @xpos \else \hldest@opt@left \fi]}%
  }%
  \def\hldest@typeh@fitb{%
    \special{\hldest@preamble /FitB]}%
  }%
  \def\hldest@typeh@fitbh{%
    \special{\hldest@preamble /FitBH
      \ifx\hldest@opt@top\empty @ypos \else \hldest@opt@top \fi]}%
  }%
  \def\hldest@typeh@fitbv{%
    \special{\hldest@preamble /FitBV
      \ifx\hldest@opt@left\empty @xpos \else \hldest@opt@left \fi]}%
  }%
  \def\hldest@typeh@fitr{%
    \special{\hldest@preamble /FitR
      \ifx\hldest@opt@left\empty @xpos\else\hldest@opt@left\fi\space
      \ifx\hldest@opt@bottom\empty @ypos\else\hldest@opt@bottom\fi\space
      \ifx\hldest@opt@right\empty @xpos\else\hldest@opt@right\fi\space
      \ifx\hldest@opt@top\empty @ypos\else\hldest@opt@top \fi]}%
  }%
  \def\hldest@typeh@xyz{%
    \begingroup
      \ifx\hldest@opt@zoom\empty
        \count1=\z@ \count2=\z@
      \else
        \count2=\hldest@opt@zoom
        \count1=\count2 \divide\count1 by 1000
        \count3=\count1 \multiply\count3 by 1000
        \advance\count2 by -\count3
      \fi
      \special{\hldest@preamble /XYZ
        \ifx\hldest@opt@left\empty @xpos\else\hldest@opt@left\fi\space
        \ifx\hldest@opt@top\empty @ypos\else\hldest@opt@top\fi\space
        \the\count1.\the\count2]}%
    \endgroup
  }%
  \def\hl@type{name}%
  \ifx\hl@type@url\empty
    \def\hl@type@url{url}%
  \fi
  \ifx\hl@type@hrefext\empty
    \def\hl@type@hrefext{url}%
  \fi
  \def\hl@opt@bstyle  {S}%
  \def\hl@opt@bwidth  {1}%
  \let\hl@opt@bcolor  \empty
  \let\hl@opt@hlight  \empty
  \let\hl@opt@bdash   \empty
  \let\hl@opt@pagefit \empty
  \let\hl@opt@cmd     \empty
  \let\hl@opt@file    \empty
  \let\hl@opt@newwin  \empty
  \def\hl@driver{%
    \ifx\@hltype\hl@raw@word
      \csname \hl@opt@cmd \endcsname
    \else
      \let\hl@BSspec\relax 
      \ifx\hl@opt@bstyle \empty
        \ifx\hl@opt@bwidth \empty
          \ifx\hl@opt@bdash \empty
            \let\hl@BSspec\empty 
          \fi
        \fi
      \fi
      \def\hlstart@preamble{%
        pdf: beginann
          <<%
            /Type/Annot%
            /Subtype/Link%
            \ifx\hl@opt@bcolor\empty\else /C[\hl@opt@bcolor]\fi
            \ifx\hl@opt@hlight\empty\else /H/\hl@opt@hlight\fi
            \ifx\hl@BSspec\relax
              /BS<<%
                /Type/Border%
                \ifx\hl@opt@bstyle\empty\else /S/\hl@opt@bstyle\fi
                \ifx\hl@opt@bwidth\empty\else /W \hl@opt@bwidth\fi
                \ifx\hl@opt@bdash\empty \else /D[\hl@opt@bdash]\fi
              >>%
            \fi
      }%
      \csname hl@typeh@\@hltype\endcsname
    \fi
  }%
  \let\hl@typeh@raw\empty
  \def\hl@typeh@name{\special{\hlstart@preamble /Dest(\@hllabel)>>}}%
  \def\hl@typeh@page{%
    \count@=\@hllabel
    \advance\count@ by-1
    \special{%
      \hlstart@preamble
      /Dest[\the\count@
            \ifx\hl@opt@pagefit\empty/Fit\else\hl@opt@pagefit\fi]%
     >>%
    }%
  }%
  \def\hl@typeh@filename{\hl@file{(\@hllabel)}}%
  \def\hl@typeh@filepage{%
    \count@=\@hllabel
    \advance\count@ by-1
    \hl@file{%
      [\the\count@ \ifx\hl@opt@pagefit\empty/Fit\else\hl@opt@pagefit\fi]%
    }%
  }%
  \def\hl@file##1{%
    \special{%
      \hlstart@preamble
      /A<<%
        /Type/Action%
        /S/GoToR%
        /D##1%
        /F(\hl@opt@file)%
        \ifx\hl@opt@newwin\empty \else
          /NewWindow \ifcase\hl@opt@newwin false\else true\fi
        \fi
      >>%
     >>%
    }%
  }%
  \def\hl@typeh@url{%
    \special{%
      \hlstart@preamble
      /A<<%
        /Type/Action%
        /S/URI%
        /URI(\@hllabel)%
      >>%
     >>%
    }%
  }%
  \def\@hlend{\endgroup}
}%
\def\href{%
  \bgroup
    \uncatcodespecials
    \catcode`\{=1 \catcode`\}=2
    \@href
}%
\def\@href#1{
  \egroup
  \edef\@hreftmp{\ifempty{#1}{}\fi}
  \expandafter\@@href\@hreftmp#1\@@
}%
\def\href@end@int{\hlend@impl{hrefint}}%
\def\href@end@ext{\hlend@impl{hrefext}}%
\def\@@href#1#2\@@{%
  \def\@hreftmp{#1}%
  \ifx\@hreftmp\hlhash
    \let\href@end\href@end@int
    \hlstart@impl{hrefint}{#2}%
  \else
    \let\href@end\href@end@ext
    \hlstart@impl{hrefext}{#1#2}%
  \fi
  \@@@href
}%
\def\@@@href{%
  \futurelet\@hreftmp\href@
}%
\def\href@{%
  \ifcat\bgroup\noexpand\@hreftmp
    \let\@hreftmp\href@@
  \else
    \let\@hreftmp\href@@@
  \fi
  \@hreftmp
}%
\def\href@@{\bgroup\aftergroup\href@end \let\@hreftmp}%
\def\href@@@#1{#1\href@end}%
\def\hldeston{\errmessage{Please enable hyperlinks with
  \string\enablehyperlinks\space before using hyperlink commands
  (consider selecting the `nolinks' driver to ignore all hyperlink
  commands in your document)}}%
\let\hldestoff\hldeston \let\hlon\hldeston \let\hloff\hldeston
\let\hlstart\hldeston \let\hlend\hldeston \let\hldest\hldeston
\let\hl@setparam\hldeston
\@hloff[*]\@hldestoff[*]%
\newif\ifusepkg@miniltx@loaded
\newcount\usepkg@recursion@level
\def\usepkg@rcrs{\the\usepkg@recursion@level}%
\let\usepkg@at@begin@document\empty
\let\usepkg@at@end@of@package\empty
\def\usepkg@word@autopict{autopict}%
\def\usepkg@word@psfrag{psfrag}%
\long\def\beginpackages#1\endpackages{%
  \let\usepackage\real@usepackage
  \let\DoNotLoadEpstopdf=t
  \let\eplaininput=\input
  #1%
  \usepkg@at@begin@document
  \global\let\usepkg@at@begin@document\empty
  \global\let\usepackage\fake@usepackage
  \let\packageinput=\input
  \let\input=\eplaininput
  \ifx\resetatcatcode\@undefined \else\resetatcatcode \fi
}%
\def\fake@usepackage{\errmessage{You should not use \string\usepackage\space outside of^^J
  \@spaces\string\beginpackages...\string\endpackages\space environment}%
}%
\def\eplain@RequirePackage{%
  \global\ece\let{usepkg@save@pkg\usepkg@rcrs}\usepkg@pkg
  \global\ece\let{usepkg@save@options\usepkg@rcrs}\usepkg@options
  \global\ece\let{usepkg@save@date\usepkg@rcrs}\usepkg@date
  \global\ece\let{usepkg@at@end@of@package\usepkg@rcrs}\usepkg@at@end@of@package
  \global\advance\usepkg@recursion@level by\@ne
  \real@usepackage
}%
\let\usepackage\fake@usepackage
\def\real@usepackage{\@getoptionalarg\@finusepackage}%
\def\@finusepackage#1{%
  \let\usepkg@options\@optionalarg
  \ifempty{#1}%
    \errmessage{No packages specified}%
  \fi
  \ifusepkg@miniltx@loaded \else
    \testfileexistence[miniltx]{tex}%
    \if@fileexists
      \global\usepkg@miniltx@loadedtrue
      \global\let\RequirePackage\eplain@RequirePackage
      \global\let\DeclareOption\eplain@DeclareOption
      \global\let\PassOptionsToPackage\eplain@PassOptionsToPackage
      \global\let\ExecuteOptions\eplain@ExecuteOptions
      \gdef\ProcessOptions{\@ifstar\eplain@ProcessOptions
                                   \eplain@ProcessOptions}%
      \global\let\AtBeginDocument\eplain@AtBeginDocument
      \global\let\AtEndOfPackage\eplain@AtEndOfPackage
      \global\let\ProvidesFile\eplain@ProvidesFile
      \global\let\ProvidesPackage\eplain@ProvidesPackage
    \else
      \errmessage{miniltx.tex not found, cannot load LaTeX packages}%
    \fi
  \fi
  \@ifnextchar[
    {\@finfinusepackage{#1}}%
    {\@finfinusepackage{#1}[]}%
}%
\def\@finfinusepackage#1[#2]{%
  \edef\usepkg@date{#2}%
  \let\usepkg@load@list\empty
  \for\usepkg@pkg:=#1\do{%
    \toks@=\expandafter{\usepkg@load@list}%
    \edef\usepkg@load@list{%
      \the\toks@
      \noexpand\usepkg@load@pkg{\usepkg@pkg}%
    }%
  }%
  \usepkg@load@list
  \ifnum\usepkg@recursion@level>0
    \global\advance\usepkg@recursion@level by\m@ne
    \expandafter\let\expandafter\usepkg@pkg\csname usepkg@save@pkg\usepkg@rcrs\endcsname
    \expandafter\let\expandafter\usepkg@options\csname usepkg@save@options\usepkg@rcrs\endcsname
    \expandafter\let\expandafter\usepkg@date\csname usepkg@save@date\usepkg@rcrs\endcsname
    \expandafter\let\expandafter\usepkg@at@end@of@package\csname usepkg@at@end@of@package\usepkg@rcrs\endcsname
    \global\ece\let{usepkg@save@pkg\usepkg@rcrs}\undefined
    \global\ece\let{usepkg@save@options\usepkg@rcrs}\undefined
    \global\ece\let{usepkg@save@date\usepkg@rcrs}\undefined
    \global\ece\let{usepkg@at@end@of@package\usepkg@rcrs}\undefined
  \fi
}%
\def\usepkg@load@pkg#1{%
  \def\usepkg@pkg{#1}%
  \ifx\usepkg@pkg\usepkg@word@autopict
    \testfileexistence[picture]{tex}%
    \if@fileexists \else
      \errmessage{Loader `picture.tex' for package `\usepkg@pkg' not found}%
    \fi
  \else
    \ifx\usepkg@pkg\usepkg@word@psfrag
      \testfileexistence[psfrag]{tex}%
      \if@fileexists \else
        \errmessage{Loader `psfrag.tex' for package `\usepkg@pkg' not found}%
      \fi
    \fi
  \fi
  \ifundefined{ver@\usepkg@pkg.sty}%
    \expandafter\@finusepkg@load@pkg
  \else
    \immediate\write-1{^^J\linenumber Eplain: package `\usepkg@pkg' already
             loaded, skipping reloading}%
  \fi
}%
\def\@finusepkg@load@pkg{%
  \testfileexistence[\usepkg@pkg]{sty}%
  \if@fileexists \else
    \errmessage{Package `\usepkg@pkg' not found}%
  \fi
  \expandafter\let\expandafter\temp\csname usepkg@options@\usepkg@pkg\endcsname
  \ifx\temp\relax
    \let\temp\empty
  \fi
  \ifx\temp\empty
    \let\temp\usepkg@options
  \else
    \ifx\usepkg@options\empty \else
      \edef\temp{\temp,\usepkg@options}%
    \fi
  \fi
  \global\ece\let{usepkg@options@\usepkg@pkg}\temp
  \let\usepackage\eplain@RequirePackage
  \global\let\usepkg@at@end@of@package\empty
  \ifx\usepkg@pkg\usepkg@word@autopict
    \input picture.tex
  \else
    \ifx\usepkg@pkg\usepkg@word@psfrag
      \input \usepkg@pkg.tex
    \else
      \input \usepkg@pkg.sty
    \fi
  \fi
  \usepkg@at@end@of@package
  \global\let\usepkg@at@end@of@package\empty
  \let\usepackage\real@usepackage
  \global\ece\let{usepkg@options@\usepkg@pkg}\undefined
  \def\Url@HyperHook##1{\hlstart@impl{url}{\Url@String}##1\hlend@impl{url}}%
}%
\def\eplain@DeclareOption#1#2{%
  \toks@{#2}%
  \expandafter\xdef\csname usepkg@option@\usepkg@pkg @#1\endcsname{\the\toks@}%
}%
\def\eplain@PassOptionsToPackage#1#2{%
  \ifempty{#1}\else
    \for\usepkg@temp:=#2\do{%
      \expandafter\let\expandafter\temp\csname usepkg@options@\usepkg@temp\endcsname
      \ifx\temp\relax
        \let\temp\empty
      \fi
      \ifx\temp\empty
        \edef\temp{#1}%
      \else
        \edef\temp{\temp,#1}%
      \fi
      \global\ece\let{usepkg@options@\usepkg@temp}\temp
    }%
  \fi
}%
\def\usepkg@exec@options#1{%
  \for\CurrentOption:=#1\do{%
    \expandafter\let\expandafter\usepkg@temp
      \csname usepkg@option@\usepkg@pkg @\CurrentOption\endcsname
    \ifx\usepkg@temp\relax
      \expandafter\let\expandafter\temp\csname usepkg@option@\usepkg@pkg @*\endcsname
      \ifx\temp\relax
        \errmessage{Unknown option `\CurrentOption' to package `\usepkg@pkg'}%
      \else
        \temp
      \fi
    \else
      \usepkg@temp
    \fi
  }%
}%
\let\eplain@ExecuteOptions\usepkg@exec@options
\def\eplain@ProcessOptions{%
  \expandafter\usepkg@exec@options\csname usepkg@options@\usepkg@pkg\endcsname
}%
\def\usepkg@accumulate@cmds#1#2{%
  \toks@=\expandafter{#1}%
  \toks@ii={#2}%
  \xdef#1{\the\toks@\the\toks@ii}%
}%
\def\eplain@AtBeginDocument{\usepkg@accumulate@cmds\usepkg@at@begin@document}%
\def\eplain@AtEndOfPackage{\usepkg@accumulate@cmds\usepkg@at@end@of@package}%
\def\eplain@ProvidesPackage#1{%
  \@ifnextchar[
    {\eplain@pr@videpackage{#1.sty}}{\eplain@pr@videpackage#1[]}%
}%
\def\eplain@pr@videpackage#1[#2]{%
  \wlog{#1: #2}%
  \expandafter\xdef\csname ver@#1\endcsname{#2}%
  \@ifl@t@r{#2}\usepkg@date{}%
    {\message{Warning: you have requested package `\usepkg@pkg', version \usepkg@date,^^J
       \@spaces but only version `\csname ver@#1\endcsname' is available.}}%
}%
\def\eplain@ProvidesFile#1{%
  \@ifnextchar[
    {\eplain@pr@videfile{#1}}{\eplain@pr@videfile#1[]}%
}%
\def\eplain@pr@videfile#1[#2]{%
  \wlog{#1: #2}%
  \expandafter\xdef\csname ver@#1\endcsname{#2}%
}%
\def\@ifl@ter#1#2{%
  \expandafter\@ifl@t@r
    \csname ver@#2.#1\endcsname
}%
\def\@ifl@t@r#1#2{%
  \ifnum\expandafter\@parse@version#1//00\@nil<%
        \expandafter\@parse@version#2//00\@nil
    \expandafter\@secondoftwo
  \else
    \expandafter\@firstoftwo
  \fi
}%
\def\@parse@version#1/#2/#3#4#5\@nil{#1#2#3#4 }%

\def\strip@prefix#1>{}%
\def\@ifpackageloaded#1{%
  \expandafter\ifx\csname ver@#1.sty\endcsname\relax
    \expandafter\@secondoftwo
  \else
    \expandafter\@firstoftwo
  \fi
}%
\long\def\g@addto@macro#1#2{%
  \begingroup
    \toks@\expandafter{#1#2}%
    \xdef#1{\the\toks@}%
  \endgroup
}%
\def\PackageWarning#1#2{%
  \begingroup
    \newlinechar=10 %
    \def\MessageBreak{%
      ^^J(#1)\@spaces\@spaces\@spaces\@spaces
    }%
    \immediate\write16{^^JPackage #1 Warning: #2\on@line.^^J}%
  \endgroup
}%
\def\PackageWarningNoLine#1#2{%
  \PackageWarning{#1}{#2\@gobble}%
}%
\def\on@line{ on input line \the\inputlineno}%
\def\@spaces{\space\space\space\space}%
\def\@inmatherr#1{%
   \relax
   \ifmmode
     \errmessage{The command is invalid in math mode}%
   \fi
}%
\let\protected@edef\edef
\let\wlog = \@plainwlog
\catcode`@ = \@eplainoldatcode
\def\eplain{t}%
{\edef\plainversion{\fmtversion}%
 \xdef\fmtversion{3.8: 12 May 2016 (and plain \plainversion)}%
}%
\def\Foot#1{\numberedfootnote{\ #1}}%
\vrule width-67pt%
\beginpackages\vrule width-3pt\usepackage{graphicx}\usepackage{url}\usepackage[dvipsnames]{color}\endpackages
\def\grf{eps}%
\def\Foot#1{\numberedfootnote{\ #1}}%
\vrule width-67pt%
\beginpackages\vrule width-3pt\usepackage{graphicx}\usepackage{url}\usepackage[dvipsnames]{color}\endpackages
                                                                                                          \catcode`\@=9
                                                                    \newcount\Count\newcount\CountArg\newcount\CountVar\newcount\CountEnd\newcount\CountReturn%
  \newcount\CountO\newcount\CountI\newcount\CountII\newcount\CountIII\newcount\CountIV\newcount\CountV\newcount\CountVI\newcount\CountVII\newcount\CountVIII\newcount\CountIX\newcount\CountX\newcount\CountXI\newcount\CountXII%
       \def\UpCon#1{\Count=#1\advance\Count by 1\xdef#1{\the\Count}}%
\def\DownCon#1{\Count=#1\advance\Count by -1\xdef#1{\the\Count}}%
                                                                                              \newdimen\DimenHI\newdimen\DimenHII%
                                                                    \newdimen\Dimen\newdimen\DimenArg\newdimen\DimenVar\newdimen\DimenEnd\newdimen\DimenReturn%
  \newdimen\DimenO\newdimen\DimenI\newdimen\DimenII\newdimen\DimenIII\newdimen\DimenIV\newdimen\DimenV\newdimen\DimenVI\newdimen\DimenVII\newdimen\DimenVIII\newdimen\DimenIX\newdimen\DimenX\newdimen\DimenXI\newdimen\DimenXII%
                                                                                                   \def\Load#1{#1=\DimenReturn}%
                                                                                           \newskip\SkipI\newskip\SkipII\newskip\SkipIII%
                                                                                         \newtoks\Toks\newtoks\ToksEnd\newtoks\ToksReturn%
               \newtoks\ToksO\newtoks\ToksI\newtoks\ToksII\newtoks\ToksIII\newtoks\ToksIV\newtoks\ToksV\newtoks\ToksVI\newtoks\ToksVII\newtoks\ToksVIII\newtoks\ToksIX\newtoks\ToksX\newtoks\ToksXI\newtoks\ToksXII%
                                                                                         \def\Concatenate#1#2{\ToksReturn={\the#1\the#2}}%
                                                                                               \def\Split#1#2\end{\ToksReturn={#2}\ToksEnd={#1}}%
                                                                                        \def\FirstChar#1{\expandafter\Split\string#1\end}%
                                                                                 \def\IsControl#1{\FirstChar{\global}%
                                                                                                  \ToksO=\ToksEnd%
                                                                                                  \FirstChar{#1}%
                                                                                                  \if\the\ToksO\the\ToksEnd \CountReturn=1%
                                                                                                                     \else \CountReturn=0\fi%
                                                                                                  \ifx#1\Infty\CountReturn=2\fi%
                                                                                                  \ifx#1\One\CountReturn=3\fi%
                                                                                                  \ifx#1\Zero\CountReturn=4\fi}%
									                                                      \def\testVoid#1#2#3{\def\test{#1}%
									                                                                          \ifx\test\empty #2%
									                                                                                    \else #3\fi}%
                                                                                            \def\addII#1#2{\DimenReturn=0pt\DimenReturn=#1\advance\DimenReturn by#2}%
                                                                                            \def\addIII#1#2#3{\addII{#1}{#2}\advance\DimenReturn by #3}%
                                                                                            \def\addIV#1#2#3#4{\addIII{#1}{#2}{#3}\advance\DimenReturn by #4}%
                                                                                            \def\addV#1#2#3#4#5{\addIV{#1}{#2}{#3}{#4}\advance\DimenReturn by #5}%
                                                                                            \def\addVI#1#2#3#4#5#6{\addV{#1}{#2}{#3}{#4}{#5}\advance\DimenReturn by #6}%
                                                                                            \def\addVII#1#2#3#4#5#6#7{\addVI{#1}{#2}{#3}{#4}{#5}{#6}\advance\DimenReturn by #7}%
                                                                                            \def\addVIII#1#2#3#4#5#6#7#8{\addVII{#1}{#2}{#3}{#4}{#5}{#6}{#7}\advance\DimenReturn by #8}%
                                                                                            \def\addIX#1#2#3#4#5#6#7#8#9{\addVIII{#1}{#2}{#3}{#4}{#5}{#6}{#7}{#8}\advance\DimenReturn by #9}%
                                                                                            %
                                                                                            \def\divideII#1#2{\DimenReturn=#1%
                                                                                                              \divide\DimenReturn by#2}%
                                                                                            \def\kernII#1#2{\addII{\Dimen}{#1}{#2}\kern\Dimen}%
                                                                                            \def\kernIII#1#2#3{\addIII{\Dimen}{#1}{#2}{#3}\kern\Dimen}%
                                                                                            \def\kernIV#1#2#3#4{\addIV{\Dimen}{#1}{#2}{#3}{#4}\kern\Dimen}%
                                                                                            \def\kernV#1#2#3#4#5{\addV{\Dimen}{#1}{#2}{#3}{#4}{#5}\kern\Dimen}%
       												                                        \def\kn#1{\kern#1pt}%
                                                                                            %
                                                                                            %
                                                                                            %
                                                                                            
                                                                              \newbox\Box\newbox\BoxArg\newbox\BoxVar\newbox\BoxEnd\newbox\BoxReturn%
                            \newbox\BoxO\newbox\BoxI\newbox\BoxII\newbox\BoxIII\newbox\BoxIV\newbox\BoxV\newbox\BoxVI\newbox\BoxVII\newbox\BoxVIII\newbox\BoxIX\newbox\BoxX\newbox\BoxXI\newbox\BoxXII%
                                                                                                        \def\Vbox#1{\OffVskip\vbox{#1}\OnVskip}%
                                                                                                        \def\Vtop#1{\OffVskip\vtop{#1}\OnVskip}%
                                                                                            \def\hboxI#1#2{\hbox to#1{\hskip0pt plus.001pt minus.001pt #2}}%
                                                                                            \def\hboxII#1#2#3{\addII{#1}{#2}\hboxI{\DimenReturn}{#3}}%
                                                                                            \def\hboxi#1#2{\hbox to#1pt{\hskip0pt plus.001pt minus.001pt #2}}%
                                                                                            %
                                                                                            %
                                                                                            %
                                                                                            %
                                                                                            %
                                                                                            %
                                                                                            %
                                                                                            %
                                                                                            %
                                                                                            %
                                                                                            %
                                                                                            %
                                                                                            %
                                                                                               \def\Lower#1#2{\lower#1\hbox{#2}}%
                                                                                               \def\Raise#1#2{\raise#1\hbox{#2}}%
                                                                            \def\Stack#1#2#3{\setbox\Box\hbox{#1}%
                                                                                             \Dimen=\wd\Box%
                                                                                             \divide\Dimen by 2%
                                                                                             \setbox\BoxArg\hbox{#2}%
                                                                                             \DimenArg=\wd\BoxArg%
                                                                                             \divide\DimenArg by 2%
                                                                                             \ifdim\wd\Box<\wd\BoxArg \advance\DimenArg by-\Dimen%
                                                                                                                      \setbox\BoxVar\hbox{\kern\DimenArg\box\Box}%
                                                                                                                      \Vbox{\box\BoxVar\kern#3\box\BoxArg}%
                                                                                                                \else \advance\Dimen by-\DimenArg%
                                                                                                                      \setbox\BoxVar\hbox{\kern\Dimen\box\BoxArg}%
                                                                                                                      \Vbox{\box\Box\kern#3\box\BoxVar}\fi}%
                                                                                \def\StackReturn#1#2#3{\setbox\BoxReturn\vbox{\Stack{#1}{#2}{#3}}}%
                                                                                               \def\RaiseBox#1{\setbox\Box\hbox{#1}%
                                                                                                               \hbox{\raise\dp\Box\copy\Box}}%
                                                                        \def\RaiseBoxReturn#1{\setbox\Box\hbox{#1}%
                                                                                              \setbox\BoxReturn\hbox{\raise\dp\Box\copy\Box}}%
                                                                             \def\FrameBox#1#2#3{\setbox\Box\hbox{\vrule width#2%
                                                                                                                  \Vbox{\hrule height#2%
                                                                                                                        \kern#3%
                                                                                                                        \hbox{\kern#3%
                                                                                                                              \hbox{#1}%
                                                                                                                              \kern#3}%
                                                                                                                        \kern#3%
                                                                                                                        \hrule height#2}%
                                                                                                                  \vrule width#2}%
                                                                                                 \Dimen=#3\advance\Dimen by#2%
                                                                                                 \hbox{\lower\Dimen\hbox{\box\Box}}}%
                                                                                 \def\FrameBoxReturn#1#2#3{\setbox\BoxReturn\hbox{\FrameBox{#1}{#2}{#3}}}%
                                                 \def\CenterPlaceBox#1#2#3{\RaiseBox{\hbox{#1}}%
                                                                           \setbox\Box\hbox{\copy\BoxReturn}%
                                                                           \let\Ht=\Dimen \Ht=#2\advance\Ht by-\ht\Box\divide\Ht by2%
                                                                           \let\Wd=\DimenEnd\Wd=#3\advance\Wd by-\wd\Box\divide\Wd by2%
                                                                           \hbox{\kern\Wd\Vbox{\kern\Ht\copy\Box\kern\Ht}\kern\Wd}}%
                                                                       \def\CenterPlaceBoxReturn#1#2#3{\setbox\BoxReturn\hbox{\CenterPlaceBox{#1}{#2}{#3}}}%
                                                        \def\MaxII#1#2{\Dimen=#1\DimenEnd=#2\ifdim\Dimen>\DimenEnd \DimenReturn=\Dimen \else \DimenReturn=\DimenEnd\fi}%
                                                        \def\MaxIII#1#2#3{\MaxII{#1}{#2}\Dimen=\DimenReturn\MaxII{\Dimen}{#3}}%
                                                        \def\MaxIV#1#2#3#4{\MaxIII{#1}{#2}{#3}\Dimen=\DimenReturn\MaxII{\Dimen}{#4}}%
                                                        %
                                                        %
                                                        %
                                                        %
                                                        %
                                                            %
                                                            %
															%
                                                            %
                                                            \def\HBOX#1#2#3#4#5#6#7#8#9{\hbox{\vrule width#1%
                                                                                              \Vbox{\hrule height#3%
                                                                                                    \vskip#7%
                                                                                          			\hbox{\hskip#5%
                                                                                                 		  #9%
                                                                                                 		  \hskip#6}%
                                                                                           			\vskip#8%
                                                                                           		    \hrule height#4}%
                                                                                     		  \vrule width#2}}%
                                                            \def\HBox#1#2#3{\HBOX{#1}{#1}{#1}{#1}{#2}{#2}{#2}{#2}{#3}}%
                                                            \def\BlueHBOX#1#2#3#4#5#6#7#8#9{\hbox{\Blue{\vrule width#1}%
                                                                                                  \Vbox{\Blue{\hrule height#3}%
                                                                                                          	  \vskip#7%
                                                                                           					  \hbox{\hskip#5%
                                                                                                 				    #9%
                                                                                                 					\hskip#6}%
                                                                                           					  \vskip#8%
                                                                                                              \Blue{\hrule height#4}}%
                                                                                     					\Blue{\vrule width#2}}}%
                                                            \def\BlueHBox#1#2#3{\BlueHBOX{#1}{#1}{#1}{#1}{#2}{#2}{#2}{#2}{#3}}%
                                                            \def\RedOrangeHBOX#1#2#3#4#5#6#7#8#9{\hbox{\RedOrange{\vrule width#1}%
                                                                                                       \Vbox{\RedOrange{\hrule height#3}%
                                                                                                          	 \vskip#7%
                                                                                           					 \hbox{\hskip#5%
                                                                                                 				   #9%
                                                                                                 				   \hskip#6}%
                                                                                           					 \vskip#8%
                                                                                                             \RedOrange{\hrule height#4}}%
                                                                                     					\RedOrange{\vrule width#2}}}%
                                                            \def\RedOrangeHBox#1#2#3{\RedOrangeHBOX{#1}{#1}{#1}{#1}{#2}{#2}{#2}{#2}{#3}}%
                                                            \def\VioletHBOX#1#2#3#4#5#6#7#8#9{\hbox{\Violet{\vrule width#1}%
                                                                                                    \Vbox{\Violet{\hrule height#3}%
                                                                                                          \vskip#7%
                                                                                           				  \hbox{\hskip#5%
                                                                                                 				#9%
                                                                                                 			    \hskip#6}%
                                                                                           				  \vskip#8%
                                                                                                          \Violet{\hrule height#4}}%
                                                                                     				\Violet{\vrule width#2}}}%
                                                            \def\VioletHBox#1#2#3{\VioletHBOX{#1}{#1}{#1}{#1}{#2}{#2}{#2}{#2}{#3}}%
                                                            \def\BrownHBOX#1#2#3#4#5#6#7#8#9{\hbox{\Brown{\vrule width#1}%
                                                                                                    \Vbox{\Brown{\hrule height#3}%
                                                                                                          \vskip#7%
                                                                                           				  \hbox{\hskip#5%
                                                                                                 				#9%
                                                                                                 			    \hskip#6}%
                                                                                           				  \vskip#8%
                                                                                                          \Brown{\hrule height#4}}%
                                                                                     				\Brown{\vrule width#2}}}%
                                                            \def\BrownHBox#1#2#3{\BrownHBOX{#1}{#1}{#1}{#1}{#2}{#2}{#2}{#2}{#3}}%
                                                                                                        \hsize=469.75499pt%
                                                                                                        \baselineskip=12pt
                                                                                                        \lineskiplimit=0pt%
                                                                                                          \lineskip=1pt%
                                                                                                         \parindent=20pt%
                                                                                                   \parfillskip=0pt plus 1fil%
                                                                                                \def\OffVskip{\SkipI=\baselineskip
                                                                                                              \baselineskip=-1000pt%
                                                                                                              \SkipII=\lineskiplimit%
                                                                                                              \lineskiplimit=16383pt%
                                                                                                              \SkipIII=\lineskip%
                                                                                                              \lineskip=0pt}%
                                                                                                \def\OnVskip{\baselineskip=\SkipI%
                                                                                                             \lineskiplimit=\SkipII%
                                                                                                             \lineskip=\SkipIII}%
                                                                                               %
                                                                                                        %
                                       \def\HeightCenterII#1#2{\setbox\BoxO\hbox{\RaiseBox{#1}}%
                                                               \setbox\BoxArg\hbox{\RaiseBox{#2}}%
                                                               \MaxII{\ht\BoxO}{\ht\BoxArg}%
                                                               \Dimen=\DimenReturn%
                                                               \ifdim\ht\BoxO<\Dimen \addII{\Dimen}{-\ht\BoxO}%
                                                                                     \divide\DimenReturn by 2%
                                                                                     \setbox\Box\hbox{\Vbox{\kern\DimenReturn%
                                                                                                            \box\BoxO%
                                                                                                            \kern\DimenReturn}}%
                                                                                     \setbox\BoxO\hbox{\box\Box}\fi%
                                                               \ifdim\ht\BoxArg<\Dimen \addII{\Dimen}{-\ht\BoxArg}%
                                                                                      \divide\DimenReturn by 2%
                                                                                      \setbox\Box\hbox{\Vbox{\kern\DimenReturn%
                                                                                                             \box\BoxArg%
                                                                                                             \kern\DimenReturn}}%
                                                                                      \setbox\BoxArg\hbox{\box\Box}\fi}%
                                    \def\HeightCenterIII#1#2#3{\setbox\BoxO\hbox{\RaiseBox{#1}}%
                                                               \setbox\BoxArg\hbox{\RaiseBox{#2}}%
                                                               \setbox\BoxVar\hbox{\RaiseBox{#3}}%
                                                               \MaxIII{\ht\BoxO}{\ht\BoxArg}{\ht\BoxVar}%
                                                               \Dimen=\DimenReturn%
                                                               \ifdim\ht\BoxO<\Dimen \addII{\Dimen}{-\ht\BoxO}%
                                                                                     \divide\DimenReturn by 2%
                                                                                     \setbox\Box\hbox{\Vbox{\kern\DimenReturn%
                                                                                                            \box\BoxO%
                                                                                                            \kern\DimenReturn}}%
                                                                                     \setbox\BoxO\hbox{\box\Box}\fi%
                                                               \ifdim\ht\BoxArg<\Dimen \addII{\Dimen}{-\ht\BoxArg}%
                                                                                      \divide\DimenReturn by 2%
                                                                                      \setbox\Box\hbox{\Vbox{\kern\DimenReturn%
                                                                                                             \box\BoxArg%
                                                                                                             \kern\DimenReturn}}%
                                                                                      \setbox\BoxArg\hbox{\box\Box}\fi%
                                                               \ifdim\ht\BoxVar<\Dimen \addII{\Dimen}{-\ht\BoxVar}%
                                                                                     \divide\DimenReturn by 2%
                                                                                     \setbox\Box\hbox{\Vbox{\kern\DimenReturn%
                                                                                                            \box\BoxVar%
                                                                                                            \kern\DimenReturn}}%
                                                                                     \setbox\BoxVar\hbox{\box\Box}\fi}%
                                   \def\HeightCenterIV#1#2#3#4{\setbox\BoxO\hbox{\RaiseBox{#1}}%
                                                               \setbox\BoxArg\hbox{\RaiseBox{#2}}%
                                                               \setbox\BoxVar\hbox{\RaiseBox{#3}}%
                                                               \setbox\BoxEnd\hbox{\RaiseBox{#4}}%
                                                               \MaxIV{\ht\BoxO}{\ht\BoxArg}{\ht\BoxVar}{\ht\BoxEnd}%
                                                               \Dimen=\DimenReturn%
                                                               \ifdim\ht\BoxO<\Dimen \addII{\Dimen}{-\ht\BoxO}%
                                                                                     \divide\DimenReturn by 2%
                                                                                     \setbox\Box\hbox{\Vbox{\kern\DimenReturn%
                                                                                                            \box\BoxO%
                                                                                                            \kern\DimenReturn}}%
                                                                                     \setbox\BoxO\hbox{\box\Box}\fi%
                                                               \ifdim\ht\BoxArg<\Dimen \addII{\Dimen}{-\ht\BoxArg}%
                                                                                      \divide\DimenReturn by 2%
                                                                                      \setbox\Box\hbox{\Vbox{\kern\DimenReturn%
                                                                                                             \box\BoxArg%
                                                                                                             \kern\DimenReturn}}%
                                                                                      \setbox\BoxArg\hbox{\box\Box}\fi%
                                                               \ifdim\ht\BoxVar<\Dimen \addII{\Dimen}{-\ht\BoxVar}%
                                                                                     \divide\DimenReturn by 2%
                                                                                     \setbox\Box\hbox{\Vbox{\kern\DimenReturn%
                                                                                                            \box\BoxVar%
                                                                                                            \kern\DimenReturn}}%
                                                                                     \setbox\BoxVar\hbox{\box\Box}\fi%
                                                               \ifdim\ht\BoxEnd<\Dimen \addII{\Dimen}{-\ht\BoxEnd}%
                                                                                     \divide\DimenReturn by 2%
                                                                                     \setbox\Box\hbox{\Vbox{\kern\DimenReturn%
                                                                                                            \box\BoxEnd%
                                                                                                            \kern\DimenReturn}}%
                                                                                     \setbox\BoxEnd\hbox{\box\Box}\fi}%
                                        \def\WidthCenterII#1#2{\setbox\BoxO\hbox{\RaiseBox{#1}}%
                                                               \setbox\BoxArg\hbox{\RaiseBox{#2}}%
                                                               \MaxII{\wd\BoxO}{\wd\BoxArg}%
                                                               \Dimen=\DimenReturn%
                                                               \ifdim\wd\BoxO<\Dimen \addII{\Dimen}{-\wd\BoxO}%
                                                                                     \divide\DimenReturn by 2%
                                                                                     \setbox\Box\hbox{\kern\DimenReturn%
                                                                                                      \box\BoxO%
                                                                                                      \kern\DimenReturn}%
                                                                                     \setbox\BoxO\hbox{\box\Box}\fi%
                                                               \ifdim\wd\BoxArg<\Dimen \addII{\Dimen}{-\wd\BoxArg}%
                                                                                      \divide\DimenReturn by 2%
                                                                                      \setbox\Box\hbox{\kern\DimenReturn%
                                                                                                       \box\BoxArg%
                                                                                                       \kern\DimenReturn}%
                                                                                      \setbox\BoxArg\hbox{\box\Box}\fi}%
                                     \def\WidthCenterIII#1#2#3{\setbox\BoxO\hbox{\RaiseBox{#1}}%
                                                               \setbox\BoxArg\hbox{\RaiseBox{#2}}%
                                                               \setbox\BoxVar\hbox{\RaiseBox{#3}}%
                                                               \MaxIII{\wd\BoxO}{\wd\BoxArg}{\wd\BoxVar}%
                                                               \Dimen=\DimenReturn%
                                                               \ifdim\wd\BoxO<\Dimen \addII{\Dimen}{-\wd\BoxO}%
                                                                                     \divide\DimenReturn by 2%
                                                                                     \setbox\Box\hbox{\kern\DimenReturn%
                                                                                                      \box\BoxO%
                                                                                                      \kern\DimenReturn}%
                                                                                     \setbox\BoxO\hbox{\box\Box}\fi%
                                                               \ifdim\wd\BoxArg<\Dimen \addII{\Dimen}{-\wd\BoxArg}%
                                                                                      \divide\DimenReturn by 2%
                                                                                      \setbox\Box\hbox{\kern\DimenReturn%
                                                                                                       \box\BoxArg%
                                                                                                       \kern\DimenReturn}%
                                                                                      \setbox\BoxArg\hbox{\box\Box}\fi%
                                                               \ifdim\wd\BoxVar<\Dimen \addII{\Dimen}{-\wd\BoxVar}%
                                                                                     \divide\DimenReturn by 2%
                                                                                     \setbox\Box\hbox{\kern\DimenReturn%
                                                                                                      \box\BoxVar%
                                                                                                      \kern\DimenReturn}%
                                                                                     \setbox\BoxVar\hbox{\box\Box}\fi}%
                                    \def\WidthCenterIV#1#2#3#4{\setbox\BoxO\hbox{\RaiseBox{#1}}%
                                                               \setbox\BoxArg\hbox{\RaiseBox{#2}}%
                                                               \setbox\BoxVar\hbox{\RaiseBox{#3}}%
                                                               \setbox\BoxEnd\hbox{\RaiseBox{#4}}%
                                                               \MaxIV{\wd\BoxO}{\wd\BoxArg}{\wd\BoxVar}{\wd\BoxEnd}%
                                                               \Dimen=\DimenReturn%
                                                               \ifdim\wd\BoxO<\Dimen \addII{\Dimen}{-\wd\BoxO}%
                                                                                     \divide\DimenReturn by 2%
                                                                                     \setbox\Box\hbox{\kern\DimenReturn%
                                                                                                      \box\BoxO%
                                                                                                      \kern\DimenReturn}%
                                                                                     \setbox\BoxO\hbox{\box\Box}\fi%
                                                               \ifdim\wd\BoxArg<\Dimen \addII{\Dimen}{-\wd\BoxO}%
                                                                                      \divide\DimenReturn by 2%
                                                                                      \setbox\Box\hbox{\kern\DimenReturn%
                                                                                                       \box\BoxArg%
                                                                                                       \kern\DimenReturn}%
                                                                                      \setbox\BoxArg\hbox{\box\Box}\fi%
                                                               \ifdim\wd\BoxVar<\Dimen \addII{\Dimen}{-\wd\BoxVar}%
                                                                                     \divide\DimenReturn by 2%
                                                                                     \setbox\Box\hbox{\kern\DimenReturn%
                                                                                                      \box\BoxVar%
                                                                                                      \kern\DimenReturn}%
                                                                                     \setbox\BoxVar\hbox{\box\Box}\fi%
                                                               \ifdim\wd\BoxEnd<\Dimen \addII{\Dimen}{-\wd\BoxEnd}%
                                                                                     \divide\DimenReturn by 2%
                                                                                     \setbox\Box\hbox{\kern\DimenReturn%
                                                                                                      \box\BoxEnd%
                                                                                                      \kern\DimenReturn}%
                                                                                     \setbox\BoxEnd\hbox{\box\Box}\fi}%
                                       \def\ArrayIIxII#1#2#3#4{\ifvoid\BoxX\DimenHI=2pt\DimenHII=2pt\fi%
                                                               \HeightCenterII{#3}{#4}%
                                                               \setbox\BoxIII\box\BoxO%
                                                               \setbox\BoxIV\box\BoxArg%
                                                               \HeightCenterII{#1}{#2}%
                                                               \setbox\BoxI\box\BoxO%
                                                               \setbox\BoxII\box\BoxArg%
                                                               \WidthCenterII{\box\BoxI}{\box\BoxIII}%
                                                               \setbox\BoxI\box\BoxO%
                                                               \setbox\BoxIII\box\BoxArg%
                                                               \WidthCenterII{\box\BoxII}{\box\BoxIV}%
                                                               \setbox\BoxII\box\BoxO%
                                                               \setbox\BoxIV\box\BoxArg%
                                                               \hbox{\Vbox{\box\BoxI\kern\DimenHI\box\BoxIII}%
                                                                           \kern\DimenHII%
                                                                     \Vbox{\box\BoxII\kern\DimenHI\box\BoxIV}}}%
                                  \def\ArrayIIxIII#1#2#3#4#5#6{\ifvoid\BoxX\DimenHI=2pt\DimenHII=2pt\fi%
                                                               \HeightCenterIII{#1}{#2}{#3}%
                                                               \setbox\BoxI\box\BoxO%
                                                               \setbox\BoxII\box\BoxArg%
                                                               \setbox\BoxIII\box\BoxVar%
                                                               \HeightCenterIII{#4}{#5}{#6}%
                                                               \setbox\BoxIV\box\BoxO%
                                                               \setbox\BoxV\box\BoxArg%
                                                               \setbox\BoxVI\box\BoxVar%
                                                               \WidthCenterII{\box\BoxI}{\box\BoxIV}%
                                                               \setbox\BoxI\box\BoxO%
                                                               \setbox\BoxIV\box\BoxArg%
                                                               \WidthCenterII{\box\BoxII}{\box\BoxV}%
                                                               \setbox\BoxII\box\BoxO%
                                                               \setbox\BoxV\box\BoxArg%
                                                               \WidthCenterII{\box\BoxIII}{\box\BoxVI}%
                                                               \setbox\BoxIII\box\BoxO%
                                                               \setbox\BoxVI\box\BoxArg%
                                                               \hbox{\Vbox{\box\BoxI\kern\DimenHI\box\BoxIV}%
                                                                     \kern\DimenHII%
                                                                     \Vbox{\box\BoxII\kern\DimenHI\box\BoxV}%
                                                                     \kern\DimenHII%
                                                                     \Vbox{\box\BoxIII\kern\DimenHI\box\BoxVI}}}%
                           \def\ArrayIIIxIII#1#2#3#4#5#6#7#8#9{\ifvoid\BoxX\DimenHI=2pt\DimenHII=2pt\fi%
                                                               \HeightCenterIII{#1}{#2}{#3}%
                                                               \setbox\BoxI\box\BoxO%
                                                               \setbox\BoxII\box\BoxArg%
                                                               \setbox\BoxIII\box\BoxVar%
                                                               \HeightCenterIII{#4}{#5}{#6}%
                                                               \setbox\BoxIV\box\BoxO%
                                                               \setbox\BoxV\box\BoxArg%
                                                               \setbox\BoxVI\box\BoxVar%
                                                               \HeightCenterIII{#7}{#8}{#9}%
                                                               \setbox\BoxVII\box\BoxO%
                                                               \setbox\BoxVIII\box\BoxArg%
                                                               \setbox\BoxIX\box\BoxVar%
                                                               \WidthCenterIII{\box\BoxI}{\box\BoxIV}{\box\BoxVII}%
                                                               \setbox\BoxI\box\BoxO%
                                                               \setbox\BoxIV\box\BoxArg%
                                                               \setbox\BoxVII\box\BoxVar%
                                                               \WidthCenterIII{\box\BoxII}{\box\BoxV}{\box\BoxVIII}%
                                                               \setbox\BoxII\box\BoxO%
                                                               \setbox\BoxV\box\BoxArg%
                                                               \setbox\BoxVIII\box\BoxVar%
                                                               \WidthCenterIII{\box\BoxIII}{\box\BoxVI}{\box\BoxIX}%
                                                               \setbox\BoxIII\box\BoxO%
                                                               \setbox\BoxVI\box\BoxArg%
                                                               \setbox\BoxIX\box\BoxVar%
                                                               \hbox{\Vbox{\box\BoxI\kern\DimenHI\box\BoxIV\kern\DimenHI\box\BoxVII}%
                                                                     \kern\DimenHII%
                                                                     \Vbox{\box\BoxII\kern\DimenHI\box\BoxV\kern\DimenHI\box\BoxVIII}%
                                                                     \kern\DimenHII%
                                                                     \Vbox{\box\BoxIII\kern\DimenHI\box\BoxVI\kern\DimenHI\box\BoxIX}}}%
                               \def\ArrayIIxIV#1#2#3#4#5#6#7#8{\ifvoid\BoxX\DimenHI=2pt\DimenHII=2pt\fi%
                                                               \HeightCenterIV{#1}{#2}{#3}{#4}%
                                                               \setbox\BoxI\box\BoxO%
                                                               \setbox\BoxII\box\BoxArg%
                                                               \setbox\BoxIII\box\BoxVar%
                                                               \setbox\BoxIV\box\BoxEnd%
                                                               \HeightCenterIV{#5}{#6}{#7}{#8}%
                                                               \setbox\BoxV\box\BoxO%
                                                               \setbox\BoxVI\box\BoxArg%
                                                               \setbox\BoxVII\box\BoxVar%
                                                               \setbox\BoxVIII\box\BoxEnd%
                                                               \WidthCenterIV{\box\BoxI}{\box\BoxII}{\box\BoxIII}{\box\BoxIV}%
                                                               \setbox\BoxI\box\BoxO%
                                                               \setbox\BoxII\box\BoxArg%
                                                               \setbox\BoxIII\box\BoxVar%
                                                               \setbox\BoxIV\box\BoxEnd%
                                                               \WidthCenterIV{\box\BoxV}{\box\BoxVI}{\box\BoxVII}{\box\BoxVIII}%
                                                               \setbox\BoxV\box\BoxO%
                                                               \setbox\BoxVI\box\BoxArg%
                                                               \setbox\BoxVII\box\BoxVar%
                                                               \setbox\BoxVIII\box\BoxEnd%
                                                               \hbox{\Vbox{\box\BoxI\kern\DimenHI\box\BoxV}%
                                                                     \kern\DimenHII%
                                                                     \Vbox{\box\BoxII\kern\DimenHI\box\BoxVI}%
                                                                     \kern\DimenHII%
                                                                     \Vbox{\box\BoxIII\kern\DimenHI\box\BoxVII}%
                                                                     \kern\DimenHII%
                                                                     \Vbox{\box\BoxIV\kern\DimenHI\box\BoxVIII}}}%
                                                                                        \def\Row#1#2#3{\setbox\Box\hbox{\RaiseBox{#1}}%
                                                                                                       \setbox\BoxArg\hbox{\RaiseBox{#2}}%
                                                                                                       \let\FillerI=\Dimen\let\FillerII=\DimenArg%
                                                                                                       \ifdim\ht\Box<\ht\BoxArg \FillerI=\ht\BoxArg\FillerII=0pt%
                                                                                                                                \advance\FillerI by -\ht\Box%
                                                                                                                                \divide\FillerI by 2%
                                                                                                                          \else \FillerII=\ht\BoxArg\FillerI=0pt%
                                                                                                                                \advance\FillerII by -\ht\BoxArg%
                                                                                                                                \divide\FillerII by 2\fi%
                                                                                                       \hbox{\Vbox{\kern\FillerI \box\Box \kern\FillerI}%
                                                                                                             \kern#3%
                                                                                                             \Vbox{\kern\FillerII \box\BoxArg \kern\FillerII}}}%
                                                                                  \def\RowReturn#1#2#3{\setbox\BoxReturn\hbox{\Row{#1}{#2}{#3}}}%
                                                                                                         \headline={\hfil}%
                                                                                                        \vsize=643.20255pt%
                                                                                                           \hsize=6.5in%
                                                                                                           \hoffset=0pt%
                                                                                                           \voffset=0pt%
                                                                                    %
                                                                  				    %
                                                                               		%
%
                                                                                    \def\Paragraph#1#2{\Dimen=\parindent%
                                                                                                       \parindent=#2%
                                                                                                       \par
                                                                                                       \vskip#1%
                                                                                                       \indent%
                                                                                                       \parindent=\Dimen}%
                                                                                    \def\Par{\Paragraph{0pt}{0pt}}%
                                                                                    \def\PAr{\Paragraph{\baselineskip}{0pt}}%
                                                                                    \def\PaR{\Paragraph{0pt}{20pt}}%
                                                                                    \def\PAR{\Paragraph{\baselineskip}{20pt}}%
														                                      \parindent=0pt%
                                                                                                 \def\Quad{\hskip7.5pt\relax}%
                                                                                                 \def\QUad{\hskip5pt\relax}%
                                                                            \font\RomanFive=cmr5
                                                                            \font\RomanSix=cmr6
                                                                            \font\RomanSeven=cmr7
                                                                            \font\RomanEight=cmr8
                                                                            \font\RomanNine=cmr9
                                                                            \font\RomanTen=cmr10
                                                                            \font\RomanTwelve=cmr12
                                                                            \font\RomanSeventeen=cmr17
                                                                            \font\BoldfaceFive=cmbx5
                                                                            \font\BoldfaceSix=cmbx6
                                                                            \font\BoldfaceSeven=cmbx7
                                                                            \font\BoldfaceEight=cmbx8
                                                                            \font\BoldfaceNine=cmbx9
                                                                            \font\BoldfaceTen=cmbx10
                                                                            \font\BoldfaceTwelve=cmbx12
                                                                            \font\BoldfaceSeventeen=cmbx17
                                                                            \font\SlantNine=cmsl9
                                                                            \font\SlantTen=cmsl10
                                                                            \font\ItalicEight=cmti8
                                                                            \font\ItalicNine=cmti9
                                                                            \font\ItalicTen=cmti10
                                                                            \font\TypewriterEight=cmtt8
                                                                            \font\TypewriterNine=cmtt9
                                                                            \font\TypewriterTen=cmtt10
                                                                            \font\MathIFive=cmmi5
                                                                            \font\MathISix=cmmi6
                                                                            \font\MathISeven=cmmi7
                                                                            \font\MathIEight=cmmi8
                                                                            \font\MathINine=cmmi9
                                                                            \font\MathITen=cmmi10
                                                                            \font\MathITwelve=cmmi12
                                                                            \font\MathIISix=cmsy6
                                                                            \font\MathIIEight=cmsy8
                                                                            \font\MathIINine=cmsy9
                                                                            \font\MathIITen=cmsy10
                                                                            \font\MathIIIEight=cmex8
                                                                            \font\MathIIINine=cmex9
                                                                            \font\MathIIITen=cmex10
                                                                            \font\EuropeanFive=eurm5
                                                                            \font\EuropeanSix=eurm6
                                                                            \font\EuropeanSeven=eurm7
                                                                            \font\EuropeanEight=eurm8
                                                                            \font\EuropeanNine=eurm9
                                                                            \font\EuropeanTen=eurm10
                                                                            \font\EuropeanBoldfaceFive=eurb5
                                                                            \font\EuropeanBoldfaceSix=eurb6
                                                                            \font\EuropeanBoldfaceSeven=eurb7
                                                                            \font\EuropeanBoldfaceEight=eurb8
                                                                            \font\EuropeanBoldfaceNine=eurb9
                                                                            \font\EuropeanBoldfaceTen=eurb10
                                                                            \font\GermanFive=eufm5
                                                                            \font\GermanSix=eufm6
                                                                            \font\GermanSeven=eufm7
                                                                            \font\GermanEight=eufm8
                                                                            \font\GermanNine=eufm9
                                                                            \font\GermanTen=eufm10
                                                                            \font\GermanBoldfaceFive=eufb5
                                                                            \font\GermanBoldfaceSix=eufb6
                                                                            \font\GermanBoldfaceSeven=eufb7
                                                                            \font\GermanBoldfaceEight=eufb8
                                                                            \font\GermanBoldfaceNine=eufb9
                                                                            \font\GermanBoldfaceTen=eufb10
                                                                            \font\RussianFive=wncyr5
                                                                            \font\RussianSix=wncyr6
                                                                            \font\RussianSeven=wncyr7
                                                                            \font\RussianEight=wncyr8
                                                                            \font\RussianNine=wncyr9
                                                                            \font\RussianTen=wncyr10
                                                                            \font\RussianBoldfaceFive=wncyb5
                                                                            \font\RussianBoldfaceSix=wncyb6
                                                                            \font\RussianBoldfaceSeven=wncyb7
                                                                            \font\RussianBoldfaceEight=wncyb8
                                                                            \font\RussianBoldfaceNine=wncyb9
                                                                            \font\RussianBoldfaceTen=wncyb10
                                                                            \font\RussianItalicFive=wncyi5
                                                                            \font\RussianItalicSix=wncyi6
                                                                            \font\RussianItalicSeven=wncyi7
                                                                            \font\RussianItalicEight=wncyi8
                                                                            \font\RussianItalicNine=wncyi9
                                                                            \font\RussianItalicTen=wncyi10
                                                                            \font\AMSMathIFive=msam5
                                                                            \font\AMSMathISix=msam6
                                                                            \font\AMSMathISeven=msam7
                                                                            \font\AMSMathIEight=msam8
                                                                            \font\AMSMathINine=msam9
                                                                            \font\AMSMathITen=msam10
                                                                            \font\AMSMathIIFive=msbm5
                                                                            \font\AMSMathIISix=msbm6
                                                                            \font\AMSMathIISeven=msbm7
                                                                            \font\AMSMathIIEight=msbm8
                                                                            \font\AMSMathIINine=msbm9
                                                                            \font\AMSMathIITen=msbm10
                                                                            \font\EulerScriptFive=eusm5
                                                                            \font\EulerScriptSix=eusm6
                                                                            \font\EulerScriptSeven=eusm7
                                                                            \font\EulerScriptEight=eusm8
                                                                            \font\EulerScriptNine=eusm9
                                                                            \font\EulerScriptTen=eusm10
                                                                            \font\BoldSymbolFive=cmbsy5
                                                                            \font\BoldSymbolSix=cmbsy6
                                                                            \font\BoldSymbolSeven=cmbsy7
                                                                            \font\BoldSymbolEight=cmbsy8
                                                                            \font\BoldSymbolNine=cmbsy9
                                                                            \font\BoldSymbolTen=cmbsy10
%
														                                          \newcount\FontNumber%
														                                          \newcount\FontSize%
                                                                                                  \font\rm=cmr10
               \def\SetFont#1#2{\FontNumber=#1\FontSize=#2\Count=0%
                                \ifnum\FontNumber=1 \ifnum\FontSize=5 \Count=1 \RomanFive \fi%
                                                    \ifnum\FontSize=6 \Count=1 \RomanSix \fi%
                                                    \ifnum\FontSize=7 \Count=1 \RomanSeven \fi%
                                                    \ifnum\FontSize=8 \Count=1 \RomanEight \fi%
                                                    \ifnum\FontSize=9 \Count=1 \RomanNine \fi%
                                                    \ifnum\FontSize=10\Count=1\RomanTen\fi%
                                                    \ifnum\FontSize=12 \Count=1 \RomanTwelve \fi%
                                                    \ifnum\FontSize=17 \Count=1 \RomanSeventeen \fi%
                                                    \hbox to 0pt{}\ifnum\Count=0 \font\Temp=cmr10 at \FontSize pt \Temp \fi \fi%
                                \ifnum\FontNumber=2 \ifnum\FontSize=5 \Count=1 \BoldfaceFive \fi%
                                                    \ifnum\FontSize=6 \Count=1 \BoldfaceSix \fi%
                                                    \ifnum\FontSize=7 \Count=1 \BoldfaceSeven \fi%
                                                    \ifnum\FontSize=8 \Count=1 \BoldfaceEight \fi%
                                                    \ifnum\FontSize=9 \Count=1 \BoldfaceNine \fi%
                                                    \ifnum\FontSize=10 \Count=1 \BoldfaceTen \fi%
                                                    \ifnum\FontSize=12 \Count=1 \BoldfaceTwelve \fi%
                                                    \ifnum\FontSize=17 \Count=1 \BoldfaceSeventeen \fi
                                                    \hbox to 0pt{}\ifnum\Count=0 \font\Temp=cmbx10 at \FontSize pt \Temp \fi\fi%
                                \ifnum\FontNumber=3 \ifnum\FontSize=9 \Count=1 \SlantNine \fi%
                                                    \ifnum\FontSize=10 \Count=1 \SlantTen \fi%
                                                    \hbox to 0pt{}\ifnum\Count=0 \font\Temp=cmsl10 at \FontSize pt \Temp \fi \fi%
                                \ifnum\FontNumber=4 \ifnum\FontSize=8 \Count=1 \ItalicEight \fi%
                                                    \ifnum\FontSize=9 \Count=1 \ItalicNine \fi%
                                                    \ifnum\FontSize=10 \Count=1 \ItalicTen \fi%
                                                    \hbox to 0pt{}\ifnum\Count=0 \font\Temp=cmti10 at \FontSize pt \Temp \fi \fi%
                                \ifnum\FontNumber=5 \ifnum\FontSize=8 \Count=1 \TypewriterEight \fi%
                                                    \ifnum\FontSize=9 \Count=1 \TypewriterNine \fi%
                                                    \ifnum\FontSize=10 \Count=1 \TypewriterTen \fi%
                                                    \hbox to 0pt{}\ifnum\Count=0 \font\Temp=cmtt10 at \FontSize pt \Temp \fi\fi%
                                \ifnum\FontNumber=6 \ifnum\FontSize=5 \Count=1 \MathIFive \fi%
                                                    \ifnum\FontSize=6 \Count=1 \MathISix \fi%
                                                    \ifnum\FontSize=7 \Count=1 \MathISeven \fi%
                                                    \ifnum\FontSize=8 \Count=1 \MathIEight \fi%
                                                    \ifnum\FontSize=9 \Count=1 \MathINine \fi%
                                                    \ifnum\FontSize=10 \Count=1 \MathITen \fi%
                                                    \ifnum\FontSize=12 \Count=1 \MathITwelve \fi%
                                                    \hbox to 0pt{}\ifnum\Count=0 \font\Temp=cmmi10 at \FontSize pt \Temp \fi \fi%
                                \ifnum\FontNumber=7 \ifnum\FontSize=6 \Count=1 \MathIISix \fi%
                                                    \ifnum\FontSize=8 \count255y=1 \MathIIEight \fi%
                                                    \ifnum\FontSize=9 \Count=1 \MathIINine \fi%
                                                    \ifnum\FontSize=10 \Count=1 \MathIITen \fi%
                                                    \hbox to 0pt{}\ifnum\Count=0 \font\Temp=cmsy10 at \FontSize pt \Temp \fi \fi%
                                \ifnum\FontNumber=8 \ifnum\FontSize=8 \Count=1 \MathIIIEight \fi%
                                                    \ifnum\FontSize=9 \Count=1 \MathIIINine \fi%
                                                    \ifnum\FontSize=10 \Count=1 \MathIIITen \fi%
                                                    \hbox to 0pt{}\ifnum\Count=0 \font\Temp=cmex10 at \FontSize pt \Temp \fi \fi
                                \ifnum\FontNumber=9 \ifnum\FontSize=5 \count255y=1 \GermanFive \fi%
                                                    \ifnum\FontSize=6 \Count=1 \GermanSix \fi%
                                                    \ifnum\FontSize=7 \Count=1 \GermanSeven \fi%
                                                    \ifnum\FontSize=8 \Count=1 \GermanEight \fi%
                                                    \ifnum\FontSize=9 \Count=1 \GermanNine \fi%
                                                    \ifnum\FontSize=10\Count=1\GermanTen\fi%
                                                    \hbox to 0pt{}\ifnum\Count=0 \font\Temp=eufm10 at \FontSize pt \Temp \fi \fi%
                                \ifnum\FontNumber=10\ifnum\FontSize=5 \Count=1 \GermanBoldfaceFive \fi%
                                                    \ifnum\FontSize=6 \Count=1 \GermanBoldfaceSix \fi%
                                                    \ifnum\FontSize=7 \Count=1 \GermanBoldfaceSeven \fi%
                                                    \ifnum\FontSize=8 \Count=1 \GermanBoldfaceEight \fi%
                                                    \ifnum\FontSize=9 \Count=1 \GermanBoldfaceNine \fi%
                                                    \ifnum\FontSize=10\Count=1\GermanBoldfaceTen\fi%
                                                    \hbox to 0pt{}\ifnum\Count=0 \font\Temp=eufb10 at \FontSize pt \Temp \fi \fi%
                                \ifnum\FontNumber=11\ifnum\FontSize=5 \Count=1 \RussianFive \fi%
                                                    \ifnum\FontSize=6 \Count=1 \RussianSix \fi%
                                                    \ifnum\FontSize=7 \Count=1 \RussianSeven \fi%
                                                    \ifnum\FontSize=8 \Count=1 \RussianEight \fi%
                                                    \ifnum\FontSize=9 \Count=1 \RussianNine \fi%
                                                    \ifnum\FontSize=10\Count=1\RussianTen\fi%
                                                    \hbox to 0pt{}\ifnum\Count=0 \font\Temp=wyncr10 at \FontSize pt \Temp \fi \fi%
                                \ifnum\FontNumber=12\ifnum\FontSize=5 \Count=1 \RussianBoldfaceFive \fi%
                                                    \ifnum\FontSize=6 \Count=1 \RussianBoldfaceSix \fi%
                                                    \ifnum\FontSize=7 \Count=1 \RussianBoldfaceSeven \fi%
                                                    \ifnum\FontSize=8 \Count=1 \RussianBoldfaceEight \fi%
                                                    \ifnum\FontSize=9 \Count=1 \RussianBoldfaceNine \fi%
                                                    \ifnum\FontSize=10\Count=1\RussianBoldfaceTen\fi%
                                                    \hbox to 0pt{}\ifnum\Count=0 \font\Temp=wyncb10 at \FontSize pt \Temp \fi \fi%
                                \ifnum\FontNumber=13\ifnum\FontSize=5 \Count=1 \RussianItalicFive \fi%
                                                    \ifnum\FontSize=6 \Count=1 \RussianItalicSix \fi%
                                                    \ifnum\FontSize=7 \Count=1 \RussianItalicSeven \fi%
                                                    \ifnum\FontSize=8 \Count=1 \RussianItalicEight \fi%
                                                    \ifnum\FontSize=9 \Count=1 \RussianItalicNine \fi%
                                                    \ifnum\FontSize=10\Count=1\RussianItalicTen\fi%
                                                    \hbox to 0pt{}\ifnum\Count=0 \font\Temp=wynci10 at \FontSize pt \Temp \fi \fi%
                                \ifnum\FontNumber=14\ifnum\FontSize=5 \Count=1 \EuropeanFive \fi%
                                                    \ifnum\FontSize=6 \Count=1 \EuropeanSix \fi%
                                                    \ifnum\FontSize=7 \Count=1 \EuropeanSeven \fi%
                                                    \ifnum\FontSize=8 \Count=1 \EuropeanEight \fi%
                                                    \ifnum\FontSize=9 \Count=1 \EuropeanNine \fi%
                                                    \ifnum\FontSize=10\Count=1\EuropeanTen\fi%
                                                    \hbox to 0pt{}\ifnum\Count=0 \font\Temp=eurm10 at \FontSize pt \Temp \fi \fi%
                                \ifnum\FontNumber=15\ifnum\FontSize=5 \Count=1 \EuropeanBoldfaceFive \fi%
                                                    \ifnum\FontSize=6 \Count=1 \EuropeanBoldfaceSix \fi%
                                                    \ifnum\FontSize=7 \Count=1 \EuropeanBoldfaceSeven \fi%
                                                    \ifnum\FontSize=8 \Count=1 \EuropeanBoldfaceEight \fi%
                                                    \ifnum\FontSize=9 \Count=1 \EuropeanBoldfaceNine \fi%
                                                    \ifnum\FontSize=10\Count=1\EuropeanBoldfaceTen\fi%
                                                    \hbox to 0pt{}\ifnum\Count=0 \font\Temp=eurb10 at \FontSize pt \Temp \fi \fi
                                \ifnum\FontNumber=16\ifnum\FontSize=5 \Count=1 \AMSMathIFive \fi%
                                                    \ifnum\FontSize=6 \Count=1 \AMSMathISix \fi%
                                                    \ifnum\FontSize=7 \Count=1 \AMSMathISeven \fi%
                                                    \ifnum\FontSize=8 \Count=1 \AMSMathIEight \fi%
                                                    \ifnum\FontSize=9 \Count=1 \AMSMathINine \fi%
                                                    \ifnum\FontSize=10\Count=1\AMSMathITen\fi%
                                                    \hbox to 0pt{}\ifnum\Count=0 \font\Temp=msam10 at \FontSize pt \Temp \fi \fi
                                \ifnum\FontNumber=17\ifnum\FontSize=5 \Count=1 \AMSMathIIFive \fi%
                                                    \ifnum\FontSize=6 \Count=1 \AMSMathIISix \fi%
                                                    \ifnum\FontSize=7 \Count=1 \AMSMathIISeven \fi%
                                                    \ifnum\FontSize=8 \Count=1 \AMSMathIIEight \fi%
                                                    \ifnum\FontSize=9 \Count=1 \AMSMathIINine \fi%
                                                    \ifnum\FontSize=10\Count=1\AMSMathIITen\fi%
                                                    \hbox to 0pt{}\ifnum\Count=0 \font\Temp=msbm10 at \FontSize pt \Temp \fi \fi
                                \ifnum\FontNumber=18\ifnum\FontSize=5 \Count=1 \EulerScriptFive \fi%
                                                    \ifnum\FontSize=6 \Count=1 \EulerScriptSix \fi%
                                                    \ifnum\FontSize=7 \Count=1 \EulerScriptSeven \fi%
                                                    \ifnum\FontSize=8 \Count=1 \EulerScriptEight \fi%
                                                    \ifnum\FontSize=9 \Count=1 \EulerScriptNine \fi%
                                                    \ifnum\FontSize=10\Count=1\EulerScriptTen\fi%
                                                    \hbox to 0pt{}\ifnum\Count=0 \font\Temp=eusm10 at \FontSize pt \Temp \fi \fi
                                \ifnum\FontNumber=19\ifnum\FontSize=5 \Count=1 \BoldSymbolFive \fi%
                                                    \ifnum\FontSize=6 \Count=1 \BoldSymbolSix \fi%
                                                    \ifnum\FontSize=7 \Count=1 \BoldSymbolSeven \fi%
                                                    \ifnum\FontSize=8 \Count=1 \BoldSymbolEight \fi%
                                                    \ifnum\FontSize=9 \Count=1 \BoldSymbolNine \fi%
                                                    \ifnum\FontSize=10\Count=1\BoldSymbolTen\fi%
                                                    \hbox to 0pt{}\ifnum\Count=0 \font\Temp=cmbsy10 at \FontSize pt \Temp \fi \fi}%
                                                                     \FontNumber=1\FontSize=10%
%
                                                                     %
                                                                     %
                                                                     %
                                                                     %
                                                                     %
                                                                     %
                                                                     %
                                                                     %
                                                                     %
                                                                     %
                                                                     %
                                                                     %
                                                                     %
                                                                     %
                                                                     %
                                                                     %
                                                                     %
                                                                     %
                                                                     %
%
                                                                     %
                                                                     %
                                                                     %
                                                                     %
                                                                     %
                                                                     %
                                                                     %
                                                                     %
                                                                     %
                                                                     %
                                                                     %
                                                                     %
                                                                     %
                                                                     %
                                                                     %
                                                                     %
                                                                     %
                                                                     %
                                                                     %
%
                                                                  \def\SF#1#2#3{\def\OldFontNumber{\count100}\OldFontNumber=\FontNumber%
                                                                                \def\OldFontSize{\count99}\OldFontSize=\FontSize%
                                                                                \SetFont{#1}{#2}%
                                                                                #3%
                                                                                \SetFont{\number\OldFontNumber}{\number\OldFontSize}}%
%
                   \def\Rm#1{\SF{1}{\number\FontSize}{#1}}@@@\def\RmBx#1{\hbox{\Rm{#1}}}@@@@@@@@@@\def\RmChBx#1{\hbox{\Rm{\char'#1}}}@@@@%
                   \def\Bf#1{\SF{2}{\number\FontSize}{#1}}@@@\def\BfBx#1{\hbox{\Bf{#1}}}@@@@@@@@@@@@@@%
                   @@@@@@@@@@@@@@@@@%
                   \def\It#1{\SF{4}{\number\FontSize}{#1}}@@@\def\ItBx#1{\hbox{\It{#1}}}@@@@@@@@@@\def\ItChBx#1{\hbox{\It{\char'#1}}}@@@@%
                   \def\Tw#1{\SF{5}{\number\FontSize}{#1}}@@@@@@@@@@@@@\def\TwChBx#1{\hbox{\Tw{\char'#1}}}@@@@%
                   \def\Mi#1{\SF{6}{\number\FontSize}{#1}}@@@@@@@@\def\MiCh#1{\Mi{\char'#1}}@@@@@\def\MiChBx#1{\hbox{\Mi{\char'#1}}}@@@@%
                   \def\Mii#1{\SF{7}{\number\FontSize}{#1}}@@\def\MiiBx#1{\hbox{\Mii{#1}}}@@@@@@\def\MiiChBx#1{\hbox{\Mii{\char'#1}}}@@%
                   \def\Miii#1{\SF{8}{\number\FontSize}{#1}}@@@\def\MiiiChBx#1{\hbox{\Miii{\char'#1}}}%
                   \def\Ft#1{\SF{9}{\number\FontSize}{#1}}@@@@@@@@@@@@@@@@@%
                   @@@@@@@@@@@@@@@@%
                   \def\Ru#1{\SF{11}{\number\FontSize}{#1}}@@\def\RuBx#1{\hbox{\Ru{#1}}}@@@@@@@@@@\def\RuChBx#1{\hbox{\Ru{\char'#1}}}@@@@%
                   \def\Rb#1{\SF{12}{\number\FontSize}{#1}}@@@@@@@@@@@@@@@@%
                   @@@@@@@@@@@@@@@@%
                   \def\Eu#1{\SF{14}{\number\FontSize}{#1}}@@\def\EuBx#1{\hbox{\Eu{#1}}}@@@@@@@@@@@@@@%
                   \def\Eb#1{\SF{15}{\number\FontSize}{#1}}@@@@@@@@@@@@@@@@%
                   \def\Ai#1{\SF{16}{\number\FontSize}{#1}}@@@@@@@@@@@@\def\AiChBx#1{\hbox{\Ai{\char'#1}}}@@@@%
                   \def\Aii#1{\SF{17}{\number\FontSize}{#1}}@\def\AiiBx#1{\hbox{\Aii{#1}}}@@@\def\AiiCh#1{\Aii{\char'#1}}@@@@@%
                   \def\Es#1{\SF{18}{\number\FontSize}{#1}}@@@@@@@@@@@@@@@@%
                   @@@@@@@@@@@@@@@@%
                   \def\RM#1#2{\SF{1}{#1}{#2}}@@@\def\RMBx#1#2{\hbox{\RM{#1}{#2}}}@@@@@@@@\def\RMChBx#1#2{\SF{1}{#1}{\hbox{\char'#2}}}@@@@%
                   \def\BF#1#2{\SF{2}{#1}{#2}}@@@\def\BFBx#1#2{\hbox{\BF{#1}{#2}}}@@@@@@@@@@@@%
                   @@@@@@@@@@@@@@@%
                   @@@@@@@@@@@@@@@%
                   \def\TW#1#2{\SF{5}{#1}{#2}}@@@\def\TWBx#1#2{\hbox{\TW{#1}{#2}}}@@@@@\def\TWCh#1#2{\SF{5}{#1}{\char'#2}}@@@\def\TWChBx#1#2{\SF{5}{#1}{\hbox{\char'#2}}}@@@@%
                   @@@@@@@@@@@\def\MIChBx#1#2{\SF{6}{#1}{\hbox{\char'#2}}}@@@@%
                   \def\MII#1#2{\SF{7}{#1}{#2}}@@@@@\def\MIICh#1#2{\SF{7}{#1}{\char'#2}}@@\def\MIIChBx#1#2{\SF{7}{#1}{\hbox{\char'#2}}}@@@%
                   @@\def\MIIICh#1#2{\SF{8}{#1}{\char'#2}}@\def\MIIIChBx#1#2{\SF{8}{#1}{\hbox{\char'#2}}}@@%
                   \def\FT#1#2{\SF{9}{#1}{#2}}@@@\def\FTBx#1#2{\hbox{\FT{#1}{#2}}}@@@@@@@@\def\FTChBx#1#2{\SF{9}{#1}{\hbox{\char'#2}}}@@@@%
                   @@@@@@@@@@@@%
                   @@@@@@@@@@@@%
                   @@@@@@@@@@@@%
                   @@@@@@@@@@@@%
                   \def\EU#1#2{\SF{14}{#1}{#2}}@@\def\EUBx#1#2{\hbox{\EU{#1}{#2}}}@@@@@@@\def\EUChBx#1#2{\SF{14}{#1}{\hbox{\char'#2}}}@@@%
                   \def\EB#1#2{\SF{15}{#1}{#2}}@@@@@@@@@@@@%
                   @@@@@@@\def\AICh#1#2{\SF{16}{#1}{\char'#2}}@@\def\AIChBx#1#2{\SF{16}{#1}{\hbox{\char'#2}}}@@@%
                   \def\AII#1#2{\SF{17}{#1}{#2}}@\def\AIIBx#1#2{\hbox{\AII{#1}{#2}}}@@@\def\AIICh#1#2{\SF{17}{#1}{\char'#2}}@\def\AIIChBx#1#2{\SF{17}{#1}{\hbox{\char'#2}}}@@%
                   \def\ES#1#2{\SF{18}{#1}{#2}}@@\def\ESBx#1#2{\hbox{\ES{#1}{#2}}}@@@@@@@@@@%
                   @@@@@@@@@\def\BSChBx#1#2{\SF{19}{#1}{\hbox{\char'#2}}}@@@%
    \def\GR#1#2{\def\Aachen{#2}\def\Bachen{ps}\ifx\Aachen\Bachen\MIChBx{#1}{40}\else\def\Bachen{xi}\ifx\Aachen\Bachen\MIChBx{#1}{30}\else\def\Bachen{th}\ifx\Aachen\Bachen\MIChBx{#1}{22}\else%
                \if#2a\MIChBx{#1}{13}\else\if#2b\MIChBx{#1}{14}\else\if#2c\MIChBx{#1}{15}\else\if#2d\MIChBx{#1}{16}\else\if#2e\MIChBx{#1}{17}\else\if#2z\MIChBx{#1}{20}\else\if#2h\MIChBx{#1}{21}%
                \else\if#2i\MIChBx{#1}{23}\else\if#2k\MIChBx{#1}{24}\else\if#2l\MIChBx{#1}{25}\else\if#2m\MIChBx{#1}{26}\else\if#2n\MIChBx{#1}{27}\else\if#2o\RMBx{#1}{o}\else\if#2p\MIChBx{#1}{31}%
                \else\if#2r\MIChBx{#1}{32}\else\if#2s\MIChBx{#1}{33}\else\if#2t\MIChBx{#1}{34}\else\if#2u\MIChBx{#1}{35}\else\if#2f\MIChBx{#1}{36}\else\if#2x\MIChBx{#1}{37}\else\if#2w\MIChBx{#1}{41}%
                \fi\fi\fi\fi\fi\fi\fi\fi\fi\fi\fi\fi\fi\fi\fi\fi\fi\fi\fi\fi\fi\fi\fi\fi
                \def\Bachen{PS}\ifx\Aachen\Bachen\MIChBx{#1}{11}\else\def\Bachen{XI}\ifx\Aachen\Bachen\MIChBx{#1}{4}\else\def\Bachen{TH}\ifx\Aachen\Bachen\MIChBx{#1}{2}\else%
                \if#2A\MIChBx{#1}{101}\else\if#2B\MIChBx{#1}{102}\else\if#2C\MIChBx{#1}{0}\else\if#2D\MIChBx{#1}{1}\else\if#2E\MIChBx{#1}{105}\else\if#2Z\MIChBx{#1}{132}\else\if#2H\MIChBx{#1}{110}%
                \else\if#2I\MIChBx{#1}{111}\else\if#2K\MIChBx{#1}{113}\else\if#2L\MIChBx{#1}{3}\else\if#2M\MIChBx{#1}{115}\else\if#2N\MIChBx{#1}{116}\else\if#2O{117}\else\if#2P\MIChBx{#1}{5}%
                \else\if#2R\MIChBx{#1}{122}\else\if#2S\MIChBx{#1}{6}\else\if#2T\MIChBx{#1}{124}\else\if#2U\MIChBx{#1}{7}\else\if#2F\MIChBx{#1}{10}\else\if#2X\MIChBx{#1}{130}\else\if#2W\MIChBx{#1}{12}%
                \fi\fi\fi\fi\fi\fi\fi\fi\fi\fi\fi\fi\fi\fi\fi\fi\fi\fi\fi\fi\fi\fi\fi\fi}%
    \def\Gr#1{\GR{\FontSize}{#1}}%
    \def\gR#1{\GR{7}{#1}}%
\def\MMiiBx#1{\def\A{#1}\def\B{0}%
              \if A\A \hbox{\kn {.1}\MiiBx{A}\kn{.3}}\def\B{1}\fi
              \if B\A \hbox{\MiiBx{B}\kn{.5}}\def\B{1}\fi
              \if C\A \hbox{\kn{.2}\MiiBx{C}\kn{.5}}\def\B{1}\fi
              \if D\A \hbox{\kn{.1}\MiiBx{D}\kn{.5}}\def\B{1}\fi
              \if E\A \hbox{\kn{.1}\MiiBx{E}\kn{.7}}\def\B{1}\fi
              \if F\A \hbox{\kn{.1}\MiiBx{F}\kn{1.4}}\def\B{1}\fi
              \if G\A \hbox{\kn{0}\MiiBx{G}\kn{.3}}\def\B{1}\fi
              \if I\A \hbox{\kn{.5}\MiiBx{I}\kn{1.2}}\def\B{1}\fi
              \if J\A \hbox{\kn{0}\MiiBx{J}\kn{2.1}}\def\B{1}\fi
              \if K\A \hbox{\kn{.1}\MiiBx{K}\kn{.1}}\def\B{1}\fi
              \if L\A \hbox{\kn{.1}\MiiBx{L}\kn{.1}}\def\B{1}\fi
              \if M\A \hbox{\kn{.1}\MiiBx{M}\kn{0}}\def\B{1}\fi
              \if N\A \hbox{\kn{.5}\MiiBx{N}\kn{2.1}}\def\B{1}\fi
              \if O\A \hbox{\kn{0}\MiiBx{O}\kn{.2}}\def\B{1}\fi
              \if P\A \hbox{\kn{.2}\MiiBx{P}\kn{.8}}\def\B{1}\fi
              \if Q\A \hbox{\kn{-.4}\MiiBx{Q}\kn{0}}\def\B{1}\fi
              \if R\A \hbox{\kn{.3}\MiiBx{R}\kn{.5}}\def\B{1}\fi
              \if S\A \hbox{\kn{.3}\MiiBx{S}\kn{1}}\def\B{1}\fi
              \if T\A \hbox{\kn{.1}\MiiBx{T}\kn{3.5}}\def\B{1}\fi
              \if U\A \hbox{\kn{.6}\MiiBx{U}\kn{1.1}}\def\B{1}\fi
              \if V\A \hbox{\kn{0}\MiiBx{V}\kn{1}}\def\B{1}\fi
              \if W\A \hbox{\kn{0}\MiiBx{W}\kn{.8}}\def\B{1}\fi
              \if X\A \hbox{\kn{.1}\MiiBx{X}\kn{1.4}}\def\B{1}\fi
              \if Y\A \hbox{\kn{.1}\MiiBx{Y}\kn{1.3}}\def\B{1}\fi
              \if Z\A \hbox{\kn{.1}\MiiBx{Z}\kn{1}}\def\B{1}\fi
              \if 1\B \else \MiiBx{#1}\fi}%
\def\IItBx#1{\def\A{#1}\def\B{0}%
              \if A\A \hbox{\kn{.1}\ItBx{A}\kn{.3}}\def\B{1}\fi
              \if B\A \hbox{\ItBx{B}\kn{.8}}\def\B{1}\fi
              \if C\A \hbox{\kn{-.3}\ItBx{C}\kn{1.3}}\def\B{1}\fi
              \if D\A \hbox{\kn{.1}\ItBx{D}\kn{.8}}\def\B{1}\fi
              \if E\A \hbox{\kn{.1}\ItBx{E}\kn{1}}\def\B{1}\fi
              \if F\A \hbox{\kn{.1}\ItBx{F}\kn{1.4}}\def\B{1}\fi
              \if G\A \hbox{\kn{-.3}\ItBx{G}\kn{.9}}\def\B{1}\fi
              \if H\A \hbox{\kn{-.1}\ItBx{H}\kn{1.5}}\def\B{1}\fi
              \if I\A \hbox{\kn{.2}\ItBx{I}\kn{1.6}}\def\B{1}\fi
              \if J\A \hbox{\kn{0}\ItBx{J}\kn{1.9}}\def\B{1}\fi
              \if K\A \hbox{\kn{.1}\ItBx{K}\kn{1.4}}\def\B{1}\fi
              \if L\A \hbox{\kn{.1}\ItBx{L}\kn{.7}}\def\B{1}\fi
              \if M\A \hbox{\kn{.1}\ItBx{M}\kn{1.6}}\def\B{1}\fi
              \if N\A \hbox{\kn{.5}\ItBx{N}\kn{2.1}}\def\B{1}\fi
              \if O\A \hbox{\kn{-.2}\ItBx{O}\kn{1}}\def\B{1}\fi
              \if P\A \hbox{\kn{.2}\ItBx{P}\kn{1.1}}\def\B{1}\fi
              \if Q\A \hbox{\kn{-.4}\ItBx{Q}\kn{1}}\def\B{1}\fi
              \if R\A \hbox{\kn{.3}\ItBx{R}\kn{.7}}\def\B{1}\fi
              \if S\A \hbox{\kn{.3}\ItBx{S}\kn{1.3}}\def\B{1}\fi
              \if T\A \hbox{\kn{-.1}\ItBx{T}\kn{2.5}}\def\B{1}\fi
              \if U\A \hbox{\kn{-.3}\ItBx{U}\kn{2}}\def\B{1}\fi
              \if V\A \hbox{\kn{-.8}\ItBx{V}\kn{1.9}}\def\B{1}\fi
              \if W\A \hbox{\kn{-.4}\ItBx{W}\kn{1.8}}\def\B{1}\fi
              \if X\A \hbox{\kn{.1}\ItBx{X}\kn{1.6}}\def\B{1}\fi
              \if Y\A \hbox{\kn{-.5}\ItBx{Y}\kn{2}}\def\B{1}\fi
              \if Z\A \hbox{\kn{.1}\ItBx{Z}\kn{1.7}}\def\B{1}\fi
              \if 1\B \else \ItBx{#1}\fi}%
            \def\rm#1{$\RmBx{#1}$}%
            \def\bf#1{$\BfBx{#1}$}%
            \def\it#1{$\ItBx{#1}$}%
            %
            %
            %
            %
            %
            %
            %
            %
            %
            %
            %
            %
            %
            %
            %
                                                            %
                                                            \def\vruleHWD#1#2#3{\vrule height#1 depth#3 width#2}%
                                                            \def\vruleWHD#1#2#3{\vrule height#2 depth#3 width#1}%
                                                            \def\vruleHW#1#2{\vrule height#1 width#2 depth0pt}%
                                                            \def\vruleWH#1#2{\vrule height#2 width#1 depth0pt}%
                                                            \def\hruleHW#1#2{\hrule height#1 width#2 depth0pt}%
                                                            \def\hruleWH#1#2{\hrule height#2 width#1 depth0pt}%
                                                            \def\vrulehdw#1#2#3{\vrule height#1pt depth#2pt width#3pt}%
                                                            \def\vrulehwd#1#2#3{\vrule height#1pt depth#3pt width#2pt}%
                                                            \def\vrulewhd#1#2#3{\vrule height#2pt depth#3pt width#1pt}%
                                                            \def\vrulehw#1#2{\vrule height#1pt width#2pt depth0pt}%
                                                            \def\vrulewh#1#2{\vrule height#2pt width#1pt depth0pt}%
                                                            \def\vrulewd#1#2{\vrule depth#2pt width#1pt height0pt}%
                                                            \def\vruleh#1{\vrule height#1pt depth0pt width0pt}%
                                                            \def\vruled#1{\vrule depth#1pt height0pt width0pt}%
                                                            \def\vrulew#1{\vrule width#1pt height0pt depth0pt}%
                                                            \def\hrulehw#1#2{\hrule height#1pt width#2pt depth0pt}%
                                                            \def\hrulewh#1#2{\hrule height#2pt width#1pt depth0pt}%
                                                            \def\hruleh#1{\hrule height#1pt depth0pt width0pt}%
                                                                                                  \def\Mark#1{\hldest{xyz}{}{#1}}%
                                                                                              \def\MarkWithText#1#2{\definexref{#1}{#2}{}}%
                                                                                            \def\Link#1#2{\hlstart{name}{}{#1}#2\hlend}%
                                                                                    %
                                                                                           %
                                                                                   %
                                                                                                    \def\LinkText#1{\refn{#1}}%
                                                                                                     \enablehyperlinks[pdftex]%
                                                                                                         \hlopts{bwidth=0}%
                                                                                             \def\TheCount#1{\global\advance#1 by 1%
                                                                                                             \count255=#1%
                                                                                                             \the\count255}%
                                                                                             \newcount\SectionCount \SectionCount=0%
                                                                                             \newcount\ItemCount    \ItemCount=0%
                                                                                             \newcount\EqCount      \EqCount=0%
                                                            \def\SectionLabel#1{XsectionX#1}
                                                                            \def\Section#1#2{\PAr\noindent\hfil\BF{14}{\TheCount{\SectionCount}. #2}%
                                                                                             \MarkWithText{\SectionLabel{#1}}{\the\SectionCount}
                                                                                             \gdef\SectionPreLabel{#1}
                                                                                             \ItemCount=0\vskip4pt}%
                                                            \def\LinkSection#1#2{\Link{\SectionLabel{#1}}{#2}}
                                                            \def\ItemLabel#1#2{XsectionX#1XitemX#2}
                                                                \def\Item#1#2{\Paragraph{1\baselineskip}{30pt}\Bf{(\the\SectionCount.\TheCount{\ItemCount}) #2\ \ }%
                                                                              \gdef\ItemPreLabel{#1}%
                                                                              \MarkWithText{\ItemLabel{\SectionPreLabel}{#1}}{\the\ItemCount}%
                                                                              \EqCount=0\ }%
                                                            \def\LinkItem#1#2#3{\Link{\ItemLabel{#1}{#2}}{#3}}
                                                            \def\LinkLocalItemText#1{\LinkText{\ItemLabel{\SectionPreLabel}{#1}}}
                                                            \def\pLinkLocalItemText#1{(\LinkLocalItemText{#1})}
                                                            \def\LinkItemText#1#2{\LinkText{\SectionLabel{#1}}.\LinkText{\ItemLabel{#1}{#2}}}
                                                            \def\pLinkItemText#1#2{(\LinkItemText{#1}{#2})}
                                                            \def\hLinkItemText#1#2{\hbox{\vrulew{2}\RM{7}{by (\LinkItemText{#1}{#2})}\vrulew{2}}}
                                                            \def\DisplayLabel#1#2#3{XsectionX#1XitemX#2XdisplayX#3}%
                                                            %
                                                            %
                                                            \def\LinkLocalDisplayText#1{\LinkText{\DisplayLabel{\SectionPreLabel}{\ItemPreLabel}{#1}}}%
                                                            \def\pLinkLocalDisplayText#1{(\LinkLocalDisplayText{#1})}%
                                                            \def\hLinkLocalDisplayText#1{\hbox{\vrulew{2}\RM{7}{by (\LinkLocalDisplayText{#1})}\vrulew{2}}}%
                                                            \def\LinkDisplayText#1#2#3{\LinkText{\SectionLabel{#1}}.\LinkText{\ItemLabel{#1}{#2}}.\LinkText{\DisplayLabel{#1}{#2}{#3}}}%
                                                            \def\pLinkDisplayText#1#2#3{(\LinkDisplayText{#1}{#2}{#3})}%
                                                            \def\hLinkDisplayText#1#2#3{\hbox{\vrulew{2}\RM{7}{by (\LinkDisplayText{#1}{#2}{#3})}\vrulew{2}}}%
                                                            %
                                                            %
                                                            \def\d#1{\endgraf\line{\hfil#1\hfil}}%
                                                            \def\D#1#2#3{\par\Vbox{\kn{#1}%
                                                                                   \line{\hfil#3\hfil}%
                                                                                   \kn{#2}}}%
                                                            \def\Dc#1{\D{8}{8}{#1}\vskip-5pt}%
													        \def\DL#1#2#3#4{\par\setbox0\hbox{#4}\setbox1\hbox to\hsize{\hfil#3\hfil}\setbox2\hbox{#4\kern-\wd0\box1}\Vbox{\kn{#1}\box2\kn{#2}}}
													        \def\DLc#1#2{\par\setbox0\hbox{#2}\setbox1\hbox to\hsize{\hfil#1\hfil}\setbox2\hbox{#2\kern-\wd0\box1}\Vbox{\kn{8}\box2\kn{8}}\vskip-5pt}
                                                            \def\DI#1#2#3#4{\par\Vbox{\kn{#1}
                                                                                      \line{\hfil%
                                                                                            #3%
                                                                                            \hfil%
                                                                                            (\TheCount{\EqCount})%
                                                                                            \MarkWithText{\DisplayLabel{\SectionPreLabel}{\ItemPreLabel}{#4}}%
                                                                                            {\the\EqCount}}%
                                                                                      \kn{#2}}}%
                                                            \def\DIc#1#2{\DI{10}{10}{#1}{#2}\vskip-5pt}
                                                            \def\DII#1#2#3#4#5{\par\Vbox{\kn{#1}
                                                                                         \line{#4%
                                                                                              \hfil%
                                                                                              #3%
                                                                                              \hfil%
                                                                                              (\TheCount{\EqCount})%
                                                                                              \MarkWithText{\DisplayLabel{\SectionPreLabel}{\ItemPreLabel}{#5}}{\the\EqCount}}%
                                                                                         \kn{#2}}%
                                                                               \Par}%
                                                                                                        \vbadness=10000%
                                                                            \newcount\FigCount \FigCount=0%
                                                                            \def\FigureLabel#1{XfigureX#1}%
                                                                            %
                                                                            %
                                                                            %
                                                                            \def\Figure#1#2#3#4#5#6{\ifx\undefined\MoreFigure \par\fi%
							                                                                        \vskip 0pt plus1in%
                                                                                                    \global\advance\FigCount by 1%
                                                                                                    \MarkWithText{\FigureLabel{#6}}{\the\FigCount}%
                                                                              						\setbox\Box=\hbox{\includegraphics[width=#3]{#2.\grf}}%
                                                                              						\hbox{\Vbox{\kern#1%
                                                                                                   	            \line{\hfill\box\Box\hfill}%
                                                                                                   	            \kern#4%
                                                                                                   	            \line{\hfill Fig. \the\FigCount: #5\hfill}
												     											                \kn{5}}}}%
                                    										\def\Figures#1#2#3#4#5#6#7#8#9{\par%
												 										                   \global\advance\FigCount by 1%
                                                                                                           \Link{\FigureLabel{#9}}{\the\FigCount}%
                                                                   						                   \setbox\Box=\hbox{\includegraphics[width=#3]{#2}}%
                                                                   										   \setbox1=\hbox{\includegraphics[width=#6]{#5}}%
                                                                   										   \hbox{\Vbox{\kern#1%
                                                                                   						   \line{\hfill%
													        										             \box\Box%
																							                     \kern#4%
																							                     \box1%
																							                     \hfill}%
                                                                                   						   \kern#7%
                                                                                   						   \line{\hfill Fig. \the\FigCount: #8\hfill}
													       									               \kn{5}}}}%
                                                            \def\Black#1{\color{black}#1}%
                                                            \def\Blue#1{\color{blue}#1\color{black}}%
                                                            \def\Green#1{\color{green}#1\color{black}}%
                                                                \definecolor{gruen}{rgb}{.0424,.53,.1074}%
                                                            \def\Gruen#1{{\color{gruen}#1}\color{black}}%
                                                            \def\Red#1{\color{red}#1\color{black}}%
                                                                \definecolor{purple}{rgb}{.4,0,.4}%
                                                            \def\Purple#1{\color{purple}#1\color{black}}%
                                                                \definecolor{violet}{rgb}{.85,0,1}%
                                                            \def\Violet#1{\color{violet}#1\color{black}}%
                                                                \definecolor{redorange}{rgb}{1,.283,0}%
                                                            \def\RedOrange#1{\color{redorange}#1\color{black}}%
                                                                \definecolor{brown}{rgb}{1,.5,0}%
                                                            \def\Brown#1{\color{brown}#1\color{black}}%
                                                                \definecolor{gray}{gray}{0.7}%
                                                            \def\Gray#1{\color{gray}#1\color{black}}%
                                                                \definecolor{pink}{rgb}{1,.72,.76}%
                                                            %
                                                           	%
															%
															%
                                                           	%
															%
															%
                                                           	%
															%
															%
                                                           	%
															%
															%
                                                           	%
															%
															%
                                                           	%
															%
															%
                                                            %
                                                           	\def\LPAREN#1{\RMChBx{#1}{50}}%
                                                           	\def\RPAREN#1{\RMChBx{#1}{51}}%
															%
															%
                                                            \def\BRAK#1{\TWChBx{#1}{74}}%
                                                            \def\Brak{\hbox{\lower1pt\hbox{\BRAK{12}}}}%
                                                            \def\KET#1{\TWChBx{#1}{76}}%
                                                            \def\Ket{\hbox{\lower1pt\hbox{\KET{12}}}}%
                                                            %
                                                            %
                                                            \def\BBRAK#1{\MIIChBx{#1}{34}}%
                                                            \def\BBrak{\hbox{\lower1pt\hbox{\BRAK{12}}}}%
                                                            \def\KKET#1{\MIIChBx{#1}{35}}%
                                                            \def\KKet{\hbox{\lower1pt\hbox{\KET{12}}}}%
                                                            %
                                                            %
                                                 						\def\OVERBAR#1#2#3{\Vbox{\hrule height#1pt%
						                   								                         \kn{#2}%
														   		                                 \hbox{#3}}}%
                                                 						\def\UNDERBAR#1#2{\setbox\BoxArg\hbox{#2}%
                                                                                          \setbox\Box\hbox{\Vbox{\copy\BoxArg%
						 											                                             \kn{1}%
																	                                             \hrule height#1pt}}%
                                                                   					      \addIII{1pt}{#1pt}{\dp\BoxArg}%
															                              \hbox{\lower\DimenReturn\box\Box}}%
                                                                        \def\Underbar#1{\UNDERBAR{.5}{#1}}%
                                                            \def\GrowingLBrace#1{\RaiseBoxReturn{#1}%
                                                                                 \DimenI=\ht\BoxReturn%
                                                                                 \advance\DimenI by 2pt%
                                                                                 \setbox\BoxI\box\BoxReturn%
                                                                                 \setbox\BoxArg\hbox{\Vbox{\MiiiChBx{70}%
                                                                                                           \MiiiChBx{74}%
                                                                                                           \MiiiChBx{72}}}%
                                                                                 \RaiseBoxReturn{\box\BoxArg}%
                                                                                 \setbox\BoxII\box\BoxReturn%
                                                                                 \DimenII=\ht\BoxII%
                                                                                 \ifdim\DimenII>\DimenI \Dimen=\DimenII%
                                                                                                        \advance\Dimen by-\DimenI%
                                                                                                        \divide\Dimen by2%
                                                                                                        \advance\Dimen by 1pt%
                                                                                                        \setbox\BoxVar\hbox{\box\BoxII\kn{5}\Vbox{\kern\Dimen%
                                                                                                                                                   \box\BoxI%
                                                                                                                                                   \kern\Dimen}}%
                                                                                                         \hbox{\lower14.5pt\box\BoxVar}%
                                                                                                   \else \Dimen=\DimenI%
                                                                                                         \advance\Dimen by-\DimenII%
                                                                                                         \divide\Dimen by2%
                                                                                                         \setbox\BoxII\hbox{\Vbox{\MiiiChBx{70}%
                                                                                                                                  \hbox{\kn{4}\vrule height\Dimen width1.2pt}%
                                                                                                                                  \MiiiChBx{74}%
                                                                                                                                  \hbox{\kn{4}\vrule height\Dimen width1.2pt}%
                                                                                                                                  \MiiiChBx{72}}}%
                                                                                                         \RaiseBoxReturn{\box\BoxII}%
                                                                                                         \setbox\BoxII\box\BoxReturn%
                                                                                                         \setbox\BoxVar\hbox{\lower1pt\box\BoxII\kn{5}\box\BoxI}%
                                                                                                         \Dimen=\ht\BoxVar%
                                                                                                         \divide\Dimen by2%
                                                                                                         \advance\Dimen by-4pt%
                                                                                                         \hbox{\lower\Dimen\box\BoxVar}%
                                                                                                   \fi}%
                                                                                        %
                                                                                        %
                                                                                        \def\ArrowDn#1{\varrowlength=#1pt\hbox{$\mapdown$}}%
                                                                                        \def\ArrowRt#1{\harrowlength=#1pt\hbox{$\mapright$}}%
                                                                                        \def\ArrowLt#1{\harrowlength=#1pt\hbox{$\mapleft$}}%
                                                                                                          \catcode`\@=12%
 \def\druleTWxy#1#2#3#4{\let\SlopeHeight=\CountVar \SlopeHeight=#4%
                        \let\SlopeWidth=\CountArg \SlopeWidth=#3%
                        \ifnum\SlopeWidth<0 \multiply\SlopeHeight by-1 \multiply\SlopeWidth by-1\fi%
                        \let\Width=\DimenVar \Width=#2%
                        \let\Height=\DimenEnd \Height=\Width \divide\Height by\SlopeWidth \multiply\Height by\SlopeHeight%
                        \let\Scale=\DimenO \Scale=.1pt%
                        \let\ScaleBox=\BoxArg \setbox\ScaleBox\hbox{\vruleHW{\Scale}{\Scale}}%
                        \let\Sum=\DimenI%
                        \let\Inc=\DimenII \Inc=\Scale%
                        \let\Line=\BoxVar%
                        \let\Sliver=\BoxI%
                        \let\LineWidth=\DimenIII \LineWidth=#1%
                        \CountEnd=\Scale%
                        \let\LineTimes=\CountI \LineTimes=\LineWidth \divide\LineTimes by\CountEnd%
                        \let\Counter=\CountO \Counter=0%
                        \def\IFi{\ifdim\Sum<\Width \advance\Sum by\Inc%
                                                   \Vbox{\copy\ScaleBox\kern\Sum}%
                                                   \IFi\fi}%
                        \def\IFii{\ifdim\Sum>0pt \hbox{\kern\Sum\copy\ScaleBox}%
                                  \advance\Sum by-\Inc%
                                  \IFii\fi}%
                        \def\IFiii{\ifdim\Sum>0pt \advance\Sum by-\Inc%
                                                  \Vbox{\copy\ScaleBox\kern\Sum}%
                                                  \IFiii\fi}%
                        \def\IFiv{\ifdim\Sum<\Width \hbox{\kern\Sum\copy\ScaleBox}%
                                                    \advance\Sum by\Inc%
                                                    \IFiv\fi}%
                        \def\TIMES{\ifnum\Counter<\LineTimes \advance\Counter by1%
                                                             \hbox{\kern\Scale\kern-\wd\Sliver\copy\Sliver}%
                                                             \TIMES\fi}%
                        \ifnum\SlopeHeight>-1 \ifnum\SlopeHeight<\SlopeWidth \advance\Width by-\Scale%
                                                                             \multiply\Inc by\SlopeHeight \divide\Inc by\SlopeWidth%
                                                                             \Sum=0pt%
                                                                             \setbox\Sliver\hbox{\IFi}%
                                                                             \setbox\Line\hbox{\kern\wd\Sliver\TIMES}%
                                                                       \else \advance\Height by-\Scale%
                                                                             \multiply\Inc by\SlopeWidth \divide\Inc by\SlopeHeight%
                                                                             \Sum=\Width%
                                                                             \setbox\Sliver\hbox{\Vbox{\IFii}}%
                                                                             \setbox\Line\hbox{\kern\wd\Sliver\TIMES}\fi%
                                        \else \ifnum-\SlopeHeight<\SlopeWidth \advance\Width by-\Scale%
                                                                              \multiply\Inc by-\SlopeHeight \divide\Inc by\SlopeWidth%
                                                                              \Sum=\Width%
                                                                              \setbox\Sliver\hbox{\IFiii}%
                                                                              \setbox\Line\hbox{\kern\wd\Sliver\TIMES}%
                                                                        \else \advance\Height by-\Scale%
                                                                              \multiply\Inc by-\SlopeWidth \divide\Inc by\SlopeHeight%
                                                                              \Sum=0pt%
                                                                              \setbox\Sliver\hbox{\Vbox{\IFiv}}%
                                                                              \setbox\Line\hbox{\kern\wd\Sliver\TIMES}\fi\fi%
                        \box\Line}%
 %
 %
 %
 %
 %
 %
 %
                                                            \def\ForAll#1{(\MiiChBx{70}\kn{2}#1)}%
                                                            \def\ForAllSuch#1#2{\ForAll{#1\kn{1}:\kn{4}#2}}%
                                                            \def\ThereIs#1{(\MiiChBx{71}\kn{2}#1)}%
                                                            \def\ThereIsShriek#1{(\MiiChBx{71}\kn{1}\RmBx{!}\kn{2}#1)}%
															\def\Implies{\kn{3}\hbox{\leaders\hbox to 4pt{\hss\Rm{=}\hss}\hskip12pt\kn{-2}\MiiChBx{51}}}%
                                                            \def\Imp#1{\setbox\Box=\vbox{\hrulewh{5}{.27}\kn{1.085}\hrulewh{5}{.27}\kn{.94}}%
                                                                        \setbox\BoxArg=\hbox{\kn{1}\RM{7}{#1}}%
                                                                        \setbox\BoxI\hbox{\lower.5pt\hbox{\TWChBx{7}{76}}}%
                                                                        \addII{\wd\BoxArg}{0pt}%
									 								    \hbox{\lower2pt\hbox{\kn{0}\Vbox{\hbox{\kn{-.5}\box\BoxArg}%
                                                                                                         \kn{1}%
                                                                                                         \hbox{\kn{1.5}\hboxI{\DimenReturn}{\leaders\hbox to 3pt{\hss\box\Box\hss}\hfil}\kn{-2}\box\BoxI}\kn{2}}}}}%
                                                 			\def\Implied{\kn{2}\hbox{\MiiChBx{50}\kn{-4}\leaders\hbox to 5pt{\hss\Rm{=}\hss}\hskip16pt}\kn{2}}%
                                                            \def\Impd#1{\setbox1=\vbox{\hrulewh{5}{.27}\kn{1.085}\hrulewh{5}{.27}\kn{.94}}%
                                                                                       \hbox{\hskip3pt\Vbox{\kn{1}%
                                                                                   					        \setbox\Box=\hbox{\RM{7}{#1}}%
                                                                                   							\addII{\Dimen}{\wd\Box}{7pt}%
                                                                                  					        \hbox{\kn{3}\box\Box}%
                                                                                   					        \kn{1}%
                                                                                   					        \hboxI{\Dimen}{\MIIChBx{7}{50}\kn{-4}\leaders\hbox to 5pt{\hss\box1\hss}\hfil}}%
														                                     \hskip-1pt}}%
                                                            \def\Ifff#1{\setbox\Box=\vbox{\hrulewh{5}{.27}\kn{1.085}\hrulewh{5}{.27}\kn{.94}}%
                                                                        \setbox\BoxArg=\hbox{\kn{1}\RM{7}{#1}}%
                                                                        \setbox\BoxVar\hbox{\lower.5pt\hbox{\TWChBx{7}{74}}}%
                                                                        \setbox\BoxI\hbox{\lower.5pt\hbox{\TWChBx{7}{76}}}%
                                                                        \addII{\wd\BoxArg}{4pt}%
									 								    \hbox{\lower2pt\hbox{\kn{2}\Vbox{\hbox{\kn{3}\box\BoxArg}%
                                                                                          \kn{1}%
                                                                                          \hbox{\box\BoxVar\kn{-.5}\hboxI{\DimenReturn}{\leaders\hbox to 3pt{\hss\box\Box\hss}\hfil}\kn{-2}\box\BoxI}\kn{2}}}}}%
                                                            \def\Iff{\Ifff{\vrulew{10}}}%
                                                            \def\SET#1#2{$\hbox{\MIIChBx{#1}{146}#2\MIIChBx{#1}{147}}$}%
                                                            \def\Set#1{\SET{10}{#1}}%
                                                            \def\SETSuch#1#2#3{\SET{#1}{#2\thinspace:\thinspace#3}}%
                                                            \def\SetSuch#1#2{\SETSuch{10}{#1}{#2}}%
                                                            \def\CAP#1{\MIIChBx{#1}{134}}%
															\def\Cap{\CAP{10}}%
                                                            \def\CUP#1{\MIIChBx{#1}{133}}%
															\def\Cup{\CUP{10}}%
                                                            \def\COP#1{\MIIChBx{#1}{64}}%
                                                            \def\Cop{\COP{10}}%
															\def\Union#1#2{$\bigcup\limits_{\hbox{#1}}{\hbox{#2}}$}%
															\def\Intersection#1#2{$\bigcap\limits_{\hbox{#1}}{\hbox{#2}}$}%
                                                            \def\Cartesian#1{\setbox\Box\hbox{#1}%
                                                                               \Dimen=\ht\Box%
                                                                               \advance\Dimen by4pt%
                                                                               \hbox{\lower\Dimen\hbox{\Stack{\hbox{\MIIChBx{20}{2}}}{\box\Box}{.5pt}}}}%
                                                            \def\Id#1{\OpGr{i}$_{\RMBx{7}{#1}}$}%
                                                            \def\Restriction#1{\thinspace\vrulewhd{1}{10}{1}\thinspace\lower4pt\hbox{\RM{7}{#1}}}%
                                                            \def\IN#1{\kn{1}\MIIChBx{#1}{62}\kn{1}}%
															\def\In{\IN{10}}%
															\def\iN{\IN{7}}%
                                                            \def\NI#1{\kn{1}\MIIChBx{#1}{63}\kn{1}}%
															\def\Ni{\NI{10}}%
                                                            \def\NIn{\hbox{\In\kn{-7}\MIIChBx{10}{66}\kn{8}}}%
															\def\Nin{\NIn}%
                                                            \def\OUT#1{\kn{1}\MIIChBx{#1}{63}\kn{1}}%
															\def\Out{\OUT{10}}%
															%
                                                            %
															%
                                                            %
                                                 			\def\SIN#1{\kn{1}\MIIChBx{#1}{32}\kn{1}}%
															\def\Sin{\SIN{10}}%
                                                 			\def\SOUT#1{\kn{1}\MIIChBx{#1}{33}\kn{1}}%
															\def\Sout{\SOUT{10}}%
															%
                                                            \def\VEE#1{\MIIChBx{#1}{137}}%
															\def\Vee{\VEE{10}}%
                                                            \def\WEDGE#1{\MIIChBx{#1}{136}}%
															\def\Wedge{\WEDGE{10}}%
															%
       														\def\Equals{$\hbox{\thinspace\RmBx{=}\thinspace}$}%
                                                            \def\NEq{\Equals\kn{-9.5}\MiiChBx{66}\kn{9.5}}%
                                                            %
                                                            %
                                                            %
                                                            \def\Eq#1{\setbox\Box=\vbox{\hrulewh{5}{.27}\kn{1.085}\hrulewh{5}{.27}\kn{.94}}%
									 								  \hbox{\kn{3}%
                                                                            \Vbox{\setbox\BoxArg=\hbox{\kn{1}\RM{7}{#1}}%
                                                                                  \addII{\wd\BoxArg}{3pt}%
                                                                                  \box\BoxArg%
                                                                                  \kn{1}%
                                                                                  \hboxI{\DimenReturn}{\leaders\hbox to 3pt{\hss\box\Box\hss}\hfil}}}\kn{2}}%
                                                            \def\EQUIV#1{\MIIChBx{#1}{21}\kn{1}}%
															\def\Equiv{\EQUIV{10}}%
															%
 															%
                                                            %
															%
															%
															%
															%
                                                            \def\CDOT#1{\MIIChBx{#1}{1}}%
															\def\Cdot{\CDOT{10}}%
															\def\cDot{\CDOT{7}}%
                                                            \def\CROSS#1{\MIIChBx{#1}{2}}%
															\def\Cross{\CROSS{10}}%
                                                            \def\CIRC#1{\MIIChBx{#1}{16}}%
															\def\Circ{\kn{.5}\MiiChBx{16}\kn{.5}}%
															\def\cIrc{\CIRC{7}}%
                                                            %
  															\def\Null{\ItChBx{34}}%
                                                            \def\Void{\Null}
                                                            \def\Proof{\par\ \Eu{Proof: }\rm}%
                                                            \def\QED{\Eu{Q}.\Eu{E}.\Eu{D}.}%
															\def\Hookrightarrow{\hbox{\MiChBx{54}\kn{-1.5}\MiiChBx{41}}}%
                     										\def\Function#1sends#2in#3to#4in#5\end{\testVoid{#1}{}{#1%
																                                                   \kn{2}%
																			                                       \vrulehwd{10}{.75}{2}%
																			                                       \kn{2}}%
																			                       \testVoid{#3}{\testVoid{#2}{}{#2%
																			                                                     \Hookrightarrow %
																			                                                     #4}}{#3%
																										                              \testVoid{#2}{\thinspace%
																				         								                            \Hookrightarrow%
																													                                \thinspace%
																													                                #5}{\Out %
																													                                    #2%
																													                                    \thinspace%
																				         								                                \Hookrightarrow%
																													                                    \thinspace%
																													                                    #4%
																													                                    \testVoid{#5}{}{\In %
																														                                                #5}}}}%
                                                            %
                                                            %
                                                            \def\Doublehookarrow{\Vbox{\hbox{\MIChBx{7}{50}\kn{-2}\MIIChBx{7}{0}\kn{-1.5}\MIChBx{7}{55}}%
                                                                                       \kn{-3}%
                                                                                       \hbox{\lower1.5pt\hbox{\MIChBx{7}{54}}\kn{-1.5}\MIIChBx{7}{0}\kn{-2}\MIChBx{7}{53}}%
                                                                                       \kn{-.5}}}%
                     										\def\Bijection#1sends#2in#3to#4in#5\end{\testVoid{#1}{}{#1%
																                                                   \kn{2}%
																			                                       \vrulehwd{10}{.75}{2}%
																			                                       \kn{2}}%
																			                       \testVoid{#3}{\testVoid{#2}{}{#2%
																			                                                     \Doublehookarrow %
																			                                                     #4}}{#3%
																										                              \testVoid{#2}{\thinspace%
																				         								                            \Doublehookarrow%
																													                                \thinspace%
																													                                #5}{\Out %
																													                                    #2%
																													                                    \thinspace%
																				         								                                \Doublehookarrow%
																													                                    \thinspace%
																													                                    #4%
																													                                    \testVoid{#5}{}{\In %
																														                                                #5}}}}%
															\def\Card#1{\##1}%
															\def\Inv#1{\hbox{{#1}$^{-1}$}}%
															%
                                                  			\def\OVER#1#2#3#4#5{\let\Top=\Box \let\Bottom=\BoxArg \let\Frac=\BoxVar%
                                                                                \setbox\Top\hbox{#1}%
													 					        \setbox\Bottom\hbox{#2}%
                                                                                \ifdim\wd\Top<\wd\Bottom \addII{-\wd\Top}{\wd\Bottom}%
												     										             \divideII{\DimenReturn}{2}%
                                                                                     	                 \setbox\Frac\hbox{\Vbox{\hbox{\kern\DimenReturn%
															              											                   \box\Top}%
																											                     \kern#5%
																											                     #4{\hrule height#3}%
																											                     \kern#5%
																											                     \box\Bottom}}%
                                                                                		          \else \addII{-\wd\Bottom}{\wd\Top}%
															      							             \divideII{\DimenReturn}{2}%
                                                                                     		             \setbox\Frac\hbox{\Vbox{\box\Top%
															     													             \kern#5%
																												                 #4{\hrule height#3}%
																											                     \kern#5%
																												                 \hbox{\kern\DimenReturn%
																												                 \box\Bottom}}}\fi%
                                                                                \addII{\ht\Frac}{\dp\Frac}%
												    						    \divideII{\DimenReturn}{2}%
																			    \addII{\DimenReturn}{-4pt}%
                                                                                \hbox{\raise\DimenReturn\hbox{\box\Frac}}}%
                                                     		\def\Over#1#2{\OVER{#1}{#2}{.5pt}{\color{black}}{1.5pt}}%
																		\def\Op#1{\Blue{#1}}%
																		\def\OpGr#1{\Op{\Gr{#1}}}%
                                                                        %
                                                                        %
                                                                  		\def\Plus{$+$}%
                                                                        \def\Minus{$-$}%
																		\def\||{\MiiChBx{153}}%
																		%
																		%
																	    %
																  	   	\def\LL#1{\hbox{\Op{\MiiChBx{142}}#1\Op{\MiiChBx{143}}}}%
                                                                        \def\LLL#1{\hbox{\Op{\MiiChBx{142}\kern-3pt\MiiChBx{142}}#1\Op{\MiiChBx{143}\kern-3pt\MiiChBx{143}}}}%
                                                                        %
                 \def\Domain#1{\setbox\Box\hbox{#1}%
                               \Dimen\ht\Box \divide\Dimen by2%
                               \DimenArg\dp\Box \advance\DimenArg by1pt%
                               \DimenEnd=\DimenArg\advance\DimenEnd by.5pt%
                               \setbox\BoxArg\hbox{\vruleHWD{\Dimen}{1pt}{\DimenArg}\kn{.5}\box\Box\kn{.5}\vruleHWD{\Dimen}{1pt}{\DimenArg}}%
                               \setbox\Box\hbox{\kn{.5}\Vbox{\copy\BoxArg\hruleHW{1pt}{\wd\BoxArg}}\kn{.5}}%
                               \hbox{\lower\DimenEnd\box\Box}}%
                  \def\Range#1{\setbox\Box\hbox{#1}%
                               \Dimen\ht\Box \divide\Dimen by2 \advance\Dimen by.5pt%
                               \setbox\BoxArg\hbox{\raise\Dimen\hbox{\vruleHW{\Dimen}{1pt}}\kn{.5}\box\Box\kn{.5}\raise\Dimen\hbox{\vruleHW{\Dimen}{1pt}}}%
                               \setbox\Box\hbox{\kn{.5}\Vbox{\hruleHW{1pt}{\wd\BoxArg}\box\BoxArg}\kn{.5}}%
                               \box\Box}%
                  \def\World#1{\setbox\Box\hbox{#1}%
                               \DimenVar=\dp\Box\advance\DimenVar by1.5pt%
                               \Dimen\ht\Box \advance\Dimen by1pt%
                               \DimenArg\dp\Box \advance\DimenArg by1pt%
                               \setbox\BoxArg\hbox{\vruleHWD{\Dimen}{1pt}{\DimenArg}\kn{.5}\box\Box\kn{.5}\vruleHWD{\Dimen}{1pt}{\DimenArg}}%
                               \setbox\Box\hbox{\kn{.5}\Vbox{\hruleHW{1pt}{\wd\BoxArg}\copy\BoxArg\hruleHW{1pt}{\wd\BoxArg}}\kn{.5}}%
                               \lower\DimenVar\box\Box}%
                  \def\Yin#1{\setbox\Box\hbox{#1}%
                             \DimenVar=\ht\Box\advance\DimenVar by1pt%
                             \DimenArg=\wd\Box\advance\DimenArg by2pt%
                             \DimenEnd=\dp\Box%
                             \advance\DimenEnd by 1pt%
                             \DimenO=\dp\Box\advance\DimenO by 2pt%
                             \hbox{\lower\DimenO\hbox{\Vbox{\Red{\hruleWH{\DimenArg}{1pt}}%
                                                      \hbox{\Red{\vruleHWD{\DimenVar}{1pt}{\DimenEnd}}\kn{1}\box\Box}%
                                                      \Red{\hruleWH{\DimenArg}{1pt}}}}}}%
                  \def\Yang#1{\setbox\Box\hbox{#1}%
                              \DimenVar=\ht\Box\advance\DimenVar by1pt%
                              \DimenArg=\wd\Box\advance\DimenArg by2pt%
                              \DimenI=\wd\Box\advance\DimenI by1pt%
                              \DimenEnd=\dp\Box%
                              \DimenO=\dp\Box\advance\DimenO by 2pt%
                              \hbox{\lower\DimenO\hbox{\Vbox{\Red{\hruleWH{\DimenArg}{1pt}}%
                                                             \hbox{\box\Box\kn{1}\Red{\vruleHWD{\DimenVar}{1pt}{\DimenEnd}}}\kn{0}%
                                                             \hbox{\kern\DimenI\Red{\vrulehw{1}{1}}}%
                                                             \Red{\hruleWH{\DimenArg}{1pt}}}}}}%
                                                                                       %
                                                                                       %
                                                                                       %
                                                                                       \def\Infty{\MIIChBx{9}{061}}%
                                                                                       \def\iNfty{\MIIChBx{7}{061}}%
                                                                                       \def\One{\EUChBx{8}{061}}%
                                                                                       \def\oNe{\EUChBx{6}{061}}%
                                                                                       \def\Two{\EUChBx{8}{062}}%
                                                                                       \def\tWo{\EUChBx{6}{062}}%
                                                                                       \def\Three{\EUChBx{8}{063}}%
                                                                                       \def\Four{\EUChBx{8}{064}}%
                                                                                       \def\Five{\EUChBx{8}{065}}%
                                                                                       \def\Six{\EUChBx{8}{066}}%
                                                                                       \def\Seven{\EUChBx{8}{067}}%
                                                                                       \def\Eight{\EUChBx{8}{070}}%
                                                                                       \def\Zero{\EUChBx{8}{060}}%
                                                                                       \def\zEro{\EUChBx{6}{060}}%
                                                                                       %
                                                                                       %
                                                        \def\Infinity#1{\hbox{\Vbox{\Op{\Infty}%
                                                                                    \kn{1}}\kn{-1}}\Of{#1}}%
                                                            \def\Schlange{\hbox{\lower3pt\hbox{\char'176}}}%
                                                            \def\And{\hbox{ and }}%
                                                            \def\Andd{\quad\And\quad}%
                                                            \def\Or{\hbox{ or }}%
                                                            \def\Orr{\quad\Or\quad}%
                                                            \def\CLUB#1{\MIIChBx{#1}{174}}%
															\def\Club{\CLUB{10}}%
															\def\cLub{\CLUB{7}}%
                                                            \def\SPADE#1{\MIIChBx{#1}{177}}%
															\def\Spade{\SPADE{10}}%
															\def\sPade{\SPADE{7}}%
															\def\RedDiamond{\hbox{\Red{\AIChBx{10}{7}}}}%
															\def\ReddIamond{\hbox{\Red{\AIChBx{7}{7}}}}%
                                                            \def\RedHeart{\hbox{\Red{\MIIChBx{10}{176}%
                                                                                \kn{-5.3}%
                                                                                \Vbox{\hrulewh{2.8}{4}%
										                         					  \kn{1.75}}%
										                        				\kn{-5.175}%
                                                                                \font\Temp=cmsy10 at 9.75pt \Temp\char'176%
									         									\kn{-7.49}%
                                                                                \font\Temp=cmsy10 at 9.5pt \Temp\char'176%
																				\kn{-7.3}%
                                                                                \font\Temp=cmsy10 at 9.25pt \Temp\char'176%
																				\kn{-7.08}%
                                                                                \font\Temp=cmsy10 at 9pt \Temp\char'176%
																				\kn{-6.9}%
                                                                                \font\Temp=cmsy10 at 8.75pt \Temp\char'176%
																				\kn{-6.7}%
                                                                                \font\Temp=cmsy10 at 8.5pt \Temp\char'176%
																				\kn{-3.5}%
                                                                                \vbox{\hrulewh{.5}{1.1}\kn{.4}}%
																				\kn{-.95}%
                                                                                \vbox{\hrulewh{1.25}{.75}\kn{.75}}%
																				\kn{-1.725}%
                                                                                \vbox{\hrulewh{2.2}{.75}\kn{1.25}}%
																				\kn{-3}%
                                                                                \vbox{\hrulewh{.75}{3.75}\kn{2.05}}%
																				\kn{2.25}%
                                                                                \vbox{\hrulewh{.75}{3.75}\kn{2.05}}%
																				\kn{-4}%
                                                                                \vbox{\hrulewh{.75}{3.7}\kn{2.5}}%
																				\kn{2.8}%
                                                                                \vbox{\hrulewh{.75}{3.7}\kn{2.5}}%
																				\kn{-4.8}%
                                                                                \vbox{\hrulewh{.75}{3.5}\kn{3}}%
																			    \kn{3.7}%
                                                                                \vbox{\hrulewh{.75}{3.5}\kn{2.9}}}%
																				\kn{1.5}}}
                                                            \def\RedhEart{\Vbox{\hbox{\hbox{\Red{\MIIChBx{7}{176}%
																				            \kn{-3.4}%
                                                                                  			\vbox{\hrulewh{1.5}{3.4}%
																							      \kn{.95}}%
																						    \kn{-3.3}%
                                                                                  			\font\Temp=cmsy17 at 6.85pt \Temp\char'176%
																						    \kn{-5}%
                                                                                  		    \font\Temp=cmsy17 at 6.7pt \Temp\char'176%
												         									\kn{-4.8}%
                                                                                  			\font\Temp=cmsy17 at 6.65pt \Temp\char'176%
																						    \kn{-4.9}%
                                                                                  			\font\Temp=cmsy17 at 6.5pt \Temp\char'176%
																						    \kn{-4.6}%
                                                                                  			\font\Temp=cmsy17 at 6.45pt \Temp\char'176%
																						    \kn{-4.75}%
                                                                                  			\font\Temp=cmsy17 at 6.3pt \Temp\char'176%
																						    \kn{-4.1}%
                                                                                  		    \vbox{\hrulewh{.6}{2.6}\kn{1.95}}%
																						    \kn{-.25}%
                                                                                  			\vbox{\hrulewh{.6}{3}\kn{1.6}}%
																						    \kn{-.2}%
                                                                                  			\vbox{\hrulewh{.6}{3.5}\kn{1.2}}%
																						    \kn{.9}%
                                                                                  			\vbox{\hrulewh{.6}{3.5}\kn{1.2}}%
																						    \kn{-.2}%
                                                                                  			\vbox{\hrulewh{.6}{3}\kn{1.6}}%
																						    \kn{-.35}%
                                                                                  			\vbox{\hrulewh{.6}{2.6}\kn{1.95}}%
																						    \kn{-2.05}%
                                                                                  			\vbox{\hrulewh{.65}{1}\kn{.4}}}%
																						    \kn{2}}}%
																					\kn{-1.5}}}%
															%
															%
															%
															%
															%
															%
                                                            \def\IronCross{\vbox{\hbox{\Green{\AIChBx{10}{172}}}\kn{-.5}}}%
                                                            \def\iRonCross{\vbox{\hbox{\Green{\AIChBx{7}{172}}}\kn{-.5}}}%
															\def\GreenStar{\hbox{\Green{\AIChBx{10}{106}}}}%
                                                            \def\gReenStar{\hbox{\Green{\AIChBx{7}{106}}}}%
                                                           	%
															%
															%
                 \def\BlueCubE#1#2#3#4#5{\setbox\Box\hbox{#1}\setbox\BoxI\hbox{\Green{\vrule width2pt}\Vbox{\Green{\hrule height2pt}\kn{1}\hbox{\kn{1}\box\Box\kn{1}}\kn{1}\Green{\hrule height2pt}}\Green{\vrule width2pt}}%
                                          \setbox\Box\hbox{#2}\setbox\BoxII\hbox{\Red{\vrule width1pt}\Vbox{\Red{\hrule height1pt}\kn{1}\hbox{\kn{1}\box\Box\kn{1}}\kn{1}\Red{\hrule height1pt}}\Red{\vrule width1pt}}%
                                          \setbox\Box\hbox{#3}\setbox\BoxIII\hbox{\Red{\vrule width1.5pt}\Vbox{\Red{\hrule height1.5pt}\kn{1}\hbox{\kn{1}\box\Box\kn{1}}\kn{1}\Red{\hrule height1.5pt}}\Red{\vrule width1.5pt}}%
                                          \setbox\Box\hbox{#4}\setbox\BoxIV\hbox{\Black{\vrule width1pt}\Vbox{\Black{\hrule height1pt}\kn{1}\hbox{\kn{1}\box\Box\kn{1}}\kn{1}\Black{\hrule height1pt}}\Black{\vrule width1pt}}%
                                          \setbox\Box\hbox{#5}\setbox\BoxV\hbox{\Black{\vrule width1.5pt}\Vbox{\Black{\hrule height1.5pt}\kn{1}\hbox{\kn{1}\box\Box\kn{1}}\kn{1}\Black{\hrule height1.5pt}}\Black{\vrule width1.5pt}}%
                                          \setbox\BoxVI\hbox{}%
                                          \setbox\Box\hbox{\kn{1}\Vbox{\copy\BoxII\kern\ht\BoxI\kern\dp\BoxI\kern\ht\BoxIV\kern\dp\BoxIV}%
                                                        \kn{1}\Vbox{\copy\BoxVI\kern\ht\BoxII\kern\dp\BoxII\kern\ht\BoxI\kern\dp\BoxI\kern\ht\BoxIV\kern\dp\BoxIV}%
                                                        \kn{1}\Vbox{\copy\BoxV\kern\ht\BoxVI\kern\dp\BoxVI\kern\ht\BoxII\kern\dp\BoxII\kern\ht\BoxI\kern\dp\BoxI\copy\BoxIV}%
                                                        \kn{1}\Vbox{\copy\BoxI\kern\ht\BoxIV\kern\dp\BoxIV}%
                                                        \kn{1}\Vbox{\copy\BoxIII\kern\ht\BoxI\kern\dp\BoxI\kern\ht\BoxIV\kern\dp\BoxIV}\kn{1}}%
                                          \setbox\BoxVII\hbox{\Blue{\vrule width1.5pt}\Vbox{\Blue{\hrule height1.5pt}\kn{1}\hbox{\kn{1}\box\Box\kn{1}}\kn{1}\Blue{\hrule height1.5pt}}\Blue{\vrule width1.5pt}}%
                                          \DimenArg=\dp\BoxVII\advance\DimenArg by\ht\BoxVII\divide\DimenArg by2\advance\DimenArg by-3pt%
                                          \lower\DimenArg\box\BoxVII}%
                  \def\DDiag{\raise1.3297pt\hbox{\kn{.15}\Vbox{\kn{.2}\hbox{\kn{-.3}\vbox{\kn{-.6}\hbox{\AIICh{7}{37}}}\kn{-1.1}}\kn{-1.25}}}}%
                  \def\Cube{\setbox\Box\hbox{\DDiag}%
                            \def\M{\Dimen}\M=\ht\Box\multiply\M by3\divide\M by2%
                            \def\N{\DimenArg}\N=\ht\Box%
                            \def\O{\dimen242}\O=\M\advance\O by-.4pt%
                            \setbox\BoxI\hbox{\kn{.1}%
                                          \Vbox{\copy\Box%
                                                \kn{-3.3}\kern\ht\Box%
                                                \copy\Box%
                                                \kn{.15}}}%
                            \setbox2\hbox{\kn{-.2}%
                                          \Vbox{\hruleHW{\M}{.4pt}\kn{.225}}%
                                          \Vbox{\hbox{\kn{-.4}%
                                                      \dimen243=\M\advance\dimen243by.5pt%
                                                      \vruleHW{.4pt}{\dimen243}}%
                                                \kern\M\kn{-.75}%
                                                \hruleHW{.4pt}{\M}%
                                                \kn{.225}}%
                                          \kn{-.2}%
                                          \Vbox{\hruleHW{\M}{.4pt}\kn{.225}}}%
                            \setbox3\hbox{\kn{-.2}%
                                          \Vbox{\copy\Box%
                                                \kn{-3.35}\kern\ht\Box%
                                                \Gray{\copy\Box}%
                                                \kn{.2}}}%
                            \setbox4\hbox{\kn{-.2}%
                                          \vruleHW{\M}{.4pt}%
                                          \Vbox{\hruleHW{.4pt}{\M}%
                                                \kern\M\kn{-.8}%
                                                \Gray{\hruleHW{.4pt}{\M}}}%
                                          \kn{-.2}%
                                          \Vbox{\hruleHW{.4pt}{.4pt}%
                                                \Gray{\hruleHW{\O}{.4pt}}}}%
                            \setbox\Box\hbox{\kn{.2}%
                                          \raise\N\box4%
                                          \box3%
                                          \kern-\M%
                                          \kern-\N%
                                          \box\BoxI%
                                          \box2}%
                            \hbox{\lower4pt\hbox{\Vbox{\box\Box\kn{-.2}}}}}%
                  \def\CubE#1#2#3#4#5#6{\setbox\Box\hbox{#1}\setbox1\hbox{\Vbox{\copy\Box\kern\dp\Box}}%
                                        \setbox\Box\hbox{#2}\setbox3\hbox{\Vbox{\copy\Box\kern\dp\Box}}%
                                        \setbox\Box\hbox{#3}\setbox2\hbox{\Vbox{\copy\Box\kern\dp\Box}}%
                                        \setbox\Box\hbox{#4}\setbox4\hbox{\Vbox{\copy\Box\kern\dp\Box}}%
                                        \setbox\Box\hbox{#5}\setbox5\hbox{\Vbox{\copy\Box\kern\dp\Box}}%
                                        \setbox\Box\hbox{#6}\setbox6\hbox{\Vbox{\copy\Box\kern\dp\Box}}%
                                        \setbox\BoxXII\hbox{\Green{\DDiag}}%
                                        \setbox\BoxVIII\hbox{\raise4pt\hbox{\Cube}}%
                                        \setbox10\hbox{\ReddIamond}%
                                        \def\RDht{\dimen240}\RDht=7.5pt\def\RDwd{\dimen241}\RDwd=.5pt\advance\RDwd by-\wd10%
                                        \def\GSht{\dimen238}\GSht=24.5pt\def\GSwd{\dimen239}\GSwd=-7pt%
                                        \def\BSht{\dimen236}\BSht=31.5pt\def\BSwd{\dimen237}\BSwd=8.2pt%
                                        \def\RHht{\dimen234}\RHht=8pt\def\RHwd{\dimen235}\RHwd=23.5pt%
                                        \def\BCht{\dimen232}\BCht=-15.2pt\def\BCwd{\dimen233}\BCwd=8.2pt%
                                        \def\ICht{\dimen230}\ICht=-6.5pt\def\ICwd{\dimen231}\ICwd=22pt%
                                        \def\Half{\dimen242}\Half=-\GSwd\advance\Half by\wd5%
                                        \def\Whole{\dimen243}\Whole=\RDwd\advance\Whole by\wd4%
                                \setbox9\hbox{\ifdim\Half>\Whole\kern\Half\else \kern\Whole\fi
                                              \copy\BoxVIII%
                                              \kern-\wd\BoxVIII%
                                              \kn{11}%
                                              \Vbox{\copy\BoxXII%
                                                    \kn{-2.5}%
                                                    \hbox{\kn{3}\iRonCross}%
                                                    \kn{-5}}%
                                              \kn{-20.9}%
                                              \kn{8}%
                                              \vrulewd{.4}{9}%
                                              \kn{-2.9}%
                                              \Vtop{\kern9pt%
                                                    \cLub}%
                                              \kn{-10.94}%
                                             \kn{16.1}%
                                             \Vtop{\kn{-8}\Red{\hrulehw{.4}{4}}}%
                                             \kn{-2}%
                                             \Vtop{\kn{-10}\RedhEart}%
                                             \kn{-22.95}%
                                             \kn{7.75}%
                                             \Vbox{\hrulewh{.4}{13}\kn{12.5}}%
                                             \kn{-2.9}%
                                             \Vbox{\sPade\kn{25}}%
                                             \kn{-10.95}%
                                             \kn{-5}%
                                             \Vbox{\copy\BoxXII\kn{15.5}}%
                                             \kn{-8}%
                                             \Vbox{\gReenStar\kn{17.75}}%
                                             \kn{.465}%
                                             \kn{-3}%
                                             \Vbox{\Red{\hrulewh{6}{.4}}\kn{7.5}}%
                                             \kn{-8}%
                                             \Vtop{\kn{-10.5}\ReddIamond}%
                                             \kn{-.34}%
                                             \Half=-\wd1\divide\Half by2%
                                             \kern\BSwd\kern\Half\Vbox{\copy1\kern\BSht}%
                                             \kern-\BSwd\kern\Half%
                                             \Half=-\ht2\divide\Half by2%
                                             \kern\RHwd\Vbox{\copy2\kern\RHht\kern\Half}%
                                             \kern-\RHwd\kern-\wd2%
                                             \Half=-\wd3\divide\Half by2%
                                             \kern\BCwd\kern\Half\Vtop{\kern-\BCht\copy3}%
                                             \kern-\BCwd\kern\Half%
                                             \Half=-\ht4\divide\Half by2%
                                             \kern\RDwd\kern-\wd4\Vbox{\copy4\kern\RDht\kern\Half}%
                                             \kern-\RDwd%
                                             \kern\GSwd\kern-\wd5\Vbox{\copy5\kern\GSht}%
                                             \kern-\GSwd\kern\ICwd\Vtop{\hbox{}\kern-\ICht\copy6}}%
                                             \Half=\wd9\divide\Half by2\advance\Half by4pt%
                                             \hbox{\lower4pt\box9}}%
\def\CubeUL{\setbox\BoxI\hbox{\sPade\RedhEart\gReenStar}%
            \setbox\BoxII\hbox{\cLub\ReddIamond\iRonCross}%
            \ifdim\wd\BoxI<\wd\BoxII\DimenO=\wd\BoxII\else\DimenO=\wd\BoxI\fi%
            \DimenI=\DimenO\advance\DimenIby2pt%
            \setbox\Box\hbox{\Vbox{\hbox to\DimenO{\hfil\box\BoxI\hfil}%
                                \kn{1}%
                                \hbox{\kn{-1}\vruleHW{1pt}{\DimenI}\kn{-1}}%
                                \kn{1}%
                                \hbox to\DimenO{\hfil\box\BoxII\hfil}}}%
            \hbox{\lower7pt\hbox{\HBox{1pt}{1pt}{\box\Box}}}}%
\def\CubeUR{\setbox\BoxI\hbox{\sPade\RedhEart\iRonCross}%
            \setbox\BoxII\hbox{\cLub\ReddIamond\gReenStar}%
            \ifdim\wd\BoxI<\wd\BoxII\DimenO=\wd\BoxII\else\DimenO=\wd\BoxI\fi%
            \DimenI=\DimenO\advance\DimenIby2pt%
            \setbox\Box\hbox{\Vbox{\hbox to\DimenO{\hfil\box\BoxI\hfil}%
                                \kn{1}%
                                \hbox{\kn{-1}\vruleHW{1pt}{\DimenI}\kn{-1}}%
                                \kn{1}%
                                \hbox to\DimenO{\hfil\box\BoxII\hfil}}}%
            \hbox{\lower7pt\hbox{\HBox{1pt}{1pt}{\box\Box}}}}%
\def\CubeLR{\setbox\BoxI\hbox{\sPade\ReddIamond\gReenStar}%
            \setbox\BoxII\hbox{\cLub\RedhEart\iRonCross}%
            \ifdim\wd\BoxI<\wd\BoxII\DimenO=\wd\BoxII\else\DimenO=\wd\BoxI\fi%
            \DimenI=\DimenO\advance\DimenIby2pt%
            \setbox\Box\hbox{\Vbox{\hbox to\DimenO{\hfil\box\BoxI\hfil}%
                                \kn{1}%
                                \hbox{\kn{-1}\vruleHW{1pt}{\DimenI}\kn{-1}}%
                                \kn{1}%
                                \hbox to\DimenO{\hfil\box\BoxII\hfil}}}%
            \hbox{\lower7pt\hbox{\HBox{1pt}{1pt}{\box\Box}}}}%
\def\CubeLL{\setbox\BoxI\hbox{\sPade\ReddIamond\iRonCross}%
            \setbox\BoxII\hbox{\cLub\RedhEart\gReenStar}%
            \ifdim\wd\BoxI<\wd\BoxII\DimenO=\wd\BoxII\else\DimenO=\wd\BoxI\fi%
            \DimenI=\DimenO\advance\DimenIby2pt%
            \setbox\Box\hbox{\Vbox{\hbox to\DimenO{\hfil\box\BoxI\hfil}%
                                \kn{1}%
                                \hbox{\kn{-1}\vruleHW{1pt}{\DimenI}\kn{-1}}%
                                \kn{1}%
                                \hbox to\DimenO{\hfil\box\BoxII\hfil}}}%
            \hbox{\lower7pt\hbox{\HBox{1pt}{1pt}{\box\Box}}}}%
\def\CubeDia{\setbox\Box\hbox{\Vbox{\kn{-.5}\hbox{\AiiCh{36}\kn{-9}\AiiCh{37}}}}%
             \setbox\BoxI\hbox to \wd\Box{\Vbox{\hruleWH{\wd\Box}{.4pt}\kn{2.25}}}%
             \setbox\BoxII\hbox{\raise1pt\hbox{\copy\Box\kern-\wd\Box\box\BoxI\kn{-4.75}\vruleWHD{.4pt}{\ht\Box}{\dp\Box}}}%
             \hbox{\box\BoxII\kn{4}}}%
                            \def\Tet{\setbox\Box\hbox{\kn{-.4}\AIICh{5}{36}\kn{-1.75}}%
                                     \setbox\BoxI\hbox{\raise\dp\Box\box\Box}%
                                     \setbox\Box\hbox{\kn{-.4}\Vbox{\kn{-.35}\AIIChBx{5}{37}}\kn{-1.75}}%
                                     \setbox\BoxII\box\Box%
                                     \lower1pt\hbox{\kn{.5}\hbox{\Vbox{\hbox{\copy\BoxI\copy\BoxII}\kn{-.2}\hbox{\copy\BoxII\copy\BoxI}}\kn{.5}}}}%
                            \def\TetXII{\hbox{\kn{1}\Tet\raise4pt\hbox{\kn{-6.25}\Red{\MIIChBx{7}{17}}}\kn{3}}}%
                            \def\TetVI{\hbox{\kn{1}\Tet\lower.6pt\hbox{\kn{-6.25}\Red{\MIIChBx{7}{17}}}\kn{3}}}%
                            \def\TetIII{\hbox{\kn{1}\Tet\raise1.6pt\hbox{\kn{-4.25}\Blue{\MIIChBx{7}{17}}\kn{1}}}}%
                            \def\TetIX{\hbox{\kn{1}\Tet\raise1.6pt\hbox{\kn{-8.5}\Blue{\MIIChBx{7}{17}}}\kn{5}}}%
                            \def\TetIX{\Spade}%
                            \def\TetXII{\RedHeart}%
                            \def\TetIII{\Club}%
                            \def\TetVI{\RedDiamond}%
                             \def\TeT#1#2#3#4{\setbox\Box\hbox{#1}\setbox\BoxI\hbox{\Vbox{\copy\Box\kern-\dp\Box}}%
                                              \DimenI=\ht\BoxI\divide\DimenI by2%
                                              \setbox\Box\hbox{#2}\setbox\BoxII\hbox{\Vbox{\copy\Box\kern-\dp\Box}}%
                                              \DimenII=\wd2\divide\DimenII by2%
                                              \setbox\Box\hbox{#3}\setbox\BoxIII\hbox{\Vbox{\copy\Box\kern-\dp\Box}}%
                                              \DimenIII=\ht\BoxI\divide\DimenIIIby2%
                                              \setbox\Box\hbox{#4}\setbox\BoxIV\hbox{\Vbox{\copy\Box\kern-\dp\Box}}%
                                              \DimenV=\wd\BoxIV\divide\DimenVby2%
                                              \DimenIV=\ht\BoxIV\divide\DimenIVby2%
                                              \setbox\Box\hbox{\Tet}%
                                              \DimenO=\ht\Box\divide\DimenO by2%
                                              \Dimen=\wd\Box\divide\Dimen by2%
                                              \setbox\BoxV\hbox{\Vbox{\copy\BoxI%
                                                                  \kern-\DimenI%
                                                                  \kern\DimenO%
                                                                  \kn{1}%
                                                                  \kern\ht\BoxIV}%
                                                            \Vbox{\hbox{\kern\Dimen\kern-\DimenII\copy\BoxII}%
                                                                  \kn{1}%
                                                                  \copy\Box%
                                                                  \kn{1}%
                                                                  \hbox{\kern\Dimen\kern-\DimenV\copy\BoxIV}}%
                                                            \Vbox{\copy\BoxIII%
                                                                  \kern-\DimenIII%
                                                                  \kern\DimenO%
                                                                  \kn{1}%
                                                                  \kern\ht\BoxIV}}%
                                              \DimenVI=\ht\BoxV\divide\DimenVIby2\advance\DimenVIby-2pt%
                                              \lower\DimenVI\box\BoxV}%
                             \def\MerInv#1#2#3#4{\setbox\Box\hbox{#1}\setbox\BoxI\hbox{\Vbox{\copy\Box\kern\dp\Box}}%
                                                 \setbox\Box\hbox{#2}\setbox\BoxII\hbox{\Vbox{\copy\Box\kern\dp\Box}}%
                                                 \setbox\Box\hbox{#3}\setbox\BoxIII\hbox{\Vbox{\copy\Box\kern\dp\Box}}%
                                                 \setbox\Box\hbox{#4}\setbox\BoxIV\hbox{\Vbox{\copy\Box\kern\dp\Box}}%
                                                 \def\HI{\DimenO}\HI=\ht\BoxI\divide\HI by2%
                                                 \def\WI{\Dimen}\WI=\wd\BoxI\divide\WI by2%
                                                 \def\HII{\DimenI}\HII=\ht\BoxII\divide\HII by2%
                                                 \def\WII{\DimenArg}\WII=\wd\BoxII\divide\WII by2%
                                                 \def\HIII{\DimenII}\HIII=\ht\BoxIII\divide\HIII by2%
                                                 \def\WIII{\DimenVar}\WIII=\wd\BoxIII\divide\WIII by2%
                                                 \def\HIV{\DimenIII}\HIV=\ht\BoxIV\divide\HIV by2%
                                                 \def\WIV{\DimenEnd}\WIV=\wd\BoxIV\divide\WIV by2%
                                                 \def\H{\DimenIV}%
                                                 \set\Box\hbox{\ifdim\WIII>\WI\H=\WIII\advance\H by-\WI\kern\H\fi%
                                                                \Vbox{\ifdim\HII>\HI\H=\HII\advance\H by-\HI\kern\H\fi%
                                                                      \copy\BoxI%
                                                                      \ifdim\HII>\HI\H=\HII\advance\H by-\HI\kern\H\fi}%
                                                                \ifdim\WIII>\WI\H=\WIII\advance\H by-\WI\kern\H\fi}%
                                                 \setbox\BoxArg\hbox{\ifdim\WIV>\WII\H=\WIV\advance\H by-\WII\kern\H\fi%
                                                                \Vbox{\ifdim\HII<\HI\H=\HI\advance\H by-\HII\kern\H\fi%
                                                                      \copy\BoxII%
                                                                      \ifdim\HII<\HI\H=\HI\advance\H by-\HII\kern\H\fi}%
                                                                \ifdim\WIV>\WII\H=\WIV\advance\H by-\WII\kern\H\fi}%
                                                 \setbox\BoxIII\hbox{\ifdim\WIII<\WI\H=\WI\advance\H by-\WIII\kern\H\fi%
                                                                \Vbox{\ifdim\HIV>\HIII\H=\HIV\advance\H by-\HIII\kern\H\fi%
                                                                      \copy\BoxIII%
                                                                      \ifdim\HIV>\HIII\H=\HIV\advance\H by-\HIII\kern\H\fi}%
                                                                \ifdim\WIII<\WI\H=\WI\advance\H by-\WIII\kern\H\fi}%
                                                 \setbox\BoxIV\hbox{\ifdim\WII>\WIV\H=\WII\advance\H by-\WIV\kern\H\fi%
                                                                \Vbox{\ifdim\HIII>\HIV\H=\HIII\advance\H by-\HIV\kern\H\fi%
                                                                      \copy\BoxIV%
                                                                      \ifdim\HIII>\HIV\H=\HIII\advance\H by-\HIV\kern\H\fi}%
                                                                \ifdim\WII>\WIV\H=\WII\advance\H by-\WIV\kern\H\fi}%
                                                 \setbox\Box\hbox{\Vbox{\Blue{\hrule}\hbox{\Blue{\vrule}\kn{1}\Vbox{\kn{1}\hbox{\\Box\kn{1}\box\BoxArg}\kn{2}\hbox{\box\BoxIII\kn{1}\box\BoxIV}\kn{1}}\kn{1}\Blue{\vrule}}\Blue{\hrule}}}%
                                                 \H=\ht\Box\divide\H by 2\advance\H by-2pt%
                                                 \hbox{\lower\H\box\Box}}%
                  \def\Subsets#1{\EB{12}{2}$^{\hbox{#1}}$}%
                  \def\Constant#1#2{\setbox\Box\hbox{#1}%
                                    \setbox\BoxI\hbox{\Vbox{\kn{-1}\box\Box}}%
                                    \Dimen=\dp\BoxI\advance\Dimen by\ht\BoxI%
                                    \setbox\Box\hbox{#2}%
                                    \setbox\BoxII\hbox{\Vbox{\box\Box\kn{-1}}}%
                                    \setbox\Box\hbox{$\buildrel\mapsto\over{\box\BoxI}$}%
                                    \hbox{\kn{1}\lower\Dimen\hbox{$\buildrel{\box\BoxII}\over{\box\Box}$}}\kn{1}}%
                             \def\BasisO{\Eu{0}$_{\FTChBx{7}{142}}$}%
                             \def\BasisI{\Eu{1}$_{\FTChBx{7}{142}}$}%
                             \def\BasisInf{\BSChBx{10}{61}$_{\FTChBx{7}{142}}$}%
                             \def\Arc{\Of{\BasisO;\BasisI;\BasisInf}}%
                                                                        \def\Distinct{\hbox{ distinct}}%
                                                                        \def\Sim{\MiiChBx{30}}%
                                                                        \def\sIm{\MiiChBx{7}{30}}%
																		\def\sIm{\MIIChBx{7}{30}}%
                                                                        \def\OSim{\MiiChBx{30}\kn{-5.75}\Vbox{\TWBx{7}{o}\kn{1.1}}\kn{2}}%
                                                                        %
                                                                        %
                                                                        %
                                                                        \def\Star#1{\setbox\BoxArg\hbox{\RedOrangeHBox{1pt}{1pt}{#1}}%
                                                                                    \setbox\Box\hbox{\AICh{5}{106}}%
                                                                                    \Stack{\box\Box}{\box\BoxArg}{1pt}}%
                                                                        \def\sTAR#1{\setbox\Box\hbox{\RedOrangeHBox{1pt}{1pt}{\RMBx{7}{#1}}}%
                                                                                    \setbox\BoxArg\hbox{\AICh{4}{106}}%
                                                                                    \Stack{\box\BoxArg}{\box\Box}{1pt}}%
                          	                                            \def\Flower#1{\setbox\Box\hbox{\RedOrangeHBox{1pt}{1pt}{#1}}%
                                                                                      \setbox\BoxArg\hbox{\AICh{5}{147}}%
                                                                                      \StackedVbox{\box\BoxArg}{\box\Box}{1pt}}%
        				                                                \def\fLOWER#1{\setbox\Box\hbox{\RedOrangeHBox{1pt}{1pt}{\RMBx{7}{#1}}}%
                                                                                      \setbox\BoxArg\hbox{\AICh{4}{147}}%
                                                                                      \StackedVbox{\box\BoxArg}{\box\Box}{1pt}}%
                                                                        \def\Bump#1{\setbox\Box\hbox{\VioletHBox{1pt}{1pt}{#1}}%
                                                                                    \setbox\BoxArg\hbox{\AICh{5}{155}}%
                                                                                    \StackedVbox{\box\BoxArg}{\box\Box}{1pt}}%
                                                                        \def\bUMP#1{\setbox\Box\hbox{\VioletHBox{1pt}{1pt}{\RMBx{7}{#1}}}%
                                                                                    \setbox\BoxArg\hbox{\AICh{4}{155}}%
                                                                                    \StackedVbox{\box\BoxArg}{\box\Box}{1pt}}%
                                                                        \def\\{{\TwChBx{134}}}%
                                                                        \def\?{\thinspace}%
                                                                        \def\Prime{\hbox{\hskip-2pt\RmChBx{23}}}%
                                                                        %
                           								                \def\Quotes#1{\Rm{\char'140}\nobreak\hskip-.8pt\nobreak\Rm{\char'140}\nobreak#1\nobreak\Rm{\char'047}\nobreak\hskip-.8pt\nobreak\Rm{\char'047}}%
   											                            \def\VECTOR#1{\setbox\Box\hbox{#1}%
														                              \addII{\wd\Box}{1pt}%
                  												                      \setbox\BoxArg\hbox{\Vbox{\Blue{\MIChBx{7}{76}}%
												  			                                              \kn{-2.11}}%
														                                            \kn{-1}%
														                                            \Blue{\vruleWHD{\DimenReturn}{.35pt}{0pt}}%
														                                            \kn{-5}%
														                                            \Vbox{\Blue{\MIChBx{7}{76}}%
														                                            \kn{-2.11}}}%
                  												                                    \Stack{\box\BoxArg}{\box\Box}{1pt}}%
											                             \def\Line#1{\setbox\BoxArg\hbox{#1}%
											                                         \setbox\Box\hbox{\MIIChBx{7}{44}}%
            												                         \ifdim\wd\BoxArg<\wd\Box \Stack{\copy\Box}{\copy\BoxArg}{1pt}%
                                                                                                        \else \ifdim\wd\BoxArg<13pt \addII{-13pt}{\wd\BoxArg}%
                                                                                                                                   \setbox\Box\hbox{\MIIChBx{7}{40}%
                                                                                                                                                     \kern\DimenReturn%
                                                                                                                                                     \MIIChBx{7}{41}}%
                                                                                                                             \else \addII{\wd\BoxArg}{-14pt}%
                                                                                                                                   \setbox\Box\hbox{\MIIChBx{7}{40}%
                                                                                                                                                     \kn{-.5}%
                                                                                                                                                     \Vbox{\hruleWH{\DimenReturn}{.25pt}%
                                                                                                                                                           \kn{1.6}}%
                                                                                                                                                     \kn{-.5}%
                                                                                                                                                     \MIIChBx{7}{41}}\fi\fi%
                                                                                    \kn{1}%
													                                \Stack{\copy\Box}{\copy\BoxArg}{.75pt}%
													                                \kn{1}}%
                                                                         \def\FilledRightHalfArrow{\setbox\BoxArg\vbox{\hbox{\Vbox{\hbox{\MiCh{52}}%
                                                                                                                        \vbox  to-5pt{}%
                                                                                                                        \hbox{\kn{6.9}\vrulewhd{.4}{2.6}{0}}}}}%
                                                                                                   \setbox13\vbox{\box\BoxArg%
                                                                                                                  \kn{2.4}}%
                                                                                                   \box13}%
                                                                         \def\OverFilledRightHalfArrow{\setbox\BoxArg\vbox{\hbox{\Vbox{\hbox{\kn{6.9}\vrulewh{2.6}{.4}}%
                                                                                                                                  \kn{1.25}%
                                                                                                                                  \hbox{\MiCh{52}}%
                                                                                                                                  \vbox to-5pt{}%
                                                                                                                                  \hbox{\kn{6.9}\vrulewhd{.4}{2.6}{0}}}}}%
                                                                                                        \setbox13\vbox{\box\BoxArg%
                                                                                                                       \kn{2.4}}%
                                                                                                        \box13}%
                                                                          \def\FilledLeftHalfArrow{\setbox\BoxArg\vbox{\hbox{\Vbox{\hbox{\MiCh{50}}%
                                                                                                                        \vbox  to-5pt{}%
                                                                                                                        \hbox{\kn{2.675}\vrulewhd{.4}{2.6}{0}}}}}%
                                                                                                   \setbox13\vbox{\box\BoxArg%
                                                                                                                  \kn{2.4}}%
                                                                                                   \box13}%
                                                                          \def\OverFilledLeftHalfArrow{\setbox\BoxArg\vbox{\hbox{\Vbox{\hbox{\kn{.5}\vrulewh{2.6}{.4}}%
                                                                                                                                  \kn{1.25}%
                                                                                                                                  \hbox{\MiCh{50}}%
                                                                                                                                  \vbox  to-5pt{}%
                                                                                                                                  \hbox{\kn{2.675}\vrulewhd{.4}{2.6}{0}}}}}%
                                                                                                       \setbox13\vbox{\box\BoxArg%
                                                                                                                      \kn{2.4}}%
                                                                                                       \box13}%
                                                                          \def\setRel#1{\let\TopBox=\Box \let\BaseBox=\BoxArg%
                                                                                        \setbox\BaseBox\hbox{#1}%
                                                                                        \setbox\TopBox\hbox{\MiCh{52}}%
                                                                                        \ifdim\wd\BaseBox<\wd\TopBox \Stack{\box\TopBox}{\box\BaseBox}{-1.5pt}%
                                                                                                      \else \addIII{\wd\BaseBox}{-\wd\TopBox}{2pt}%
                                                                                                            \Vbox{\hbox{\vruleWH{\DimenReturn}{.4pt}%
                                                                                                                        \kn{-2}%
                                                                                                                        \lower2.3pt\box\TopBox}%
                                                                                                                        \kn{-1.5}%
                                                                                                                        \box\BaseBox}\fi}%
                                                                     \def\setRelation#1{\let\TopBox=\Box \let\BaseBox=\BoxArg%
                                                                                        \setbox\BaseBox\hbox{#1}%
                                                                                        \setbox\TopBox\hbox{\FilledRightHalfArrow}%
                                                                                        \ifdim\wd\BaseBox<\wd\TopBox \Stack{\box\TopBox}{\box\BaseBox}{-1pt}%
                                                                                                      \else \addIII{\wd\BaseBox}{-\wd\TopBox}{2pt}%
                                                                                                            \Vbox{\hbox{\vruleWH{\DimenReturn}{.4pt}%
                                                                                                                        \kn{-2}%
                                                                                                                        \lower2.3pt\box\TopBox}%
                                                                                                                        \kn{-1}%
                                                                                                                        \box\BaseBox}\fi}%
                                                                     \def\SetRelation#1{\let\TopBox=\Box \let\BaseBox=\BoxArg%
                                                                                        \setbox\BaseBox\hbox{#1}%
                                                                                        \setbox\TopBox\hbox{\OverFilledRightHalfArrow}%
                                                                                        \ifdim\wd\BaseBox<\wd\TopBox \Stack{\box\TopBox}{\box\BaseBox}{-1pt}%
                                                                                                      \else \addIII{\wd\BaseBox}{-\wd\TopBox}{2pt}%
                                                                                                            \Vbox{\hbox{\vruleWH{\DimenReturn}{.4pt}%
                                                                                                                        \kn{-2}%
                                                                                                                        \lower2.3pt\box\TopBox}%
                                                                                                                        \kn{-1}%
                                                                                                                        \box\BaseBox}\fi}%
                                                                 \def\setRelBack#1{\let\TopBox=\Box \let\BaseBox=\BoxArg%
                                                                                        \setbox\BaseBox\hbox{#1}%
                                                                                        \setbox\TopBox\hbox{\MiCh{50}}%
                                                                                        \ifdim\wd\BaseBox<\wd\TopBox \Stack{\box\TopBox}{\box\BaseBox}{-1.5pt}%
                                                                                                               \else \addIII{\wd\BaseBox}{-\wd\TopBox}{2pt}%
                                                                                                                     \Vbox{\hbox{\lower2.3pt\box\TopBox%
                                                                                                                           \kn{-2}%
                                                                                                                           \vruleWH{\DimenReturn}{.4pt}}%
                                                                                                                     \kn{-1.5}%
                                                                                                                     \box\BaseBox}\fi}%
                                                                 \def\SetRelationBack#1{\let\TopBox=\Box \let\BaseBox=\BoxArg%
                                                                                        \setbox\BaseBox\hbox{#1}%
                                                                                        \setbox\TopBox\hbox{\OverFilledLeftHalfArrow}%
                                                                                        \ifdim\wd\BaseBox<\wd\TopBox \Stack{\box\TopBox}{\box\BaseBox}{-1pt}%
                                                                                                               \else \addIII{\wd\BaseBox}{-\wd\TopBox}{2pt}%
                                                                                                                     \Vbox{\hbox{\lower2.3pt\box\TopBox%
                                                                                                                           \kn{-2}%
                                                                                                                           \vruleWH{\DimenReturn}{.4pt}}%
                                                                                                                     \kn{-1}%
                                                                                                                     \box\BaseBox}\fi}%
                                                                 \def\setRelationBack#1{\let\TopBox=\Box \let\BaseBox=\BoxArg%
                                                                                        \setbox\BaseBox\hbox{#1}%
                                                                                        \setbox\TopBox\hbox{\FilledLeftHalfArrow}%
                                                                                        \ifdim\wd\BaseBox<\wd\TopBox \Stack{\box\TopBox}{\box\BaseBox}{-1pt}%
                                                                                                               \else \addIII{\wd\BaseBox}{-\wd\TopBox}{2pt}%
                                                                                                                     \Vbox{\hbox{\lower2.3pt\box\TopBox%
                                                                                                                           \kn{-2}%
                                                                                                                           \vruleWH{\DimenReturn}{.4pt}}%
                                                                                                                     \kn{-1}%
                                                                                                                     \box\BaseBox}\fi}%
                                                                        \def\RightPolar#1{\setbox\Box\hbox{#1}%
                                                                                          \setbox\BoxArg\hbox{\Vbox{\hruleWH{\wd\Box}{.5pt}%
                                                                                                                    \kn{2.25}}%
                                                                                                              \kn{-.75}%
                                                                                                              \MIIChBx{10}{16}}%
                                                                                                         \hbox{\Vbox{\box\BoxArg%
                                                                                                               \box\Box}}}%
                                                                         \def\LeftPolar#1{\setbox\Box\hbox{#1}%
                                                                                          \setbox\BoxArg\hbox{\MIIChBx{10}{16}%
                                                                                                              \kn{-.75}%
                                                                                                              \Vbox{\hruleWH{\wd\Box}{.5pt}%
                                                                                                                    \kn{2.25}}}%
                                                                                                 \hbox{\Vbox{\box\BoxArg%
                                                                                                             \hbox{\kn{1.5}\box\Box}}}}%
                                                                         %
                                                                         %
                                                                    \def\UnderLine#1{\setbox\Box\hbox{#1}%
											                                         \setbox\BoxArg\hbox{\MIIChBx{7}{44}}%
            												                         \ifdim\wd\Box<\wd\BoxArg \Stack{\box\Box}{\box\BoxArg}{1pt}%
                                                                                                 \else \ifdim\wd\Box<13pt \addII{\Dimen}{-13pt}{\wd\Box}%
                                                                                                                        \setbox\BoxArg\hbox{\MIIChBx{7}{40}%
                                                                                                                                       \kern\Dimen%
                                                                                                                                       \MIIChBx{7}{41}}%
                                                                                                                  \else \addII{\Dimen}{\wd\Box}{-14pt}%
                                                                                                                        \setbox\BoxArg\hbox{\MIIChBx{7}{40}%
                                                                                                                                       \kn{-.5}%
                                                                                                                                       \Vbox{\hruleWH{\Dimen}{.25pt}%
                                                                                                                                             \kn{1.6}}%
                                                                                                                                       \kn{-.5}%
                                                                                                                                       \MIIChBx{7}{41}}\fi\fi%
                         											                \kn{1}%
													                                \lower4.25pt\hbox{\Stack{\box\Box}{\box\BoxArg}{1.75pt}}%
													                                \kn{1}}%
    											                         \def\Radior#1{\setbox\Box\hbox{#1}%
											                                           \addII{\Dimen}{\wd\Box}{1pt}%
                  											                           \setbox1\hbox{\Vbox{\MIIChBx{7}{140}%
											                                                               \kn{-2.2}}%
											                                                         \kn{-3.6}%
											                                                         \vruleWHD{\Dimen}{.35pt}{0pt}%
											                                                         \kn{-3.6}%
											                                                         \Vbox{\MIIChBx{7}{141}%
											                                                         \kn{-2.2}}}%
                  												                       \Stack{\box1}{\box\Box}{1pt}}%
    											                         \def\RADIOR#1{\setbox\Box\hbox{#1}%
											                                           \addII{\Dimen}{\wd\Box}{1pt}%
                  											                           \setbox1\hbox{\Vbox{\Blue{\MIIChBx{7}{140}}%
											                                                               \kn{-2.2}}%
											                                                         \kn{-3.6}%
											                                                         \Blue{\vruleWHD{\Dimen}{.35pt}{0pt}}%
											                                                         \kn{-3.6}%
											                                                         \Vbox{\Blue{\MIIChBx{7}{141}}%
											                                                         \kn{-2.2}}}%
                  												                       \Stack{\box1}{\box\Box}{1pt}}%
                          								                \def\Restrict#1#2{\setbox\Box\hbox{#2}%
								                                                          \addIII{\dimen249}{-8pt}{\ht\Box}{\dp\Box}%
								                                                          \divideII{\dimen249}{\dimen249}{2}%
                                            									          \hbox{#1%
									    						                                \kn{1}%
															                                    \vrulewhd{1}{9}{2}%
															                                    \kn{1}%
															                                    \Vtop{\hbox{}%
															     	                                  \kern\dimen249%
															                                          \box\Box}}}%
											%
											%
											%
											%
											%
											%
											%
											\def\WideHat#1{$\widehat{\hbox{#1}}$}%
                                            \def\Of#1{\hbox{\Vbox{\Brown{\MIIIChBx{5}{070}\MIIIChBx{5}{072}}\kn{-2}}%
                                                            \Vbox{\hbox{#1}}%
                                                            \Vbox{\Brown{\MIIIChBx{5}{071}\MIIIChBx{5}{073}}\kn{-2}}}}%
                                            \def\oF#1{\hbox{\Vbox{\MIIIChBx{5}{070}\kn{-3}\MIIIChBx{5}{072}\kn{-1}}%
                                                            \Vbox{\hbox{#1}}%
                                            				\Vbox{\MIIIChBx{5}{071}\kn{-3}\MIIIChBx{5}{073}\kn{-1}}}}%
                                                              			\def\Pr#1{\hbox{\kn{1}
                                                                                        \Vbox{\hbox{\kn{-1.4}\TWChBx{7}{023}}%
                                                                                			  \kn{-3.4}%
                                                                                			  \hbox{\vruleWH{.25pt}{7.2pt}}%
                                                                                              \kn{-.3}%
                                                                                              \hbox{\kn{-.8}\TWChBx{7}{022}}%
                                                                                              \kn{-5}}%
                                                                          				\kn{-1.2}%
                                                                          				\Vbox{\hbox{\kn{1}#1}}%
                                                                          				\kn{1.5}%
                                                                          				\Vbox{\hbox{\kn{-.8}%
											 									                    \TWChBx{7}{022}}%
                                                                                		      \kn{-3.4}%
                                                                               				  \hbox{\kn{1.235}%
										 										                    \vruleWH{.25pt}{7.2pt}}%
                                                                                			  \kn{-.3}%
                                                                                			  \hbox{\kn{-1.4}%
																				                    \TWChBx{7}{023}}%
                                                                                		\kn{-5}}}}%
                                                                        \def\Nr#1{\hbox{\kn{1}\Vbox{\hbox{\kn{-1.4}\Violet{\TWChBx{7}{023}}}%
                                                                                				    \kn{-3.4}%
                                                                                				    \hbox{\vruleWH{.25pt}{7.2pt}}%
                                                                                                    \kn{-.3}%
                                                                                                    \hbox{\kn{-.8}\Violet{\TWChBx{7}{022}}}%
                                                                                                    \kn{-5}}%
                                                                          				\kn{-1.2}%
                                                                          				\Vbox{\hbox{\kn{1}#1}}%
                                                                          					  \kn{1.5}%
                                                                          					  \Vbox{\hbox{\kn{-.8}%
											 									                          \Violet{\TWChBx{7}{022}}}%
                                                                                					\kn{-3.4}%
                                                                               						\hbox{\kn{1.235}%
										 										                          \vruleWH{.25pt}{7.2pt}}%
                                                                                					\kn{-.3}%
                                                                                					\hbox{\kn{-1.4}%
																				                          \Violet{\TWChBx{7}{023}}}%
                                                                                					\kn{-5}}}}%
                                                              			\def\pR#1{\hbox{\kn{1}%
                                                                                        \Vbox{\hbox{\kn{-1}\TWChBx{5}{023}}%
                                                                                			  \kn{-2.4}%
                                                                                			  \hbox{\vruleWH{.25pt}{5.5pt}}%
                                                                                              \kn{-.2}%
                                                                                              \hbox{\kn{-.4}\TWChBx{5}{022}}%
                                                                                              \kn{-3.7}}%
                                                                          				\kn{-2}%
                                                                          				\Vbox{\hbox{\kn{1}#1}}%
                                                                          				\kn{1}%
                                                                                        \Vbox{\hbox{\kn{-1.6}\TWChBx{5}{022}}%
                                                                                			  \kn{-2.4}%
                                                                                			  \hbox{\vruleWH{.25pt}{5.5pt}}%
                                                                                              \kn{-.2}%
                                                                                              \hbox{\kn{-2}\TWChBx{5}{023}}%
                                                                                		      \kn{-3.7}}}}%
                                                                        %
                  %
                    \def\CUT{\lower0pt\hbox{\kn{.5}\AiChBx{145}\kn{.5}}}
                    \def\NCUT{\lower0pt\hbox{\kn{.5}\AiChBx{164}\kn{.5}}}%
                    \def\Def#1{$\hbox{\Gruen{\Bf{#1}}}$}%
                    \def\Ton#1{\lower1pt\hbox{\BrownHBox{.5pt}{.5pt}{#1}}}%
                    \def\Sing#1{\Lbrace#1\Rbrace}
                    \def\Id#1{\OpGr{i}$_{\RMBx{7}{#1}}$}
                    \def\Order{\Blue{\MiiChBx{36}}}\def\oRder{\Blue{\MIIChBx{7}{36}}}\def\NOrder{\Order\hskip-6pt/}
                    \def\OrderShriek{\hbox{\Stack{\BF{7}{!}}{\Order}{-1pt}}}
                    \def\ORder{\Stack{\Blue{\Prime}}{\Order}{-5pt}}
                    \def\Restriction#1{\thinspace\vrulewhd{1}{10}{1}\thinspace\lower4pt\hbox{\RM{7}{#1}}}
                    \def\Period{\thinspace.}%
                    \def\Semicolon{\thinspace;}%
                    \def\Comma{\thinspace,}%
                    \def\Colon{\thinspace:}%
\def\CInt#1#2{\setbox3\hbox{\Vbox{\Brown{\hrulehw{1}{1}\kn{2}\hrulehw{1}{1}\kn{2}\hrulehw{1}{1}}}}%
              \setbox1\hbox{#1}\setbox2\hbox{#2}\dimen241\wd1\dimen242\wd2\dimen243\wd3%
              \divide\dimen241by2\divide\dimen242by2\divide\dimen243by2%
              \ifdim\dimen242>\dimen241 \dimen240=\dimen242\advance\dimen240by-\dimen241\dimen241=\dimen240 \dimen240=\dimen242\advance\dimen240by-\dimen243\dimen243=\dimen240\dimen242=0pt%
                                  \else \dimen240=\dimen241\advance\dimen240by-\dimen242\dimen242=\dimen240 \dimen240=\dimen241\advance\dimen240by-\dimen243\dimen243=\dimen240\dimen241=0pt\fi%
              \setbox1\hbox{\kern\dimen241\box1\kern\dimen241}\setbox2\hbox{\kern\dimen242\box2\kern\dimen242}\setbox3\hbox{\kern\dimen243\box3\kern\dimen243}%
              \setbox4\hbox{\Vbox{\box2\kn{1}\box3\kn{1}\box1}}\dimen240\ht4\advance\dimen240by-6pt\divide\dimen240by2\dimen241\ht4\advance\dimen241by2pt%
              \baselineskip=12pt\lineskip=1pt\lineskiplimit=0pt%
              \setbox5\vbox{\hbox{\vrule height1pt width3pt}\kern-12pt\kern\dimen241\hbox{\vrule depth1pt width3pt}\kern-2pt}%
              \setbox5\hbox{\box5}%
              \setbox4\hbox{\vrule height\dimen241 depth2pt width1pt\copy5\box4\box5\vrule height\dimen241 depth2pt width1pt}%
              \setbox4\hbox{\Stack{\oRder}{\box4}{1pt}}%
              \setbox4\hbox{\kn{1}\lower\dimen240\box4\kn{1}}%
              \box4}%
\def\OInt#1#2{\setbox3\hbox{\Vbox{\Brown{\hrulehw{1}{1}\kn{2}\hrulehw{1}{1}\kn{2}\hrulehw{1}{1}}}}%
              \setbox1\hbox{#1}\setbox2\hbox{#2}\dimen241\wd1\dimen242\wd2\dimen243\wd3%
              \divide\dimen241by2\divide\dimen242by2\divide\dimen243by2%
              \ifdim\dimen242>\dimen241 \dimen240=\dimen242\advance\dimen240by-\dimen241\dimen241=\dimen240 \dimen240=\dimen242\advance\dimen240by-\dimen243\dimen243=\dimen240\dimen242=0pt%
                                  \else \dimen240=\dimen241\advance\dimen240by-\dimen242\dimen242=\dimen240 \dimen240=\dimen241\advance\dimen240by-\dimen243\dimen243=\dimen240\dimen241=0pt\fi%
              \setbox1\hbox{\kern\dimen241\box1\kern\dimen241}\setbox2\hbox{\kern\dimen242\box2\kern\dimen242}\setbox3\hbox{\kern\dimen243\box3\kern\dimen243}%
              \setbox4\hbox{\Vbox{\box2\kn{1}\box3\kn{1}\box1}}\dimen240\ht4\advance\dimen240by-6pt\divide\dimen240by2\dimen241\ht4\advance\dimen241by2pt%
              \baselineskip=12pt\lineskip=1pt\lineskiplimit=0pt%
              \setbox5\vbox{\hbox{\vrule height1pt width3pt}\kern-12pt\kern\dimen241\hbox{\vrule depth1pt width3pt}\kern-2pt}%
              \setbox4\hbox{\copy5\vrule height\dimen241 depth2pt width1pt\kn{1}\box4\kn{1}\vrule height\dimen241 depth2pt width1pt\box5}%
              \setbox4\hbox{\Stack{\oRder}{\box4}{1pt}}%
              \setbox4\hbox{\kn{1}\lower\dimen240\box4\kn{1}}%
              \box4}%
\def\Mik#1#2{\setbox3\hbox{\Vbox{\Brown{\hrulehw{6}{2}}}}%
              \setbox1\hbox{#1}\setbox2\hbox{#2}\dimen241\wd1\dimen242\wd2\dimen243\wd3%
              \divide\dimen241by2\divide\dimen242by2\divide\dimen243by2%
              \ifdim\dimen242>\dimen241 \dimen244=\dimen242\dimen240=\dimen242\advance\dimen240by-\dimen241\dimen241=\dimen240 \dimen240=\dimen242\advance\dimen240by-\dimen243\dimen243=\dimen240\dimen242=0pt%
                                  \else \dimen244=\dimen241\dimen240=\dimen241\advance\dimen240by-\dimen242\dimen242=\dimen240 \dimen240=\dimen241\advance\dimen240by-\dimen243\dimen243=\dimen240\dimen241=0pt\fi%
              \setbox1\hbox{\kern\dimen241\box1\kern\dimen241}\setbox2\hbox{\kern\dimen242\box2\kern\dimen242}\setbox3\hbox{\kern\dimen243\box3\kern\dimen243}%
              \multiply\dimen244by2\advance\dimen244by4pt%
              \setbox4\hbox{\Vbox{\box2%
                                  \kn{1}%
                                  \box3%
                                  \hbox{\kn{-2}\Brown{\vruleHW{2pt}{\dimen244}}\kn{-2}}%
                                  \kn{1}%
                                  \box1}}%
              \dimen240\ht4\advance\dimen240by-6pt\divide\dimen240by2\dimen241\ht4\advance\dimen241by2pt%
              \baselineskip=12pt\lineskip=1pt\lineskiplimit=0pt%
              \setbox5\vbox{\hbox{\vrule height1pt width3pt}\kern-12pt\kern\dimen241\hbox{\vrule depth1pt width3pt}\kern-2pt}%
              \setbox4\hbox{\vrule height\dimen241 depth2pt width1pt\copy5\box4\box5\vrule height\dimen241 depth2pt width1pt}%
              \setbox4\hbox{\Stack{\oRder}{\box4}{1pt}}%
              \setbox4\hbox{\kn{1}\lower\dimen240\box4\kn{1}}%
              \box4}%
\def\Mak#1#2{\setbox3\hbox{\Vbox{\Brown{\hrulehw{6}{2}}}}%
              \setbox1\hbox{#1}\setbox2\hbox{#2}\dimen241\wd1\dimen242\wd2\dimen243\wd3%
              \divide\dimen241by2\divide\dimen242by2\divide\dimen243by2%
              \ifdim\dimen242>\dimen241 \dimen244=\dimen242\dimen240=\dimen242\advance\dimen240by-\dimen241\dimen241=\dimen240 \dimen240=\dimen242\advance\dimen240by-\dimen243\dimen243=\dimen240\dimen242=0pt%
                                  \else \dimen244=\dimen241\dimen240=\dimen241\advance\dimen240by-\dimen242\dimen242=\dimen240 \dimen240=\dimen241\advance\dimen240by-\dimen243\dimen243=\dimen240\dimen241=0pt\fi%
              \setbox1\hbox{\kern\dimen241\box1\kern\dimen241}\setbox2\hbox{\kern\dimen242\box2\kern\dimen242}\setbox3\hbox{\kern\dimen243\box3\kern\dimen243}%
              \advance\dimen244by8pt%
              \setbox4\hbox{\Vbox{\box2%
                                  \kn{1}%
                                  \hbox{\kn{-2.5}\Brown{\vruleHW{2pt}{\dimen244}}\kn{-3}}%
                                  \box3%
                                  \kn{1}%
                                  \box1}}%
              \dimen240\ht4\advance\dimen240by-6pt\divide\dimen240by2\dimen241\ht4\advance\dimen241by2pt%
              \baselineskip=12pt\lineskip=1pt\lineskiplimit=0pt%
              \setbox5\vbox{\hbox{\vrule height1pt width3pt}\kern-12pt\kern\dimen241\hbox{\vrule depth1pt width3pt}\kern-2pt}%
              \setbox4\hbox{\vrule height\dimen241 depth2pt width1pt\copy5\box4\box5\vrule height\dimen241 depth2pt width1pt}%
              \setbox4\hbox{\Stack{\oRder}{\box4}{1pt}}%
              \setbox4\hbox{\kn{1}\lower\dimen240\box4\kn{1}}%
              \box4}%
                    \def\Restrict#1#2{\hbox{#1\kn{1}\vrulehwd{10}{.75}{3}\kn{1}\lower3pt\hbox{#2}}}%
                    \def\Period{\thinspace.}\def\Semicolon{\thinspace;}\def\Comma{\thinspace,}\def\Colon{\thinspace:}%
                    \def\BlackSquare{\AIChBx{12}{4}}\def\bLackSquare{\AIChBx{5}{4}}%
                    \def\GreenSquare{\Green{\AIChBx{12}{4}}}\def\gReenSquare{\Green{\AIChBx{5}{4}}}%
                    \def\RedSquare{\Red{\AIChBx{12}{4}}}\def\rEdSquare{\Red{\AIChBx{5}{4}}}%
                    \def\Hex{\MIICh{12}{65}\kn{-10.75}\MIICh{12}{64}}%
                    \def\HexII{\hbox{\Hex\kn{-5.8}\raise3.1pt\hbox{\Blue{\AIChBx{5}{116}}}\kn{-2.9}\raise3.1pt\hbox{\Blue{\AIChBx{5}{110}}}}}%
                    \def\HexIV{\hbox{\Hex\kn{-5.8}\raise.2pt\hbox{\Red{\AIChBx{5}{110}}}\kn{-2.9}\raise.2pt\hbox{\Red{\AIChBx{5}{116}}}\kn{4.8}}}%
                    \def\HexVIII{\hbox{\Hex\kn{-10.65}\raise.2pt\hbox{\Red{\AIChBx{5}{116}}}\kn{-2.9}\raise.2pt\hbox{\Red{\AIChBx{5}{110}}}\kn{4.8}}}%
                    \def\HexXII{\hbox{\Hex\kn{-7.45}\raise3pt\hbox{\Red{\Vbox{\AIChBx{5}{116}%
                                                                                \AIChBx{5}{110}}}}\kn{3}}}%
                    \def\HexX{\hbox{\Hex\kn{-10.5}\raise3.1pt\hbox{\Blue{\AIChBx{5}{110}}}\kn{-2.9}\raise3pt\hbox{\Blue{\AIChBx{5}{116}}}\kn{4.6}}}%
                    \def\HexVI{\hbox{\Hex\kn{-7.45}\lower2.5pt\hbox{\Blue{\Vbox{\AIChBx{5}{116}%
                                                                               \AIChBx{5}{110}}}}\kn{3}}}%
                    \def\MidBox#1#2{\Dimen=#1\setbox\Box\hbox{#2}\DimenArg=-\wd\Box\dimen242=-\ht\Box\advance\dimen242by-\dp\Box%
                                    \advance\DimenArg by\Dimen\divide\DimenArg by2\advance\DimenVarby\Dimen\divide\DimenVarby2%
                                    \hbox{\kern\DimenArg\Vbox{\kern\DimenVar\box\Box\kern\DimenVar}\kern\DimenArg}}%
                    \def\HeX#1#2#3#4#5#6{\setbox\BoxI\hbox{\MidBox{10pt}{#1}}%
                                         \setbox\BoxII\hbox{\MidBox{10pt}{#2}}%
                                         \setbox\BoxIII\hbox{\MidBox{10pt}{#3}}%
                                         \setbox\BoxIV\hbox{\MidBox{10pt}{#4}}%
                                         \setbox\BoxV\hbox{\MidBox{10pt}{#5}}%
                                         \setbox\BoxVI\hbox{\MidBox{10pt}{#6}}%
                                                 \lower13pt\hbox{\Vbox{\hbox{\kn{10}\box\BoxVI}%
                                                             \kn{-3}%
                                                             \hbox{\box\BoxI\kn{10}\box\BoxV}%
                                                             \hbox{\box\BoxII\kn{10}\box\BoxIV}%
                                                             \kn{-3}%
                                                             \hbox{\kn{10}\box\BoxIII}}%
                                                             \kn{-20}\raise13pt\hbox{\Hex}\kn{10}}}%
																		\def\Span#1{\hbox{\setbox\Box=\hbox{#1}%
                                                                                   		  \Vbox{\hrule height1pt%
                                                                                         		\hboxII{3pt}{\wd\Box}{\vrulewd{1}{2}\hfil\vrulewd{1}{2}}%
                                                                                         		\kn{1}%
                                                                                         	    \hboxII{3pt}{\wd\Box}{\kn{1.5}\box\Box\kn{1.5}}}}}%
  \def\Vol#1{\lower4pt\hbox{\kn{1}\Vbox{\kn{1}\BlueHBox{1pt}{1pt}{\vruleh{7}#1}\kn{1}}\kn{.5}}}%
  \def\vOl#1{\lower4pt\hbox{\kn{1}\BlueHBox{1pt}{1pt}{\vruleh{5}#1}\kn{.5}}}%
  \def\Hol#1#2#3#4#5#6{\setbox\BoxEnd\BlueHBox{.5pt}{1pt}{\ArrayIIxIII{\RM{7}{#1}}{\RM{7}{#2}}{\RM{7}{#3}}{\RM{7}{#4}}{\RM{7}{#5}}{\RM{7}{#6}}}%
                       \setbox\Box\hbox{\kn{1}\box\BoxEnd\kn{1}}%
                       \Dimen=\ht\Box\divide\Dimen by2\advance\Dimen by -2.5pt%
                       \hbox{\lower\Dimen\copy\Box}}%
  \def\PB#1{\Purple{\Bf{#1}}}\def\OB#1{\RedOrange{\Bf{#1}}}%
  \def\pB#1{\Purple{\BFBx{7}{#1}}}\def\oB#1{\RedOrange{\BFBx{7}{#1}}}%
  \def\GacTurn{\OpGr{r}\kn{0}$_{\Green{\BFBx{7}{(}}\pB{a}\oB{c}\Green{\BFBx{7}{)}}}$}\def\GACTURN{\hbox{\Stack{\sIm}{\OpGr{r}}{1pt}$_{\Green{\BFBx{7}{(}}\pB{a}\oB{c}\Green{\BFBx{7}{)}}}$}}%
  \def\GadTurn{\OpGr{r}\kn{0}$_{\Green{\BFBx{7}{(}}\pB{a}\oB{d}\Green{\BFBx{7}{)}}}$}%
  \def\RabTurn{\OpGr{r}\kn{0}$_{\Red{\BFBx{7}{(}}\pB{a}\oB{b}\Red{\BFBx{7}{)}}}$}\def\RABTURN{\hbox{\Stack{\sIm}{\OpGr{r}}{1pt}$_{\Red{\BFBx{7}{(}}\pB{a}\oB{b}\Red{\BFBx{7}{)}}}$}}%
  \def\RacTurn{\OpGr{r}\kn{0}$_{\Red{\BFBx{7}{(}}\pB{a}\oB{c}\Red{\BFBx{7}{)}}}$}%
  \def\BadTurn{\OpGr{r}\kn{0}$_{\Black{\BFBx{7}{(}}\pB{a}\oB{d}\Black{\BFBx{7}{)}}}$}\def\BADTURN{\hbox{\Stack{\sIm}{\OpGr{r}}{1pt}$_{\Black{\BFBx{7}{(}}\pB{a}\oB{d}\Black{\BFBx{7}{)}}}$}}%
  \def\BabTurn{\OpGr{r}\kn{0}$_{\Black{\BFBx{7}{(}}\pB{a}\oB{b}\Black{\BFBx{7}{)}}}$}%
  \def\GInv{\OpGr{s}$_{\gReenSquare}$}\def\GINV{\hbox{\Stack{\sIm}{\OpGr{s}}{1pt}$_{\gReenSquare}$}}%
  \def\BInv{\OpGr{s}$_{\bLackSquare}$}\def\BINV{\hbox{\Stack{\sIm}{\OpGr{s}}{1pt}$_{\bLackSquare}$}}%
  \def\RInv{\OpGr{s}$_{\rEdSquare}$}\def\RINV{\hbox{\Stack{\sIm}{\OpGr{s}}{1pt}$_{\rEdSquare}$}}%
  \def\ACInv{\OpGr{e}\kn{0}$_{\Green{\BFBx{7}{(}}\pB{a}\oB{c}\Red{\BFBx{7}{)}}}$}%
  \def\ADInv{\OpGr{e}\kn{0}$_{\Green{\BFBx{7}{(}}\pB{a}\oB{d}\Black{\BFBx{7}{)}}}$}%
  \def\ABInv{\OpGr{e}\kn{0}$_{\Red{\BFBx{7}{(}}\pB{a}\oB{b}\Black{\BFBx{7}{)}}}$}%
  \def\BCInv{\OpGr{e}\kn{0}$_{\Green{\BFBx{7}{(}}\pB{b}\oB{c}\Black{\BFBx{7}{)}}}$}%
  \def\BDInv{\OpGr{e}\kn{0}$_{\Red{\BFBx{7}{(}}\pB{b}\oB{d}\Green{\BFBx{7}{)}}}$}%
  \def\CDInv{\OpGr{e}\kn{0}$_{\Red{\BFBx{7}{(}}\pB{c}\oB{d}\Black{\BFBx{7}{)}}}$}%
  \def\AVInv{\hbox{\OVERBAR{.5}{1}{\PB{a},\OB{a}}}}\def\AvInv{\hbox{\OVERBAR{.5}{1}{\pB{a},\oB{a}}}}%
  \def\BVInv{\hbox{\OVERBAR{.5}{1}{\PB{b},\OB{b}}}}\def\BvInv{\hbox{\OVERBAR{.5}{1}{\pB{b},\oB{b}}}}%
  \def\CVInv{\hbox{\OVERBAR{.5}{1}{\PB{c},\OB{c}}}}\def\CvInv{\hbox{\OVERBAR{.5}{1}{\pB{c},\oB{c}}}}%
  \def\DVInv{\hbox{\OVERBAR{.5}{1}{\PB{d},\OB{d}}}}%
  \def\AVbcInv{\hbox{\OpGr{th}$_{\BFBx{7}{(}\pB{a}\BFBx{7}{;}\pB{b}\pB{c}\BFBx{7}{)}}$}}\def\AVbcINV{\hbox{\Stack{\sIm}{\OpGr{th}}{1pt}$_{\BFBx{7}{(}\pB{a}\BFBx{7}{;}\pB{b}\pB{c}\BFBx{7}{)}}$}}%
  \def\AVcbInv{\hbox{\OpGr{th}$_{\BFBx{7}{(}\pB{a}\BFBx{7}{;}\pB{c}\pB{b}\BFBx{7}{)}}$}}\def\AVcbINV{\hbox{\Stack{\sIm}{\OpGr{th}}{1pt}$_{\BFBx{7}{(}\pB{a}\BFBx{7}{;}\pB{c}\pB{b}\BFBx{7}{)}}$}}%
  \def\BVcdInv{\hbox{\OpGr{th}$_{\BFBx{7}{(}\pB{b}\BFBx{7}{;}\pB{c}\pB{d}\BFBx{7}{)}}$}}%
  \def\BVdcInv{\hbox{\OpGr{th}$_{\BFBx{7}{(}\pB{b}\BFBx{7}{;}\pB{d}\pB{c}\BFBx{7}{)}}$}}%
  \def\CVdaInv{\hbox{\OpGr{th}$_{\BFBx{7}{(}\pB{c}\BFBx{7}{;}\pB{d}\pB{a}\BFBx{7}{)}}$}}%
  \def\CVadInv{\hbox{\OpGr{th}$_{\BFBx{7}{(}\pB{c}\BFBx{7}{;}\pB{a}\pB{d}\BFBx{7}{)}}$}}%
  \def\DVabInv{\hbox{\OpGr{th}$_{\BFBx{7}{(}\pB{d}\BFBx{7}{;}\pB{a}\pB{b}\BFBx{7}{)}}$}}%
  \def\DVbaInv{\hbox{\OpGr{th}$_{\BFBx{7}{(}\pB{d}\BFBx{7}{;}\pB{b}\pB{a}\BFBx{7}{)}}$}}%
  \def\GSquare{\hbox{\lower2pt\hbox{\Stack{\sIm}{\GreenSquare}{1pt}}}}%
  \def\BSquare{\hbox{\lower2pt\hbox{\Stack{\sIm}{\BlackSquare}{1pt}}}}%
  \def\RSquare{\hbox{\lower2pt\hbox{\Stack{\sIm}{\RedSquare}{1pt}}}}%
  \def\wurf#1#2#3#4#5#6#7#8{\setbox\BoxEnd\hbox{\raise12pt\hbox{\Op{\MIIIChBx{7}{104}}}{\ArrayIIxIV{\RM{7}{#1}}{\RM{7}{#2}}{\RM{7}{#3}}{\RM{7}{#4}}{\RM{7}{#5}}{\RM{7}{#6}}{\RM{7}{#7}}{\RM{7}{#8}}}%
                                                \raise12pt\hbox{\Op{\MIIIChBx{7}{105}}}}%
                            \setbox\Box\hbox{\kn{1}\box\BoxEnd\kn{1}}%
                            \Dimen=\ht\Box\divide\Dimen by2%
                            \hbox{\lower\Dimen\copy\Box}}%
  \def\Wurf#1{\Blue{\BRAK{10}}#1\Blue{\KET{10}}}%
  \def\Wuerf#1#2#3{\setbox\BoxX\hbox{1}\DimenHII=4pt\DimenHI=2pt%
                   \setbox\Box\hbox{\Stack{\Op{\Sim}}{\BlueHBox{.5pt}{1pt}{\ArrayIIxIII{\RMBx{7}{#1}}{\RMBx{7}{#2}}{\RMBx{7}{#3}}{\AvInv}{\BvInv}{\CvInv}}}{1pt}}%
                   \Dimen=\ht\Box\divide\Dimen by2\advance\Dimen by -5.5pt%
                   \lower\Dimen\copy\Box}%
  \def\InvWuerf#1#2#3{\setbox\Box\hbox{\Wuerf{#1}{#2}{#3}}\hbox{\box\Box\kn{2}\Vbox{\RMBx{7}{-1}\kn{12}}}}%
  \def\GTriangle{\Stack{\sIm}{\hbox{\lower2pt\hbox{\Green{\AiChBx{110}}}}}{1pt}}%
  \def\RTriangle{\Stack{\sIm}{\hbox{\lower2pt\hbox{\Red{\AiChBx{110}}}}}{1pt}}%
  \def\BTriangle{\Stack{\sIm}{\hbox{\lower2pt\hbox{\AiChBx{110}}}}{1pt}}%
  \def\LibraP#1#2{\hbox{$_{\RMBx{7}{#1}}$\OpGr{p}$_{\RMBx{7}{#2}}$}}%
  \def\LibraL#1#2{\hbox{$_{\RMBx{7}{#1}}$\OpGr{l}$_{\RMBx{7}{#2}}$}}%
  \def\LibraR#1#2{\hbox{$_{\RMBx{7}{#1}}$\OpGr{r}$_{\RMBx{7}{#2}}$}}%
  \def\LibraF#1#2{\hbox{$_{\RMBx{7}{#1}}$\OpGr{f}$_{\RMBx{7}{#2}}$}}%
  \def\LibraTH#1#2{\hbox{$_{\RMBx{7}{#1}}$\OpGr{th}$_{\RMBx{7}{#2}}$}}%
  \def\GIn{\OpGr{p}$_{\gReenSquare}$}%
  \def\BIn{\OpGr{p}$_{\bLackSquare}$}%
  \def\RIn{\OpGr{p}$_{\rEdSquare}$}%
  \def\Load{load}\def\Loads{M$^{\ESBx{7}{F}}$}%
  \def\WWurf#1{\Blue{\BBRAK{10}}#1\Blue{\KKET{10}}}%
  \def\Quin#1#2#3#4#5{\setbox\Box\hbox{\kn{1}\BlueHBox{1pt}{2pt}{\ArrayIIIxIII{}{#3}{}{#1}{#5}{#2}{}{#4}{}}\kn{1}}%
                      \Dimen=\ht\Box%
                      \divide\Dimen by 2%
                      \advance\Dimen by -1pt%
                      \lower\Dimen\box\Box}%
  \def\BLoads{\RuChBx{27}}%
  \def\MLL#1#2#3{\hbox{\LL{#3}\lower5pt\hbox{\Vbox{\RMBx{7}{#1}\kn{7}\RMBx{7}{#2}}}}}%
  \def\MAB#1#2{\Stack{\hbox{\RMBx{7}{#1}\MIIChBx{7}{44}\RMBx{7}{#2}}}{\MiiBx{M}}{1pt}}%
  \def\Mab#1#2{\hbox{\RmBx{M}$_{\RMBx{7}{(#1,#2})}$}}%
  \def\MABF#1#2{\setbox\BoxArg\hbox{\MiiBx{M}\vrulehdw{1}{3}{4}}\Stack{\hbox{\RMBx{7}{#1}\MIIChBx{7}{44}\RMBx{7}{#2}}}{\box\BoxArg}{1pt}}%
  \def\MNF{\MiiBx{M}$_{\AIChBx{5}{3}}$\kn{1}}%
  \def\Mabf#1#2#3{\hbox{\Rm{M}\lower2pt\hbox{\pR{\RMBx{7}{#1},\RMBx{7}{#2},\RMBx{7}{#3}}}}}%
  \def\COreo{\AiChBx{10}}\def\CCOreo{\AiChBx{11}}%
  \def\Arc#1{\setbox\Box\hbox{#1}%
             \Dimen=\wd\Box\advance\Dimen by-6pt%
             \hbox{\vruleh{11}\kn{2}\Vbox{\hbox{\MiiiChBx{172}\Vbox{\hrule width\Dimen height 1.2pt\kn{0}}\MiiiChBx{173}}\kn{2}\hbox{\kn{2}\box\Box}}\kn{2}}}%
  \def\CArc#1{\setbox\Box\hbox{#1}%
             \Dimen=\wd\Box\advance\Dimen by2pt%
             \hbox{\vruleh{11}\Vbox{\hbox{\MiiiChBx{76}\kn{-5}\Vbox{\hrule width\Dimen height 1.2pt\kn{0}}\kn{-5}\MiiiChBx{76}}\kn{-1}\hbox{\kn{5}\box\Box}}\kn{-2}}}%
 \def\BZero#1{\Stack{\hbox{\Op{\RM{5}{#1}}}}{\EUBx{8}{0}}{1pt}}%
 \def\BOne#1{\Stack{\hbox{\Op{\RM{5}{#1}}}}{\EUBx{8}{1}}{1pt}}%
 \def\BInfty#1{\Stack{\hbox{\Op{\RM{5}{#1}}}}{\Infty}{1pt}}%
 \def\BPlus#1{\Stack{\hbox{\Op{\RM{5}{#1}}}}{\RmBx{+}}{1pt}}%
 \def\BMinus#1{\kn{1}\Stack{\hbox{\Op{\RM{5}{#1}}}}{\RmBx{--}}{1pt}\kn{1}}%
 \def\BTimes#1{\kn{1}\Stack{\hbox{\Op{\RM{5}{#1}}}}{\Cdot}{1pt}\kn{1}}%
 \def\BOVER#1#2#3#4#5#6{\let\Top=\Box \let\Bottom=\BoxArg \let\Frac=\BoxVar%
                        \setbox\Top\hbox{#1}%
					    \setbox\Bottom\hbox{#2}%
                        \ifdim\wd\Top<\wd\Bottom \addII{-\wd\Top}{\wd\Bottom}%
												 \divideII{\DimenReturn}{2}%
                                                 \setbox\Frac\hbox{\Vbox{\hbox{\kern\DimenReturn%
															              	   \box\Top}%
																		 \kern#5%
																		 #4{\hrule height#3}%
																	     \kern#5%
																		 \copy\Bottom}%
                                                                    \kn{1}%
                                                                    \Vbox{\hbox{\Op{\RMBx{5}{#6}}}%
                                                                          \kn{-2}%
                                                                          \kern#5%
                                                                          \kern\ht\Bottom}}%
                                           \else \addII{-\wd\Bottom}{\wd\Top}%
												 \divideII{\DimenReturn}{2}%
                                                 \setbox\Frac\hbox{\Vbox{\box\Top%
															     		 \kern#5%
																		 #4{\hrule height#3}%
																		 \kern#5%
																		 \hbox{\kern\DimenReturn%
																			   \copy\Bottom}}%
                                                                    \kn{1}%
                                                                    \Vbox{\hbox{\Op{\RMBx{5}{#6}}}%
                                                                          \kn{-2}%
                                                                          \kern#5%
                                                                          \kern\DimenReturn%
                                                                          \kern\ht\Bottom}}\fi%
                        \addII{\ht\Frac}{\dp\Frac}%
						\divideII{\DimenReturn}{2}%
						\addII{\DimenReturn}{-4pt}%
                        \raise\DimenReturn\hbox{\box\Frac}}%
 \def\BOver#1#2#3{\BOVER{#1}{#2}{.5pt}{\color{black}}{1.5pt}{#3}}%
 \def\BDet#1#2#3#4#5{\Stack{\hbox{\Op{\RMBx{5}{#5}}}}{\Ft{det}}{1pt}$\pmatrix{#1&#2\cr#3&#4\cr}$}%
 \def\BField{F$_{\hbox{\Op{\RMBx{7}{b}}}}$}%
 \def\Rus#1#2{\RuBx{b}$_{\RMBx{7}{(#1;#2)}}$}%
 \def\Russ#1#2#3{\RuBx{b}$_{\RMBx{7}{(#1,#2;#3)}}$}%
 \def\Russs#1#2#3#4{\RuBx{b}$_{\RMBx{7}{(#1,#2,#3;#4)}}$}%
 \def\Conjugate#1{\setbox\Box\hbox{#1}\hbox{\copy\Box\Vbox{\RMBx{7}{*}\kern\ht\Box\kn{-3}}}}%
 \def\oRderIff{\hbox{\kn{2}\MIIChBx{7}{36}\kn{-7}\MIIChBx{7}{37}\kn{2}}}%
 \def\OrdDolar#1{\setbox\Box\hbox{#1}\hbox{\copy\Box\Vbox{\oRderIff\kern\ht\Box}}}%
 \def\Polar#1{\setbox\Box\hbox{#1}\hbox{\copy\Box\Vbox{\RMBx{7}{o}\kern\ht\Box}}}%
 \def\Pollar#1{\setbox\Box\hbox{#1}\hbox{\copy\Box\Vbox{\RMBx{7}{oo}\kern\ht\Box}}}%
 \def\Polllar#1{\setbox\Box\hbox{#1}\hbox{\copy\Box\Vbox{\RMBx{7}{ooo}\kern\ht\Box}}}%
 \def\Dolar#1{\setbox\Box\hbox{#1}\hbox{\copy\Box\Vbox{\BSChBx{6}{5}\kern\ht\Box}}}%
 \def\Dollar#1{\setbox\Box\hbox{#1}\hbox{\copy\Box\Vbox{\hbox{\BSChBx{6}{5}\BSChBx{6}{5}}\kern\ht\Box}}}%
 \def\Dolllar#1{\setbox\Box\hbox{#1}\hbox{\copy\Box\Vbox{\hbox{\BSChBx{6}{5}\BSChBx{6}{5}\BSChBx{6}{5}}\kern\ht\Box}}}%
 \def\Functions#1#2{\hbox{\hbox{\hbox{#2}\kn{.5}{\Vbox{\RMBx{7}{#1}\kn{5}}}}}}%
 \def\Bijections#1#2{\hbox{\Functions{#1}{#2}\kn{-2}\RM{12}{!}}}%
 \def\DTran#1{\setbox\Box\hbox{#1}%
              \Dimen=\ht\Box\advance\Dimen by\dp\Box\advance\Dimen by6pt%
              \DimenVar=\wd\Box\advance\DimenVar by2pt%
              \DimenO=\dp\Box\advance\DimenO by3pt%
              \setbox\Box\hbox{\kn{1}\Red{\vruleHW{\Dimen}{2pt}}%
                               \Vbox{\Red{\hruleHW{2pt}{\DimenVar}}%
                                     \kn{1}%
                                     \hbox{\kn{1}#1\kn{1}}%
                                     \kn{1}%
                                     \Red{\hruleHW{2pt}{\DimenVar}}}%
                               \Red{\vruleHW{\Dimen}{2pt}}\kn{1}}%
              \hbox{\lower\DimenO\box\Box}}%
 \def\LTran#1{\setbox\Box\hbox{#1}%
              \Dimen=\ht\Box\advance\Dimen by\dp\Box\advance\Dimen by6pt%
              \DimenVar=\wd\Box\advance\DimenVar by2pt\divide\DimenVar by2%
              \DimenO=\dp\Box\advance\DimenO by3pt%
              \setbox\Box\hbox{\kn{1}\Red{\vruleHW{\Dimen}{2pt}}%
                               \Vbox{\hbox{\Red{\vruleHW{2pt}{\DimenVar}}\Gray{\vruleHW{2pt}{\DimenVar}}}%
                                     \kn{1}%
                                     \hbox{\kn{1}#1\kn{1}}%
                                     \kn{1}%
                                     \hbox{\Red{\vruleHW{2pt}{\DimenVar}}\Gray{\vruleHW{2pt}{\DimenVar}}}}%
                               \Gray{\vruleHW{\Dimen}{2pt}}\kn{1}}%
              \hbox{\lower\DimenO\box\Box}}%
 \def\RTran#1{\setbox\Box\hbox{#1}%
              \Dimen=\ht\Box\advance\Dimen by\dp\Box\advance\Dimen by6pt%
              \DimenVar=\wd\Box\advance\DimenVar by2pt\divide\DimenVar by2%
              \DimenO=\dp\Box\advance\DimenO by3pt%
              \setbox\Box\hbox{\kn{1}\Gray{\vruleHW{\Dimen}{2pt}}%
                               \Vbox{\hbox{\Gray{\vruleHW{2pt}{\DimenVar}}\Red{\vruleHW{2pt}{\DimenVar}}}%
                                     \kn{1}%
                                     \hbox{\kn{1}#1\kn{1}}%
                                     \kn{1}%
                                     \hbox{\Gray{\vruleHW{2pt}{\DimenVar}}\Red{\vruleHW{2pt}{\DimenVar}}}}%
                               \Red{\vruleHW{\Dimen}{2pt}}\kn{1}}%
              \hbox{\lower\DimenO\box\Box}}%
 \def\LRTran#1{\setbox\Box\hbox{#1}%
              \Dimen=\ht\Box\advance\Dimen by\dp\Box\advance\Dimen by6pt%
              \DimenVar=\wd\Box\advance\DimenVar by2pt\divide\DimenVar by3%
              \DimenO=\dp\Box\advance\DimenO by3pt%
              \setbox\Box\hbox{\kn{1}\Red{\vruleHW{\Dimen}{2pt}}%
                               \Vbox{\hbox{\Red{\vruleHW{2pt}{\DimenVar}}\Gray{\vruleHW{2pt}{\DimenVar}}\Red{\vruleHW{2pt}{\DimenVar}}}%
                                     \kn{1}%
                                     \hbox{\kn{1}#1\kn{1}}%
                                     \kn{1}%
                                     \hbox{\Red{\vruleHW{2pt}{\DimenVar}}\Gray{\vruleHW{2pt}{\DimenVar}}\Red{\vruleHW{2pt}{\DimenVar}}}}%
                               \Red{\vruleHW{\Dimen}{2pt}}\kn{1}}%
              \hbox{\lower\DimenO\box\Box}}%
 \def\DiagTran#1{\setbox\Box\hbox{#1}%
                 \setbox\BoxArg\hbox{\raise\dp\Box\box\Box}%
                 \let\Height=\DimenO \Height=\ht\BoxArg \advance\Height by6pt%
                 \let\Numerator=\CountO \Numerator=3%
                 \let\Denominator=\CountI \Denominator=2%
                 \let\SmallTriBase=\DimenVar \SmallTriBase=1.2pt%
                 \let\HeightPlus=\DimenI \HeightPlus=\Height \advance\HeightPlus by1.7pt%
                 \let\Base=\DimenII \Base=\HeightPlus \divide\Base by\Numerator \multiply\Base by\Denominator%
                 \let\SideRule=\DimenIII \SideRule=\Base \advance\Base by-\SmallTriBase%
                 \let\SideFrontBox=\BoxO \setbox\SideFrontBox\hbox{\Red{\druleTWxy{2.26pt}{\Base}{2}{-3}}}%
                 \let\BoxFill=\DimenIV \BoxFill=\ht\BoxArg \advance\BoxFill by2pt%
                 \let\FudgeI=\DimenV \FudgeI=\SideRule \advance\FudgeI by-.5pt%
                 \let\FudgeII=\DimenVI \FudgeII=\SmallTriBase \advance\FudgeII by.5pt%
                 \let\LeftBackBox=\BoxI \setbox\LeftBackBox\hbox{\Vbox{\hbox{\kern\FudgeII\Gray{\vruleHW{2pt}{\FudgeI}}}%
                                                                       \kern\BoxFill%
                                                                       \hbox{\kern\FudgeI\Gray{\vruleHW{2pt}{\FudgeII}}}}}%
                 \advance\FudgeI by-.5pt%
                 \advance\FudgeII by.5pt%
                 \let\RightBackBox=\BoxII \setbox\RightBackBox\hbox{\Vbox{\hbox{\Gray{\vruleHW{2pt}{\FudgeII}}\kern\FudgeI}%
                                                                          \kern\BoxFill%
                                                                          \hbox{\Gray{\vruleHW{2pt}{\FudgeI}}\kern\FudgeII}}}%
                 \let\MiddleBox=\BoxIII \setbox\MiddleBox\hbox{\Vbox{\hbox{\Gray{\vruleHW{2pt}{\wd\BoxArg}}}%
                                                                     \kn{1}%
                                                                     \copy\BoxArg%
                                                                     \kn{1}%
                                                                     \hbox{\Gray{\vruleHW{2pt}{\wd\BoxArg}}}}}
                 \setbox\Box\hbox{\box\LeftBackBox\kern-\wd\SideFrontBox\copy\SideFrontBox\kn{-.35}\box\MiddleBox\kn{-.35}\box\RightBackBox\kern-\wd\SideFrontBox\box\SideFrontBox}%
                 \hbox{\lower3pt\box\Box}}%
 \def\Su#1{\setbox\BoxO\hbox{\RM{7}{#1}}%
           \hbox{\lower2pt\box\BoxO}}%
 \def\Schlange#1{\setbox\BoxVar\hbox{#1}\setbox\BoxArg\hbox{\TWCh{14}{176}}%
                 \DimenVar=\wd\BoxVar\DimenArg=\wd\BoxArg\DimenI=0pt\DimenII=0pt%
                 \ifdim\DimenVar<\DimenArg \DimenII=\DimenArg\advance\DimenII by-\DimenVar\divide\DimenI by2\fi%
                 \ifdim\DimenArg<\DimenVar \DimenI=\DimenVar\advance\DimenI by-\DimenArg\divide\DimenI by2\fi%
                 \hbox{\Vbox{\hbox{\kern\DimenI\box\BoxArg}\kn{-6}\hbox{\kern\DimenII\box\BoxVar}}}}%
 \def\Power#1#2{\setbox\BoxArg\hbox{#1}%
                   \Dimen=\ht\BoxArg%
                   \advance\Dimen by-1.82pt%
                   \hbox{\box\BoxArg%
                         \kn{.5}%
                         \Vbox{\RMBx{7}{#2}%
                               \kern\Dimen}%
                         \kn{1}}}%
 \def\SEQUENCE#1#2{\setbox\BoxArg\hbox{#1}%
                   \Dimen=\ht\BoxArg%
                   \advance\Dimen by-1.82pt%
                   \hbox{\box\BoxArg%
                         \kn{.5}%
                         \Vbox{\hbox{\Underbar{\RM{7}{#2}}}%
                               \kern\Dimen}}}%
 \def\Sequence#1{\setbox\BoxArg\hbox{#1}%
                 \Dimen=\ht\BoxArg%
                 \advance\Dimen by-1.82pt%
                 \hbox{\box\BoxArg%
                       \kn{.5}%
                       \Vbox{\AIIBx{7}{N}%
                             \kern\Dimen}}}%
\vskip2in
                                                                                      \D{8}{0}{\BF{16}{THE PROJECTIVE LINE AS A MERIDIAN}}
                                                                                      \D{16}{-4}{(by Kelly McKennon, May 2017)}
\vfill\eject
                                                                                       \Dc{\BFBx{14}{Prologue}}
\PAR This paper is written in \TeX\ \Tw{eplain}, for compilation using \Tw{pdf.tex}. This was done so that links and colors would be available without the restrictions of LaTex. Unfortunately, \Bf{arXiv}\ does not compile submissions in \Tw{pdf.tex}, and the author was unable to make the inner \Tw{eplain}\ links work using the \TeX\ compiler employed by \Bf{arXiv}. There is however a \Tw{pdf}\ version available on the \Bf{world wide web} where the links do work: specifically at www.iis.sinica.edu.tw/page/library/TechReport/tr2017/tr17002.pdf\thinspace.
\PaR There are some symbols used here which are not universally standard, all of which should be covered in the appendix. The most common are the following:
                                                                            \Dc{X\Cop Y}
\Par for the set complement of a subset Y of a set X;
                                                     \Dc{\Domain{\OpGr{f}}\Andd\Range{\OpGr{f}}\thinspace, respectively,}
\Par for the domain and range, respectively, of a function \OpGr{f};
                                             \Dc{\setRelation{\OpGr{f}}\Of{S}\ \Equiv\ \SetSuch{\OpGr{f}\Of{x}}{x\In S}}
\Par for the image of a subset S of the domain of a function \OpGr{f};
                    \Dc{\Underbar{n}\ \Equiv\ \Set{\One,\Two,\dots \Eu{n}}}
\Par for the set of the first n positive integers and
                                                                          \Dc{\Power{X}{Y}}
\Par for the family of all functions from a set Y to a set X.
\PaR Some of the material in these papers appeared originally in http://vixra.org/abs/1306.0233 , and several mistakes in that paper have been rectified here.
\PaR The author can be reached at KellyMack@Proton.Com .

\vfil\eject
                                                                                                      \def\xrefpageword{Page }%
                                                                                                 \d{\BF{14}{Table of Contents}}
  \PAr \line{\hfill\hbox to 5.5in{\LinkSection{P}{Introduction}\dotfill\xref{pageP}\ \thinspace}\hfill}
  \PAr \line{\hfill\hbox to 5.5in{\LinkSection{ED}{Erlanger Definition I}\dotfill\xref{pageED}\ \thinspace}\hfill}
  \PAr \line{\hfill\hbox to 5.5in{\LinkSection{E}{Erlanger Definition II}\dotfill\xref{pageE}}\hfill}
  \PAr \line{\hfill\hbox to 5.5in{\LinkSection{LL}{Involution Libras}\dotfill\xref{pageLL}}\hfill}
  \PAr \line{\hfill\hbox to 5.5in{\LinkSection{CL}{Meridians on the Cube}\dotfill\xref{pageCL}}\hfill}
  \PAr \line{\hfill\hbox to 5.5in{\LinkSection{W}{Wurfs and the Cross Ratio}\dotfill\xref{pageW}}\hfill}
  \PAr \line{\hfill\hbox to 5.5in{\LinkSection{EA}{Meridian Exponentials and Arcs}\dotfill\xref{pageEA}}\hfill}
  \PAr \line{\hfill\hbox to 5.5in{\LinkSection{SetTheory}{Appendix I: Mathematical Notation and Terminology}\dotfill\xref{pageSetTheory}}\hfill}
  \PAr \line{\hfill\hbox to 5.5in{\LinkSection{Graphs}{Appendix II: Graphs}\dotfill\xref{pageGraphs}}\hfill}
  \PAr \line{\hfill\hbox to 5.5in{\LinkSection{Topology}{Appendix III: Topology}\dotfill\xref{pageTopology}}\hfill}
  \PAr \line{\hfill\hbox to 5.5in{\LinkSection{Notation}{Notation}\dotfill\xref{pageNotation}}\hfill}
  \PAr \line{\hfill\hbox to 5.5in{\LinkSection{Index}{Index}\dotfill\xref{pageIndex}}\hfill}
  \PAr \line{\hfill\hbox to 5.5in{\LinkSection{Eponymy}{Eponymy}\dotfill\xref{pageEponymy}}\hfill}
  \PAr \line{\hfill\hbox to 5.5in{\LinkSection{Acknowledgement}{Acknowledgement}\dotfill\xref{pageAcknowledgement}}\hfill}
  \PAr \line{\hfill\hbox to 5.5in{\LinkSection{Bibliography}{Bibligography}\dotfill\xref{pageBibliography}}\hfill}
  \vfill\eject
                                                                                                   \Section{P}{Introduction}\xrdef{pageP}
\Item{P}{Purpose} The purpose of the present article is to examine the essence of what has commonly been described as a \Quotes{projective line}, but which is here named a \Quotes{meridian}. This shall be done in several papers: this first paper devoted to the meridian itself, the second to the character and form of the family of projective isomorphisms of one meridian onto another and the third to some connections between meridians and higher dimensional projective space.
\PaR In this first paper we shall view the meridian from various points of view:
\Paragraph{0pt}{10pt}(1) as a set acted upon by a family of involutions;
\Paragraph{0pt}{10pt}(2) as a set acted upon by a 3-transitive group of permutations;
\Paragraph{0pt}{10pt}(3) as a set with a quinary operator;
\Paragraph{0pt}{10pt}(4) as an equivalence class of quadruples, relating to the cross ratio.
\PaR In the final section of this first paper we shall investigate how the existence of a certain single-valued exponential on a meridian is characteristic of the meridian corresponding to the field of real numbers.
\PaR Most of the terminology applied here is standard, but not all. Both standard and non-standard terminology is detailed in the \LinkSection{SetTheory}{appendix}. This appendix is somewhat more encompassing than necessary for what is required here, so that it can serve for the sequel as well. There is also an \LinkSection{Notation}{index for notation} as well as the \LinkSection{Index}{index for terminology}.
\Item{H}{Historical} The first considerations of perspective are perhaps coeval with the development of man's eyesight, and records of such are traceable to antiquity.\Foot{A famous projective theorem concerning hexagons was published by Pappus of Alexandria in the first half of the fourth century AD.}
\PaR Some of the Renaissance painters and architects tried to organize what was then known about perspective to aid them in their art: \It{exempli gratia}, Filipo Brunelleshi (1337-1446), Leone Battista Alberti (1404-1472), Piero della Francesca (1410-1492), Leonardo Da Vinci (1452-1519) and Albrecht D\"urer (1471-1528).
\PaR A crucial intellectual step in understanding perspective is the conception of points at \Quotes{infinity}. This surfaced in the early seventeenth century, and manifested itself, probably independently, through the minds of Johannes Kepler (1571-1630) and G\'erard Desargues (1591-1661). The work of the latter inspired the young Blaise Pascal (1623-1669) to write in 1639 a significant treatise on projective geometry, the mathematical formalization of the ideas of perspective.
\PaR It was in the early nineteenth century, along with the rapid and accelerating progress of nearly all science at the time, that interest and understanding in projective geometry again surfaced: \It{exempli gratia}, through Gaspard Monge (1746-1818) at the turn of the century, Jean-Victor Poncelet (1788-1867) in 1822 and Jakob Steiner (1796-1863) in 1833.
\PaR The interest in perspective, by its very nature, is geometrical. Another mathematical interest, of which the published roots are traceable even further back into antiquity than the study of perspective, is the solution of algebraic equations.\Foot{Evidence of such interest is present in Babylonian clay tablets before 2000BC.} Perhaps the most brilliant and influential work on the subject was done by \"Evariste Galois (1811-1832) during the last three years of his short life, when he developed and published what afterwards came to be known as \Quotes{Galois theory}.\Foot{Galois died in a dual and, foreseeing his possible death, wrote a famous letter to Auguste Chevalier describing his ideas for solving equations. The celebrated mathematician Hermann Weyl once said, \Quotes{This letter, if judged by the novelty and profundity of ideas it contains, is perhaps the most substantial piece of writing in the whole history of mankind.}} Several concepts which later were to loom large over the landscape of mathematics, were inherent in Galois' work: the concepts of a \Quotes{group} and a \Quotes{field}. Niels Henrik Abel (1802-1829), a contemporary of Galois doing important work on the solution of equations, also implicitly used the concept of a field in his work.\Foot{This is another example of an idea whose time had come. However it was to be in use for about 60 years before it was finally formalized by Heinrich Martin Weber (1842-1913) in 1893. The English term \Quotes{field} was coined in the same year by Eliakim Hastings Moore (1863-1932). The German word \It{K\"orper} (meaning body or corpus), probably more apt, was introduced by Richard Dedekind (1831-1916) in 1871, as a common term for the two fields of real and complex numbers.}
 \PaR The algebraic concept of a field is intimately related to the geometric concept of a projective space. This was brought to light in a 1857 paper by Karl Georg Christian von Staudt (1798-1867). In his landmark book \LinkText{BibSTAUDT} published in 1847, von Staudt had already laid down the first rigorous axiom system for projective geometry, stripping away the superfluous notions of length and angle, and drawing attention to the fundamental notions of harmonic conjugates and polarity which animate the symmetry forming the heart of the subject. In his 1857 paper he introduced the concept of a \It{Wurf}\kn{2}\Foot{or \Quotes{throw} in English.} of which we shall have more to say \LinkSection{W}{\It{infra}}. These \It{Wurf}\thinspace s\Foot{or \Quotes{W\"urfe} in German.} lead to the construction of a field and, conversely, every field arises in such a way from some projective space. This foreshadowed a future wherein the study of projective spaces was to proceed along two parallel paths, one employing the manipulative tools of algebra, and the other the visual figures of geometry.
 \Item{FT1}{Definition of Projective Space} The algebraic path has as its foundation (along with the concept of a field) the construct of a \Quotes{vector space}. Although the inherent ideas had been around since the early eighteenth century, the formal definition as it is today was given by Giuseppe Peano (1858-1932) in 1888. Vector spaces permeate much of present day mathematics, as well as physics and engineering. Because of their ubiquity and familiarity, the temptation to adopt them as the vehicle of projective geometry is rather strong.
 \PaR Here is one of several possible equivalent definitions. We consider a set X and a group
                                                                                          \DIc{\Ft{Homograph}\Of{X}}{1}
 \Par of permutations of X such that there exists a vector space V over a field F and a bijection \OpGr{b}\ of X onto the family of lines through the origin o of V such that
     \DIc{\Ft{Homograph}\Of{X}\Equals\SetSuch{\Function\Inv{\OpGr{b}}\Circ\setRelation{\OpGr{f}}\Circ\OpGr{b}sendsxinXto\Inv{\OpGr{b}}\Of{\SetSuch{\OpGr{f}\Of{t}}{t\In\OpGr{b}\Of{x}}}inX\end}{\OpGr{f}\ is a
                                                                                       linear automorphism of V\Period}}{2}
 \Par Thus the permutations in \Ft{Homograph}\Of{X}\ are induced by the vector space automorphisms of V applied to the lines through the origin of V\Period\ We shall say that X is a \Def{projective} \Def{space} \Def{with} \Def{defining} \Def{family} \Def{of} \Def{homographies} \Ft{Homograph}\Of{X}\Period\ The function \OpGr{b}\ of \pLinkLocalDisplayText{2} will be said to be a \Def{vector} \Def{representation} of the projective space X.
 \Item{FT2}{Lines and Collineations} Suppose that S is a subset of a projective space X\Period\ It is easy to show that if \OpGr{b}\ and \OpGr{c}\ are two vector representations of X, then \setRelation{\OpGr{b}}\Of{S}\ is the set of lines of a \Two-dimensional subspace of the vector space if, and only if, \setRelation{\OpGr{c}}\Of{S}\ is as well.\Foot{By definition \setRelation{\OpGr{b}}\Of{S}\Equals\SetSuch{\OpGr{b}\Of{x}}{x\In S}\Period\ \It{Cf.} \pLinkDisplayText{Graphs}{FDefs}{19}.} Such sets S are said to be \Def{lines} of X. A bijection from one projective space X onto another Y is said to be a \Def{collineation} if it sends the lines of X to the lines of Y\Period\ If a projective space is itself a line, it is said to be \One-dimensional.
 \Item{FT3}{Projective Isomorphisms} If X and Y are two projective spaces, we say that a bijection \Function\OpGr{f}sendsinXtoinY\end\ is a \Def{projective} \Def{isomorphism} if
                                          \DIc{\Ft{Homograph}\Of{Y}\Equals\SetSuch{\OpGr{f}\Circ\OpGr{h}\Circ\Inv{\OpGr{f}}}{\OpGr{h}\In\Ft{Homograph}\Of{X}}\Period}{1}
 \Par When X\Equals Y, a projective isomorphism is sometimes called a \Def{projective} \Def{automorphism}. It is easy to see and to show that
                                                    \DIc{\ForAll{\OpGr{h}\In\Ft{Homograph}\Of{X}}\quad\OpGr{h}\ is a projective automorphism of X\Period}{2}
 \Par Sometimes the converse of \pLinkLocalDisplayText{2} holds and sometimes it doesn't.
 \PaR It is evident that projective isomorphisms are collineations. If a finite dimensional projective space X is not \One-dimensional, it can be shown that collineations are projective isomorphisms.
 \Item{FT4}{Projective Automorphisms from Homogeneous Coordinates} If \OpGr{b}\ is a vector representation of a projective space X, we shall say that the field of the vector space is a \Def{representation} \Def{field} of X. One can show that any two representation fields of a projective space are isomorphic as fields. There is a description of the general form of a projective automorphism of a projective space X in terms of the automorphisms of the representation fields of X. To present this description, we shall introduce the idea of \Quotes{homogeneous coordinates} for projective spaces.
 \PaR Let F be a field (of characteristic different from 2) and \Eu{n}\ a natural number. We shall say that a subset S of $\overbrace{F\Cross\dots\Cross F}^{\RMBx{7}{n times}}$ is \Eu{n}-\Def{homogeneous} if no element of S has each coordinate zero,
                                         \DI{8}{0}{\ForAll{\Pr{s$_1$,\dots,s$_{\EUBx{7}{n}}$}\In S}\ForAll{k\In F}\quad\Pr{k\Cdot s$_1$,\dots,k\Cdot s$_{\EUBx{7}{n}}$}\In S}{1}
                                                \DII{8}{8}{\ForAll{\Set{\Pr{s$_1$,\dots,s$_{\EUBx{7}{n}}$},\Pr{t$_1$,\dots,t$_{\EUBx{7}{n}}$}\Sin S}}\ThereIs{k\In F}%
                                                      \quad\Pr{k\Cdot s$_1$,\dots,k\Cdot s$_{\EUBx{7}{n}}$}\Equals\Pr{t$_1$,\dots,t$_{\EUBx{7}{n}}$}\Period}{and}{2}
 \Par We shall write the family of all \Eu{n}-homogeneous subsets of $\overbrace{F\Cross\dots\Cross F}^{\RMBx{7}{n times}}$ by
                                                                                        \DIc{F$^{(\EUBx{7}{n})}$\Period}{3}
 \PaR Now let V be a vector space over the field F of dimension \Eu{n}\Period\ We shall denote the family of all lines through the origin o of V by
                                                                                       \DIc{V$_{\FTBx{7}{proj}}$\Period}{4}
 \Par Let \Pr{b$_1$,\dots,b$_n$}\ is a basis for the vector space. Then there are \Eu{n}\ unique functions \Function\OpGr{l}$_i$sendsinVtoinF\end\ such that, for each x\In V,
                                                  \DIc{x\ \Equals\ \OpGr{l}$_{\oNe}$\Of{x}\Cdot b$_1$+\dots+\OpGr{l}$_{\EUBx{7}{n}}$\Of{x}\Cdot b$_n$\Period}{5}
 \Par We define the bijection
         \DIc{\Function\Op{\Gr{n}$_{\EUBx{7}{n}}$}sendsLinV$_{\FTBx{7}{proj}}$to\SetSuch{\Pr{\OpGr{l}$_{\oNe}$\Of{x},\dots,\OpGr{l}$_{\EUBx{7}{n}}$\Of{x}}}{x\In L}in F$^{(\EUBx{7}{n})}$\end\Period}{6}
 \Par Now let \OpGr{a}\ be a field automorphism of F. We define \Function\OpGr{a}$^{(\EUBx{7}{n})}$sendsinF$^{(\EUBx{7}{n})}$toinF$^{(\EUBx{7}{n})}$\end\ by
  \DIc{\ForAll{S\In F$^{(\EUBx{7}{n})}$}\quad\OpGr{a}$^{(\EUBx{7}{n})}$\Of{S}\ \Equiv\ \SetSuch{\Pr{\OpGr{a}\Of{s$_1$},\dots,\OpGr{a}\Of{s$_{\EUBx{7}{n}}$}}}{\Pr{s$_1$,\dots,s$_{\EUBx{7}{n}}$}\In S}}{7}
 \Par Thus, relative to the given basis, each field automorphism \OpGr{a}\ acts on the family of lines through the origin by
                \DIc{\Function\Inv{(\Op{\Gr{n}$_{\EUBx{7}{n}}$})}\Circ\OpGr{a}$^{(\EUBx{7}{n})}$\Circ\Op{\Gr{n}$_{\EUBx{7}{n}}$}sendsinV$_{\FTBx{7}{proj}}$toinV$_{\FTBx{7}{proj}}$\end\Period}{8}
 \Par Let \OpGr{b}\ be a vector representation of a projective space X on V. For each x\In X, the set \Op{\Gr{n}$_{\EUBx{7}{n}}$}\Circ\OpGr{b}\Of{x}\ is said to be a \Def{set} \Def{of} \Def{homogeneous} \Def{coordinates} \Def{for} x\Period\ When \OpGr{a}\ is a field automorphism of the field F of the vector space V, the mapping
                                   \DIc{\Inv{\OpGr{b}}\Circ\Inv{(\Op{\Gr{n}$_{\EUBx{7}{n}}$})}\Circ\OpGr{a}$^{(\EUBx{7}{n})}$\Circ\Op{\Gr{n}$_{\EUBx{7}{n}}$}\Circ\OpGr{b}}{9}
 \Par can without difficulty be shown to be a projective automorphism. We shall call it the \Def{projective} \break\Def{automorphism} \Def{associated} \Def{with} \Def{field} \Def{automorphism} \OpGr{a}\ \Def{and} \Def{the} \Def{vector} \Def{representation} \OpGr{b}\Period\ In the case in which X is not \One-dimensional, this projective automorphism is frequently called an \Def{automorphic} \Def{collineation} \Def{associated} \Def{with} \Def{field} \Def{automorphism} \OpGr{a}\Period
 \Item{FT5}{Anatomy of a Projective Automorphism} It can be proved that every projective automorphism of a projective space X is the composition of a homography and a projective automorphism associated with a field automorphism. This fact is sometimes viewed in the literature as being a part of the \Def{extended} \Def{fundamental} \Def{theorem} \Def{of} \Def{projective} \Def{geometry}. The \Def{first} \Def{form} \Def{of} \Def{the} \Def{fundamental} \Def{theorem} \Def{of} \Def{projective} \Def{geometry}, or what we in these papers call the \Def{fundamental} \Def{theorem} \Def{of} \Def{projective} \Def{geometry} is as follows:
 \Item{FT6}{The Fundamental Theorem} Let X be a projective space. For any two distinct points x and y in X, there is exactly one line
                                                                                                \DIc{\Line{x,y}}{1}
 \Par containing those two points. A subset S of X is called a \Def{projective} \Def{subspace} \Def{of} X if, for each \Set{x,y}\Sin X\Comma\ either x\Equals y or the line \Line{x,y}\ is a a subset of S\Period\ The intersection of all subspaces containing a subset A\Sin X will be written
                                                                                                 \DIc{A$^{oo}$}{2}
 and is called the \Def{subspace} X \Def{spanned} \Def{by} A or, more simply, the \Def{span} \Def{of} A. A subset S of X is said to be \Def{independent} if no subset of S has the same span as any of its proper subsets.
 A maximal independent subset is called a \Def{simplex}. A \Def{basis} for X is a subset B of X such that the complement in B of every singleton is a simplex. The \Def{dimension} of a projective space is \One\ less than the cardinality of a simplex. Thus the dimension of a line is \One\Comma\ and so \Three\ is the cardinality of a basis for a line. What follows is the \Def{fundamental} \Def{theorem}:
 \PaR Let \Set{a$_1$,\dots,a$_m$}\ and \Set{b$_1$,\dots,b$_n$}\ be bases for X --- then
     \DIc{m\Equals n\Andd\ThereIsShriek{\Function\OpGr{f}sendsinStoinS\end\ a homography}\quad b$_1$\Equals\OpGr{f}\Of{a$_1$}, b$_2$\Equals\OpGr{f}\Of{a$_2$},\dots, b$_n$\Equals\OpGr{f}\Of{a$_n$}\Period}{3}
 \Item{FT7}{Perspectivities} Let X be a projective space projectively isomorphic with a projective subspace\Foot{A projective subspace S of a projective space P is actually a projective space itself, where\ \ \ \ \ \ \ \ \ \ \ \ \ \ \ \ \ \ \ \ \ \ \ \ \ \ \ \ \ \ \ \ \ \ \ \ \ \ \ \ \break\hbox to1.4in{} \Ft{Homograph}\Of{S}\Equiv\SetSuch{\OpGr{f}\Restriction{S}}{\OpGr{f}\In\Ft{Homograph}\Of{P}\And\setRelation{\OpGr{f}}\Of{S}\Equals S}\Period} of a projective space Y of one dimension greater than X. Let \OpGr{b}\ and \OpGr{c}\ be projective isomorphisms of X onto distinct subspaces of Y. Let p be any element of Y which lies on neither of the ranges \Range{\OpGr{b}}\ nor \Range{\OpGr{c}}\Period\ Given any point x\In X, the line \Line{p,\OpGr{b}\Of{x}}\ intersects \Range{\OpGr{c}}\ at exactly one point \Line{p,\OpGr{b}\Of{x}}\Wedge\Range{\OpGr{c}}\Period\ The function
                                                          \DIc{\Function sendsxinXto\Inv{\OpGr{c}}\Of{\Line{p,\OpGr{b}\Of{x}}\Wedge\Range{\OpGr{c}}}inX\end}{1}
 \Par is called a \Def{perspectivity} \Def{of} X. A composition of perspecitivies is called a \Def{projectivity}. It can be shown that
                                             \DIc{\Ft{Homograph}\Of{X}\ \Equals\ \SetSuch{\Function\OpGr{f}sendsinXtoinX\end}{\OpGr{f}\ is a projectivity}\Period}{2}
 \Par This fact is also considered a part of the \Quotes{extended} fundamental theorem of projective geometry, and is integral in showing the equivalence of the definition of projective space given here to the various synthetic definitions.
 \Item{MP}{Meridians} A projective space of dimension \One \ is often called a \Quotes{projective line}. This is because the most well-known example may be viewed as a line, corresponding to the field of real numbers. While suggestive, such a denomination is not optimal -- for several reasons. There are other common examples, such as finite projective spaces, and the projective space corresponding to the field of complex numbers, which are not lines. Another reason is that even the one dimensional projective space corresponding to the field of real numbers can be viewed as a circle instead of a line, and in some ways is more suggestively viewed in that manner. Consequently we here adopt a different term for a one dimensional projective space: a \Def{meridian}.\Foot{The word \Quotes{meridian} literally means \Quotes{middle of the day}. If one were on the plane of the ecliptic at the middle of the day, the sun would be directly overhead. A circle passing twice through the polar axis of the earth and the line directly overhead, when viewed from below, appears as a line overhead. This line, located on the plane \Quotes{at infinity}, is called a \Quotes{celestial meridian}. It of course is relative to its observer, since it always lies directly above the observer. It is a physical example of what is here defined as a meridian.}
 \PaR From the point of view of the definition of projective space we have just given, we may regard a meridian as the set of lines passing through the origin of a two dimensional vector space. The synthetic definition is not quite so obvious, as synthetic methods using lines and planes require spaces of dimension larger than \One . We shall review and introduce several alternative methods \It{infra}\Period
 \vfill\eject
                                                                                              \Section{ED}{Erlanger Definition I}\xrdef{pageED}
 \Item{EP}{Erlanger Programm} Felix Klein (1840-1925), while working at the University of Erlangen-N\"urnberg, proposed that geometries be characterized in terms of groups \Mii{G}\ of permutations. His idea was that each permutation of a set X amounted to changing that set so as to be viewed from a new perspective. Thus the properties of X\Comma\ which did not change after applications of these permutations, were the essential and intrinsic properties of X (relative to \Mii{G}). This procedure became known as the Erlanger Programm and has been quite successful in a number of various contexts, some not \It{per} \It{se} involving geometry. We shall utilize this method for our initial definitions of a meridian.
 \Item{L}{Funtion Libras}\setbox\Box\hbox{\quad}\kern-\wd\Box\Foot{This is as special case of a more general \Quotes{libra} which will be defined \LinkSection{LL}{\It{infra}}. The choice of the word \Quotes{libra} is due to the fact that it to some degree is an avatar of a balance or set of scales. We shall have more to say on this in a later article.}\quad  Let X and Y be two sets. For any three bijections \OpGr{a}, \OpGr{b}\And \OpGr{c}\ of X onto Y, we define
                                                               \DIc{\LLL{\OpGr{a},\OpGr{b},\OpGr{c}}\ \Equiv\ \OpGr{a}\Circ\Inv{\OpGr{b}}\Circ\OpGr{c}\Period}{1}
 A family \Mii{L}\ of bijective functions from X onto Y will be called a \Def{function} \Def{libra} if
                                                       \DIc{\ForAll{\Set{\OpGr{a},\OpGr{b},\OpGr{c}}\Sin\Mii{L}}\quad \LLL{\OpGr{a},\OpGr{b},\OpGr{c}}\In\Mii{L}\Period}{2}
 \Par For any family  \Mii{F}\ of bijective functions from X onto Y, we shall write \LLL{\Mii{F}}\ for the intersection of all function libras from X to Y which contain \Mii{F}\ as a subfamily. We say that \LLL{\Mii{F}}\ is the \Def{function} \Def{libra} \Def{generated} \Def{by} \Mii{F}\Period\
 \PaR If X\Equals Y and the elements of \Mii{F}\ are involutions (self-inverse bijections), then
          \DIc{\LLL{\Mii{F}}\ \Equals\ \SetSuch{\OpGr{f}$_{\RM{7}{1}}$\Circ\dots\Circ\OpGr{f}$_{\RM{7}{2n-1}}$}{\Eu{n}\In\Aii{N}\And\Set{\OpGr{f}$_{\RM{7}{1}}$,\dots,\OpGr{f}$_{\RM{7}{2n-1}}$}\Sin\Mii{F}}\Period}{3}
 In this case \LLL{\Mii{F}}\ may be a group, but not necessarily.
 \Item{MD}{Definitions} Let M be a set with at least four elements. Let \Mii{M}\ be a family of non-trivial involutions\Foot{Self-inverse permutations of M not equal to \Id{M}.} of M such that
                                                         \DI{8}{0}{\ForAll{\OpGr{a},\OpGr{b}\In\Mii{M}}\quad\OpGr{a}\Circ\OpGr{b}\Circ\OpGr{a}\In\Mii{M}\Semicolon}{1}
                     \DI{8}{0}{\ForAllSuch{\Set{a,b,c,d}\Sin M}{\Set{a,c}\Cap\Set{b,d}\Equals\Void}\ThereIsShriek{\OpGr{f}\In\Mii{M}}\quad\OpGr{f}\Of{a}\Equals c\Andd\OpGr{f}\Of{b}\Equals d\Semicolon}{2}
                                           \DII{8}{8}{\ForAll{a,b\In M}\quad\SetSuch{\OpGr{f}\In\Mii{M}}{\OpGr{f}\Of{a}\Equals b}\ is a function libra of permutations of M.}{and}{3}
  We shall call such a family \Mii{M}\ a \Def{meridian} \Def{family} \Def{of} \Def{involutions} \Def{of} M, and M will be said to be a \Def{meridian} \Def{relative} \Def{to} \Mii{M}.
 \PaR We note that, because the elements of \Mii{M}\ are self-inverse, \pLinkLocalDisplayText{3} has the following consequence
                                                                                                      \DIc{\ForAll{t\In %
        M}\ForAllSuch{\Set{\OpGr{a},\OpGr{b},\OpGr{c}}\Sin\Mii{M}}{\OpGr{a}\Of{t}\Equals\OpGr{b}\Of{t}\Equals\OpGr{c}\Of{t}}\quad\OpGr{a}\Circ\OpGr{b}\Circ\OpGr{c}\Equals\OpGr{c}\Circ\OpGr{b}\Circ\OpGr{a}\Period}{4}
 \PaR For a, b, c and d as in \pLinkLocalDisplayText{2} we write \Vol{a,c;b,d}\ for the function \OpGr{f}\In\Mii{M}\ satisfying
                                                     \DIc{\ForAll{\OpGr{a},\OpGr{b}\In\Mii{M}}\quad\Vol{a,c;b,d}\Of{a}\Equals c\Andd\Vol{a,c;b,d}\Of{b}\Equals d\Period}{5}
 \PaR Two meridians M$_1$ and M$_2$, respectively, relative to meridian families \Mii{M}$_1$ and \Mii{M}$_2$, respectively, of involutions, are said to be \Def{isomorphic} \Def{as} \Def{meridians} if there is a bijection \Function\OpGr{c} sendsinM$_1$toinM$_2$\end\ such that
                                                             \DIc{\SetSuch{\Inv{\OpGr{c}}\Circ\OpGr{f}\Circ\OpGr{c}}{\OpGr{f}\In\Mii{M}$_2$}\ \Equals\ \Mii{M}$_1$.}{6}
 We say that \LLL{\Mii{M}}\ (\It{cf.} \pLinkDisplayText{ED}{L}{3}) is the group of \Def{homographies} of M.
 \PaR By a \Def{meridian} \Def{basis} we shall mean any subset of a meridian of cardinality \Three . By an \Def{ordered} \Def{meridian} \Def{basis} we mean an ordered triple \Pr{a,b,c}\ where \Set{a,b,c}\ has cardinality \Three .
 \Item{EI}{Example: Field Meridian} We recall that a group is a set G with a binary operation \Function sends\Pr{x,y}inG\Cross Gtox\Cdot yinG\end\ for which there exists e\In G such that
                 \DIc{\ForAll{\Set{x,y,z}\Sin G}\ThereIsShriek{m\In G}\quad e\Cdot x\Equals x\Cdot e\Equals x,\quad x\Cdot m\Equals m\Cdot x\Equals e\Andd x\Cdot (y\Cdot z)\Equals(x\Cdot y)\Cdot z\Period}{1}
 The element e is called the \Def{identity} of \Cdot\ and the element m the \Def{inverse} of x relative to \Cdot\Period\ A group and its binary operation are called \Def{abelian} if x\Cdot y\Equals y\Cdot x for all \Set{x,y}\Sin G.
 \PaR We further recall that a field is a set F with one abelian group binary operation + and another binary operation \thinspace\Cdot\ \thinspace for which there  exists an element \One\ of F distinct from the identity \Zero\ of the group operation + such that \thinspace\Cdot\ \thinspace is a abelian group operation when restricted to the cartesian product of the complement of \Set{\Zero}\ in F with itself, and such that\Foot{It is common to define x\Cdot y to be \Zero\ if either x or y is \Zero\Period\ With this definition condition \pLinkDisplayText{ED}{EI}{2} becomes trivial when any of x, y or z equals \Zero\Period}
                                                                  \DIc{\ForAll{\Set{x,y,z}\Sin F}\quad x\Cdot(y+z)\ \Equals\ (x\Cdot y)+(x\Cdot z)\Period}{2}
 With a field F we write -x for the inverse of an element x relative to the binary operation +, \Inv{x}\ for the inverse of x relative to the operation \Cdot\Comma\ and \Over{x}{y}\ for x\Cdot\Inv{y}\Period
 \PaR Let F be a field such that \One$+$\One\NEq\Zero.\Foot{Thus, F is a field with characteristic different than \Two.} We shall say that F is a \Def{meridian} \Def{field}.
 \PaR Let \Infty\ be any object not in F and let M be the union of F with the singleton\Foot{A \Def{singleton} is any set containing a single element.} \Set{\Infty}. For \Set{a,b,c,d}\Sin F such that a\Cdot d\NEq b\Cdot c, the \Def{homography} \Function\Vol{a,b,c,d}sendsinMtoinM\end\ is defined by \setbox0\hbox{\Over{a\Cdot x+b}{c\Cdot x+d}\quad\hfil if x\In F and c\Cdot x+d\NEq\Zero;}
                             \DIc{\Vol{a,b,c,d}\Of{x}\ \Equiv\ \lower22pt\hbox{\raise51pt\hbox{\MIIICh{18}{50}}\Vbox{\copy0%
                                                                                         \kn{4}%
                                                                                         \hbox to\wd0{\kn{10}\Over{b}{d}\hfil if x\Equals\Infty\And d\NEq\Zero;}%
                                                                                         \kn{4}%
                                                                                         \hbox to\wd0{\kn{9}\Infty\hfil otherwise.}}}}{3}
 \PaR Let
                                                            \DIc{\Mii{M}\ \Equiv\ \SetSuch{\Vol{a,b,c,-a}}{\Set{a,b,c}\Sin F\And b\Cdot c+a\Cdot a\NEq\Zero}\Period}{4}
 It is a pedestrian exercise to show that \Mii{M}\ is a meridian family of involutions of M
 \Item{T}{Theorem} Let \Mii{M}\ be a meridian family of involutions of a set M, and let \Set{\Zero,\One,\Infty} be a basis for M. Let F be the complement in M of the singleton \Set{\Infty}. We define
                                                                       \DIc{\ForAll{\Set{x,y}\Sin F}\quad x+y\ \Equiv\ \Vol{\Infty,\Infty;x,y}\Of{\Zero}}{1}
 and
                                                         \DIc{\ForAllSuch{\Set{x,y}\Sin F}{\Zero\NIn\Set{x,y}}\quad x\Cdot y\ \Equiv\ \Vol{\Infty,\Zero;x,y}\Of{\One}\Period}{2}
 Then F is a field relative to the operations + and \thinspace\Cdot\thinspace, and it is not of characteristic \Two\Period
 \Proof Let x and y be in F. Since \Vol{\Infty,\Infty;x+y,\Zero}\ leaves \Infty\ fixed and sends \Zero\ to x+y, if follows from the uniqueness part of \pLinkDisplayText{ED}{MD}{2} that
                                                                       \DIc{\Vol{\Infty,\Infty;x+y,\Zero}\Equals\Vol{\Infty,\Infty;x,y}\Period}{2.5}
 \Par Since all the constituents of \Vol{\Infty,\Infty;x,\Zero}\Circ\Vol{\Infty,\Infty;\Zero,\Zero}\Circ\Vol{\Infty,\Infty;x,\Zero}\ leave \Infty\ fixed, it follows from \pLinkDisplayText{ED}{MD}{3} that it is in \Mii{M}\Period\ We have
                 \Dc{\Vol{\Infty,\Infty;x,\Zero}\Circ\Vol{\Infty,\Infty;\Zero,\Zero}\Circ\Vol{\Infty,\Infty;y,\Zero}\Of{y}\Equals\Vol{\Infty,\Infty;x,\Zero}\Circ\Vol{\Infty,\Infty;\Zero,\Zero}\Of{\Zero}%
                                                                 \Equals\Vol{\Infty,\Infty;x,\Zero}\Of{\Zero}\Equals x\Equals\Vol{\Infty,\Infty;x,y}\Of{y}}
  and so from the uniqueness part of \pLinkDisplayText{ED}{MD}{2} follows
                                      \DIc{\Vol{\Infty,\Infty;x,\Zero}\Circ\Vol{\Infty,\Infty;\Zero,\Zero}\Circ\Vol{\Infty,\Infty;y,\Zero}\Equals\Vol{\Infty,\Infty;x,y}\Period}{2.6}
  \Par We have
                                             \DIc{\Vol{\Infty,\Infty;x+y,\Zero}\ \Eq{\hLinkLocalDisplayText{2.5}}\ \Vol{\Infty,\Infty;x,y}\ \Eq{\hLinkLocalDisplayText{2.6}}\
                                                         \Vol{\Infty,\Infty;x,\Zero}\Circ\Vol{\Infty,\Infty;\Zero,\Zero}\Circ\Vol{\Infty,\Infty;y,\Zero}\Period}{3}
 Similar reasoning shows that, if \Zero\NIn\Set{x,y}, then
                    \DIc{\Vol{\Infty,\Zero;x\Cdot y,\One}\ \Equals\ \Vol{\Infty,\Zero;x,y}\ \Equals\ \Vol{\Infty,\Zero;x,\One}\Circ\Vol{\Infty,\Zero;\One,\One}\Circ\Vol{\Infty,\Zero;y,\One}\Period}{4}
 \PaR For \Set{x,u}\Sin F such that \Zero\NEq u, we define
                                                             \DIc{-x\ \Equiv\ \Vol{\Infty,\Infty;\Zero,\Zero}\Of{x}\Andd\Inv{u}\ \Equiv\ \Vol{\Infty,\Zero;\One,\One}\Of{u}\Colon}{5}
 we have
                                      \DI{8}{0}{x+\Zero\Equals\Vol{\Infty,\Infty;x,\Zero}\Of{\Zero}\Equals x,\quad u\Cdot \One\Equals\Vol{\Zero,\Infty;u,\One}\Of{\One}\Equals u,}{6}
             \D{8}{0}{(-x)+x\Equals\Vol{\Infty,\Infty;-x,x}\Of{\Zero}\Eq{\hLinkLocalDisplayText{3}}\Vol{\Infty,\Infty;-x,\Zero}\Circ\Vol{\Infty,\Infty;\Zero,\Zero}\Circ\Vol{\Infty,\Infty;x,\Zero}\Of{\Zero}}%
                                                                                                         \DI{-4}{-4}{}{7}
                                               \D{0}{8}{\Equals\Vol{\Infty,\Infty;-x,\Zero}\Circ\Vol{\Infty,\Infty;\Zero,\Zero}\Of{x}\Eq{\hLinkLocalDisplayText{5}}\Vol{\Infty,\Infty;-x,\Zero}\Of{-x}\Equals \Zero,}
           \D{8}{0}{\Inv{u}\Cdot u\Equals\Vol{\Infty,\Zero;\Inv{u},u}\Of{\One}\Eq{\hLinkLocalDisplayText{4}}\Vol{\Infty,\Zero;\Inv{u},\One}\Circ\Vol{\Infty,\Zero;\One,\One}\Circ\Vol{\Infty,\Zero;u,\One}\Of{\One}}%
                                                                                                        \DI{-4}{-4}{}{8}
                      \D{0}{8}{\Equals\Vol{\Infty,\Zero;\Inv{u},\One}\Circ\Vol{\Infty,\Zero;\One,\One}\Of{u}\Eq{\hLinkLocalDisplayText{5}}\Vol{\Infty,\Zero;\Inv{u},\One}\Of{\Inv{u}}\Equals \One\Comma}
 and furthermore, for \Set{y,z,v,w}\Sin F such that \Zero\NIn\Set{v,w}
    \DIc{x+y\Equals\Vol{\Infty,\Infty;x,y}\Of{\Zero}\Equals\Vol{\Infty,\Infty;y,x}\Of{\Zero}\Equals y+x,\quad u\Cdot v\Equals\Vol{\Infty,\Zero;u,v}\Of{\One}\Equals\Vol{\Infty,\Zero;v,u}\Of{\One}\Equals v\Cdot u,}{9}
\D{8}{0}{x+(y+z)\Equals\Vol{\Infty,\Infty;x,y+z}\Of{\Zero}\Eq{\hLinkLocalDisplayText{3}}\Vol{\Infty,\Infty;x,\Zero}\Circ\Vol{\Infty,\Infty;\Zero,\Zero}\Circ\Vol{\Infty,\Infty;y+z,\Zero}\Of{\Zero}\Eq{\hLinkLocalDisplayText{3}}}
 \DI{4}{0}{\Vol{\Infty,\Infty;x,\Zero}\Circ\Vol{\Infty,\Infty;\Zero,\Zero}\Circ\Vol{\Infty,\Infty;y,\Zero}\Circ\Vol{\Infty,\Infty;\Zero,\Zero}\Circ\Vol{\Infty,\Infty;z,\Zero}\Of{\Zero}\Eq{\hLinkLocalDisplayText{3}}}{10}
           \D{4}{8}{\Vol{\Infty,\Infty;x+y,\Zero}\Circ\Vol{\Infty,\Infty;\Zero,\Zero}\Circ\Vol{\Infty,\Infty;z,\Zero}\Of{\Zero}\Eq{\hLinkLocalDisplayText{3}}\Vol{\Infty,\Infty;x+y,z}\Of{\Zero}\Equals(x+y)+z,}
               \DL{8}{0}{u\Cdot (v\Cdot w)\Equals\Vol{\Infty,\Zero;u,v\Cdot w}\Of{\One}\Eq{\hLinkLocalDisplayText{4}}\Vol{\Infty,\Zero;u,\One}\Circ\Vol{\Infty,\Zero;\One,\One}\Circ\Vol{\Infty,\Zero;v\Cdot
                                                                                    w,\One}\Of{\One}\Eq{\hLinkLocalDisplayText{4}}}{and}
       \DI{4}{0}{\Vol{\Infty,\Zero;u,\One}\Circ\Vol{\Infty,\Zero;\One,\One}\Circ\Vol{\Infty,\Zero;v,\One}\Circ\Vol{\Infty,\Zero;\One,\One}\Circ\Vol{\Infty,\Zero;w,\One}\Of{\One}\Eq{\hLinkLocalDisplayText{4}}}{11}
       \D{4}{8}{\Vol{\Infty,\Zero;u\Cdot v,\One}\Circ\Vol{\Infty,\Zero;\One,\One}\Circ\Vol{\Infty,\Zero;w,\One}\Of{\One}\Eq{\hLinkLocalDisplayText{4}}\Vol{\Infty,\Zero;u\Cdot v,w}\Of{\One}\Equals(u\Cdot v)\Cdot w}
 \PaR It follows from \pLinkLocalDisplayText{6} through \pLinkLocalDisplayText{11} that + and \Cdot\ are abelian binary group operations. It remains to show the distributive property \pLinkDisplayText{ED}{EI}{2}. Let then x, y and z be elements of F such that x and y+z are distinct from \Zero. We define
                 \DIc{\OpGr{th}\ \Equiv\ \Vol{\Infty,\Zero;x,\One}\Circ\Vol{\Infty,\Zero;\One,\One}\Circ\Vol{\Infty,\Infty;y,z}\Circ\Vol{\Infty,\Zero;\One,\One}\Circ\Vol{\Infty,\Zero;x,\One}\Period}{12}
 It follows from \pLinkDisplayText{ED}{MD}{1} that \OpGr{th}\ is in \Mii{M}. Direct calculation shows that
                                                                                   \DIc{\OpGr{th}\Of{\Infty}\ \Equals\ \Infty\Period}{13}
 We have
          \D{8}{0}{\OpGr{th}\Of{\Zero}\Equals\Vol{\Infty,\Zero;x,l}\Circ\Vol{\Infty,\Zero;\One,\One}\Circ\Vol{\Infty,\Infty;y,z}\Circ\Vol{\Infty,\Zero;\One,\One}\Circ\Vol{\Infty,\Zero;x,\One}\Of{\Zero}\Equals%
                                         \Vol{\Infty,\Zero;x,\One}\Circ\Vol{\Infty,\Zero;\One,\One}\Circ\Vol{\Infty,\Infty;y,z}\Circ\Vol{\Infty,\Zero;\One,\One}\Of{\Infty}\Equals}
                 \D{4}{0}{\Vol{\Infty,\Zero;x,\One}\Circ\Vol{\Infty,\Zero;\One,\One}\Circ\Vol{\Infty,\Infty;y,z}\Of{\Zero}\Equals\Vol{\Infty,\Zero;x,\One}\Circ\Vol{\Infty,\Zero;\One,\One}\Of{y+z}\Equals}
                                                                                                        \DI{-2}{-2}{}{14}
            \D{0}{4}{\Vol{\Infty,\Zero;x,\One}\Circ\Vol{\Infty,\Zero;\One,\One}\Circ\Vol{\Infty,\Zero;\One,y+z}\Circ\Vol{\Infty,\Zero;\One,y+z}\Of{y+z}\Equals\Vol{\Infty,\Zero;x,\One}%
                                                         \Circ\Vol{\Infty,\Zero;\One,\One}\Circ\Vol{\Infty,\Zero;\One,y+z}\Of{\One}\Eq{\hLinkLocalDisplayText{4}}}
                                                                           \D{0}{8}{\Vol{\Infty,\Zero;x,y+z}\Of{\One}\Equals x\Cdot(y+z)\Period}
 \PaR Furthermore, since \pLinkLocalDisplayText{4} implies that \Vol{\Infty,\Zero;x,y}\Equals\Vol{\Infty,\Zero;x,\One}\Circ\Vol{\Infty,\Zero;\One,\One}\Circ\Vol{\Infty,\Zero;y,\One}, we have
                                                                                          \D{8}{0}{\OpGr{th}\Of{x\Cdot y}\Equals}
                 \D{8}{0}{\Vol{\Infty,\Zero;x,\One}\Circ\Vol{\Infty,\Zero;\One,\One}\Circ\Vol{\Infty,\Infty;y,z}\Circ\Vol{\Infty,\Zero;\One,\One}\Circ\Vol{\Infty,\Zero;x,\One}\Circ\Vol{\Infty,\Zero;x,y}%
                                             \Circ\Vol{\Infty,\Zero;x,y}\Circ\Vol{\Infty,\Zero;x,\One}\Circ\Vol{\Infty,\Zero;\One,\One}\Circ\Vol{\Infty,\Zero;y,\One}\Of{\One}}
\DI{8}{0}{\Equals\Vol{\Infty,\Zero;x,\One}\Circ\Vol{\Infty,\Zero;\One,\One}\Circ\Vol{\Infty,\Infty;y,z}\Circ\Vol{\Infty,\Zero;y,\One}\Of{\One}\Equals\Vol{\Infty,\Zero;x,\One}\Circ\Vol{\Infty,\Zero;\One,\One}\Of{z}\Equals}{15}
              \D{8}{8}{\Vol{\Infty,\Zero;x,\One}\Circ\Vol{\Infty,\Zero;\One,\One}\Circ\Vol{\Infty,\Zero;z,\One}\Of{\One}\Eq{\hLinkLocalDisplayText{4}}\Vol{\Infty,\Zero;x,z}\Of{\One}\Equals x\Cdot z\Period}
 It follows from \pLinkLocalDisplayText{13} and \pLinkLocalDisplayText{15} that
                                                                                    \Dc{\OpGr{th}\ \Equals\ \Vol{\Infty,\Infty;x\Cdot y,x\Cdot z}}
 which implies
                                                          \DIc{\OpGr{th}\Of{\Zero}\Equals\Vol{\Infty,\Infty;x\Cdot y,x\Cdot z}\Of{\Zero}\Equals x\Cdot y+x\Cdot z\Period}{16}
 From \pLinkLocalDisplayText{14} and \pLinkLocalDisplayText{16} follows that x\Cdot(y+z)\Equals x\Cdot y+x\Cdot z\Period
 \PaR If F were of characteristic \Two\ and x any element of F, then
                                                  \Dc{\Vol{\Infty,\Infty;x,x}\Of{\Zero}\Equals x+x\Equals\Zero\Equals\Vol{\Infty,\Infty;\Zero,\Zero}\Of{\Zero}}
 \Par and so \Vol{\Infty,\Infty;x,x}\ and \Vol{\Infty,\Infty;\Zero,\Zero}\ agree would agree on the two distinct points \Infty\ and \Zero\Comma\ which by \pLinkDisplayText{ED}{MD}{2} would imply that \Vol{\Infty,\Infty;x,x}\Equals\Vol{\Infty,\Infty;\Zero,\Zero}\Period\ But then \Vol{\Infty,\Infty;\Zero,\Zero}\ would leave each element of M fixed, which would be absurd, since \Id{M}\ is not an element of \Mii{M}\Period\ \QED
 \Item{Circle}{Example: Circle Meridian} In \pLinkItemText{ED}{EI} we extended a field by one point to construct a meridian. One can do a similar thing with a plane to construct a \Quotes{projective plane}.
 \PaR Let \Eb{P}\ be a euclidean plane. For each line L in \Eb{P}\ we shall add a point \Infinity{L}\ to L not in \Eb{P}. We do this in such a way that points \Infinity{L$_1$} and \Infinity{L$_2$} are equal if, and only if, the lines L$_1$ and L$_2$ are parallel. The set \Infinity{\Eb{P}}\ of all these \Quotes{points at infinity} is called the \Def{line} \Def{at} \Def{infinity} and its union with \Eb{P}\ will be denoted by
                                                                                                     \DIc{\Aii{P}\Period}{1}
 \PaR Let C be a circle in \Eb{P}\Period\ Let x and y be elements of C. By \Line{x,y}\ we shall mean the line through x and y. If p is any point in \Aii{P}\Cop C,\Foot{By \Aii{P}\Cop C we mean \SetSuch{x\In\Aii{P}}{x\NIn C}\Period\ \It{Cf.} \pLinkDisplayText{SetTheory}{Sets}{14}.} and x is any point in C, then the line \Line{p,x}, unless it is tangent to C, intersects C at exactly one other point: we write
                                                                                             \DIc{\Op{p\Lower{4pt}{\RMBx{7}{C}}}\Of{x}}{2}
 for this other point. In the case of tangency, we define the value of \pLinkLocalDisplayText{2}\ to be just x. It is evident that each such function \Op{p\Lower{4pt}{\RMBx{7}{C}}}\ thus defined is an involution on C. In fact, it can be shown that the family
                                                                   \DIc{\Mii{M}\Of{C}\ \Equiv\ \SetSuch{\Op{p\Lower{4pt}{\RMBx{7}{C}}}}{p\In\Aii{P}\Cop C}}{3}
 is a meridian family of involutions. Relative to this family, C is said to be a \Def{circle} \Def{meridian}.
                                                                             \Figure{-2pt}{CircleInv}{2.5in}{5pt}{Circle Meridian Involution}{FI}
 \PaR Each circle meridian is isomorphic to every other circle meridian as a meridian.\Foot{The same constructions which have here been adopted for a circle can be adduced for an ellipse, or even for a hyperbola or parabola in \Aii{P}\Period\ Each meridian obtained in this manner is isomorphic to a circle meridian.}
 \Item{LineM}{Example: Line Meridian} Let S\Equiv\Set{a,b,c,d}\ be a subset of \Aii{P}\ of cardinality 4. Then the set \Set{\Line{a,b},\Line{a,c},\Line{a,d},\Line{b,c},\Line{b,d},\Line{c,d}}\ is by definition a \Def{complete} \Def{quadrilateral}. If two lines of this complete quadrilateral intersect S in a common point, they will be said to be \Def{adjacent}: otherwise \Def{opposite}. For instance, the lines \Line{a,b}\ and \Line{a,c}\ are adjacent and the lines \Line{a,b}\ and \Line{c,d}\ are opposite.
 \PaR For two distinct lines L and M in \Aii{P}\ we shall denote their point of intersection by L\Wedge M: thus
                                                                                        \DIc{L\Cap M\ \Equals\ \Set{L\Wedge M}\Period}{1}
                                                                             \Figure{5pt}{CompleteQuadrilateral}{5in}{5pt}{Cubic Triple of Pairs on a Line}{FII}
  \PaR Let L be a line in \Aii{P}\ and S as above. If we set r\Equiv\Line{a,b}\Wedge L, r\Prime\Equiv\Line{c,d}\Wedge L, s\Equiv\Line{a,c}\Wedge L, s\Prime\Equiv\Line{b,d}\Wedge L, t\Equiv\Line{a,d}\Wedge L and t\Prime\Equiv\Line{b,c}\Wedge L, then the set
                                                                  \DIc{\Set{\Set{r,r\Prime},\Set{s,s\Prime},\Set{t,t\Prime}}}{2}
 of pairs of points on L will be called a \Def{cubic} \Def{triple} \Def{of} \Def{pairs} \Def{of} \Def{points} \Def{on} \Def{L}.  It can be shown that cubic triples which agree on five of the six involved points, must agree on all six points. Thus, if \Set{a,b}\ and \Set{c,d}\ are disjoint subsets of a line L, then
                               \DIc{\ForAll{t\In C}\ThereIsShriek{\Vol{a,b;c,d}\Of{t}\In L}\quad\Set{\Set{a,b},\Set{c,d},\Set{t,\thinspace\Vol{a,b;c,d}\Of{t}}}\ is a cubic triple.}{3}
 \Par It can be shown that the collection of functions \Vol{a,b;c,d}\ thus defined form a meridian family \Mii{M}\Of{L} of involutions of the line L. Relative to this family \Mii{M}\Of{L}, the line L is called a \Def{line} \Def{meridian}. Every line meridian in \Aii{P}\ is isomorphic to every other line meridian.
 \PaR Furthermore, each line meridian is isomorphic to each circle meridian. If the circle C is tangent to the line L, we can obtain an isomorphism as follows: if p is the point of tangency, choose q in C distinct from p and take the function
                                                                                   \DIc{\Function sendsxinCto(\Line{q,x}\Wedge L)inL\end}{4}
 where  \Line{q,x}\Wedge L\Equals\Infinity{L}\ when x\Equals q.
                                                                    \Figure{-10pt}{stereo}{2.5in}{5pt}{Meridian Isomorphism from a Circle to a Line}{FCL}
 \Item{SphereM}{Example: Sphere Meridian} Let \Eb{E}\ denote three dimensional Euclidian space. As in \Eb{P}\ we associate to each line L a point \Infinity{L}\ distinct from \Eb{E}\ in such a way that, for two lines L$_1$ and L$_2$, the points \Infinity{L$_1$}\ and \Infinity{L$_2$}\ are equal if, and only if, the two lines are parallel. The set
                                                                      \DIc{\Infinity{\Eb{E}}\ \Equiv\ \SetSuch{\Infinity{L}}{L is a line in \Eb{E}}}{1}
 is called the \Def{plane} \Def{at} \Def{infinity}. We denote
                                                                                     \DIc{\Aii{E}\ \Equiv\ \Eb{E}\Cup\Infinity{\Eb{E}}}{2}
 and call \Aii{E}\ a \Def{three} \Def{dimensional} \Def{real} \Def{projective} \Def{space}. A subset X of \Infinity{\Eb{E}}\ is a \Def{line} \Def{at} \Def{infinity} if
                                                              \DIc{\ThereIs{P a plane in \Eb{E}}\quad X\Equals\SetSuch{\Infinity{L}}{L is a line in P}\Period}{3}
 \PaR Let now S be a sphere in \Eb{E}. Let L be any line in \Aii{E} not tangent to S, which intersects S. Then there are two intersection points p and q. If P is the tangent plane at p and Q is the tangent plane at q, then P\Cap Q is called the \Def{line} \Def{dual} \Def{to} \Def{L} \Def{relative} \Def{to} \Def{the} \Def{sphere} S. Given any two lines L and M in \Aii{E}\ which do not intersect, and any point x\NIn(L\Cup M), there exists exactly one line
                                                                                                    \DIc{\Line{(L,M;x)}}{4}
 which passes through x and intersects both L and M. If L is as above, and x is any point on S, the line \Line{(L,P\Cap Q;x)}\ intersects S at another point
                                                                                                    \DIc{\Vol{S;L}\Of{x}}{5}
 (where \Vol{S;L}\Of{x}\ is defined to be x if x is on L). The family
                                                            \DIc{\Mii{S}\ \Equiv\ \SetSuch{\Vol{S;L}}{L a line in \Aii{P}\ with \Card{(L\Cap S)}\Equals\Two}}{6}
 is a meridian family of involutions of S.\Foot{The expression \Card{(L\Cap S)}\ means the cardinality of the set L\Cap S\Period\ \It{Cf.} \pLinkDisplayText{Graphs}{Cardinality}{2}.} The meridian S, relative to this family, will be called a \Def{sphere} \Def{meridian}.\Foot{When an ordered basis \Pr{\Zero,\Infty,\One}\ for S is chosen with \Zero\ and \Infty\ at opposite ends of a diameter of S and \One\ equidistant from \Zero\ and \Infty\Comma\ S is sometimes called a \Def{Riemann} \Def{Sphere}.}
 \PaR The field \Aii{C}\ of complex numbers, as is well-known, may be associated with the plane \Eb{P}, which is then called the \Def{gaussian} \Def{plane} \Eb{G}\Period\ If we add one point \Infinity{\Eb{G}}\Comma\ to \Eb{G}\Comma\ we shall denote it by \Aii{G}\Period\ This set \Aii{G}\ may be viewed as a meridian in consequence of Theorem \pLinkItemText{ED}{T}. This meridian is isomorphic to the sphere meridian. There is a classical isomorphism between these two meridians called the \Def{stereographic} \Def{projection} which is described by the following figure (where p is the point on top of the sphere):
 \Par
                                                       \Figure{-10pt}{Stereograph}{4in}{5pt}{Stereographic Projection of a Sphere onto the Complex Plane}{StC}
 \Par Here the complex plane runs through the \Quotes{equator} of the sphere. There is another projection, analogous to the projection of the circle given in \pLinkItemText{ED}{LineM}, where the plane is tangent to the sphere.
 \Item{FT}{Fundamental Theorem}  Let \Mii{M}\ be a meridian family of involutions of a set M. Let \Set{a,b,c}\ and \Set{u,v,w}\ be subsets of M, each of cardinality \Three. Then
                                               \DIc{\ThereIsShriek{\OpGr{f}\In\LLL{\Mii{M}}}\quad\OpGr{f}\Of{a}\Equals u, \ \ \OpGr{f}\Of{b}\Equals v\Andd\OpGr{f}\Of{c}\Equals w\Period}{1}
 \Proof We shall break the existence part of the problem down into the various possible special cases, and verify that \pLinkLocalDisplayText{1} holds for each case.
 \PaR [Case \One : \Set{a,b,c}\Cap\Set{u,v,w}\Equals\Void\ ]\quad If d\Equiv\Vol{a,u;b,v}\Of{c}\ equals w, we let \OpGr{f}\Equiv\Vol{a,u;b,v}\Period\ Otherwise we let
                                                                      \DIc{\OpGr{f}\ \Equiv\ \Vol{u,v;w,w}\Circ\Vol{u,v;d,w}\Circ\Vol{a,u;b,v}\Period}{2}
 \PaR [Case \Two : a\Equals u\And\Set{b,c}\Cap\Set{v,w}\Equals\Void\ ]\quad Same proof as in Case \One .
 \PaR [Case \Three  : a\Equals u\And b\Equals v ]\quad Same proof as in Case \One .
 \PaR [Case \Four : a\Equals v\And b\Equals u ]\quad Same proof as in Case \One .
 \PaR [Case \Five : a\Equals v\And \Set{b,c}\Cap\Set{u,w}\Equals\Void\ ]\quad Let d\Equiv\Vol{a,a;c,u}\Of{b}\ and then \OpGr{f}\Equiv\Vol{u,u;a,w}\Circ\Vol{u,a;d,w}\Circ\Vol{a,a;c,u}\Period
 \PaR [Case \Six : a\Equals v, b\Equals w ]\quad Let \OpGr{f}\Equiv\Vol{a,a;b,u}\Circ\Vol{a,b;c,u}\Period
 \PaR [Case \Seven : a\Equals v\And c\Equals w ]\quad Same proof as in Case \Five  .
 All other cases can be subsumed by one of these seven by permuting \Set{a,b,c} (and \Set{u,v,w} accordingly). Thus existence is shown.
 \PaR Suppose that there are two functions \OpGr{f}\And\OpGr{th}\ sending a to u, b to v and c to w. Then \OpGr{f}\Circ\Inv{\OpGr{th}}\ fixes all three points a, b and c. Let us rename the ordered basis \Pr{a,b,c}\ as \Pr{o,1,\Infty}\ and let F, + and \Cdot\ be as in \pLinkItemText{ED}{T}. Then \LLL{\Mii{M}}\ is a family of homographies. Thus there exists \Set{p,q,r,s}\Sin F such that \Vol{p,q,r,s}\Equals\OpGr{f}\Circ\Inv{\OpGr{th}}\Period\ We have
                                                \D{8}{0}{o\Equals\OpGr{f}\Circ\Inv{\OpGr{th}}\Of{o}\Equals\Over{p\Cdot o+q}{r\Cdot o+s}\quad\Implies\quad q\Equals o\Comma}
                                                       \D{8}{0}{\Infty\Equals\OpGr{f}\Circ\Inv{\OpGr{th}}\Of{\Infty}\Equals\Over{p}{r}\quad\Implies\quad r\Equals o}
                               \DL{8}{8}{\One\Equals\OpGr{f}\Circ\Inv{\OpGr{th}}\Of{\One}\Equals\Over{p\Cdot\One+q}{r\Cdot\One+s}\Equals\Over{p}{s}\quad\Implies\quad p\Equals s\Period}{and}
 Consequently \OpGr{f}\Circ\Inv{\OpGr{th}}\ is the identity function and so \OpGr{f}\Equals\OpGr{th}\Period\ \QED
 \Item{FTDef}{Notation} We shall find it useful in the sequel to denote the function \OpGr{f}\ of \pLinkDisplayText{ED}{FT}{1} by
                                                                                          \DIc{\Hol{a}{b}{c}{u}{v}{w}\Period}{1}
 \Item{TTP}{Corollary} Let \Mii{M}\ be a meridian family of involutions of a set M. Let a and b be distinct points of M. Then
                                                  \DIc{\ForAllSuch{\OpGr{f}\In\LLL{\Mii{M}}}{\OpGr{f}\Of{a}\Equals b\And\OpGr{f}\Of{b}\Equals a}\quad\OpGr{f}\In\Mii{M}\Period}{1}
  \Proof Let c\In M be distinct from a and b. By the uniqueness part of \pLinkItemText{ED}{FT}, the functions \OpGr{f}\ and \Vol{a,b;c,\OpGr{f}\Of{c}}\ are identical. \QED
 \Item{TT}{Theorem} Let \Mii{M}\ be a meridian family of involutions of a set M and let \OpGr{f}\ be a non-involution element of \LLL{\Mii{M}}. Then there exists \Set{\OpGr{a},\OpGr{b}}\Sin\Mii{M}\ such that \DIc{\OpGr{f}\Equals\OpGr{a}\Circ\OpGr{b}\Period}{1}
 \Proof By hypothesis there exists m\In M such that m, \OpGr{f}\Of{m} and \OpGr{f}\Circ\OpGr{f}\Of{m}\ are distinct. We have
                                            \Dc{\OpGr{f}\Circ\Vol{\OpGr{f}\Of{m},\OpGr{f}\Of{m};m,\OpGr{f}\Circ\OpGr{f}\Of{m}}\Of{\OpGr{f}\Circ\OpGr{f}\Of{m}}\Equals\OpGr{f}\Of{m}}
 and
                                        \Dc{\OpGr{f}\Circ\Vol{\OpGr{f}\Of{m},\OpGr{f}\Of{m};m,\OpGr{f}\Circ\OpGr{f}\Of{m}}\Of{\OpGr{f}\Of{m}}\Equals\OpGr{f}\Circ\OpGr{f}\Of{m}\Period}
 It follows from \pLinkItemText{ED}{TTP} that \OpGr{f}\Circ\Vol{\OpGr{f}\Of{m},\OpGr{f}\Of{m};m,\OpGr{f}\Circ\OpGr{f}}\ is in is in \Mii{M}\Period\ Letting \OpGr{b}\Equiv\Vol{\OpGr{f}\Of{m},\OpGr{f}\Of{m};m,\OpGr{f}\Circ\OpGr{f}\Of{m}}\ and \break\OpGr{a}\Equiv\OpGr{f}\Circ\Vol{\OpGr{f}\Of{m},\OpGr{f}\Of{m};m,\OpGr{f}\Circ\OpGr{f}\Of{m}}\Comma\ we obtain \pLinkLocalDisplayText{1}. \QED
 \Item{FTC}{Theorem}\kn{-10}\Foot{This is the meridian form of the result discussed in \pLinkItemText{P}{FT5}.} Let M be a meridian with meridian family \Mii{M}\ of involutions. Let \Pr{\Zero,\One,\Infty}\ be an ordered basis for M. Let F be the field associated with this basis as in \pLinkItemText{ED}{T}. Let \OpGr{a}\ be any field automorphism of F and extend \OpGr{a}\ to all of M by defining \OpGr{a}\Of{\Infty}\Equiv\Infty\Period\ Then
                                                                                \DIc{\OpGr{a}\ is a meridian automorphism of M\Period}{1}
 \PaR Furthermore, if \OpGr{f}\ is any meridian automorphism, there exists a basis and field automorphism \OpGr{a}\ as in \pLinkLocalDisplayText{1} and an element \OpGr{h}\In\LLL{\Mii{M}}\ such that
                                                                                   \DIc{\OpGr{f}\Equals\OpGr{a}\Circ\OpGr{h}\Period}{2}
 \Proof Let \OpGr{h}\ be any element of \LLL{\Mii{M}}\ and choose \Set{a,b,c,d}\Sin F such that \OpGr{f}\Equals\Vol{a,b,c,d}\ as in \pLinkDisplayText{ED}{EI}{3}. If d\Equals\Zero, then
  \DIc{\OpGr{a}\Circ\OpGr{f}\Circ\Inv{\OpGr{a}}\Of{\Infty}\Equals\OpGr{a}\Circ\OpGr{f}\Of{\Infty}\Equals\OpGr{a}\Of{\Infty}\Equals\Infty\Equals\Vol{\OpGr{a}\Of{a},\OpGr{a}\Of{b},\OpGr{a}\Of{c},\Zero}\Of{\Infty}%
                                                              \Equals\Vol{\OpGr{a}\Of{a},\OpGr{a}\Of{b},\OpGr{a}\Of{c},\OpGr{a}\Of{d}}\Of{\Infty}\Period}{3}
  \Par If d\NEq\Zero\Comma\ then
                        \DIc{\OpGr{a}\Circ\OpGr{f}\Circ\Inv{\OpGr{a}}\Of{\Infty}\Equals\OpGr{a}\Circ\OpGr{f}\Of{\Infty}\Equals\OpGr{a}\Of{\Over{b}{d}}\Equals\Over{\OpGr{a}\Of{b}}{\OpGr{a}\Of{d}}%
                                                              \Equals\Vol{\OpGr{a}\Of{a},\OpGr{a}\Of{b},\OpGr{a}\Of{c},\OpGr{a}\Of{d}}\Of{\Infty}\Period}{4}
  \Par If c\Cdot x+d\NEq\Zero, then
           \DIc{\OpGr{a}\Circ\OpGr{f}\Circ\Inv{\OpGr{a}}\Of{x}\Equals\OpGr{a}\Circ\OpGr{f}\Of{\Inv{\OpGr{a}}\Of{x}}\Equals\OpGr{a}\Of{\Over{a\Cdot\Inv{\OpGr{a}}\Of{x}+b}{c\Cdot\Inv{\OpGr{a}}\Of{x}+d}}\Equals
                        \Over{\OpGr{a}\Of{a}\Cdot x+\OpGr{a}\Of{b}}{\OpGr{a}\Of{c}\Cdot x+\OpGr{a}\Of{d}}\Equals\Vol{\OpGr{a}\Of{a},\OpGr{a}\Of{b},\OpGr{a}\Of{c},\OpGr{a}\Of{d}}\Of{x}\Period}{5}
  \Par Assertion \pLinkLocalDisplayText{1} now follows from \pLinkLocalDisplayText{3}, \pLinkLocalDisplayText{4} and \pLinkLocalDisplayText{5}.
  \PaR Let \OpGr{f}\ be any meridian automorphism. Let q\Equiv\OpGr{f}\Of{\Zero}\Comma\ r\Equiv\OpGr{f}\Of{\One}\ and s\Equiv\OpGr{f}\Of{\Infty}\Period\ Let \OpGr{a}\Equiv\OpGr{f}\Circ\Inv{\Hol{\zEro}{\oNe}{\iNfty}{q}{r}{s}}\Period\ Then \OpGr{a}\ is a meridian automorphism leaving \Zero, \One\ and \Infty\ fixed. Noting that, for \Set{x,y}\Sin F\Comma
                      \Dc{\OpGr{a}\Circ\Vol{\Infty,\Infty;x,y}\Circ\Inv{\OpGr{a}}\Of{\Infty}\Equals\Infty\Andd\OpGr{a}\Circ\Vol{\Infty,\Infty;x,y}\Circ\Inv{\OpGr{a}}\Of{\OpGr{a}\Of{x}}\Equals
                                                                         \OpGr{a}\Circ\Vol{\Infty,\Infty;x,y}\Of{x}\Equals\OpGr{a}\Of{y}\Comma}
 \Par we see that
                                                 \Dc{\OpGr{a}\Circ\Vol{\Infty,\Infty;x,y}\Circ\Inv{\OpGr{ps}y}\Equals\Vol{\Infty,\Infty;\OpGr{a}\Of{x},\OpGr{ps}y\Of{y}}}
  whence follows that
                          \DIc{\OpGr{a}\Of{x+y}\Equals\OpGr{a}\Circ\Vol{\Infty,\Infty;x,y}\Of{\Zero}\Equals\Vol{\Infty,\Infty;\OpGr{ps}yyy\Of{x},\OpGr{ps}y\Of{y}}\Circ\OpGr{a}\Of{\Zero}\Equals
                                                         \Vol{\Infty,\Infty;\OpGr{a}\Of{x},\OpGr{a}\Of{y}}\Of{\Zero}\Equals\OpGr{ps}yy\Of{x}+\OpGr{a}\Of{y}}{6}
  \Par Similarly we have
  \Dc{\OpGr{a}\Circ\Vol{\Zero,\Infty;x,y}\Circ\Inv{\OpGr{a}}\Of{\Infty}\Equals\Zero\Andd\OpGr{ps}y\Circ\Vol{\Zero,\Infty;x,y}\Circ\Inv{\OpGr{a}}\Of{\OpGr{a}\Of{x}}\Equals
                                                                         \OpGr{a}\Circ\Vol{\Zero,\Infty;x,y}\Of{x}\Equals\OpGr{a}\Of{y}\Comma}
  and so
                                                  \Dc{\OpGr{a}\Circ\Vol{\Zero,\Infty;x,y}\Circ\Inv{\OpGr{a}}\Equals\Vol{\Zero,\Infty;\OpGr{a}\Of{x},\OpGr{a}\Of{y}}}
  whence follows that
                          \DIc{\OpGr{a}\Of{x\Cdot y}\Equals\OpGr{a}\Circ\Vol{\Zero,\Infty;x,y}\Of{\One}\Equals\Vol{\Infty,\Zero;\OpGr{a}\Of{x},\OpGr{a}\Of{y}}\Circ\OpGr{a}\Of{\One}\Equals
                                                       \Vol{\Zero,\Infty;\OpGr{a}\Of{x},\OpGr{a}\Of{y}}\Of{\One}\Equals\OpGr{a}\Of{x}\Cdot\OpGr{a}\Of{y}}{7}
  \Par From \pLinkLocalDisplayText{6} and \pLinkLocalDisplayText{7} follows that \OpGr{a}\Restriction{F}\ is a field automorphism. Letting \OpGr{h}\Equiv\Hol{\zEro}{\oNe}{\iNfty}{q}{r}{s}\Comma\ we obtain \pLinkLocalDisplayText{2}. \QED
 \Item{AOS}{Projective Automorphisms of the Sphere Meridian} In contrast to the case of the circle meridian, the sphere meridian has projective automorphisms which are not homographies. This is connected with the fact that there is a hausdorff topology inherent in the meridian structure of the circle meridian\Foot{\It{Cf.} \pLinkItemText{EA}{D13}.}, while there is none inherent in the meridian structure of the sphere meridian. There is a topology on S inherited from the metric on euclidian space \Eu{E}\Comma\ but this topology cannot be obtained just by examining \Ft{Homograph}\Of{S}\Period\ However if we do apply the particular metric inherited from \Eu{E}\ to S\Comma\ then we can describe exactly the meridian automorphisms continuous relative to that topology.\Foot{There are however an uncountable number of non-continuous projective automorphisms of S.}
 \PaR We recall from \pLinkItemText{ED}{SphereM} that S induces a duality of lines in \Aii{E}\Period\ There is a analogous duality of points with planes as well. Let H be any plane in \Aii{E}\ which intersects S in more than a single point. The points of intersection comprise a circle. At each point of that circle there is a plane in \Aii{E}\ tangent to S at that point. The intersection \Set{h}\ of all these tangent planes is a singleton, the point of which h is called the \Def{point} \Def{dual} \Def{to} H \Def{relative} \Def{to} S\Period\ This \Def{dual} \Def{pair} \Pr{H,h} induces an involution \Vol{H;h}\ of which the fixed points are the points on H\Cap S and of which the values on other points of S are defined as follows
                                                         \DIc{\ForAll{x\In S\Cop(H\Cap S)}\quad\Set{\Vol{H;h}\Of{x}}\Equals\Line{x,h}\Cap(S\Cop\Set{x})\Period}{1}
 \PaR Thus \Vol{H;h}\Of{x}\ is the point other than x on S through which the line through x and h pierces. The function \Vol{H;h}\ is evidently continuous and it is not difficult to show that it is a projective automorphism of S. In fact it can be shown that each projective automorphism of S which is continuous relative to the metric inherited from \Eu{E}\  is of the form \OpGr{f}\Circ\Vol{H;h}, where \OpGr{f}\ is a homography and H is a plane which intersects S at more than one point.\Foot{To this purpose we note that if the automorphism \OpGr{f}\ in \pLinkDisplayText{ED}{FTC}{2} is continuous, then so is the automorphism \OpGr{a}\ which comes from a field automorphism. We choose an ordered basis \Pr{\Zero,\One,\Infty}\ on S so that the field determined by that basis corresponds to the field isomorphism. It is known that continuous field automorphisms of the complex field are the identity function and complex conjugation. If H is the plane through S containing the three basis elements \Zero, \One\ and \Infty\Comma\ it is not difficult to see that the automorphism \Vol{H;h}\ is just complex conjugation.}
 \Par
                                                             \Figure{-10pt}{PlanePointDual}{5in}{5pt}{Continuous Involutive Projective Automorphism}{CPPA}
 \Item{DF}{Theorem} Let \Mii{M}\ be a meridian family of involutions of a set M. Let \OpGr{f}\ be in \Mii{M}\ and suppose that \OpGr{f}\Of{m}\Equals m for some m\In M. Then there is exactly one other element n of M distinct from m such that \OpGr{f}\Of{n}\Equals n.
 \Proof There can be no more than one such element n since otherwise it would follow from the fundamental theorem that \OpGr{f}\ would be the identity function on M. Let a be any element of M distinct from m and let b\Equiv\OpGr{f}\Of{a}. If b\Equals a, we may let n\Equiv a. Else let n\Equiv\Vol{a,a;b,b}\Of{m}. We have
                                                                 \Dc{\Vol{a,a;b,b}\Circ\OpGr{f}\Of{a}\Equals b\Andd \Vol{a,a;b,b}\Circ\OpGr{f}\Of{b}\Equals a}
 \Par whence from \pLinkItemText{ED}{TTP} follows that \Vol{a,a;b,b}\Circ\OpGr{f}\ is an involution. Consequently
                                                                        \DIc{\Vol{a,a;b,b}\Circ\OpGr{f}\Equals\OpGr{f}\Circ\Vol{a,a;b,b}\Period}{1}
  Letting n\Equiv\Vol{a,a;b,b}\Of{m}, we have
                                \Dc{\OpGr{f}\Of{n}\Equals\OpGr{f}\Circ\Vol{a,a;b,b}\Of{m}\Eq{\hLinkLocalDisplayText{1}}\Vol{a,a;b,b}\Circ\OpGr{f}\Of{m}\Equals\Vol{a,a;b,b}\Of{m}\Equals n.}
 \QED
 \Item{MO}{Definitions} Let \Mii{M}\ be a meridian family of involutions of a set M. A \Def{meridian} \Def{orbit} is a family \Set{\Set{a,b},\Set{b,c},\Set{c,d},\Set{d,a}} of pairs for which there exists an element \OpGr{p}\In\LLL{\Mii{M}}\ such that
                                                    \DIc{\OpGr{p}\Of{a}\Equals b,\quad\OpGr{p}\Of{b}\Equals c,\quad\OpGr{p}\Of{c}\Equals d\Andd\OpGr{p}\Of{d}\Equals a.}{1}
  Such a function \OpGr{p}\ will be called a \Def{meridian} \Def{cycle}.
  \PaR For any \Eu{n}\In\Aii{N}\Comma\ we define
                                               \DIc{\LLL{\Mi{M}}$_{\EUBx{7}{n}}$\ \Equiv\ \SetSuch{\OpGr{f}\In\LLL{\Mii{M}}}{\OpGr{f}\ has exactly \Eu{n}\ fixed points.}}{2}
  \Item{HI}{Theorem} Let \Mii{M}\ be a meridian family of involutions of a set M. Then
                                                    \D{8}{0}{\ForAllSuch{\Set{\Set{a,b},\Set{b,c}}}{\Set{a,b,c}\Sin M\And\Card{\Set{a,b,c}}\Equals\Three}\ThereIsShriek{d\In M}}
                                                                                                        \DI{-2}{-2}{}{2}
                                                                         \D{0}{8}{\Set{\Set{a,b},\Set{b,c},\Set{c,d},\Set{d,a}}\ is a meridian orbit.}
 \Proof Let d\Equiv\Vol{a,a;c,c}\Of{b}. Letting \OpGr{p}\Equiv\Vol{a,b;c,d}\Circ\Vol{a,a;c,c}, direct computation shows that \pLinkDisplayText{ED}{MO}{1} holds. This proves existence.
 \PaR For any other \OpGr{p}\ satisfying \pLinkDisplayText{ED}{MO}{1}, \Vol{a,c;b,b}\Circ\Vol{a,b;c,d}\ agrees with \OpGr{p}\ at a, b and d -- and so must agree everywhere by the fundamental theorem. It follows that
                                                          \DIc{d\Equals\OpGr{p}\Of{c}\Equals\Vol{a,c;b,b}\Circ\Vol{a,b;c,d}\Of{c}\Equals\Vol{a,c;b,b}\Of{d}\Period}{3}
 It follows from \pLinkItemText{ED}{DF} that b and d are the only fixed points of \Vol{a,c;b,b}. Thus d is the unique element of M for which \pLinkDisplayText{ED}{MO}{1} holds. \QED
 \Item{HII}{Theorem} Let \Mii{M}\ be a meridian family of involutions of a set M. Then
            \DIc{\ForAllSuch{\Set{\Set{a,b},\Set{b,d}}}{\Set{a,b,d}\Sin M\And\Card{\Set{a,b,d}}\Equals\Three}\ThereIsShriek{\OpGr{t}\In\LLL{\Mii{M}}$_1$}\quad\OpGr{t}\Of{d}\Equals d\Andd\OpGr{t}\Of{a}\Equals b\Period}{1}
  \PaR In addition,
                                                             \DIc{\Set{\Set{a,b},\Set{b,\OpGr{t}\Of{b}},\Set{\OpGr{t}\Of{b},d},\Set{d,a}}\ is a meridian orbit.}{2}
  \Proof Since \Vol{d,d;b,b}\ and \Vol{d,d;a,b}\ have only the one fixed point d in common, it follows from the fundamental theorem that if they agreed on any other point, they would be equal. It follows that \Vol{d,d;b,b}\Circ\Vol{d,d;a,b} is in \LLL{\Mii{M}}$_1$, fixes d and sends a to b. This proves the existence of \pLinkLocalDisplayText{1}.
  \PaR Suppose that \OpGr{t}\ satisfies \pLinkLocalDisplayText{1}. Let c\Equiv\OpGr{t}\Of{b}\Period\ Since\ \Vol{a,c;b,b}\Circ\OpGr{t}\kn{8} interchanges a and b, it follows from \pLinkItemText{ED}{TTP} that it is an involution. Thus
                                      \D{8}{0}{d\Equals\Vol{a,c;b,b}\Circ\OpGr{t}\Circ\Vol{a,c;b,b}\Circ\OpGr{t}\Of{d}\Equals\Vol{a,c;b,b}\Circ\OpGr{t}\Circ\Vol{a,c;b,b}\Of{d}\ \Implies}
                                                      \D{8}{8}{\OpGr{t}\Of{\Vol{a,c;b,b}\Of{d}}\Equals\Vol{a,c;b,b}\Of{d}\ \Implies\ \Vol{a,c;b,b}\Of{d}\Equals }
  \Par since d is the sole fixed point of \OpGr{t}\Period\ Letting \OpGr{p}\Equiv\Vol{a,c;b,b}\Circ\Vol{a,b;c,d}, direct calculation shows that the set \Set{\Set{a,b},\Set{b,c},\Set{c,d},\Set{d,a}}\ is a meridian orbit with meridian cycle \OpGr{p}. By \pLinkDisplayText{ED}{HI}{2}, it follows that c is unique with this property. Since \Vol{d,d;b,b}\Circ\Vol{d,d;a,b}\ can also serve for \OpGr{t}\ in the foregoing, we have that
                                                                        \Dc{\OpGr{t}\Of{b}\Equals c\Equals\Vol{d,d;b,b}\Circ\Vol{d,d;a,b}\Of{b}\Period}
  It follows from the fundamental theorem that \OpGr{t}\Equals\Vol{d,d;b,b}\Circ\Vol{d,d;a,b}\Period\ \QED
 \Item{Har}{Meridian Harmony} Let \Mii{M}\ be a meridian family of involutions of a set M and suppose that \Set{\Set{a,b},\Set{b,c},\Set{c,d},\Set{d,a}} is a meridian orbit. Then the pair of pairs \Set{\Set{a,c},\Set{b,d}}\ is said to be a \Def{harmonic} \Def{pair}. The following are equivalent assertions about \Set{a,b,c,d}\Sin M:
                                                                                  \DI{8}{0}{\Set{\Set{a,c},\Set{b,d}}\ is a harmonic pair;}{1}
                                         \DI{8}{0}{\ThereIs{\OpGr{a}\In\Mii{M}$_2$}\quad\OpGr{a}\Of{a}\Equals c,\quad\OpGr{a}\Of{b}\Equals b\Andd\OpGr{a}\Of{d}\Equals d\Semicolon}{2}
                                          \DI{8}{8}{\ThereIs{\OpGr{t}\In\LLL{\Mii{M}}$_1$}\quad \OpGr{t}\Of{a}\Equals b,\quad\OpGr{t}\Of{b}\Equals c\Andd\OpGr{t}\Of{d}\Equals d\Period}{3}
 \Proof [(1)\Implies(2)]\quad Suppose that \pLinkLocalDisplayText{1} holds and let \OpGr{p}\ be as in \pLinkDisplayText{ED}{MD}{1}. Let \OpGr{ps}\Equiv\Vol{a,c;b,b}\Circ\OpGr{p}\Circ\Vol{a,c;b,b}\ and r\Equiv\Vol{a,c;b,b}\Of{d}\Period\ We have
                                                     \Dc{\OpGr{ps}\Of{c}\Equals r,\quad\OpGr{ps}\Of{b}\Equals c,\quad\OpGr{ps}\Of{a}\Equals b\Andd\OpGr{ps}\Of{r}\Equals a}
 whence follows that \Set{\Set{a,b},\Set{b,c},\Set{c,r},\Set{r,a}}\ is a meridian orbit. This, with \pLinkItemText{ED}{HI}, implies that r\Equals d. Thus, if we let \OpGr{a}\Equiv\Vol{a,c;b,b}, \pLinkLocalDisplayText{2} holds.
 \PaR[(2)\Implies(3)]\quad Suppose now that \pLinkLocalDisplayText{2} holds. By \pLinkItemText{ED}{HII} there exists \OpGr{t}\In\LLL{\Mii{M}}$_1$\ such that \pLinkDisplayText{ED}{HII}{1} and \pLinkDisplayText{ED}{HII}{2} hold. It follows from \pLinkItemText{ED}{HI} that \OpGr{t}\Of{b}\Equals c. Thus \pLinkLocalDisplayText{3} holds.
 \PaR[(3)\Implies(1)]\quad This follows from \pLinkDisplayText{ED}{HII}{1} and \pLinkDisplayText{ED}{HII}{2}. \QED
 \PAR The following figure illustrates harmony within the context of a circle meridian. The set \Set{\Set{a,c},\Set{b,d}} is a pair of harmonic pairs and \Set{\Set{a,b},\Set{b,c},\Set{c,d},\Set{d,a}} is a meridian orbit:
                                                 \Figure{0pt}{Harmony}{3in}{5pt}{\hbox{Harmony on a Circle Meridian}}{FIII}
 \Item{MC}{Mea Culpa} So far as I know \LinkText{BibTITS} was the first who was aware of the characterization of a meridian presented in the present section.
 \PaR Since I have not done an exhaustive search of literature on the subject, I cannot claim originality for any of the contents contained in the present article . They merely represent my perception of the beauty of the subject entertained.
 \vfil\eject
                                                                                              \Section{E}{Erlanger Definition II}\xrdef{pageE}
\Item{I}{Introduction} Another way of expressing the \LinkItem{ED}{FT}{fundamental theorem} is to say that \LLL{\Mii{M}}\ is \Def{\Three-transitive} as a group acting on M. This begs the question: \Quotes{how far does a \Three-transitive group \Mii{G}\ acting on a set X go towards introducing a meridian structure on X?} It is the business of the present section to supply an answer, but first we formalize the meaning of \Three-transitive.
 \Item{D}{Definitions} Let X be a set with cardinality at least 4. By a \Def{basis} \Def{for} X we shall mean a subset of X of cardinality \Three .
  \PaR A group \Mii{G}\ of permutations of X will be said to be \Def{\Three-transitive} if, for each pair \Set{\Set{a,b,c},\Set{u,v,w}} of bases of X,
                                               \DIc{\ThereIsShriek{\OpGr{f}\In\Mii{G}}\quad\OpGr{f}\Of{a}\Equals u, \OpGr{f}\Of{b}\Equals v\And\OpGr{f}\Of{c}\Equals w\Period}{1}
 We shall find it helpful at times in the sequel to denote the function \OpGr{f}\ as
                                                                                             \DIc{\Hol{a}{b}{c}{u}{v}{w}\Period}{2}
  \PaR  \PaR The \Def{order} of an element \OpGr{f}\ of a group \Mii{G}\ of permutations is the least \Eu{n}\In\Aii{N}\ such that
                                                           \DIc{$\overbrace{\hbox{\OpGr{f}\Circ\dots\Circ\OpGr{f}}}^{\RMBx{7}{n times}}$\ \Equals\ \Id{X}\Period}{3}
 For a group \Mii{G}\ of permutations, we introduce the notation
                                          \DI{8}{0}{\ForAll{\Eu{n}\In\Aii{N}}\quad\Mii{G}$_{\RMBx{7}{n}}$\ \Equiv\ \SetSuch{\OpGr{f}\In\Mii{G}}{\OpGr{f}\ has order \Eu{n}}\Period}{4}
\Item{MD}{Definition} We shall say that a \Three-transitive group \Mii{G}\ of permutations of a set X is a \Def{meridian} \Def{group} \Def{of} \Def{permutations} \Def{of} X if the following two additional conditions hold:
   \DI{8}{0}{\ForAllSuch{\Set{a,b,c}\Sin M}{\Card{(\Set{a,b,c})}\Equals\Three}\ThereIsShriek{\OpGr{f}\In\Mii{G}$_4$}\quad\OpGr{f}\Of{a}\Equals b,\quad\OpGr{f}\Of{b}\Equals c\Andd\Inv{\OpGr{f}}\Of{a}\Equals\OpGr{f}\Of{c}}{1}
 \Par and
                                            \DI{0}{8}{\ForAllSuch{\OpGr{f}\In\Mii{G}}{\ThereIs{x\In X}\ \OpGr{f}\Circ\OpGr{f}\Of{x}\Equals x}\quad\OpGr{f}\In\Mii{G}$_2$\Period}{2}
 \PaR We note that the statement \Inv{\OpGr{f}}\Of{a}\Equals\OpGr{f}\Of{c}\ in (1) is superfluous, since it follows from the fact that \OpGr{f}\ has order 4. We include it because is provides clarity in some of the proofs \It{infra}.
\Item{TI}{Theorem} Let \Mii{M}\ be a meridian family of involutions of a set M. Then \LLL{\Mii{M}}\ is a meridian group of permutations.
\Proof That \pLinkDisplayText{E}{MD}{1} holds follows from \pLinkDisplayText{ED}{HI}{2}. That \pLinkDisplayText{E}{MD}{2} holds follows from Corollary \pLinkItemText{ED}{TTP}. \QED
 \Item{L3}{Lemma} Let \Mii{G}\ be a meridian group of permutations of a set X with at least \Four\ elements. Let p and q be distinct points of X. Then
                                                         \DIc{\ThereIsShriek{\OpGr{f}\In\Mii{G}$_2$}\quad\OpGr{f}\Of{p}\Equals p\Andd\OpGr{f}\Of{q}\Equals q\Period}{1}
 \Proof Let x\In X be distinct from p and q. By \pLinkDisplayText{E}{MD}{1} there exists a unique \OpGr{r}\In\Mii{G}$_4$ such that
                                                       \DIc{\OpGr{r}\Of{p}\Equals x,\quad\OpGr{r}\Of{x}\Equals q\Andd\OpGr{r}\Of{q}\Equals\Inv{\OpGr{r}}\Of{p}\Period}{2}
 Let y\Equiv\OpGr{r}\Of{q}. Direct calculation with \pLinkLocalDisplayText{2} shows that
                          \Dc{\Hol{x}{q}{y}{p}{y}{q}\Circ\OpGr{r}\Of{x}\Equals y,\quad\Hol{x}{q}{y}{p}{y}{q}\Circ\OpGr{r}\Of{y}\Equals x,\quad\Hol{x}{q}{y}{p}{y}{q}\Circ\OpGr{r}\Of{p}\Equals p\Andd%
                                                                                      \Hol{x}{q}{y}{p}{y}{q}\Circ\OpGr{r}\Of{q}\Equals q,}
which, along with \pLinkDisplayText{E}{MD}{2} establishes the existence of \OpGr{f}\Equiv\Hol{x}{q}{y}{p}{y}{q}\Circ\OpGr{r}\ in \Mii{G}$_2$.
\PaR Suppose that \OpGr{th}\In\Mii{G}$_2$ satisfies
                                                                               \Dc{\OpGr{th}\Of{p}\Equals p\Andd\OpGr{th}\Of{q}\Equals q\Period}
 \Par Let z\Equiv\OpGr{th}\Of{x}\Period\ We have
                                                     \Dc{\Hol{p}{q}{z}{x}{z}{q}\Circ\OpGr{th}\Of{p}\Equals x,\Andd\Hol{p}{q}{z}{x}{z}{q}\Circ\OpGr{th}\Of{x}\Equals q}
 \Par which by the uniqueness part of \pLinkDisplayText{E}{MD}{2} implies that \Hol{x}{q}{z}{p}{z}{q}\Circ\OpGr{th}\Equals\OpGr{r}\ and z\Equals y. Consequently
                                                         \Dc{\OpGr{th}\Of{q}\Equals\Hol{x}{q}{z}{p}{z}{q}\Circ\OpGr{r}\Of{q}\Equals z\Equals\OpGr{f}\Of{x}\Period}
 \Par Since \OpGr{th}\ and \OpGr{f}\ agree at p and q as well, it follows that they are identical. \QED
 \Item{L4}{Lemma} Let \OpGr{a}\And \OpGr{b}\ be distinct elements of \Mii{G}$_2$ which fix a common point p. Then
                                                                   \DIc{\ForAllSuch{x\In M}{\OpGr{a}\Circ\OpGr{b}\Of{x}\Equals x}\quad x\Equals p\Period}{1}
 \Proof Let x be as in \pLinkLocalDisplayText{1} and assume that x\NEq p. Then \OpGr{a}\Of{x}\Equals\OpGr{b}\Of{x}\Period\ Since \OpGr{a}\ and \OpGr{b}\ are distinct, it follows from \pLinkDisplayText{E}{MD}{2} that \OpGr{a}\Of{x}\NEq x. Thus, if y\Equiv\OpGr{a}\Of{x}, then \OpGr{a}\Of{y}\Equals\OpGr{b}\Of{y}\ and so \OpGr{a}\Equals\OpGr{b}\Period\ \QED
 \Item{L5}{Lemma} Let \Mii{G}\ be a meridian group of permutations of a set X with at least \Four\ elements. Let p, x and y be distinct points of X. Then
                                                        \DIc{\ThereIsShriek{\OpGr{t}\ a translation}\quad\OpGr{t}\Of{p}\Equals p\Andd\OpGr{t}\Of{x}\Equals y\Period}{1}
 \Proof From \pLinkItemText{E}{L3} follows that there exists \OpGr{a}\In\Mii{G}$_2$ such that
                                                                               \DIc{\OpGr{a}\Of{p}\Equals p\Andd\OpGr{a}\Of{y}\Equals y\Period}{2}
 Let \OpGr{t}\Equiv\OpGr{a}\Circ\Hol{p}{x}{y}{p}{y}{x}\Period\ It follows from Lemma \pLinkLocalItemText{L4} that \OpGr{t}\ is a translation, which establishes the existence part of \pLinkLocalDisplayText{1}.
 \PaR Let z\Equiv\OpGr{t}\Of{y}\Period\ We assert that
                                 \DIc{\ThereIs{\OpGr{r}\In\Mii{G}$_4$}\quad \OpGr{r}\Of{p}\Equals z,\quad\OpGr{r}\Of{z}\Equals y,\quad\OpGr{r}\Of{y}\Equals x\Andd\OpGr{r}\Of{x}\Equals p.}{3}
 \Par First note that, in view of \pLinkDisplayText{E}{MD}{2}, both \Hol{x}{y}{z}{z}{y}{x}\ and \Hol{x}{y}{z}{y}{x}{p}\ are in \Mii{G}$_2$. Since \Hol{x}{y}{z}{z}{y}{x}\Circ\OpGr{t}\ interchanges x and y, it also is an involution.
 Thus
                       \Dc{p\Equals\Hol{x}{y}{z}{z}{y}{x}\Circ\OpGr{t}\Circ\Hol{x}{y}{z}{z}{y}{x}\Circ\OpGr{t}\Of{p}\Equals\Hol{x}{y}{z}{z}{y}{x}\Circ\OpGr{t}\Circ\Hol{x}{y}{z}{z}{y}{x}\Of{p}\quad\Implies\quad%
                                               \OpGr{t}\Of{\Hol{x}{y}{z}{z}{y}{x}\Of{p}}\Equals\Hol{x}{y}{z}{z}{y}{x}\Of{p}\quad\Implies\quad\Hol{x}{y}{z}{z}{y}{x}\Of{p}\Equals p\Period}
 \Par Letting \OpGr{r}\Equiv\Hol{x}{y}{z}{y}{x}{p}\Circ\Hol{x}{y}{z}{z}{y}{x}, direct calculation shows that \pLinkLocalDisplayText{3} holds.
 \PaR Now suppose that \OpGr{th}\ is another translation such that \OpGr{th}\Of{p}\Equals p and \OpGr{th}\Of{x}\Equals y. Let w\Equiv\OpGr{th}\Of{y}. Proceeding as above we obtain an element \OpGr{s}\ of \Mii{G}$_2$ such that
                                                    \DIc{\OpGr{s}\Of{p}\Equals w,\quad\OpGr{s}\Of{w}\Equals y,\quad\OpGr{s}\Of{y}\Equals x\Andd\OpGr{s}\Of{x}\Equals p.}{4}
 \Par From \pLinkLocalDisplayText{3} follows
                                                            \Dc{\OpGr{r}\Of{y}\Equals x,\quad\OpGr{r}\Of{x}\Equals p\Andd\OpGr{r}\Of{p}\Equals\Inv{\OpGr{r}}\Of{y}}
 \Par and from \pLinkLocalDisplayText{4} follows
                                                         \Dc{\OpGr{s}\Of{y}\Equals x,\quad\OpGr{s}\Of{x}\Equals p\Andd\OpGr{s}\Of{p}\Equals\Inv{\OpGr{s}}\Of{y}\Period}
 \Par From \pLinkDisplayText{E}{MD}{1} follows that \OpGr{s}\Equals\OpGr{r}, which implies that w\Equals z. Thus \OpGr{t}\ and \OpGr{s}\ agree on three points, and so are equal. \QED
 \Item{T2}{Theorem} Let \Mii{G}\ be a meridian group of permutations of a set X with at least \Four\ elements. Then \Mii{G}$_2$ is a meridian family of permutations of X and \Mii{G}\Equals\LLL{\Mii{G}$_2$}\Period
 \Proof That \pLinkDisplayText{ED}{MD}{1} holds is evident, so  we proceed to establishing \pLinkDisplayText{ED}{MD}{2}. Let \Set{a,b,c,d}\Sin X be such that a\NEq b\NEq c\NEq d\NEq a\Period\ We consider two separate cases.
 \Par\Quad [Case \One : a\Equals c\And b\Equals d]\Quad That there exists a unique \OpGr{f}\In\Mii{G}$_2$ such that \pLinkDisplayText{ED}{MD}{2} holds in this case is a consequence of  \pLinkItemText{E}{L3}.
 \Par\Quad [Case \Two : a\NEq c or b\NEq d]\Quad If a\NEq c we let \OpGr{f}\Equiv\Hol{a}{b}{c}{c}{d}{a}\ --- otherwise we let \OpGr{f}\Equiv\Hol{a}{b}{d}{c}{d}{b}\Period\ It follows from \pLinkDisplayText{E}{MD}{2} that \OpGr{f}\ is an involution. That the equation in \pLinkDisplayText{ED}{MD}{2} holds is now trivial.
 \PaR Now we shall establish \pLinkDisplayText{ED}{MD}{3}. Let a and b be in X and let \OpGr{a}, \OpGr{b}\ and \OpGr{c}\ be elements of \Mii{G}$_2$, each of which sends a to b. Again we shall treat two separate cases.
 \Par\Quad [Case \One : a\NEq b]\Quad That both \OpGr{a}\Circ\OpGr{b}\Circ\OpGr{c}\Of{a}\Equals b and \OpGr{a}\Circ\OpGr{b}\Circ\OpGr{c}\Of{b}\Equals a is evident. It follows from \pLinkDisplayText{E}{MD}{2} that \OpGr{a}\Circ\OpGr{b}\Circ\OpGr{c}\Of{b}\ is an involution.
 \Par\Quad [Case \Two : a\Equals b]\Quad That \OpGr{a}\Circ\OpGr{b}\Circ\OpGr{c}\Of{a}\Equals b is trivial. We must show that \OpGr{a}\Circ\OpGr{b}\Circ\OpGr{c}\ is an involution. If \OpGr{a}\Equals\OpGr{b}\Comma\ this is trivial so we may and shall presume that \OpGr{a}\NEq\OpGr{b}\Period\  Towards our purpose we let x be any element of X distinct from a an define
                                                      \Dc{y\Equiv\OpGr{a}\Circ\OpGr{b}\Of{x},\quad z\Equiv\OpGr{c}\Of{x}\Andd\OpGr{d}\Equiv\Hol{a}{y}{z}{a}{z}{y}\Period}
 \Par By \pLinkItemText{E}{L4} we have y\NEq x\Period\ Thus \OpGr{d}\NEq\OpGr{c}\Period\ That \OpGr{a}\Circ\OpGr{b}\ and \OpGr{d}\Circ\OpGr{c}\ are in \Mii{G}$_1$ follows from Lemma \pLinkItemText{E}{L4}. We have
                                 \Dc{\OpGr{a}\Circ\OpGr{b}\Of{a}\Equals a\Equals\OpGr{d}\Circ\OpGr{c}\Of{a}\Andd\OpGr{a}\Circ\OpGr{b}\Of{x}\Equals y\Equals\OpGr{d}\Circ\OpGr{c}\Of{x}\Period}
 \Par It follows from \pLinkItemText{E}{L5} that
                                                                                 \Dc{\OpGr{a}\Circ\OpGr{b}\Equals\OpGr{d}\Circ\OpGr{c}\Period}
 \Par Consequently \OpGr{a}\Circ\OpGr{b}\Circ\OpGr{c}\Equals\OpGr{d}, which is an involution.
 \PaR It remains to verify that that \Mii{G}\Sin\LLL{\Mii{G}$_2$}. Let \OpGr{f}\ be in \Mii{G}\Period\ We may and shall presume that \OpGr{f}\NEq\Id{X}\Period\ Let w be any element of X such that \OpGr{f}\Of{w}\NEq w. Let x\Equiv\OpGr{f}\Of{w}, y\Equiv\OpGr{f}\Of{x}\ and z\Equiv\OpGr{f}\Of y.
 \Par\Quad[Case \One : y\Equals w]\Quad It follows from \pLinkDisplayText{E}{MD}{2} that \OpGr{f}\ is in \Mii{G}$_2$.
 \Par\Quad [Case \Two : y\NEq w and z\Equals w]\Quad Then \OpGr{f}\Circ\Hol{w}{x}{y}{w}{y}{x}\ sends w to x and x to w which, by \pLinkDisplayText{E}{MD}{2}, implies that \OpGr{a}\Equiv\OpGr{f}\Circ\Hol{w}{x}{y}{w}{y}{x}\ is in \Mii{G}$_2$. Consequently \OpGr{f}\Equals \OpGr{a}\Circ\Hol{w}{x}{y}{w}{y}{x}\Comma\ which is in \LLL{\Mii{G}$_2$}\Period
 \Par\Quad [Case \Three  : \Card{(\Set{w,x,y,z})}\Equals\Four]\Quad Then \OpGr{f}\ agrees with \Hol{x}{y}{z}{z}{y}{x}\Circ\Hol{x}{y}{z}{y}{x}{w}\ on w, x, and y. Thus \OpGr{f}\Equals\Hol{x}{y}{z}{z}{y}{x}\Circ\Hol{x}{y}{x}{y}{z}{w}\Comma\ and so is in \LLL{\Mii{G}$_2$}\Period\ \QED
  \Item{DS}{Remark} It is a consequence of Theorems \pLinkItemText{E}{TI} and \pLinkItemText{E}{T2} that a meridian could have been defined as a set X of at least cardinality 4 having a meridian family of permutations. This choice would have been more strictly consistent with the Erlanger Programm, which requires a group of permutations.
  \Item{Trans}{Terminology} Let \Mii{G}\ be a meridian group of permutations on a set M with at least cardinality 4. We have already appropriated the term \Def{involution} for a self-inverse element of \Mii{G}. The general term for any element of \Mii{G}\ with two fixed points is \Def{dilation}. An element of \Mii{G}\ with one fixed point is called a \Def{translation}. We shall call elements of \Mii{G}\ which are neither involutions, dilations nor translations, \Def{pure} \Def{rotations}. The following characterizations follow from the results of the present section: for \OpGr{f}\In\Mii{G}\ such that \OpGr{f}\NEq\Inv{\OpGr{f}},
     \DI{8}{0}{\OpGr{f}\ is a pure rotation\quad\Iff\quad\ThereIs{\Set{\OpGr{a},\OpGr{b}}\Sin\Mii{G}$_2$}\quad\OpGr{f}\Equals\OpGr{a}\Circ\OpGr{b},\quad\OpGr{a}\Circ\OpGr{b}\NEq\OpGr{b}\Circ\OpGr{a}\Andd\ForAll{m\In M}\
                                                                                           \OpGr{b}\Of{m}\NEq\OpGr{a}\Of{m}\Comma}{5}
                 \DI{8}{4}{\OpGr{f}\ is a translation\quad\Iff\quad\ThereIs{\Set{\OpGr{a},\OpGr{b}}\Sin\Mii{G}$_2$}\ThereIsShriek{m\In M}\quad\OpGr{f}\Equals\OpGr{a}\Circ\OpGr{b},
                                                  \quad\OpGr{a}\Circ\OpGr{b}\NEq\OpGr{b}\Circ\OpGr{a} \Andd\OpGr{a}\Of{m}\Equals\OpGr{b}\Of{m}\Equals m}{6}
  and
                                                                                           \D{4}{0}{\OpGr{f}\ is a dilation\quad\Iff}
                                                                                                        \DI{-4}{-4}{}{7}
                      \D{0}{8}{\ThereIs{m\In M\And\Set{\OpGr{a},\OpGr{b}}\Sin\Mii{G}$_2$\Distinct}\quad\OpGr{f}\Equals\OpGr{a}\Circ\OpGr{b},\quad\OpGr{a}\Circ\OpGr{b}\NEq\OpGr{b}\Circ\OpGr{a}
                                                                        \Andd\OpGr{a}\Of{m}\Equals\OpGr{b}\Of{m}\NEq m\Period}
 \Proof [\pLinkLocalDisplayText{5} holds]\Quad Let \Op{f}\ be a pure rotation. By \pLinkDisplayText{E}{MD}{2} there exists \Set{\OpGr{a},\OpGr{b}}\Sin\Mii{G}$_2$ such that \OpGr{f}\Equals\OpGr{a}\Circ\OpGr{b}\Period\ If \OpGr{b}\ and \OpGr{a}\ agreed at any point m\In M, then \OpGr{f}\ would have a fixed point. If \OpGr{a}\ and \OpGr{b}\ commuted with one another, then \OpGr{f}\ would be an involution.
\PaR If the right side of \pLinkLocalDisplayText{5} holds, than it is trivial that \OpGr{f}\ has no fixed point.
\PaR\Quad[\pLinkLocalDisplayText{6} holds]\Quad Let \Op{f}\ be a translation. By \pLinkItemText{E}{T2} we know that M is a meridian. Let m be the fixed point of \OpGr{f}, let n\In M be distinct from m, let a\Equiv\Inv{\OpGr{f}}\Of{n}\ and let b\Equiv\OpGr{f}\Of{n}\Period\ It follows from \pLinkDisplayText{ED}{Har}{3} that \Set{\Set{m,n},\Set{a,b}}\ is an harmonic pair and so \pLinkDisplayText{ED}{Har}{2} implies that there exists \OpGr{b}\In\Mii{M}\ such that \OpGr{b}\Of{a}\Equals b, \OpGr{b}\Of{m}\Equals n and \OpGr{b}\Of{n}\Equals n. Defining \OpGr{a}\Equiv\OpGr{f}\Circ\OpGr{b}, we have
                    \Dc{\OpGr{a}\Of{n}\Equals\OpGr{f}\Circ\OpGr{b}\Of{n}\Equals\OpGr{f}\Of{n}\Equals b\Andd\OpGr{a}\Of{b}\Equals\OpGr{f}\Circ\OpGr{b}\Of{b}\Equals\OpGr{f}\Of{a}\Equals n
                                                                        \quad\Imp{\hLinkDisplayText{E}{MD}{2}}\quad\OpGr{a}\In\Mii{M}}
\Par Evidently \OpGr{a}\Of{m}\Equals m and, since \OpGr{f}\ is not an involution, we have \OpGr{a}\Circ\OpGr{b}\NEq\OpGr{b}\Circ\OpGr{a}\Period
\PaR Now suppose that \OpGr{f}\Equals\OpGr{a}\Circ\OpGr{b}\ for \Set{\OpGr{a},\OpGr{b}}\Sin\Mii{G}$_2$, that \OpGr{a}\Circ\OpGr{b}\NEq\OpGr{b}\Circ\OpGr{a}\ and that \OpGr{a}\Of{m}\Equals\OpGr{b}\Of{m}\Equals m\Period\ If \OpGr{f}\Of{n}\Equals n for n\NEq m, then \OpGr{a}\ would equal \OpGr{b}\ at n as well as m and so by the fundamental theorem would be equal. Hence \OpGr{f}\ is a translation.
\PaR\Quad[\pLinkLocalDisplayText{7} holds]\Quad Let \OpGr{f}\ be a non-involutive dilation, let \Zero\ and \Infty\ be the distinct fixed points of \OpGr{f}\ and let \One\ be any point in M distinct from \Zero\ and \Infty\Period  Let \OpGr{a}\Equiv\Vol{\Zero,\Infty;\One,\OpGr{f}\Of{\One}}\ and \OpGr{b}\Equiv\Vol{\Zero,\Infty;\One,\One}\Period\ Then
     \Dc{\OpGr{f}\Equals\Hol{\zEro}{\oNe}{\iNfty}{\zEro}{\oNe}{\Op{\gR{f}\oF{\oNe}}}\Equals\Vol{\Zero,\Infty;\One,\OpGr{f}\Of{\One}}\Circ\Vol{\Zero,\Infty;\One,\One}\Equals\OpGr{a}\Circ\OpGr{b}\Comma\
                                                                          \OpGr{a}\Of{\Zero}\Equals\Infty\Equals\OpGr{b}\Of{\Zero}}
 \Par and \OpGr{a}\Circ\OpGr{b}\NEq\OpGr{b}\Circ\OpGr{a}\ since \OpGr{f}\ is not an involution.
 \PaR Finally, let \OpGr{f}\Equals\OpGr{a}\Circ\OpGr{b}\ for \Set{\OpGr{a},\OpGr{b}}\Sin\Mii{G}$_2$, \OpGr{a}\Circ\OpGr{b}\NEq\OpGr{b}\Circ\OpGr{a}\ and \OpGr{a}\Of{m}\Equals\OpGr{b}\Of{m}\NEq m for some m\In M. Then
                                                               \Dc{\OpGr{f}\Of{m}\Equals m\NEq\OpGr{a}\Of{m}\Equals\OpGr{f}\Of{\OpGr{a}\Of{m}}}
\Par and \OpGr{f}\ is not an involution since \OpGr{a}\ an \OpGr{b}\ do not commute. \QED
 \vfill\eject
                                                                                                      \Section{LL}{Involution Libras}\xrdef{pageLL}
 \Item{I}{Introduction} Our next characterization of a meridian will be in terms of \Quotes{balanced} functions on the faces of a cube. An important part of this characterization involves a simpler sort of object which we shall term a \Quotes{libra}. The basic notion behind it is a set of scales -- hence the name. For brevity however, we shall take a short cut past the scales, leaving those for later.
 \Item{LD}{Definitions} Let L be a set and \Function\LL{,,}sendsinL\Cross L\Cross LtoinL\end\ a trinary operator on L for which the following holds:
                                                                        \DI{8}{8}{\ForAll{\Set{a,b}\Sin L}\quad\LL{a,a,b}\Equals b\Equals\LL{b,a,a}}{1}
                                                          \DII{0}{8}{\ForAll{\Set{a,b,c,d,e}\Sin L}\quad\LL{\LL{a,b,c},d,e}\Equals\LL{a,b,\LL{c,d,e}}\Period}{and}{2}
 Then \LL{,,}\ will be said to be a \Def{libra} \Def{operator} and L, relative to \LL{,,}\ a \Def{libra}. A subset B of a libra will be said to be \Def{balanced} provided \LL{a,b,c} is in B whenever \Set{a,b,c}\Sin B.
 \Item{LI}{Theorem} Let \LL{,,}\ be a libra operator on a set L. Then
                                                                \DIc{\ForAll{\Set{a,b,c,e,f}\Sin L}\quad\LL{a,\LL{d,c,b},e}\Equals\LL{\LL{a,b,c},d,e}\Period}{1}
 \Proof We have
              \D{8}{0}{a\Eq{\hLinkDisplayText{LL}{LD}{1}}\LL{a,b,b}\Eq{\hLinkDisplayText{LL}{LD}{1}}\LL{a,b,\LL{c,c,b}}\Eq{\hLinkDisplayText{LL}{LD}{1}}\LL{a,b\LL{\LL{c,d,d},c,b}}\Eq{\hLinkDisplayText{LL}{LD}{2}}}
                                                                                                        \DI{0}{0}{}{2}
                                  \D{0}{8}{\LL{\LL{a,b,\LL{c,d,d}},c,b}\Eq{\hLinkDisplayText{LL}{LD}{2}}\LL{\LL{\LL{a,b,c},d,d},c,b}\Eq{\hLinkDisplayText{LL}{LD}{2}}\LL{\LL{a,b,c},d\LL{d,c,b}}}
whence follows
                                           \D{8}{0}{\LL{a,\LL{d,c,b},e}\Eq{\hLinkLocalDisplayText{2}}\LL{\LL{\LL{a,b,c},d,\LL{d,c,b}},\LL{d,c,b},e}\Eq{\hLinkDisplayText{LL}{LD}{2}}}
                                                      \D{8}{-4}{\LL{\LL{a,b,c},d,\LL{\LL{d,c,b},\LL{d,c,b},e}}\Eq{\hLinkDisplayText{LL}{LD}{1}}\LL{\LL{a,b,c},d,e}\Period}
\par\QED
 \Item{Conv}{Convention} The various compositions of libra operators with libra operators, in view of \pLinkDisplayText{LL}{LD}{1}, \pLinkDisplayText{LL}{LD}{2} and \pLinkDisplayText{LL}{LI}{1}, may be greatly simplified: we define
                                                  \DIc{\LL{a,b,c,d,e}\ \Equiv\ \LL{\LL{a,b,c},d,e}\ \Equals\ \LL{a,\LL{d,c,b},e}\ \Equals\ \LL{a,b,\LL{c,d,e}}\Period}{1}
 Each such composition may be converted to a form
                                                          \DIc{\LL{a$_1$,a$_2$,\LL{a$_3$,a$_4$,\LL{\dots\LL{a$_{\RMBx{7}{n-2}}$,a$_{\RMBx{7}{n-1}}$,a$_n$}\dots}}}}{2}
 for n an odd positive integer. We shall at times adopt the abbreviation
                                                                                              \DIc{\LL{a$_1$,a$_2$,\dots,a$_n$}}{3}
 for \pLinkLocalDisplayText{2}.
 \Item{E1}{Example} Let A be an affine space over a field F. Then the \Def{translations} of A form a vector space over F. The translation of a point a\In A by a vector v\In V is denoted by v+a. To any two distinct points a and b in A corresponds a unique vector (which we denote by b-a) such that (b-a)+a\Equals b. Then
                                                                              \DIc{\ForAll{\Set{a,b,c}\Sin A}\quad\LL{a,b,c}\ \Equiv\ (a-b)+c}{1}
 defines a libra operator. We have d\Equals\LL{a,b,c} precisely when the points a,b, c and d describe the points of a parallelogram.\Foot{taken in clockwise, or counter-clockwise order.}
                                                                          \Figure{-2pt}{Parallelogram}{2.5in}{5pt}{Affine Libra Operator}{FIV}
 \Item{E2}{Example} We have already seen an important example in \pLinkItemText{ED}{L}: the libra of functions.
 \Item{LII}{Theorem} Let \LL{,,}\ be a libra operator on a libra L and e be an element of L. Then the binary operation
                                                              \DIc{\Function\Cdot sends\Pr{x,y}inL\Cross Ltox\Cdot y\Equiv\LL{x,e,y}inL\end}{1}
 is a group operation on L, relative to which e is the identity and
                                                                              \DIc{\ForAll{x\In L}\quad\LL{e,x,e}\ is the group inverse of e.}{2}
 \Proof For \Set{x,y,z}\Sin L
                                               \D{8}{0}{(x\Cdot y)\Cdot z\Equals\LL{\LL{x,e,y},e,z}\Eq{\hLinkDisplayText{LL}{LD}{2}}\LL{x,e\LL{y,e,z}}\Equals x\Cdot(y\Cdot z),}
                                               \D{8}{0}{x\Cdot e\Equals\LL{x,e,e}\Eq{\hLinkDisplayText{LL}{LD}{1}}x\Eq{\hLinkDisplayText{LL}{LD}{1}}\LL{e,e,x}\Equals e\Cdot x,}
                     \D{8}{0}{x\Cdot\LL{e,x,e}\Equals\LL{x,e,\LL{e,x,e}}\Eq{\hLinkDisplayText{LL}{LD}{2}}\LL{\LL{x,e,e},x,e}\Eq{\hLinkDisplayText{LL}{LD}{1}}\LL{x,x,e}\Eq{\hLinkDisplayText{LL}{LD}{1}}e}
              \DL{9}{8}{\LL{e,x,e}\Cdot x\Equals\LL{\LL{e,x,e},e,x}\Eq{\hLinkDisplayText{LL}{LD}{2}}\LL{e,x,\LL{e,e,x}}\Eq{\hLinkDisplayText{LL}{LD}{1}}\LL{e,x,x}\Eq{\hLinkDisplayText{LL}{LD}{1}}e\Period}{and}
 \Item{LIII}{Theorem} Let G be a group with binary operation \Cdot\Period\ Define the trinary operator
                                                                  \DIc{\Function\LL{,,}sends\Pr{a,b,c}inG\Cross G\Cross Gtoa\Cdot\Inv{b}\Cdot cinG\end\ .}{1}
 Then \LL{,,}\ is a libra operator.
  \par\Vbox{\hruleh{2}}
  \Proof For \Set{r,s,t,t,v}\Sin G,
                                                             \D{8}{0}{\LL{r,s,s}\Equals r\Cdot\Inv{s}\Cdot s\Equals r\Equals s\Cdot\Inv{s}\Cdot r\Equals\LL{s,s,r}}
 \DL{8}{4}{\LL{\LL{r,s,t},u,v}\Equals(r\Cdot\Inv{s}\Cdot t)\Cdot\Inv{u}\Cdot v\Equals r\Cdot\Inv{s}\Cdot(t\Cdot\Inv{u}\Cdot v)\Equals\LL{r,s,\LL{t,u,v}}\Period}{and}
 \QED
 \Item{D2}{Definition} The libra operator defined in \pLinkItemText{LL}{LIII} will be called the \Def{group} \Def{libra} \Def{operator}.
 \Item{D3}{Definitions} A function \OpGr{f}\ from one libra L$_1$ into another L$_2$ which preserves the libra operator is called a \Def{libra} \Def{homomorphism}. Thus a libra homomorphism \OpGr{f}\ is characterized by
                                               \DIc{\ForAll{\Set{a,b,c}\Sin L$_1$}\quad\LL{\OpGr{f}\Of{a},\OpGr{f}\Of{b},\OpGr{f}\Of{c}}\Equals\OpGr{f}\Of{\LL{a,b,c}}\Period}{1}
 A bijective libra homomorphism is a \Def{libra} \Def{isomorphism}.
 \Item{LIV}{Theorem} Let G and H be two groups, and let \OpGr{f}\ be a group homomorphism from G into H. Then \OpGr{f}\ is also a libra homomorphism.
 \Proof For \Set{a,b,c}\Sin G we have
            \Dc{\LL{\OpGr{f}\Of{a},\OpGr{f}\Of{b},\OpGr{f}\Of{c}}\Equals\OpGr{f}\Of{a}\Cdot\Inv{(\OpGr{f}\Of{b})}\Cdot\OpGr{f}\Of{d}\Equals\OpGr{f}\Of{a\Cdot\Inv{b}\Cdot c}\Equals \OpGr{f}\Of{\LL{a,b,c}}\Period}
 \Item{D4}{Definitions and Notation} A libra L will be called \Def{abelian} if
                                                                                \DIc{\LL{a,b,c}\Equals\LL{c,b,a}\ for all \Set{a,b,c}\Sin L.}{1}
 \Par Evidently L is abelian if and only if each of its corresponding groups is abelian.
 \PaR For a and b in a libra L, we define the functions
                       \DIc{\Function\LibraP{a}{b}sendsxinLto\LL{a,x,b}inL\end,\quad\Function\LibraL{a}{b}sendsxinLto\LL{a,b,x}inL\end\Andd\Function\LibraR{a}{b}sendsxinLto\LL{x,a,b}inL\end\Period}{2}
 The functions \LibraR{a}{b}\ and \LibraL{a}{b}, respectively, are called \Def{libra} \Def{right} \Def{translations} and \Def{libra} \Def{left} \Def{translations}, respectively. When L is an abelian libra, the function \LibraP{a}{b}\ will be called a \Def{libra} \Def{inner} \Def{involution}.
 \Item{LV}{Theorem} Let \LL{,,}\ be an abelian libra operation on a set L. Let \Gr{P}\Of{L}\ denote the set of inner involutions on L. Then
                                                                    \DI{8}{0}{\ForAll{\OpGr{f}\In\Gr{P}\Of{L}}\quad\OpGr{f}\Equals\Inv{\OpGr{f}}\Comma}{1}
                                                          \DI{8}{0}{\ForAll{\Set{a,b}\Sin L}\ThereIsShriek{\OpGr{f}\In\Gr{P}\Of{L}}\quad\OpGr{f}\Of{a}\Equals b}{2}
                                         \DII{8}{8}{\ForAll{\Set{\OpGr{a},\OpGr{b},\OpGr{c}}\Sin\Gr{P}\Of{L}}\quad\OpGr{a}\Circ\OpGr{b}\Circ\OpGr{c}\In\Gr{P}\Of{L}\Period}{and}{3}
 \Proof For \Set{r,s,t,u,v,w,x}\Sin L
                   \D{8}{0}{\LibraP{r}{s}\Circ\LibraP{r}{s}\Of{x}\Equals\LL{r,\LL{r,x,s},s}\Equals\LL{r,\LL{s,x,r},s}\Eq{\hLinkDisplayText{LL}{LI}{1}}\LL{r,r,x,s,s}\Eq{\hLinkDisplayText{LL}{LD}{1}}x\Comma}
                                                                        \D{8}{8}{\LibraP{r}{s}\Equals\LL{r,r,s}\Eq{\hLinkDisplayText{LL}{LD}{1}}s\Comma}
 and, if we let a\Equiv\LL{r,u,v} and b\Equiv\LL{w,t,s},
            \D{8}{0}{\LibraP{r}{s}\Circ\LibraP{t}{u}\Circ\LibraP{u}{w}\Of{x}\Equals\LL{r,\LL{t,\LL{v,x,w},u},s}\Eq{\hLinkDisplayText{LL}{LI}{1}}\LL{r,\LL{\LL{t,w,\LL{x,v,u}}},s}\Eq{\hLinkDisplayText{LL}{LI}{1}}}
                                                     \D{8}{8}{\LL{r,\LL{x,v,u},\LL{w,t,s}}\Equals\LL{r,u,v,x,w,t,s}\Equals\LL{\LL{r,u,v},x,\LL{w,t,s}}\Equals\LibraP{a}{b}}
 It remains only to show that \LibraP{r}{s}\ is the only element of \Gr{P}\Of{L}\ which sends r to s. Suppose that \LibraP{t}{u} is another such. Then s\Equals\LL{t,r,u}, whence, for each x\In L,
                                             \D{8}{0}{\LibraP{r}{s}\Of{x}\Equals\LL{r,x,s}\Equals\LL{r,x\LL{t,r,u}}\Eq{\hLinkDisplayText{LL}{LI}{1}}\LL{r,\LL{r,t,x},u}\Equals}
                                 \D{8}{-8}{\LL{r,\LL{x,t,r},u}\Eq{\hLinkDisplayText{LL}{LI}{1}}\LL{\LL{r,r,t},x,u}\Eq{\hLinkDisplayText{LL}{LD}{1}}\LL{t,x,u}\Equals\LibraP{t}{u}\Of{x}\Period}
 \QED
  \Item{DV}{Definition} A family \Gr{P}\ of operators of a set S satisfying conditions \pLinkDisplayText{LL}{LV}{1}, \pLinkDisplayText{LL}{LV}{2} and \pLinkDisplayText{LL}{LV}{3} of \pLinkItemText{LL}{LV} will be said to be an \Def{inner} \Def{involution} \Def{libra} \Def{on} S.
 \Item{EAT}{Example and Theorem} Let \Mii{M}\ be a meridian involution family on a meridian M. For \Set{a,b}\Sin M, the family
                                                                       \DIc{\MAB{a}{b}\ \Equiv\ \SetSuch{\OpGr{f}\In\Mii{M}}{\OpGr{f}\Of{a}\Equals b}}{1}
 is an abelian function libra by definition \pLinkDisplayText{ED}{MD}{3}. By \pLinkDisplayText{ED}{MD}{2} and by \pLinkItemText{LL}{DV}, \MAB{a}{b}\ is an inner involution libra on
                                                                                       \DIc{\Mab{a}{b}\ \Equiv\ M\Cop\Set{a,b}\Period}{2}
 \Par Let \Set{a,b,c,d}\Sin M. Then
  \DI{8}{0}{if a\NEq b\And c\NEq d, then\MAB{a}{b}\ and \MAB{c}{d}\ are isomorphic as libras,}{4}
  \DI{8}{0}{\MAB{a}{a}\ and \MAB{c}{c}\ are isomorphic as libras}{5}
  \DII{8}{8}{if a\NEq b, then \MAB{a}{b}\And\MAB{c}{c}\ are not isomorphic as libras.}{and}{6}
 \Proof Suppose first that a\NEq b\And c\NEq d. By the Fundamental Theorem there exists \OpGr{ps}\In\LLL{\Mii{M}}\ such that
                                                                               \Dc{\OpGr{ps}\Of{a}\Equals c\Andd\OpGr{ps}\Of{b}\Equals d\Period}
 \Par The function
                                                          \DIc{\Function sends\OpGr{th}in\MAB{a}{b}to\Inv{\OpGr{ps}}\Circ\OpGr{th}\Circ\OpGr{ps}in\MAB{c}{d}\end}{7}
 \Par is a libra isomorphism.
 \PaR Now suppose that a\Equals b and that c\Equals d. By the Fundamental Theorem there exists \OpGr{ps}\In\LLL{\Mii{M}}\ such that
                                                                                             \Dc{\OpGr{ps}\Of{a}\Equals c\Period}
 \Par Again, the function in \pLinkLocalDisplayText{7} is a libra isomorphism.
 \PaR Finally we suppose that a\NEq b and assume that there existed a libra isomorphism \Op{\Ft{w}}\ from \MAB{a}{b}\ onto \Mab{c}{c}\Period\ Let t\In M be distinct from both a and b. Since \Vol{a,b;t,t}\ fixes t and is an involution, it follows from Theorem \pLinkItemText{ED}{DF} that there would be exactly one other point u\In M such that \Vol{a,b;t,t}\Of{u}\Equals u. From \pLinkDisplayText{ED}{MD}{1} follows that \Vol{a,b;t,u}\Circ\Vol{a,b;t,t}\Circ\Vol{a,b;t,u}\ is in\MAB{a}{b}\Period\ Furthermore
 \Dc{\Vol{a,b;t,u}\Circ\Vol{a,b;t,t}\Circ\Vol{a,b;t,u}\ is in\MAB{a}{b}\ and \Vol{a,b;t,u}\Circ\Vol{a,b;t,t}\Circ\Vol{a,b;t,u}\Of{u}\Equals\Vol{a,b;t,u}\Circ\Vol{a,b;t,t}\Of{t}\Equals\Vol{a,b;t,u}\Of{t}\Equals u\Period}
 \Par Since \Vol{a,b;t,t}\ and \Vol{a,b;t,u}\Circ\Vol{a,b;t,t}\Circ\Vol{a,b;t,u}\ would agree on both t and u, it follows from the uniqueness part of \pLinkDisplayText{ED}{MD}{2} that
                                                                     \DIc{\Vol{a,b;t,t}\Equals\Vol{a,b;t,u}\Circ\Vol{a,b;t,t}\Circ\Vol{a,b;t,u}\Period}{8}
 \Par Let \OpGr{a}\ and \OpGr{b}\ be the elements of \MAB{c}{c}\ such that  \Op{\Ft{w}}\Of{\OpGr{a}}\Equals\Vol{a,b;t,t}\ and \Op{\Ft{w}}\Of{\OpGr{b}}\Equals\Vol{a,b;t,u}\Period\ From \pLinkLocalDisplayText{8} follows that
           \Dc{\Op{\Ft{w}}\Of{\OpGr{b}}\Circ\Op{\Ft{w}}\Of{\OpGr{a}}\Circ\Op{\Ft{w}}\Of{\OpGr{a}}\Equals \Op{\Ft{w}}\Of{\OpGr{b}}\quad\Implies\quad\OpGr{b}\Circ\OpGr{a}\Circ\OpGr{b}\Equals\OpGr{a}\Period}
 \Par Since \OpGr{b}\ would fix c, from \pLinkItemText{ED}{DF} follows that it would have another distinct fixed point p. Then
                                       \Dc{\OpGr{a}\Of{p}\Equals\OpGr{b}\Circ\OpGr{a}\Circ\OpGr{b}\Of{p}\Equals\OpGr{b}\Of{\OpGr{a}\Of{p}}\quad\Implies\quad\OpGr{a}\Of{p}\In\Set{c,p}\Period}
 \Par But \OpGr{a}\Of{c}\Equals c, and so \OpGr{a}\Of{p}\Equals p\Period\ It follows from the uniqueness part of \pLinkDisplayText{ED}{MD}{2} that \OpGr{a}\Equals\OpGr{b}\Period\ This however is absurd because \Vol{a,b;t,t}\ and \Vol{a,b;t,u}\ are distinct. \QED
 \Item{EIIL}{Example} Let C be a circle as in \pLinkItemText{ED}{Circle}. Let L be a line in the projective plane \Aii{P}. For each p\In L, the function \Op{p\Lower{4pt}{\RMBx{7}{C}}}\ is defined as in \pLinkItemText{ED}{Circle}. The family \SetSuch{\Op{p\Lower{4pt}{\RMBx{7}{C}}}}{p\In L} is an inner involution libra on C.
                                                                                 \Figure{-20pt}{CircleInv}{2.5in}{5pt}{Circle Meridian Involution}{1}
 \PaR If the line L intersects C, either at two points a and b, or at a point of tangency c, then then the corresponding inner involution meridian is either the family \MAB{a}{b}\Comma\ or the family \MAB{c}{c}\ of \pLinkItemText{LL}{EAT}\Period
 \Item{LVI}{Corollary} Relative to the trinary operator
                             \Dc{\Function sends\Pr{\OpGr{a},\OpGr{b},\OpGr{c}}in\Gr{P}\Of{L}\Cross\Gr{P}\Of{L}\Cross\Gr{P}\Of{L}to\OpGr{a}\Circ\OpGr{b}\Circ\OpGr{c}in\Gr{P}\Of{L}\end\Comma}
 \hbox{\Gr{P}\Of{L}\ is a function libra.}
 \Item{LVII}{Theorem} Let \Gr{P}\ be a family of permutations of a set S such that
                                                                       \DI{8}{0}{\ForAll{\OpGr{f}\In\Gr{P}}\quad\OpGr{f}\Equals\Inv{\OpGr{f}}\Comma}{1}
                                                          \DI{8}{0}{\ForAll{\Set{a,b}\Sin S}\ThereIsShriek{\OpGr{f}\In\Gr{P}}\quad\OpGr{f}\Of{a}\Equals b\Comma}{2}
                                               \DII{8}{8}{\ForAll{\Set{\OpGr{a},\OpGr{b},\OpGr{c}}\Sin\Gr{P}}\quad\OpGr{a}\Circ\OpGr{b}\Circ\OpGr{c}\In\Gr{P}\Period}{and}{3}
 \PaR Let \Gr{TH}\Equiv\SetSuch{\OpGr{a}\Circ\OpGr{b}}{\Set{\OpGr{a},\OpGr{b}}\Sin\Gr{P}}\Period\ Then
\DI{8}{0}{\ForAllSuch{\Set{\OpGr{a},\OpGr{b},\OpGr{c},\OpGr{d}}\Sin\Gr{P}}{\ThereIs{s\In S}\quad \OpGr{a}\Circ\OpGr{b}\Of{s}\Equals\OpGr{c}\Circ\OpGr{d}\Of{s}}\quad\OpGr{a}\Circ\OpGr{b}\Equals\OpGr{c}\Circ\OpGr{d}\Comma}{4}
                                                   \DI{8}{0}{\ForAll{\Set{a,b,c}\Sin S}\ThereIsShriek{$_a$\OpGr{th}$_b$\In\Gr{TH}}\quad$_a$\OpGr{th}$_b$\Of{a}\Equals b}{5}
                                                                             \DII{8}{8}{\Gr{TH}\ is an abelian group under composition.}{and}{6}
 \Proof If \OpGr{a}\Circ\OpGr{b}\Of{s}\Equals\OpGr{c}\Circ\OpGr{d}\Of{s}, then \OpGr{c}\Circ\OpGr{a}\Circ\OpGr{b}\Of{s}\Equals\OpGr{d}\Of{s}\ and so \pLinkLocalDisplayText{3} and \pLinkLocalDisplayText{2} imply that \OpGr{c}\Circ\OpGr{a}\Circ\OpGr{b}\Equals\OpGr{d}\Period\ This \pLinkLocalDisplayText{4} holds.
 \PaR Let \OpGr{a}\ be the function in \Gr{P}\ which leaves a fixed and \OpGr{b}\ be the one which sends a to b. Then \OpGr{b}\Circ\OpGr{a}\Of{a}\Equals b\Period\ That \OpGr{b}\Circ\OpGr{a}\ is unique with this property follows from \pLinkLocalDisplayText{4}, which proves \pLinkLocalDisplayText{5}.
 \PaR For \Set{\OpGr{a},\OpGr{b},\OpGr{c},\OpGr{d}}\Sin\Gr{P}\ we have (\OpGr{a}\Circ\OpGr{b})\Circ(\OpGr{c}\Circ\OpGr{d})\Equals(\OpGr{a}\Circ\OpGr{b}\Circ\OpGr{c})\Circ\OpGr{d}\ which is in \Gr{TH}\ by \pLinkLocalDisplayText{3}. That \OpGr{a}\Circ\OpGr{a}\Equals\Id{S}\ follows from \pLinkLocalDisplayText{1}. For \Set{\OpGr{a},\OpGr{b}}\Sin\Gr{P}, we have (\OpGr{a}\Circ\OpGr{b})\Circ(\OpGr{b}\Circ\OpGr{a})\Equals\Id{S}\ by  \pLinkLocalDisplayText{1}. Thus \Gr{TH}\ is a group.
 \PaR For \Set{\OpGr{a},\OpGr{b},\OpGr{c},\OpGr{d}}\Sin\Gr{P}\ we have \Inv{(\OpGr{a}\Circ\OpGr{b}\Circ\OpGr{c})}\Equals\OpGr{c}\Circ\OpGr{b}\Circ\OpGr{a}\Period\ By \pLinkLocalDisplayText{1} and \pLinkLocalDisplayText{2} this implies that
                                                                    \Dc{\OpGr{a}\Circ\OpGr{b}\Circ\OpGr{c}\Equals\OpGr{c}\Circ\OpGr{b}\Circ\OpGr{a}\Period}
 Consequently
  \Dc{(\OpGr{a}\Circ\OpGr{b})\Circ(\OpGr{c}\Circ\OpGr{d})\Equals(\OpGr{a}\Circ\OpGr{b}\Circ\OpGr{c})\Circ\OpGr{d}\Equals(\OpGr{c}\Circ\OpGr{b}\Circ\OpGr{a})\Circ\OpGr{d}\Equals\OpGr{c}\Circ(\OpGr{b}\Circ\OpGr{a}\Circ\OpGr{d})%
                                                      \Equals\OpGr{c}\Circ(\OpGr{d}\Circ\OpGr{a}\Circ\OpGr{b})\Equals(\OpGr{c}\Circ\OpGr{d})\Circ(\OpGr{a}\Circ\OpGr{b})}
 which proves \pLinkLocalDisplayText{6}. \QED
 \Item{LVIII}{Theorem} Let \Gr{P}\ be an inner involution libra on a set S. For all \Set{a,c}\Sin S we shall denote by \LibraF{a}{c}\ the function in \Gr{P}\ which sends a to c and define
                                                                        \DIc{\ForAll{\Set{a,b,c}\Sin S}\quad\LL{a,b,c}\ \Equiv\ \LibraF{a}{c}\Of{b}}{1}
 \PaR Then \LL{,,}\ is a libra operator on S and, for each x\In S,
                                                                   \DIc{\Function sends\OpGr{p}in\Gr{P}to\OpGr{p}\Of{x}inS\end\ is a libra isomorphism.}{2}
 \Proof Let a, b and c be generic elements of S. We have
                                                              \Dc{\LL{a,a,b}\Equals\LibraF{a}{b}\Of{a}\Equals b\Andd\LL{a,b,b}\Equals\LibraF{a}{b}\Of{b}\Equals a}
 by definition, which is just \pLinkDisplayText{LL}{LD}{1}.
 \PaR Let a, b, c, d and e be generic elements of S and let \Gr{TH}\ be as in \pLinkItemText{LL}{LVII}. For all \Set{x,y}\Sin S, let \LibraTH{x}{y}\ be as in \pLinkDisplayText{LL}{LVII}{5}. It follows from \pLinkItemText{LL}{LVII} that
                     \DI{8}{0}{\ForAll{x\In S}\quad\LibraF{a}{x}\Circ\LibraF{b}{x}\Of{b}\Equals a\quad\Implies\quad\ForAll{x\In S}\quad \LibraF{a}{x}\Circ\LibraF{b}{x}\Of{b}\Equals\LibraTH{b}{a}\Of{b}}{3}
           \DII{8}{8}{\kn{-12}\ForAll{x\In S}\quad\LibraF{x}{e}\Circ\LibraF{x}{d}\Of{b}\Equals e\quad\Implies\quad\ForAll{x\In S}\quad \LibraF{x}{e}\Circ\LibraF{x}{d}\Of{b}\Equals\LibraTH{d}{e}\Of{b}\Period}{and}{4}
 \Par Letting u\Equiv\LibraTH{b}{a}\Of{c}\ and v\Equiv\LibraTH{d}{e}\Of{c}, we have
                   \D{8}{0}{\LL{\LL{a,b,c},d,e}\Equals\LL{\LibraF{a}{c}\Of{b},d,e}\Equals\LL{\LibraF{a}{c}\Circ\LibraF{b}{c}\Of{c},d,e}\Eq{by (3)}\LL{\LibraTH{b}{a}\Of{c},d,e}\Equals\LL{u,d,e}\Equals}
  \D{8}{0}{\LibraF{a}{e}\Of{d}\Equals\LibraF{u}{e}\Circ\LibraF{u}{d}\Of{u}\Eq{by(4)}\LibraTH{d}{e}\Equals\LibraTH{d}{e}\Circ\LibraTH{b}{a}\Of{c}\Eq{\hLinkDisplayText{LL}{LVII}{6}}%
                                                 \LibraTH{b}{a}\Circ\LibraTH{d}{e}\Of{c}\Equals\LibraTH{b}{a}\Of{v}\Eq{by (3)}\LibraF{a}{v}\Circ\LibraF{b}{v}\Of{v}\Equals}
              \D{8}{8}{\LibraF{a}{v}\Of{b}\Equals\LL{a,b,v}\Equals\LL{a,b,\LibraTH{d}{e}}\Eq{by (4)}\LL{a,b,\LibraF{c}{e}\Circ\LibraF{c}{d}\Of{c}}\Equals\LL{a,b,\LibraF{c}{e}\Of{d}}\Equals\LL{a,b,\LL{c,d,e}}}
 \Par which establishes that \LL{,,}\ is a libra operator.
 \PaR Let x be in S and \OpGr{a}, \OpGr{b}\ and \OpGr{c}\ be in \Gr{P}. Let c\Equiv\OpGr{c}\Of{x}, b\Equiv\OpGr{b}\Of{c}\ and a\Equiv\OpGr{a}\Of{b}\Period\ Then
              \Dc{\OpGr{a}\Equals\LibraF{a}{b}\Comma\quad\OpGr{b}\Equals\LibraF{b}{c}\Andd\OpGr{c}\Equals\LibraF{c}{x}\Period}
 We have
                                  \D{8}{0}{\LLL{\OpGr{a},\OpGr{b},\OpGr{c}}\Of{x}\Equals\OpGr{a}\Circ\OpGr{b}\Circ\OpGr{c}\Of{x}\Equals\LibraF{a}{b}\Circ\LibraF{b}{c}\Circ\OpGr{c}\Of{x}%
                                                        \Equals\LibraF{a}{b}\Circ\LibraF{b}{c}\Circ\LibraF{c}{x}\Of{x}\Equals a\Equals\LL{a,x,b,b,x,c,c,x,x}\Equals}
                      \D{8}{8}{\LL{\LL{a,x,b},\LL{b,x,c},\LL{c,x,x}}\Equals\LL{\LibraF{a}{b}\Of{x},\LibraF{b}{c}\Of{x},\LibraF{c}{x}\Of{x}}\Equals\LL{\OpGr{a}\Of{x},\OpGr{b}\Of{x},\OpGr{c}\Of{x}}}
 \Par which establishes \pLinkDisplayText{LL}{LD}{2}
 \QED
 \Item{EE}{Example} Let C be the circle libra in the context of \pLinkItemText{ED}{Circle}. Let L be a line in the projective plane \Aii{P}. As in example \pLinkItemText{LL}{EIIL}, the points p of L correspond to involutions \Op{p\Lower{4pt}{\RMBx{7}{C}}}\ of C in such a way that \SetSuch{\Op{p\Lower{4pt}{\RMBx{7}{C}}}}{p\In L}\ comprise an inner involution libra. This induces a libra on the points of L. The first figure below illustrates the induced libra on the line and how computation of the libra operation is independent of the choice of point in C used to compute it:
                                                                       \Figure{-2pt}{LineLibra}{3in}{5pt}{Libra on the Line Induced by a Circle}{FVI}
 \PaR  In view of Theorem \pLinkItemText{LL}{LVIII}, an isomorphic libra is induced on C as well. The next figure illustrates this libra operation on C:
                                                                       \Figure{-2pt}{CircleLibra}{3in}{5pt}{Libra on the Circle Induced by a Line}{FVII}
 \vfill\eject
                                                                                               \Section{CL}{Meridians on the Cube}\xrdef{pageCL}
  \Item{I}{Introduction} Consider a cube
                                                                                               \d{\CubE{}{}{}{}{}{}\quad\Period}
   \Par We shall write \Es{F}\ for the set of faces of the cube. There are six faces \IronCross, \GreenStar, \Spade, \Club, \RedHeart\ and \RedDiamond.
  \PaR We shall denote by \Ru{K}\ the family of permutations of the faces of the cube which preserve the cube's physical structure intact.
  \PaR In addition we shall require other permutations of \Es{F}\ which are not in \Ru{K}.  We shall write the permutation of \Es{F}\ which interchanges the green faces and leaves the others fixed by \GIn\Period\ The involutions \BIn\ and \RIn\ are defined similarly. We shall write \Ru{K}$^+$\ for the smallest group of permutations of \Es{F}\ containing \GIn, \BIn, \RIn\ and all the elements of \Ru{K}.\Foot{The generating set for \Ru{K}$^+$\ may be made much smaller of course.}
 \PaR It turns out that a meridian can be characterized by functions which are \Quotes{balanced} in a certain sense --- and in such a way that any balanced function remains balanced when \Es{F}\ is permuted by any of the functions in \Ru{K}$^+$. It is the purpose of the present section to show how this is done.
 \Item{BD}{Definitions and Notation} Two distinct faces will be said to be \Def{opposing} if they are of the same color. Thus \IronCross\ and \GreenStar\ are opposing faces, \Spade\ and \Club\ are opposing faces, and \RedHeart\ and \RedDiamond\ are opposing faces. For any face \Eb{F}, we shall denote its opposing face by
                                                                                                  \DIc{\WideHat{\Eb{F}}\Period}{1}
 \PaR Two faces will be said to be \Def{adjacent} if they are distinct and not opposing.
 \PaR By a \Def{half} \Def{cube} we shall mean a subset of \Eb{F}\ consisting of three mutually adjacent faces (thus, one face of each of the three colors green, black and red).
 \PaR Let M be a set with cardinality at least 4. By a \Def{pre-\Load} we shall mean a function from any subset of \Es{F}\ into M. A pre-\Load\ is a \Def{\Load} if its domain is all of \Es{F}\Period\  We shall write
                                                                                                          \Dc{\Loads}
 for the family of all loads. A pre-\Load\ which is constant on a half cube will be said to be singular. If a pre-load is not singular, we shall say that it is \Def{regular}.
  \PaR For \Set{a,b,c,d,e,f}\Sin M, we denote the following pre-\Load s as follows:
                                                       \DI{8}{8}{\WWurf{a,b,c,d}\ \Equiv\ \Set{\Pr{\RedHeart,a},\Pr{\RedDiamond,b},\Pr{\Spade,c},\Pr{\Club,d}}\Comma}{2}
                                                \DI{8}{8}{\WWurf{a,b,c,d,e}\ \Equiv\ \Set{\Pr{\RedHeart,a},\Pr{\RedDiamond,b},\Pr{\Spade,c},\Pr{\Club,d},\Pr{\IronCross,e}}}{3}
                            \DII{8}{8}{\kn{-12}\WWurf{a,b,c,d,e,f}\ \Equiv\ \Set{\Pr{\RedHeart,a},\Pr{\RedDiamond,b},\Pr{\Spade,c},\Pr{\Club,d},\Pr{\IronCross,e},\Pr{\GreenStar,f}}\Period}{and}{4}
 \Item{PMD}{Definition and Notation} A sub-family \BLoads\ of \Loads\ will be said to be a \Def{pre-meridian} \Def{family} \Def{of} \Def{loads} \Def{for} M if the following conditions are met:\Foot{In view of the second condition, the effect of the first condition is that, if two adjacent sides carry the same value of a given load, that load is balanced if and only if the load is singular. In particular, each singular load is balanced.}
                      \DI{8}{0}{\ForAllSuch{\Set{a,b,c,d,e,f}\Sin M}{\Card{(\Set{a,b}\Cap\Set{c,d})}\Equals\One}\quad\WWurf{a,b,c,d,e,f}\In\BLoads\Iff\Set{a,b}\Cap\Set{c,d}\Cap\Set{e,f}\NEq\Void\Comma}{1}
                                                   \DI{8}{0}{\ForAll{\Op{x}\In\BLoads}\ForAll{\Op{f}\In\Ru{K}$^+$}\quad\Op{x}\In\BLoads\Iff\Op{x}\Circ\OpGr{f}\In\BLoads}{2}
                                  \DII{8}{-6}{\kn{-12}\ForAllSuch{\Set{a,b,c,d,e}\Sin M}{\WWurf{a,b,c,d,e}\ is regular}\ThereIsShriek{f\In M}\quad\WWurf{a,b,c,d,e,f}\In\BLoads\Period}{and}{3}
 \PaR For \WWurf{a,b,c,d,e}\ such that there exists a unique element of M such that \WWurf{a,b,c,d,e,x}\In\BLoads, we denote this element x as
                                                                                             \DIc{\Quin{a}{b}{c}{d}{e}\Period}{3.5}
 \Par Thus
                                                                                \DIc{\WWurf{a,b,c,d,e,\Quin{a}{b}{c}{d}{e}}\In\BLoads\Period}{4}
 \Par One can visualize the contents of \pLinkLocalDisplayText{3.5} as being attached to all the sides of the cube except the rear one: e is attached to the front \IronCross\Comma\ a is attached to the left side \RedHeart\Comma\ b is attached to the right side \RedDiamond\Comma\ c is attached to the top \Spade\ and d is attached to the bottom \Club\Period\ The operator implicitly defined in \pLinkLocalDisplayText{4} will be called a \Def{meridian} \Def{quinary} \Def{operator}.
 \PaR In view of \pLinkLocalDisplayText{2}, we have
                         \DIc{\Quin{a}{b}{c}{d}{e}\Equals\Quin{b}{a}{c}{d}{e}\Equals\Quin{a}{b}{d}{c}{e}\Equals\Quin{c}{d}{a}{b}{e}\Equals\Quin{d}{c}{a}{b}{e}\Equals\Quin{c}{d}{b}{a}{e}\Period}{4.5}
 \PaR Suppose \Set{a,b,c,d,x}\Sin M and \Set{a,b}\NEq\Set{c,d}. If \Set{a,b}\Cap\Set{c,d}\Equals\Void, then \WWurf{a,b,c,d,x}\ is regular and so \pLinkLocalDisplayText{3} implies that there exists a unique y\In M such that \WWurf{a,b,c,d,x,y}\In\BLoads\Period\ If u\In M, \Set{a,b}\Cap\Set{c,d}\Equals\Set{u} and x\NEq u, then \pLinkLocalDisplayText{1} implies that u is the unique element of M such that \WWurf{a,b,c,d,x,u}\In\BLoads\Period\ It follows that
                         \DIc{\ForAllSuch{\Set{a,b,c,d}\Sin M}{\Set{a,b}\NEq\Set{c,d}}\quad\Function sendsxin(M\Cop(\Set{a,b}\Cap\Set{c,d}))to\Vol{a,b;c,d}\Of{x}\Equiv\Quin{a}{b}{c}{d}{x}inM\end}{5}
 is well-defined.\Foot{The terminology  \Vol{a,b;c,d}\ is identical with that of \pLinkItemText{ED}{MD}. This constitutes a venial solecism however, since we shall see \It{infra}\ that they mean the same when they occur in a common context.}
 \PaR For \Set{a,b}\Sin M we denote
                                                                                      \DIc{\Mab{a}{b}\ \Equiv\ M\Cop\Set{a,b}\Period}{7}
 \par A pre-meridian family \BLoads\ of loads for M will be said to be a \Def{meridian} \Def{family} \Def{of} \Def{loads} \Def{for} M provided that
                    \D{8}{0}{\ForAll{\Set{a,b}\Sin M}\quad\Function sends\Pr{r,s,t}in\Mab{a}{b}\Cross \Mab{a}{b}\Cross \Mab{a}{b}to\MLL{a}{b}{r,t,s}\Equiv\Quin{a}{b}{r}{s}{t}in\Mab{a}{b}\end}
                                                                                                        \DI{-4}{-4}{}{8}
                                                                                \D{0}{8}{\MLL{a}{b}{,,}\ is a libra operator on \Mab{a}{b}}
 and
                 \DIc{\ForAllSuch{\Set{a,b,c,d}\Sin M}{\Set{a,b}\Cap\Set{c,d}\Equals\Void}\ForAll{\Function\Op{x}sendsin\Es{F}toinM\end}\quad\Op{x}\In\BLoads\Iff\Vol{a,b;c,d}\Circ\Op{x}\In\BLoads\Period}{9}
 \PaR We note that for a meridian family of loads, it follows from \pLinkLocalDisplayText{4.5} that the libra \Mab{a}{b} specified in \pLinkLocalDisplayText{8}
 is abelian.
 \Item{HT}{Theorem} Let \BLoads\ be a meridian family of loads for M. We define
                                                      \DIc{\Mii{M}\ \Equiv\ \SetSuch{\Vol{a,b;c,d}}{\Set{a,b,c,d}\Sin M\Andd\Set{a,b}\Cap\Set{c,d}\Equals\Void}\Period}{1}
 \PaR Then \Mii{M}\ is a meridian family of involutions of M.
 \Proof We first show that the elements of \Mii{M}\ are involutions. Let \Vol{a,b;c,d}\ be in \Mii{M}\ for \Set{a,b,c,d}\Sin M. From \pLinkDisplayText{CL}{PMD}{8}, we have
                                                       \Dc{\ForAll{x\In\Mab{a}{b}}\quad\Vol{a,b;c,d}\Circ\Vol{a,b;c,d}\Of{x}\Equals\MLL{a}{b}{c,\MLL{a}{b}{c,x,d},d}}
 \Par Since \Mab{a}{b}\ is abelian, we have \MLL{a}{b}{c,x,d}\Equals\MLL{a}{b}{d,x,c}\ and so the above becomes
         \DIc{\ForAll{x\In\Mab{a}{b}}\quad\Vol{a,b;c,d}\Circ\Vol{a,b;c,d}\Of{x}\Eq{\hLinkDisplayText{LL}{D4}{1}}\MLL{a}{b}{c,\MLL{a}{b}{d,x,c},d}\Eq{\hLinkDisplayText{LL}{LI}{1}}\MLL{a}{b}{\MLL{a}{b}{c,c,x},d,d}%
                                                                                        \Eq{\hLinkDisplayText{LL}{LD}{1}}\Equals x.}{2}
 \Par In view of \pLinkDisplayText{CL}{PMD}{2} we have
                                                                                           \DIc{\Vol{a,b;c,d}\Equals\Vol{c,d;a,b}}{3}
 and so, as above,
                                                                       \Dc{\ForAll{x\In\Mab{c}{d}}\quad\Vol{a,b;c,d}\Circ\Vol{a,b;c,d}\Of{x}\Equals x}
 which, with \pLinkLocalDisplayText{3}, implies that \Vol{a,b;c,d}\ is an involution.
 \PaR Next we prove \pLinkDisplayText{ED}{MD}{2}. Let \Set{a,b,d,e}\Sin M be such that \Set{a,e}\Cap\Set{b,d}\Equals\Void\Period\ We have
                                                 \DIc{\Vol{a,e;b,d}\Of{b}\Equals\MLL{a}{e}{b,b,d}\Equals d\Andd\Vol{a,e;b,d}\Of{d}\Equals\MLL{a}{e}{b,d,d}\Equals b\Period}{4}
 \Par From \pLinkLocalDisplayText{3} and \pLinkLocalDisplayText{4} follows that
                                                                         \DIc{\Vol{a,e;b,d}\Of{a}\Equals e\Andd\Vol{a,e;b,d}\Of{e}\Equals a\Period}{5}
 \Par Now suppose that \OpGr{f}\In\Mii{M} also satisfies
                                                 \DIc{\OpGr{f}\Of{a}\Equals e,\quad\OpGr{f}\Of{e}\Equals a,\quad\OpGr{f}\Of{b}\Equals d\Andd\OpGr{f}\Of{d}\Equals b\Period}{6}
 \Par Choose \Set{r,s,u,v}\Sin M such that \Set{r,s}\Cap\Set{u,v}\Equals\Void\ and \OpGr{f}\Equals\Vol{r,v;s,u}\Period\ We need to show that \Vol{r,v;s,u}\Equals\Vol{a,e;b,d}\ so, without loss of generality, we may suppose that s\NIn\Set{a,e,b,d}\Period\ As we are dealing with an involution, since \Vol{r,v;s,u}\Of{s}\Equals u, it follows that u\NIn\Set{a,e,b,d}\ as well. We have
                                                    \Dc{\MLL{s}{u}{a,a,e}\Equals e\Andd\MLL{s}{u}{r,a,v}\Equals\Vol{r,v;s,u}\Of{a}\Equals\OpGr{f}\Of{a}\Equals e\Period}
 \Par By \pLinkDisplayText{LL}{LV}{2} it follows that \MLL{s}{u}{a,x,e}\Equals\MLL{s}{u}{r,x,v}\ for all x\In M. This just means that \OpGr{f}\Equals\Vol{r,v;s,u}\Equals\Vol{a,e;s,u}. In similar manner we show that \Vol{a,e;b,d}\Equals\Vol{a,e;s,u}. Putting these both together we obtain that
                                                                                     \Dc{\OpGr{f}\Equals\Vol{r,v;s,u}\Equals\Vol{a,e;b,d}}
  which establishes \pLinkDisplayText{ED}{MD}{2}.
  \PaR Now we turn to \pLinkDisplayText{ED}{MD}{3}. Let \Set{a,d}\Sin M and let \OpGr{a}, \OpGr{b}\ and \OpGr{c}\ be elements of \Mii{M}, each of which sends a to d. Let b\In M be distinct from a and define
                                                           \Dc{e\ \Equiv\ \OpGr{a}\Of{b}\Comma\quad m\ \Equiv\ \OpGr{b}\Of{b}\Andd n\ \Equiv\ \OpGr{c}\Of{b}\Period}
  \Par By \pLinkDisplayText{ED}{MD}{2} we know that
                                                        \Dc{\OpGr{a}\Equals\Vol{a,d;b,e}\Comma\quad\OpGr{b}\Equals\Vol{a,d;b,m}\Andd\OpGr{c}\Equals\Vol{a,d;b,n}\Period}
 \Par For x\In M we have
                    \Dc{\OpGr{b}\Circ\OpGr{c}\Of{x}\Equals\Quin{a}{d}{b}{m}{\Quin{a}{d}{b}{n}{x}}\Equals\MLL{a}{d}{b,\MLL{a}{d}{b,x,n},m}\Equals\MLL{a}{d}{\MLL{a}{d}{b,b ,x},n,m}\Equals\MLL{a}{d}{x,n,m}}
 whence follows that
            \Dc{\OpGr{a}\Circ\OpGr{b}\Circ\OpGr{c}\Of{x}\Equals\Vol{a,d;b,e}\Of{\MLL{a}{d}{x,n,m}}\Equals\MLL{a}{d}{b,\MLL{a}{d}{x,n,m},e}\Eq{\hLinkDisplayText{LL}{LI}{1}\ and \hLinkDisplayText{LL}{LD}{1}}%
                                                                                           \MLL{a}{d}{b,x,\MLL{a}{d}{n,m,e}}\Period}
\Par If c\Equiv\MLL{a}{d}{n,m,e}, this just means that \OpGr{a}\Circ\OpGr{b}\Circ\OpGr{c}\Equals\Vol{a,d;b,c}, whence follows the conclusion of \pLinkDisplayText{ED}{MD}{3}.
 \PaR Finally we turn to \pLinkDisplayText{ED}{MD}{1}. Let \OpGr{a}\ and \OpGr{b}\ be elements of \Mii{M}\Period\ Then there exists \Set{a,b,d,e,r,s,u,v}\Sin M such that
                                                                           \Dc{\OpGr{a}\Equals\Vol{a,e;b,d}\Andd\OpGr{b}\Equals\Vol{r,v;s,u}\Period}
 For any x\In M we have
                                                \D{8}{0}{\Vol{a,e;b,d}\Circ\Vol{r,v;s,u}\Of{x}\Equals\Quin{b}{d}{a}{e}{\Quin{s}{u}{r}{v}{x}}\Eq{\hLinkDisplayText{CL}{PMD}{9}}}
                                              \D{8}{8}{\Quin{\Quin{b}{d}{a}{e}{s}}{\Quin{b}{d}{a}{e}{u}}{\Quin{b}{d}{a}{e}{r}}{\Quin{b}{d}{a}{e}{v}}{\Quin{b}{d}{a}{e}{x}}\Period}
 Setting
              \Dc{d\ \Equiv\ \Quin{b}{d}{a}{e}{r}\Comma\quad n\ \Equiv\ \Quin{b}{d}{a}{e}{s}\Comma\quad p\ \Equiv\ \Quin{b}{d}{a}{e}{u}\Comma\quad q\ \Equiv\ \Quin{b}{d}{a}{e}{v}\Andd y\ \Equiv\ \OpGr{a}\Of{x}}
 we obtain
                                                           \Dc{\OpGr{a}\Circ\OpGr{b}\Of{x}\Equals\Quin{n}{p}{m}{q}{y}\Equals\Vol{m,q;n,p}\Circ\OpGr{a}\Of{x}\Period}
 It follows that \OpGr{a}\Circ\OpGr{b}\Circ\OpGr{a}\Equals\Vol{m,q;n,p}\Comma\ which establishes \pLinkDisplayText{ED}{MD}{1}. \QED
 \Item{Ev}{Example} Let L be the example of \pLinkItemText{ED}{LineM}: a line meridian. Then a load \Function\Op{x}sendsin\Es{F}toinM\end\ is in \BLoads\ if, and only if \Set{\Set{\Op{x}\Of{\RedHeart},\Op{x}\Of{\RedDiamond}},\Set{\Op{x}\Of{\Spade},\Op{x}\Of{\Club}},\Set{\Op{x}\Of{\IronCross},\Op{x}\Of{\GreenStar}}}\ is a cubic triple of pairs of points of L.
 \vfill\eject
                                                                                             \Section{W}{Wurfs and the Cross Ratio}\xrdef{pageW}
\Item{I}{Introduction} The \Def{cross ratio} of four elements w, x, y and z of a field F is defined to be
                                                                                           \DIc{\Over{(w-y)(x-z)}{(x-y)(w-z)}\Comma}{1}
  which makes sense whenever no more than two of the elements are identical, and even when the value \Infty\ is permitted. Its use dates back to antiquity, but it seems first in 1847 (by Karl von Staudt) to be seriously considered in the context of a projective line, or what is here called a meridian. Since the cross ratio is invariant under linear fractional transformations, its avatar in a meridian must be invariant under homographies.
 Von Staudt's idea was to consider ordered quadruples of points on a line, regarding any two such quadruples as \Quotes{equal}, provided there exists a homography which maps the coordinates of one to the corresponding coordinates of another. His term for such a \Quotes{quadruple} was a \Quotes{Wurf}. This German word means a \Quotes{throw} or \Quotes{cast} in English. We shall retain the term \Quotes{Wurf} in this paper, in part because of its closer relation to a cube, which we shall find to be of some use in the exposition here -- but since this paper is written in English, we shall drop the upper case \Quotes{W} in favour of \Quotes{w}.
 \PaR In the present section we shall show how the meridian may be viewed in the context of wurfs and the cross ratio.
 \Item{C}{The Cube} The word \Quotes{Wurf} in German evokes its derivative \Quotes{W\"urfel}, which literally means a die\Foot{The singular form of the noun \Quotes{dice}.} or derivatively, a cube. In fact, a cube
                                                                                                           \Dc{\Cube}
 \Par offers a suggestive setting to manifest von Staudt's idea. We proceed to analyze a cube in detail.
 \PaR It has 8 \Def{vertices}, 12 \Def{edges} and 6 \Def{faces}. Each vertex \Eb{p}\ is connected to three other \Def{adjacent} \Def{vertices} directly via edges, and to a fourth vertex, its \Def{opposite} \Def{vertex} \Eb{o}, which is unique in the sense that \Eb{p}\ and \Eb{o}\ are on no common face. We shall call this set consisting of the three vertices adjacent to \Eb{p}\ and the vertex opposite to \Eb{p}\ the \Def{cubic} \Def{quadriad} \Def{exclusive} \Def{of} \Eb{p}. The complement of this set will be called the \Def{cubic} \Def{quadriad} \Def{of} \Eb{p}. This latter terminology is justified because, for any element \Eb{q}\ of the cubic quadriad exclusive of \Eb{p}, the cubic quadriad exclusive of \Eb{q}\ is precisely the cubic quadriad including \Eb{p}. Thus the set of vertices of a cube is a disjoint union of the sets of vertices of two cubic quadriads.
 \PaR We shall label the elements of one cubic quadriad \PB{a}, \PB{b}, \PB{c}\ and \PB{d}. The above considerations suggest that we then label the other cubic quadriad as \OB{a}, \OB{b}, \OB{c}\ and \OB{d}\ where \OB{a}\ is the opposite vertex of \PB{a}, \OB{b}\ the opposite vertex of \PB{b}\ and so on. We denote
                              \DIc{\PB{Q}\ \Equiv\ \Set{\PB{a},\PB{b},\PB{c},\PB{d}}\Comma\quad\OB{Q}\ \Equiv\ \Set{\OB{a},\OB{b},\OB{c},\OB{d}}\Andd\Bf{K}\ \Equiv\ \PB{Q}\Cup\OB{Q}\Period}{1}
 \PaR Any two elements of a cubic quadriad determine a face: we shall write \IronCross\ for the face determined such that
                                                                                          \DIc{\PB{a},\PB{b},\OB{c},\OB{d}\In\IronCross\Period}{2}
 \Par Then the face \GreenStar\ opposite to \IronCross\ satisfies
                                                                                    \DIc{\PB{c},\PB{d},\OB{a},\OB{b}\In\GreenStar\Period}{3}
 \Par We define the notation of the other four faces as follows:
          \DIc{\PB{a},\PB{c},\OB{b},\OB{d}\In\Spade,\quad\PB{b},\PB{d},\OB{a},\OB{c}\In\Club,\quad\PB{a},\PB{d},\OB{b},\OB{c}\In\RedDiamond,\quad\PB{b},\PB{c},\OB{a},\OB{d}\In\RedHeart\Period}{4}
                                                               \hskip1.8in\includegraphics[width=1.5in]{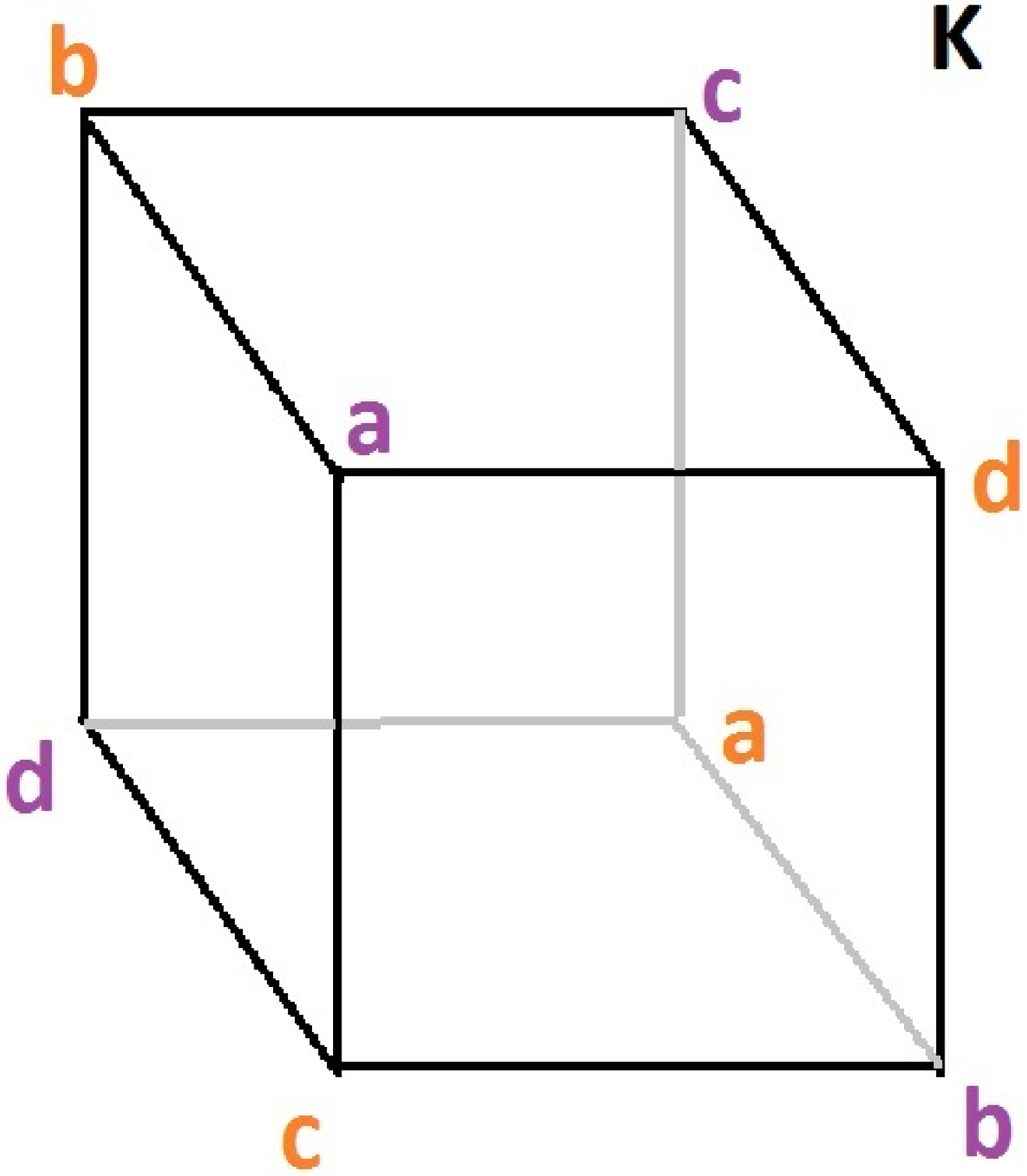}\kn{36}\hbox{\raise50pt\hbox{\CubE{}{}{}{}{}{}}}
 \PaR We now examine the permutations of the vertices of \Bf{K}\ which leave the cube intact. They can be separated into two general types: those which leave the cubic quadriads invariant and those which interchange the cubic quadriads. These permutations may all be realized by rotations around certain axes. By a \Def{vertex} \Def{axis} we mean the a line segment with opposite vertices as endpoints. By an \Def{edge} \Def{axis} we mean a line segment through the midpoints of opposite edges. By a \Def{face} \Def{axis} we mean a line segment with endpoints the centers of opposite faces.
 \PaR By a \Def{face} \Def{involution} we shall mean a rotation of 180$^o$ about a face axis.  The involutions corresponding to the face axes through \IronCross\ and \GreenStar, through \Spade\ and \Club\ and that through \RedHeart\ and \RedDiamond, respectively, will be denoted
                                                                                         \DIc{\GInv\Comma\quad\BInv\Andd\RInv\Comma}{5}
  respectively. These three involutions of the cube leave the cubic quadriads invariant and, along with the identity permutation \Id{\BF{7}{K}}\ of the cube, form a group \Ft{K}$_4$ called the \Def{Klein} \Def{4-group}. A 90$^o$ rotation around a face axis will be called a \Def{face} \Def{rotation}. The face rotation which leaves \IronCross\ and \GreenStar\ invariant and which carries \PB{a}\ to \OB{d}\ will be written
                                                                                                        \DIc{\GadTurn}{6}
  \Par The other six such rotations are defined analogously:
                                                                          \DIc{\GacTurn,\quad\RacTurn,\quad\RabTurn,\quad\BabTurn\Andd\BadTurn\Period}{7}
  \PaR These face rotations interchange cubic quadriads.
  \PaR Rotating the cube 180$^o$ about an edge axis performs what we shall call an \Def{edge} \Def{involution}. If we take the edge axis which pierces the midpoint of the edge joining \PB{a}\ to \OB{b}\ and the midpoint of the opposite edge (which joins \PB{b}\ to \OB{a}), shall write the corresponding edge involution as
                                                                                                     \DIc{\ABInv\Period}{8}
  \Par The other edge involutions are denoted analogously
                                                                             \DIc{\ACInv,\quad\ADInv,\quad\BCInv,\quad\BDInv\Andd\CDInv\Period}{9}
  \PaR We call the vertex axis through \PB{a}\ and \OB{a}
                                                                                                    \DIc{\AVInv\Period}{10}
  \Par The other vertex axes are denoted analogously
                                                                                         \DIc{\BVInv,\quad\CVInv\Andd\DVInv\Period}{11}
  \Par The 120$^o$ rotations about the vertex axes leave the cubic quadriads fixed. The one about the vertex axis \AVInv\ which sends \PB{b}\ to \PB{c}\ (and \OB{b}\ to \OB{c}) will be denoted
                                                                                                     \DIc{\AVbcInv\Period}{12}
  \Par The other such rotations are denoted analogously
                                                          \DIc{\AVcbInv,\quad\BVcdInv,\quad\BVdcInv,\quad\CVdaInv,\quad\CVadInv,\quad\DVabInv\Andd\DVbaInv\Period}{13}
  \PaR In all, there are twenty-four of these transformations of the cube. The group \Ft{K}$_4$ is normal in this larger group of rotational transformations. This means that if we compose each element of \Ft{K}$_4$ with a fixed transformation \OpGr{f}, the resulting set of four elements will be the same whether the composition takes \OpGr{f}\ before or after. Such a set is called a coset of \Ft{K}$_4$, and we shall list \Ft{K}$_4$ and the five other cosets of \Ft{K}$_4$ as a means of describing the transformations of the cube:
                                                                                      \DI{8}{0}{\Set{\OpGr{i},\RInv,\GInv,\BInv}\Comma}{14}
                                                                                 \DI{8}{0}{\Set{\AVbcInv,\BVcdInv,\CVdaInv,\DVbaInv}\Comma}{15}
                                                                                 \DI{8}{0}{\Set{\AVcbInv,\BVdcInv,\CVadInv,\DVabInv}\Comma}{16}
                                                                                    \DI{8}{0}{\Set{\RabTurn,\RacTurn,\BCInv,\ADInv}\Comma}{17}
                                                                                    \DI{8}{0}{\Set{\GacTurn,\GadTurn,\ABInv,\CDInv}\Comma}{18}
                                                                                   \DI{8}{0}{\Set{\BadTurn,\BabTurn,\ACInv,\BDInv}\Period}{19}
 \Par The permutations in each of the sets \pLinkLocalDisplayText{14}, \pLinkLocalDisplayText{15} and \pLinkLocalDisplayText{16} send each cubic quadric to itself, while permutations in the latter three sets interchange them. Thus the union of the first three is a normal subgroup. We shall denote the group of these 24 permutations of the cube by
                                                                                                    \DIc{\Ft{K}\Period}{20}
 \PaR Let \Eu{F}\ denote the set of face axes, \Eu{E}\ the set of edge axes and \Eu{V}\ the set of vertex axes. Each of these sets is left invariant by the cube transformations described above, which fact induces groups of permutations on these sets. Since \Eu{F}\ has cardinality \Three\  and thus is only subject to six permutations, the correspondence from the cube transformations is not injective. The set \Eu{E}\ has six elements and so is subject to 6! permutations: here the correspondence is not surjective. However \Eu{V}\ has cardinality 4 and the correspondence is a bijection of the group of physical transformations of the cube onto the group of permutations of \Eu{V}. For the remainder of this section we shall commit a minor abuse of notation by viewing the elements of \Ft{K}\ as permutations of the set \Eu{V}\ rather than the set \Bf{K}\ of vertices of the cube.
 \Item{W}{Definition} Let M be any set of cardinality at least \Four. By an M-\Def{quadriad}  we shall mean a function \Function\Op{x}sendsin\Eu{V}toinM\end\ such that the cardinality of the range is at least \Three. Now let \Mii{G}\ be any \Three-transitive group of permutations of M. We define a an equivalence relation \Sim\ on the set \Mii{W}\ of all M-quadriads by
                 \DIc{\ForAll{\Set{\Op{x},\Op{y}}\Sin\Mii{W}}\quad\Op{x}\Sim\Op{y}\quad\Iff\quad\ThereIs{\OpGr{s}\In\Ft{K}$_4$\And\OpGr{f}\In\Mii{G}}\Quad\Op{y}\Equals\OpGr{f}\Circ\Op{x}\Circ\OpGr{s}\Period}{1}
 We write the family of all equivalence classes relative to \Sim\ as \Ft{W}\Period\ An element of \Ft{W}\ will be said to be a \Def{wurf}. If one member of a wurf has a range of cardinality \Three , then all members of that wurf have ranges of cardinality \Three . Such wurfs will be called \Def{singular} while all other wurfs will be called \Def{regular}. There are exactly three distinct singular wurfs:
                                     \DI{8}{0}{\GSquare\ \Equiv\ \SetSuch{\Op{x}\In\Mii{W}}{either\quad\Op{x}\Of{\AVInv}\Equals\BVInv\quad or\quad\Op{x}\Of{\CVInv}\Equals\DVInv}\Comma}{2}
                                    \DI{8}{0}{\kn{-4}\BSquare\ \Equiv\ \SetSuch{\Op{x}\In\Mii{W}}{either\quad\Op{x}\Of{\AVInv}\Equals\CVInv\quad or\quad\Op{x}\Of{\BVInv}\Equals\DVInv}}{3}
                             \DII{8}{8}{\kn{-17}\RSquare\ \Equiv\ \SetSuch{\Op{x}\In\Mii{W}}{either\quad\Op{x}\Of{\AVInv}\Equals\DVInv\quad or\quad\Op{x}\Of{\BVInv}\Equals\CVInv}\Period}{and}{4}
 \PaR For \Set{a,b,c,d}\Sin M and \Set{\Rb{b},\Rb{d},\Rb{l},\Rb{q}}\Equals \Eu{V}, we shall at times use the notation
                                            \DIc{\wurf{\Rb{b}}{\Rb{d}}{\Rb{l}}{\Rb{q}}{a}{b}{c}{d}\ \Equiv\ \Set{\Pr{\Rb{b},a},\Pr{\Rb{d},b},\Pr{\Rb{l},c},\Pr{\Rb{q},d}}\Period}{8}
 and the abbreviation
                   \DIc{\Wurf{a,b,c,d}\ \Equiv\ \Set{\Pr{\AVInv\thinspace,\thinspace a},\Pr{\BVInv\thinspace,\thinspace b},\Pr{\CVInv\thinspace,\thinspace c},\Pr{\DVInv\thinspace,\thinspace d}}\Period}{9}
 When employing the above notation for a wurf, we write
                                                                                           \DIc{\Wurf{a,b,c,d}$^{\sIm}$}{10}
 for the element of \Ft{W}\ of which \Wurf{a,b,c,d}\ is a member. In terms of this notation we have
                                                                 \DI{8}{0}{\GSquare\ \Equals\ \SetSuch{\Wurf{a,b,c,d}}{a\Equals b \ or \ c\Equals d}\Comma}{11}
                                                                    \DI{8}{0}{\kn{-2}\BSquare\ \Equals\ \SetSuch{\Wurf{a,b,c,d}}{a\Equals c \ or \ b\Equals d}}{12}
                                                             \DII{8}{8}{\kn{-12}\RSquare\ \Equals\ \SetSuch{\Wurf{a,b,c,d}}{a\Equals d \ or \ b\Equals c\Period}}{and}{13}
  \Item{TI}{Theorem} Let \Set{a,b,c}\Sin M be a basis for M (no two elements of \Set{a,b,c}\ are equal). Then
                                                            \DIc{\ForAll{\Ft{x}\In\Ft{W}}\ThereIs{t\In M}\quad\Wurf{a,b,c,t}\In\Ft{x}\Period}{1}
  Furthermore,
                                     \DIc{if \Ft{x}\In\Set{\GSquare,\BSquare,\RSquare}, then the t in (1) is unique.}{2}
 \Proof [Case: \Ft{x}\In\GSquare]\quad Let t\Equiv c. If follows from \pLinkDisplayText{W}{W}{11} that \pLinkLocalDisplayText{1} holds and b is the only value for t for which it does hold.
 \PaR [Case: \Ft{x}\In\BSquare]\quad Let t\Equiv b. If follows from \pLinkDisplayText{W}{W}{12} that \pLinkLocalDisplayText{1} holds and b is the only value for t for which it does hold.
 \PaR [Case: \Ft{x}\In\RSquare]\quad Let t\Equiv a. If follows from \pLinkDisplayText{W}{W}{13} that \pLinkLocalDisplayText{1} holds and b is the only value for t for which it does hold.
 \PaR [Case: \Ft{x}\NIn\Set{\GSquare,\BSquare,\RSquare}]\quad Let \Wurf{p,q,r,s}\ be any wurf in \Ft{x}. Since \Wurf{p,q,r,s}\ is regular, \Hol{p}{q}{r}{a}{b}{c}\ is an element of \Mii{G}\Period\ We have
                                               \DIc{\Hol{p}{q}{r}{a}{b}{c}\Circ\Wurf{p,q,r,s}\In\Ft{x}\quad\Implies\Wurf{a,b,c,\Hol{p}{q}{r}{a}{b}{c}\Of{s}}\In\Ft{x}\Period}{3}
 Thus we may let t\Equiv\Hol{p}{q}{r}{a}{b}{c}\Of{s}. \QED
 \Item{DW}{Notation} Let \Mii{G}\ be a \Three-transitive group of permutations on a set M, where M has at least four elements. It follows from Theorem \pLinkItemText{W}{TI} that, for any basis \Set{x,y,z}\ for M, the function
                                                                \DIc{\Function sendstinMto\Wurf{x,y,z,t}$^{\sIm}$in\Ft{W}\end\ is surjective.}{1}
 We shall shall denote the function of \pLinkLocalDisplayText{1} by
                                                                                                \DIc{\Wuerf{x}{y}{z}\  \Period}{2}
   \Item{TII}{Theorem} Let \Mii{G}\ be a \Three-transitive group of permutations on a set M, where M has at least four elements. Then the following statements are equivalent:
                                                   \DI{8}{0}{\ForAll{\Set{a,b,c}\ a basis in M}\quad \Function\Wuerf{a}{b}{c}sendsinMtoin\Ft{W}\end\quad is a bijection,}{1}
                         \DII{8}{8}{\ForAllSuch{\OpGr{f}\In\Mii{G}}{\ThereIs{\Set{a,b}\Sin M}\ a\NEq b, \OpGr{f}\Of{a}\Equals b\And\OpGr{f}\Of{b}\Equals a}\quad\OpGr{f}\In\Mii{G}$_2$\Period}{and}{2}
 \Proof [(1)\Implied(2)]\quad Suppose first that \pLinkLocalDisplayText{2} is satisfied, and let \Set{s,t}\Sin M be such that
                                                                             \DIc{\Wuerf{a}{b}{c}\Of{s}\ \Equals\ \Wuerf{a}{b}{c}\Of{t}\Period}{2a}
 If\quad\Wuerf{a}{b}{c}\Of{t}\In\Set{\GSquare,\BSquare,\RSquare}, it follows from \pLinkDisplayText{W}{TI}{2} that s\Equals t. If\quad\Wuerf{a}{b}{c}\Of{t}\NIn\Set{\GSquare,\BSquare,\RSquare}, then \break \Set{s,t}\Cap\Set{a,b,c}\Equals\Void, and so we shall presume that this is so. There exist \OpGr{f}\In\Mii{G}\ and \OpGr{s}\In\Ft{K}$_4$ such that
                                                                         \DIc{\Wurf{a,b,c,s}\Equals\OpGr{f}\Circ\Wurf{a,b,c,t}\Circ\OpGr{s}\Period}{3}
 If \OpGr{s}\Equals\Id{\EUBx{7}{V}}, then \OpGr{f}\And\Id{M}\ agree on three distinct points, whence follows that they are identical, and so
                                                                                        \DIc{s\Equals\OpGr{f}\Of{t}\Equals t\Period}{4}
  \PaR If \OpGr{s}\Equals\GInv, then \pLinkLocalDisplayText{3} implies
                                            \DIc{\Wurf{a,b,c,s}\Equals\OpGr{f}\Circ\Wurf{b,a,t,c}\Equals\Wurf{\OpGr{f}\Of{b},\OpGr{f}\Of{a},\OpGr{f}\Of{t},\OpGr{f}\Of{c}}\Period}{5}
 It follows in particular that a\Equals\OpGr{f}\Of{b}\And b\Equals\OpGr{f}\Of{a}, which by \pLinkLocalDisplayText{2} implies that \OpGr{f}\Equals\Inv{\OpGr{f}}\Period\ Consequently from \pLinkLocalDisplayText{5} we have
                                                          \DIc{c\Equals\OpGr{f}\Of{t}\quad\Implies\quad t\Equals\OpGr{f}\Of{c}\quad\Implies\quad t\Equals s\Period}{6}
 \PaR If \OpGr{s}\Equals\BInv\ or \OpGr{s}\Equals\RInv, arguments analogous to that just used show that s\Equals t in each case. Hence, in view of \pLinkLocalDisplayText{6}, it follows that then function \Wuerf{a}{b}{c}\ is injective, and so bijective.\vskip 2pt
 \PaR [(1)\Implies(2)]\quad We now suppose that \pLinkLocalDisplayText{1} holds. Let a and b be distinct elements of M and \OpGr{f}\ an element of \Mii{G}\ such that
                                                                              \DIc{\OpGr{f}\Of{a}\Equals b\Andd\OpGr{f}\Of{b}\Equals a\Period}{7}
 Let t be any element of M not in \Set{a,b}\Period\ Let c\Equiv\OpGr{f}\Of{t} and assume that t\NEq\OpGr{f}\Of{c}\Period\ Then
                                  \D{8}{0}{\Wuerf{a}{b}{c}\Of{t}\Equals\Wurf{a,b,c,t}$^{\sIm}$\Equals(\OpGr{f}\Circ\Wurf{a,b,c,t})$^{\sIm}$\Equals\Wurf{b,a,\OpGr{f}\Of{c},c}$^{\sIm}$\Equals}
                                         \D{8}{8}{(\Wurf{b,a,\OpGr{f}\Of{c},c}\Circ\GInv)$^{\sIm}$\Equals\Wurf{a,b,c,\OpGr{f}\Of{c}}$^{\sIm}$\Equals\Wuerf{a}{b}{c}\Of{\OpGr{f}\Of{c}}}
 which violates \pLinkLocalDisplayText{1}. It follows that t\Equals\OpGr{f}\Of{c}, which implies that \OpGr{f}\In\Mii{G}$_2$. \QED
 \Item{CI}{Corollary} Let M be a set with at least four elements and let \Mii{G}\ be a \Three-transitive group of permutations of M such that \pLinkDisplayText{W}{TII}{2} holds. Define
                                                    \DIc{\Ft{G}\ \Equiv \SetSuch{\Wuerf{a}{b}{c}\Circ\InvWuerf{u}{v}{w}}{\Set{a,b,c}\And\Set{u,v,w}\ are bases for M}}{1}
 Then \Ft{G}\ is a group isomorphic with \Mii{G} and M is isomorphic with \Ft{W}\ in the sense that there exists a bijection \Function\OpGr{ps}sendsinMtoin\Ft{W}\end\ such that
                                                             \DIc{\Mii{G}\ \Equals\ \SetSuch{\Inv{\OpGr{ps}}\Circ\Ft{f}\Circ\OpGr{ps}}{\Ft{f}\In\Ft{G}}\Period}{2}
 \Proof Let a, b and c be pairwise distinct elements of M and define \OpGr{ps}\ to be \Wuerf{a}{b}{c}\Period\ \QED
  \Item{CII}{Discussion} One readily checks that, for any \OpGr{f}\ in \Ft{K},
                                                    \DIc{\ForAll{\Op{x}\And\Op{y}\ M-quadriads}\quad\OpGr{f}\Of{\Op{x}}\Sim\OpGr{f}\Of{\Op{y}}\Iff\Op{x}\Sim\Op{y}\Period}{1}
  \Par Thus, to each \OpGr{f}\In\Ft{K}\ corresponds a unique function
         \DIc{\Function\Stack{\sIm}{\OpGr{f}}{1pt}sendsin\Ft{W}toin\Ft{W}\end\quad such that \ForAll{\Ft{w}\In\Ft{W}}\ForAll{\Op{x}\In\Ft{w}}\quad\Op{x}\Circ\OpGr{f}\In\Stack{\sIm}{\OpGr{f}}{1pt}\Of{\Op{x}}\Period}{2}
  \Par It follows from the definition of the equivalence relation \Sim\ that
                                                             \DIc{\ForAll{\OpGr{s}\In\Ft{K}$_4$}\quad\Stack{\sIm}{\OpGr{s}}{1pt}\Equals\Id{\FTBx{7}{W}}\Period}{3}
  \Par Consequently, if \OpGr{a}\ and \OpGr{b}\ are elements of the same coset (where the cosets are given in \pLinkDisplayText{W}{C}{14}-\pLinkDisplayText{W}{C}{19}, then
                                                                         \DIc{\Stack{\sIm}{\OpGr{a}}{1pt}\Equals\Stack{\sIm}{\OpGr{b}}{1pt}\Period}{4}
  Thus, if we choose one element from each coset, we shall obtain a group of permutations \Ft{W}:
                                                                          \DIc{\Set{\OpGr{i},\AVbcINV,\AVcbINV,\RABTURN,\GACTURN,\BADTURN}\Period}{5}
  For reference we compute for any M-quadriad \Wurf{a,b,c,d}
                                           \DI{8}{0}{\Wurf{a,b,c,d}\Circ\AVbcInv\Equals\Wurf{a,c,d,b}\Andd\AVbcINV\Of{\Wurf{a,b,c,d}$^{\sIm}$}\Equals\Wurf{a,c,d,b}$^{\sIm}$\Comma}{6}
                                           \DI{8}{0}{\Wurf{a,b,c,d}\Circ\AVcbInv\Equals\Wurf{a,d,b,c}\Andd\AVcbINV\Of{\Wurf{a,b,c,d}$^{\sIm}$}\Equals\Wurf{a,d,b,c}$^{\sIm}$\Comma}{7}
                                           \DI{8}{0}{\Wurf{a,b,c,d}\Circ\RabTurn\Equals\Wurf{b,d,a,c}\Andd\RABTURN\Of{\Wurf{a,b,c,d}$^{\sIm}$}\Equals\Wurf{b,d,a,c}$^{\sIm}$\Comma}{8}
                                          \DI{8}{0}{\Wurf{a,b,c,d}\Circ\GacTurn\Equals\Wurf{c,d,b,a}\Andd\RABTURN\Of{\Wurf{a,b,c,d}$^{\sIm}$}\Equals\Wurf{c,d,b,a}$^{\sIm}$\Comma}{9}
                                  \DII{8}{8}{\kn{-8.5}\Wurf{a,b,c,d}\Circ\BadTurn\Equals\Wurf{d,a,b,c}\Andd\BADTURN\Of{\Wurf{a,b,c,d}$^{\sIm}$}\Equals\Wurf{d,a,b,c}$^{\sIm}$\Period}{and}{10}
  Also for reference we compute for any M-quadriad \Wurf{a,b,c,d}
                                                                               \DI{8}{0}{\Wurf{a,b,c,d}\Circ\RInv\Equals\Wurf{d,a,c,b}\Comma}{11}
                                                                                  \DI{8}{0}{\Wurf{a,b,c,d}\Circ\BInv\Equals\Wurf{c,d,a,b}}{12}
                                                                        \DII{8}{0}{\kn{-12}\Wurf{a,b,c,d}\Circ\GInv\Equals\Wurf{b,a,d,c}\Period}{and}{13}
  \Item{TV}{Theorem} Let M be a set with at least four elements and let \Mii{G}\ be a \Three-transitive group of permutations of M such that \pLinkDisplayText{W}{TII}{2} holds. Then
                                                                      \DIc{\ForAll{\OpGr{f}\In\Ft{W}}\quad\Stack{\sIm}{\OpGr{f}}{1pt}\In\Ft{G}\Period}{1}
\Proof We shall prove the theorem for \OpGr{f}\Equals\BADTURN. The proofs for the other cases are analogous.
\PaR Let \Set{a,b,c}\ be a basis for M. Let \Ft{x}\ be any element of \Ft{W}\Period\ By \pLinkDisplayText{W}{TII}{1} there exists a unique t such that \Wurf{a,b,c,t}\In\Ft{x}\Period\ We have
                              \D{8}{0}{\BADTURN\Of{\Ft{x}}\Equals\BADTURN\Of{\Wurf{a,b,c,t}$^{\sIm}$}\Eq{\hLinkDisplayText{W}{CII}{10}}\Wurf{t,a,b,c}$^{\sIm}$\Eq{\hLinkDisplayText{W}{CII}{11}}}
                             \D{8}{0}{\Wurf{t,a,b,c}$^{\sIm}$\Circ\RInv\Equals\Wurf{c,b,a,t}$^{\sIm}$\Equals(\Hol{a}{b}{c}{c}{b}{a}\Circ\Wurf{a,b,c,\Hol{a}{b}{c}{c}{b}{a}\Of{t}})$^{\sIm}$\Equals}
          \D{8}{8}{\Wuerf{c}{b}{a}\Circ\Hol{a}{b}{c}{c}{b}{a}\Circ\Inv{(\Wuerf{c}{b}{a})}\Of{\Wurf{a,b,c,t})$^{\sIm}$}\Equals\Wuerf{c}{b}{a}\Circ\Hol{a}{b}{c}{c}{b}{a}\Circ\Inv{(\Wuerf{c}{b}{a})}\Of{\Ft{x}}\Period}
Consequently \BADTURN\Equals\Wuerf{c}{b}{a}\Circ\Hol{a}{b}{c}{c}{b}{a}\Circ\Inv{(\Wuerf{c}{b}{a})}\Period\ It follows from \pLinkDisplayText{W}{CI}{2} that the latter is in \Ft{G}. \QED
  \Item{HT}{Theorem} Let M be a set with at least four elements and let \Mii{G}\ be a \Three-transitive group of permutations of M. Then necessary and sufficient conditions for \Mii{G}\ to be a meridian group of permutations are
                                \DI{8}{8}{\ForAllSuch{\OpGr{f}\In\Mii{G}}{\ThereIs{\Set{a,b}\Sin M}\ a\NEq b, \OpGr{f}\Of{a}\Equals b\And\OpGr{f}\Of{b}\Equals a}\quad\OpGr{f}\In\Mii{G}$_2$}{1}
                                      \DII{0}{8}{\ThereIsShriek{\OpGr{w}\In\Mii{G}$_4$}\quad\OpGr{w}\Of{\RSquare}\Equals\GSquare\Andd\OpGr{w}\Of{\GSquare}\Equals\BSquare\Period}{and}{2}
  \Proof Suppose first that \Mii{G}\ is a meridian group of permutations of M. Recalling that M is meridian isomorphic with \Ft{W}, we see that \pLinkLocalDisplayText{1} holds follows from \pLinkDisplayText{E}{MD}{2} and that \pLinkLocalDisplayText{2} holds follows from \pLinkDisplayText{E}{MD}{1}.
  \PaR Now suppose that conditions \pLinkLocalDisplayText{1} and \pLinkLocalDisplayText{2} hold. That \pLinkDisplayText{E}{MD}{2} holds follows from \pLinkLocalDisplayText{1}, so we need only establish \pLinkDisplayText{E}{MD}{2}. Let then a,b and c be distinct elements of M. Let \OpGr{a}\ be a meridian isomorphism from M onto \Ft{W}. Since \Mii{G}\ is \Three-transitive, there exists an element \OpGr{c}\In\Mii{G}\ such that
                        \DIc{\OpGr{c}\Of{a}\Equals\Inv{\OpGr{a}}\Of{\RSquare}\Comma\quad\OpGr{c}\Of{b}\Equals\Inv{\OpGr{a}}\Of{\GSquare}\Andd\OpGr{c}\Of{c}\Equals\Inv{\OpGr{a}}\Of{\BSquare}\Period}{3}
  \Par Letting \OpGr{f}\Equiv\Inv{\OpGr{c}}\Circ\OpGr{w}\Circ\OpGr{c}\Comma\ we see from \pLinkLocalDisplayText{3} that \pLinkDisplayText{E}{MD}{1} holds. \QED
\Item{DC}{Definition} Let M be a set with at least four elements and let \Mii{G}\ be a meridian group of permutations of M. For any basis \Set{u,v,w}\ of M the function
                                                   \DIc{\Function sends\Pr{a,b,c,d}in\Mii{W}to\InvWuerf{u}{v}{w}\Of{\Wurf{a,b,c,d}$^{\sIm}$}inM\end}{1}
 will be called a \Def{cross} \Def{ratio} \Def{on} M.
  \Item{TCR}{Theorem} Let \OpGr{k}\ be a cross ratio on M, where M is a set of cardinality at least 4 and \Mii{G} is a meridian group of permutations of M. Let \Function\OpGr{f}sendsinMtoinM\end\ be a bijection. Then the following two statements are equivalent:
                                                                                           \DI{8}{4}{\OpGr{f}\In\Mii{G}\Semicolon}{1}
                         \DI{4}{8}{\ForAllSuch{\Set{a,b,c,d}\Sin M}{\Wurf{a,b,c,d}\In\Mii{W}}\quad\OpGr{k}\Of{\Pr{a,b,c,d}}\Equals\OpGr{k}\Of{\Pr{\OpGr{f}\Of{a},\OpGr{f}\Of{b},\OpGr{f}\Of{c},\OpGr{f}\Of{d}}}\Period}{2}
 \Proof That (1) implies (2) follows directly from the definition of \Sim\Period\
 \PaR Suppose that \pLinkLocalDisplayText{2} holds. Let \Set{u,v,w} be the basis of M such that
                                                   \DIc{\OpGr{k}\Of{\Pr{a,b,c,d}}\Equals\InvWuerf{u}{v}{w}\Circ\Wurf{a,b,c,d}\ for every M-quadriad \Set{a,b,c,d}\Period}{3}
 Let t be any element of M. Then
 \D{8}{0}{\Wurf{u,v,w,t}\kn{4}\Sim\kn{4}\Hol{u}{v}{w}{\OpGr{f}\Of{u}}{\OpGr{f}\Of{v}}{\OpGr{f}\Of{w}}\Circ\Wurf{u,v,w,t}\Equals}
 \DI{8}{0}{\Wurf{\Hol{u}{v}{w}{\OpGr{f}\Of{u}}{\OpGr{f}\Of{v}}{\OpGr{f}\Of{w}}\Of{u},\Hol{u}{v}{w}{\OpGr{f}\Of{u}}{\OpGr{f}\Of{v}}{\OpGr{f}\Of{w}}\Of{v},\Hol{u}{v}{w}{\OpGr{f}\Of{u}}{\OpGr{f}\Of{v}}{\OpGr{f}\Of{w}}\Of{w},%
                                                                        \Hol{u}{v}{w}{\OpGr{f}\Of{u}}{\OpGr{f}\Of{v}}{\OpGr{f}\Of{w}}\Of{t}}\Equals}{4}
                                            \D{8}{8}{\Wurf{\OpGr{f}\Of{u},\OpGr{f}\Of{v},\OpGr{f}\Of{w},\Hol{u}{v}{w}{\OpGr{f}\Of{u}}{\OpGr{f}\Of{v}}{\OpGr{f}\Of{w}}\Of{t}}\Period}
 From \pLinkLocalDisplayText{4} follows that
                            \DIc{\Wurf{u,v,w,t}$^{\sIm}$\Equals\Wurf{\OpGr{f}\Of{u},\OpGr{f}\Of{v},\OpGr{f}\Of{w},\Hol{u}{v}{w}{\OpGr{f}\Of{u}}{\OpGr{f}\Of{v}}{\OpGr{f}\Of{w}}\Of{t}}$^{\sIm}$}{5}
 On the other hand  we have
                \DIc{\OpGr{k}\Of{\Pr{\OpGr{f}\Of{u},\OpGr{f}\Of{v},\OpGr{f}\Of{w},\OpGr{f}\Of{t}}}\Eq{\hLinkLocalDisplayText{2}}\OpGr{k}\Of{\Pr{u,v,w,t}}\quad\Implies\quad\Wurf{\OpGr{f}\Of{u},\OpGr{f}\Of{v},%
                                                                  \OpGr{f}\Of{w},\OpGr{f}\Of{t}}$^{\sIm}$\Eq{\hLinkItemText{W}{DC}}\Wurf{u,v,w,t}$^{\sIm}$}{6}
 Let a\Equiv\OpGr{f}\Of{u}, b\Equiv\OpGr{f}\Of{v}\ and c\Equiv\OpGr{f}\Of{w}\Period\ From \pLinkLocalDisplayText{5} and \pLinkLocalDisplayText{6} we have
 \D{8}{8}{\Wuerf{a}{b}{c}\Of{\Hol{u}{v}{w}{\OpGr{f}\Of{u}}{\OpGr{f}\Of{v}}{\OpGr{f}\Of{w}}\Of{t}}\Equals\Wurf{\OpGr{f}\Of{u},\OpGr{f}\Of{v},\OpGr{f}\Of{w},\Hol{u}{v}{w}{\OpGr{f}\Of{u}}{\OpGr{f}\Of{v}}{\OpGr{f}\Of{w}}\Of{t}}%
                                                                                                       $^{\sIm}$\Equals}%
                                             \D{0}{8}{\Wurf{\OpGr{f}\Of{u},\OpGr{f}\Of{v},\OpGr{f}\Of{w},\OpGr{f}\Of{t}}$^{\sIm}$\Equals\Wuerf{a}{b}{c}\Of{\OpGr{f}\Of{t}}\Period}
 Since \Wuerf{a}{b}{c}\ is bijective, it follows that \OpGr{f}\Of{t}\Equals\Hol{u}{v}{w}{\OpGr{f}\Of{u}}{\OpGr{f}\Of{v}}{\OpGr{f}\Of{w}}\Of{t}\Period\ That is to say, \OpGr{f}\ is an element of \Mii{G}\Period\ \QED
 \Item{CRT}{Theorem} Let \Mii{G}\ be a meridian family of permutations on a set M with cardinality at least 4. Let \Set{\Zero,\One,\Infty}\ be a basis for M. Then
                      \DIc{\ForAllSuch{\Set{w,x,y,z}\Sin M}{\Card{\Set{w,x,y,z}}$\ge$3}\quad \InvWuerf{\One}{\Zero}{\iNfty}\Of{\Wurf{y,x,w,z}$^{\sIm}$}\ \Equals\ \Over{(w-y)\Cdot(x-z)}{(x-y)\Cdot(w-z)}}{1}
where the binary operations are the field operations as in \pLinkItemText{ED}{T}.
 \Proof The function \Wuerf{\One}{\Zero}{\iNfty}\ is the projective mapping from M onto \Ft{W}\ satisfying
                                        \Dc{\Wuerf{\One}{\Zero}{\iNfty}\Of{\Zero}\Equals\BSquare,\quad\Wuerf{\One}{\Zero}{\iNfty}\Of{\One}\Equals\RSquare\Andd\Wuerf{\One}{\Zero}{\iNfty}\Of{\Infty}\Equals\GSquare}
 while
                                                  \Dc{\Wuerf{y}{x}{w}\Of{y}\Equals\BSquare,\quad\Wuerf{y}{x}{w}\Of{x}\Equals\RSquare\Andd\Wuerf{y}{x}{w}\Of{w}\Equals\GSquare}
 whence follows that
                                \D{8}{0}{\InvWuerf{\One}{\Zero}{\iNfty}\Of{\Wurf{y,x,w,z}$^{\sIm}$}\Of{w}\Equals \Infty,\quad\InvWuerf{\One}{\Zero}{\iNfty}\Of{\Wurf{y,x,w,z}$^{\sIm}$}\Of{x}\Equals \Zero}%
                                                                                                        \DI{-4}{-4}{}{2}
                                                                  \DL{0}{8}{\InvWuerf{\One}{\Zero}{\iNfty}\Of{\Wurf{y,x,w,z}$^{\sIm}$}\Of{y}\Equals\ \One\Period}{and}
 In terms of the field operations, the element \InvWuerf{\One}{\Zero}{\iNfty}\Of{\Wurf{y,x,w,z}$^{\sIm}$}\ becomes a linear fractional transformation. The linear fractional transformation satisfying \pLinkLocalDisplayText{2} evidently is \pLinkLocalDisplayText{1}, viewed as a function of z. \QED
 \Item{HW}{Harmonic Wurfs} Let \Mii{G}\ be a meridian family of permutations on a set M with cardinality at least 4. We know from \pLinkItemText{W}{HT} that \GINV\ is an involution which fixes the point \GSquare\Period\ From \pLinkItemText{ED}{TTP} follows that \GINV\ has a second fixed point which we shall call \GTriangle\Colon
                                                                                  \DIc{\GINV\Of{\GTriangle}\ \Equals\ \GTriangle\Period}{1}
 One can easily check that
                                                             \DIc{\GTriangle\ \Equals\ \SetSuch{\Wurf{a,b,c,d}}{\Set{\Set{a,d},\Set{b,c}}\ is a harmonic pair.}}{2}
 Furthermore, if \OpGr{w}\ is as in Theorem \pLinkItemText{W}{HT}, then
                                                               \DIc{\OpGr{w}\Of{\BSquare}\Equals\GTriangle\Andd\OpGr{w}\Of{\GTriangle}\Equals\RSquare\Period}{3}
 \PaR Similarly, there are sets \RTriangle\And\BTriangle\ such that
                 \DIc{\RINV\Of{\RTriangle}\ \Equals\ \RTriangle\ \Equals\ \SetSuch{\Wurf{a,b,c,d}}{\Set{\Set{a,c},\Set{b,d}}\ is a harmonic pair}}{4}
 and
                                            \DIc{\BINV\Of{\BTriangle}\ \Equals\ \BTriangle\ \Equals\ \SetSuch{\Wurf{a,b,c,d}}{\Set{\Set{a,b},\Set{c,d}}\ is a harmonic pair.}}{5}
\Item{Rem}{Remarks} Let M be a meridian. Then \Ft{W}\ is a meridian isomorphic with M which has three distinguished points: \GSquare\Comma\ \BSquare\And\RSquare\Period\ In other words, \Ft{W}\ comes equipped with a basis and so is a sort of intrinsic representation of M which \Quotes{almost}\ has a distinguished basis, and so \Quotes{almost}\ has a distinguished field.
\Item{Ex}{Example} The set \Set{\GTriangle,\BTriangle,\RTriangle}\ may not always have cardinality \Three . It is not hard to show that the cardinality is \One\ when it is not \Three. With the sphere, circle (and line) meridians, it has cardinality \Three. For the meridian which corresponds to the field \Set{-\One,\Zero,\One}, with ordinary multiplication and addition, it has cardinality \One.
\vfill\eject
                                                                                          \Section{EA}{Meridian Exponentials and Arcs}\xrdef{pageEA}
\Item{I}{Introduction} It could be asserted with some justification that the circle or line meridian is the most useful meridian. That being so, it is reasonable to inquire into what distinguishes it from the others.
\PaR Suppose then that M is a circle meridian, and let \Pr{\Zero,\One,\Infty}\ be any ordered basis for M. We may regard \Eu{F}\Equiv M\Cop\Set{\Infty}\ as the field of real numbers where the binary operations are defined as in \pLinkItemText{ED}{T}. There is a continuous function
                                                                                             \Dc{\Function\OpGr{w}sendsinFtoinF\end}
\Par such that
                                                        \DIc{\ForAll{\Set{\EuBx{m},\EuBx{n}}\Sin\AiiBx{N}}\quad\OpGr{w}\Of{\Over{\EUBx{7}{m}}{\EUBx{7}{n}}}\ \Equals\ \ %
      \hbox{\lower24pt\hbox{\Vbox{\hbox{$\overbrace{\Two\Cdot\dots\Cdot\Two}^{\RMBx{7}{m times}}$}\kn{6}\hbox{\kn{-2}\vrulewh{35}{1}}\kn{4}\hbox{$\overbrace{\Two\Cdot\dots\Cdot\Two}^{\RMBx{7}{n times}}$}}}}\Period}{0}
\Par This function thus has the property that
                                                            \DIc{\ForAll{\Set{x,y}\Sin F}\quad\OpGr{w}\Of{x+y}\ \Equals\ \OpGr{w}\Of{x}\Cdot\OpGr{w}\Of{y}\Comma}{1}
\Par and is commonly denoted by
                                                                           \DIc{\ForAll{x\In F}\quad\Two$^{\RMBx{7}{x}}$\ \Equiv\ \OpGr{w}\Of{x}}{2}
\PaR It turns out that that the existence of such an \Quotes{exponential} function can be described in terms of the meridian structure and that the existence of such is in a certain sense sufficient, as well as necessary, to render a meridian isomorphic to a circle meridian.
 \Item{II}{Review and Discussion} Let \Mii{M}\ be any meridian family of involutions on a set M of cardinality of at least \Four. For \Set{a,d}\Sin M, the family
                                                                       \DIc{\MAB{a}{d}\ \Equiv\ \SetSuch{\OpGr{f}\In\Mii{M}}{\OpGr{f}\Of{a}\Equals d}}{1}
 is an abelian function libra by definition \pLinkDisplayText{ED}{MD}{3} and consequence \pLinkDisplayText{ED}{MD}{4}. By \pLinkDisplayText{ED}{MD}{2} and by \pLinkItemText{LL}{DV}, \MAB{a}{d}\ is an inner involution libra on
                                                                                       \DIc{\Mab{a}{d}\ \Equiv\ M\Cop\Set{a,d}\Period}{2}
 \PaR The function libra \MAB{a}{d}\ (\It{cf}. \pLinkItemText{LL}{LVIII}) induces a libra operation on The set \Mab{a}{d}:
                                         \DIc{\ForAll{\Set{x,y,z}\Sin\Mab{a}{d}}\quad\LL{x,y,z}\ \Equiv\ \GrowingLBrace{\hbox{\Vbox{\hbox{\Vol{x,x;z,z}\Of{y}\Quad if x\NEq z;}%
                                                                                                                                       \kn{5}%
                                                                                                                                       \hbox{\kn{21}y\kn{21}\Quad if x\Equals z. }}}}}{3}
 \Par and, for each t\In\Mab{a}{d},
                                                   \DIc{\Function\Russ{a}{d}{t}sends\OpGr{f}in\MAB{a}{d}to\OpGr{f}\Of{t}in\Mab{a}{d}\end\quad is a libra isomorphism.}{4}
 \Par When a\Equals d, we shall use the abbreviation:
                                                                                    \DIc{\Rus{d}{t}\ \Equiv\ \Russ{d}{d}{t}\Period}{5}
 \PaR When a\NEq d, it may happen that an element of \MAB{a}{d}\ has no fixed points. Therefore we introduce the notation
                                                   \DIc{\MABF{a}{d}\ \Equiv\ \SetSuch{\OpGr{f}\In\MAB{a}{d}}{\ThereIs{x\In\Mab{a}{d}}\Quad\OpGr{f}\Of{x}\Equals x}\Period}{6}
 \Par It follows from Theorem \pLinkItemText{ED}{DF} that each element of\ \ \MABF{a}{d}\Quad fixes exactly two points. We recall that a \Def{basis} for M is just a subset of M of cardinality \Three, and that an \Def{ordered} \Def{basis} is an ordered triple of which the elements comprise a basis. For an ordered basis \Pr{a,l,d}, we adopt the notation
                                                                    \DIc{\Mabf{a}{l}{d}\ \Equiv\ \SetSuch{\OpGr{f}\Of{l}}{\OpGr{f}\In\MABF{a}{d}}\Period}{7}
 \Par For t\In\Mab{a}{d}\Comma\ we shall write \Russs{a}{l}{d}{t}\ for the restriction of \Russ{a}{d}{t}\ to \Mabf{a}{l}{d}\Colon
                                                         \DIc{\Function\Russs{a}{l}{d}{t}sends\OpGr{f}in\Mabf{a}{l}{d}to\OpGr{f}\Of{t}in\Mab{a}{d}\end\Period}{8}
 \Item{III}{Theorem} Let \Pr{a,l,d}\ be an ordered basis for a meridian M, suppose that \MABF{a}{d}\ is balanced in \MAB{a}{d}\ and let m be an element of \Mabf{a}{l}{d}\Period\ Then
                                                                                  \DIc{\Mabf{a}{l}{d}\ \Equals\ \Mabf{a}{m}{d}\Period}{1}
 \Proof Let k be any element of \Mabf{a}{m}{d}\ and select \OpGr{ps}\In\MABF{a}{d}\ such that \OpGr{ps}\Of{m}\Equals k\Period\ Let \OpGr{f}\In\MABF{a}{d}\QUad be such that \OpGr{f}\Of{l}\Equals m\Period\ By hypothesis we have \Vol{a,d;k,k}\Circ\OpGr{ps}\Circ\OpGr{f}\ in \MABF{a}{d}\Period\ Since
                                            \Dc{\Vol{a,d;k,k}\Circ\OpGr{ps}\Circ\OpGr{f}\Of{l}\Equals\Vol{a,d;k,k}\Circ\OpGr{ps}\Of{m}\Equals\Vol{a,d;k,k}\Of{k}\Equals k\Comma}
 \Par we have shown that
                                                                                      \DIc{\Mabf{a}{m}{d}\Sin\Mabf{a}{l}{d}\Period}{2}
 \PaR Since \OpGr{f}\Of{m}\Equals l\Comma\ we have l\In\Mabf{a}{m}{d}\ and so, interchanging the roles of l and m, \pLinkLocalDisplayText{2} becomes
                                                                                      \DIc{\Mabf{a}{l}{d}\Sin\Mabf{a}{m}{d}\Period}{3}
 \Par That \pLinkLocalDisplayText{1} holds, follows from \pLinkLocalDisplayText{2} and \pLinkLocalDisplayText{3}. \QED
 \Item{MDN}{More Discussion} Let \Pr{a,l,d}\ be an ordered basis for the meridian M and suppose that \MABF{a}{d}\ is balanced in the libra \MAB{a}{d}\Period\ From \pLinkDisplayText{EA}{II}{4} and \pLinkItemText{EA}{III} follows that \Mabf{a}{l}{d}\ is balanced in the libra \Mab{a}{d}\ and, for each t\In\Mabf{a}{l}{d}\Comma\
                                                               \DIc{\Russs{a}{l}{d}{t}\ is a libra isomorphism of \MABF{a}{d}\ onto \Mabf{a}{l}{d}\Period}{1}
 \PaR It was shown in \pLinkDisplayText{LL}{EAT}{6} that there can exist no libra isomorphism of \MAB{d}{d}\ onto \MAB{a}{d}\ if a\NEq d. Nonetheless it does make sense to consider the possibility of an isomorphism of \MAB{d}{d}\ onto \MABF{a}{d}\Period\ If we return to the example of \pLinkItemText{EA}{I}, the funtion \Function\OpGr{w}sendsxinFto\Two$^{\RMBx{7}{x}}$inF\end\ may also be expressed as
                                                            \DIc{\Function\OpGr{w}sendsxin\Mab{\iNfty}{\iNfty}to\Two$^{\RMBx{7}{x}}$in\Mabf{\zEro}{\oNe}{\iNfty}\end\Period}{2}
 \Par We define
                               \DIc{\Function\OpGr{W}sends\Vol{\Infty,\Infty;\Zero,c}in\MAB{\iNfty}{\iNfty}to\Vol{\Infty,\Zero;\One,\Two$^{\RMBx{7}{c}}$}in\MABF{\iNfty}{\zEro}\end}{3}
 \Par and observe that, for all \Set{a,b,c}\Sin\Mab{\iNfty}{\iNfty}
                                                  \D{8}{0}{\OpGr{W}\Of{\LLL{\Vol{\Infty,\Infty;\Zero,a},\Vol{\Infty,\Infty;\Zero,b},\Vol{\Infty,\Infty;\Zero,c}}}\Equals%
                                                      \OpGr{W}\Of{\Vol{\Infty,\Infty;\Zero,a}\Circ\Vol{\Infty,\Infty;\Zero,b}\Circ\Vol{\Infty,\Infty;\Zero,c}}\Equals}
                 \D{8}{0}{\OpGr{W}\Of{\Vol{\Infty,\Infty;\Zero,a-b+c}}\Equals\Vol{\Infty,\Zero;\One,\Two$^{\RMBx{7}{a-b+c}}$}\Eq{\hLinkDisplayText{ED}{T}{4}}\Vol{\Infty,\Zero;\One,\Two$^{\RMBx{7}{a}}$}%
                                                         \Circ\Vol{\Infty,\Zero;\One,\Two$^{\RMBx{7}{b}}$}\Circ\Vol{\Infty,\Zero;\One,\Two$^{\RMBx{7}{c}}$}\Equals}
                               \D{8}{8}{\LLL{\Vol{\Infty,\Zero;\One,\Two$^{\RMBx{7}{a}}$},\Vol{\Infty,\Zero;\One,\Two$^{\RMBx{7}{b}}$},\Vol{\Infty,\Zero;\One,\Two$^{\RMBx{7}{c}}$}}\Equals%
                                             \LLL{\OpGr{W}\Of{\Vol{\Infty,\Infty;\Zero,a},\OpGr{W}\Of{\Vol{\Infty,\Infty;\Zero,b}},\OpGr{W}\Of{\Vol{\Infty,\Infty;\Zero,c}}}}}
 \Par whence follows that \OpGr{W}\ is a libra isomorphism of \MAB{\iNfty}{\iNfty}\ onto \MABF{\iNfty}{\zEro}\Period\ Furthermore
                                                        \DIc{\OpGr{w}\Equals\Russs{\zEro}{\oNe}{\iNfty}{\oNe}\Circ\OpGr{W}\Circ\Inv{\Rus{\iNfty}{\zEro}}\Period}{4}
 \PaR We note that
        \DIc{\One\Equals\OpGr{w}\Of{\Zero}\quad\Implies\quad\OpGr{W}\Of{\Vol{\Infty,\Infty;\Zero,\Zero}}\Equals\Vol{\Zero,\Infty;\One,\One}\Equals\Vol{\Zero,\Infty;\OpGr{w}\Of{\Zero},\OpGr{w}\Of{\Zero}}\Comma}{5}
 \Par that
                                         \DIc{\OpGr{w}\Circ\OpGr{w}\Of{\Zero}\Equals\OpGr{w}\Of{\One}\Equals\Two\Equals\Vol{\Infty,\Infty;\Zero,\One}\Of{\One}\quad\Implies\quad%
      \OpGr{W}\Circ\OpGr{W}\Of{\Vol{\Infty,\Infty;\Zero,\Zero}}\Equals\Vol{\Infty,\Zero;\Zero,\OpGr{w}\Of{\Zero}}\Circ\OpGr{W}\Of{\Vol{\Infty,\Infty;\Zero,\Zero}}\Circ\Vol{\Infty,\Zero;\Zero,\OpGr{w}\Of{\Zero}}}{6}
 \Par and that
                                      \DIc{\setRelation{\OpGr{w}}\Of{\Mabf{\zEro}{\Op{\gR{w}\oF{\zEro}}}{\iNfty}}\Equals\setRelation{\OpGr{w}}\Of{\Mabf{\zEro}{\oNe}{\iNfty}}\Equals%
                                                     \Mabf{\oNe}{\tWo}{\iNfty}\Equals\Mabf{\Op{\gR{w}}\oF{\zEro}}{\Op{\gR{w}\cIrc\gR{w}}\oF{\zEro}}{\iNfty}\Period}{7}
 \PaR Results \pLinkLocalDisplayText{4}, \pLinkLocalDisplayText{5}, \pLinkLocalDisplayText{6} and \pLinkLocalDisplayText{7} suggest the following definitions.
 \Item{DI}{Definitions} Let M be a meridian and let \Pr{a,l,d}\ be an ordered basis for M\Period\ Suppose that \MABF{a}{d}\QUad is a balanced subset of the libra \MAB{a}{d}\ and that \OpGr{w}\ is a libra isomorphism from \Mab{d}{d}\ onto \Mabf{a}{l}{d}\Period\ If there exists a libra isomorphism \Function\OpGr{W}sendsin\MAB{d}{d}toin\MABF{a}{d}\end\QUad such that
                                                                     \DIc{\OpGr{w}\Equals\Russs{a}{l}{d}{\Op{\gR{w}}\oF{a}}\Circ\OpGr{W}\Circ\Inv{\Rus{d}{a}}\Comma}{1}
 \Par we shall say that \OpGr{w}\ is a \Def{meridian} \Def{libra} \Def{isomorphism}. If, in addition
                                                                \DI{8}{0}{\OpGr{W}\Of{\Vol{d,d;a,a}}\Equals\Vol{a,d;\OpGr{w}\Of{a},\OpGr{w}\Of{a}}\Comma}{2}
                                    \DI{8}{0}{\OpGr{W}\Circ\OpGr{W}\Of{\Vol{d,d;a,a}}\Equals\Vol{d,a;a,\OpGr{w}\Of{a}}\Circ\OpGr{W}\Of{\Vol{d,d;a,a}}\Circ\Vol{d,a;a,\OpGr{w}\Of{a}}}{3}
                                               \DII{8}{8}{\setRelation{\OpGr{w}}\Of{\Mabf{a}{l}{d}}\Equals\Mabf{\Op{\gR{w}}\oF{a}}{\Op{\gR{w}\cIrc\gR{w}}\oF{a}}{d}}{and}{4}
 \Par we shall say that \OpGr{w}\ is a \Two-\Def{exponential} \Def{isomorphism} \Def{on} \Mab{d}{d}\ \Def{originating} \Def{at} a, and we say that \OpGr{W}\ is the \Def{progenitor} \Def{of} \OpGr{w}\Period\ We note that \pLinkLocalDisplayText{1} implies
                                                                    \DIc{\OpGr{W}\Equals\Inv{\Russs{a}{l}{d}{\Op{\gR{w}}\oF{a}}}\Circ\OpGr{w}\Circ\Rus{d}{a}\Period}{5}
 \Par Furthermore \OpGr{w}\Of{a}\ is in \Mabf{a}{l}{d}\ and so Theorem \pLinkItemText{EA}{III} implies that
                                                          \DI{8}{0}{\Mabf{a}{\Op{\gR{w}\oF{a}}}{d}\Equals\Mabf{a}{l}{d}\Period}{6}
 \PaR For a \Two-exponential \OpGr{w}\ originating at a with progenitor \OpGr{W}\ as above, we shall say that d is its \Def{singular} \Def{point} and that \OpGr{w}\Of{a}\ is its \Def{issue}. We note that the issue of \OpGr{w}\ is one of the two fixed points of \OpGr{W}\Of{\Vol{d,d;a,a}}\Period
 \PaR If there exists a \Two-exponential isomorphism on the meridian M as above, we shall say that M is an \Def{exponential} \Def{meridian}.
 \Item{T1}{Theorem} Let M be an exponential meridian and let \Mii{M}\ be a meridian family of involutions on M. Then
                                           \DI{8}{8}{\ForAll{\Pr{s,t,v}\ an ordered basis for M}\ThereIs{\OpGr{ps}\ a \Two-exponential originating at s onto \Mabf{s}{t}{v}}\Period}{1}
 \Proof By hypothesis there exists an ordered basis \Pr{a,l,d}\ and a \Two-exponential isomorphism \OpGr{w}\ with progenitor \OpGr{W}\ and with base a from \Mab{d}{d}\ onto \Mabf{a}{l}{d}\ such that \pLinkDisplayText{EA}{DI}{2}, \pLinkDisplayText{EA}{DI}{3} and \pLinkDisplayText{EA}{DI}{4} hold. Thus
 \Par\halign{\kn{75}#&\kn{1}#\kn{1}&#&#&#&#&#\cr
                  &&\Mab{d}{d}&&\Mabf{a}{l}{d}&&\cr
                  x\Equals\OpGr{u}\Of{a}&\Vbox{\hbox{\kn{12}\Rus{d}{a}}\kn{2}\hbox{\ArrowLt{50}}\kn{0}}&\hbox{\Vbox{\hbox{\OpGr{u}}\kn{0}}}\ArrowDn{50}&\Vbox{\hbox{\kn{2}\OpGr{W}}\kn{5}\hbox{\kn{-15}\ArrowRt{50}}}&%
                  \kn{2}\hbox{\Vbox{\hbox{\OpGr{h}}\kn{0}}}\ArrowDn{50}&\kn{-15}\Vbox{\hbox{\kn{8}\Russs{a}{l}{d}{\OpGr{w}\oF{a}}}\kn{2}\hbox{\ArrowRt{70}}\kn{0}}&\hbox{\OpGr{h}\Of{\OpGr{w}\Of{a}}}\Equals\OpGr{w}\Of{x}\cr
                  &&\Mab{d}{d}&&\Mabf{a}{l}{d}&&\cr}
 \Par We define
                                                                                   \DIc{\OpGr{th}\Equals\Hol{a}{\Op{\gR{w}}\oF{a}}{d}{s}{t}{v}\Comma\quad
                                \Function\OpGr{PS}sends\OpGr{f}in\MAB{v}{v}to\OpGr{th}\Circ\OpGr{W}\Of{\Inv{\OpGr{th}}\Circ\OpGr{f}\Circ\OpGr{th}}\Circ\Inv{\OpGr{th}}in\MABF{s}{v}\end\Andd
                                                                           \OpGr{ps}\Equiv\Russs{s}{t}{v}{l}\Circ\OpGr{PS}\Circ\Inv{\Rus{v}{s}}\Period}{2}
 \Par Thus
 \halign{\kn{0}#&\kn{1}#\kn{1}&#&#&#&#&#&#&#&#&#\cr
 &&\Mab{v}{v}&\Vbox{\hbox{\kn{14}\OpGr{th}}\kn{4}\hbox{\ArrowLt{30}}}&\Mab{d}{d}&&\Mabf{a}{l}{d}&\Vbox{\hbox{\kn{12}\Inv{\OpGr{th}}}\kn{4}\hbox{\ArrowLt{30}}}&\Mabf{s}{t}{v}&&\cr
 y\Equals\OpGr{a}\Of{\OpGr{th}\Of{a}}\Equals\OpGr{a}\Of{s}&\Vbox{\hbox{\kn{12}\Rus{v}{s}}\kn{2}\hbox{\ArrowLt{50}}\kn{0}}&\hbox{\Vbox{\hbox{\OpGr{a}}\kn{0}}}\ArrowDn{50}&%
 &\hbox{\Vbox{\hbox{\OpGr{u}}\kn{0}}}\ArrowDn{50}&\Vbox{\hbox{\kn{-4}\OpGr{W}}\kn{5}\hbox{\kn{-15}\ArrowRt{30}}}&\kn{2}\hbox{\Vbox{\hbox{\OpGr{h}}\kn{0}}}\ArrowDn{50}%
 &&\kn{-16}\hbox{\Vbox{\hbox{\Vbox{\hbox{\OpGr{PS}\Of{\OpGr{a}}}\kn{1}\hbox{\kn{8}\vrulehw{7}{.5}\kn{1.5}\vrulehw{7}{.5}}\kn{1}\hbox{\kn{6}\OpGr{b}}\kn{-12}}}\kn{0}}}\kn{-3}\ArrowDn{50}&\kn{-20}\Vbox{\hbox{\kn{8}\Russs{s}{t}{v}{t}}\kn{2}\hbox{\ArrowRt{53}}\kn{0}}&\hbox{\OpGr{b}\Of{t}}\Equals\OpGr{b}\Of{\OpGr{th}\Of{\OpGr{w}\Of{a}}}\Equals\OpGr{ps}\Of{y}\cr
 &&\Mab{v}{v}&\Vbox{\hbox{\kn{7}\Inv{\OpGr{th}}}\kn{4}\hbox{\ArrowRt{30}}}&\Mab{d}{d}&&\Mabf{a}{l}{d}&\Vbox{\hbox{\kn{11}\OpGr{th}}\kn{4}\hbox{\ArrowRt{30}}}&\Mabf{s}{t}{v}&&\cr}
 \vskip3pt
 \Par Direct computations show that
       \DIc{\MAB{v}{v}\Equals\SetSuch{\OpGr{th}\Circ\OpGr{f}\Circ\Inv{\OpGr{th}}}{\OpGr{f}\In\MAB{d}{d}}\Andd\MAB{s}{v}\Equals\SetSuch{\OpGr{th}\Circ\OpGr{f}\Circ\Inv{\OpGr{th}}}{\OpGr{f}\In\MAB{a}{d}}\Period}{3}
 \PaR Direct computation, along with \pLinkDisplayText{EA}{DI}{2}, \pLinkDisplayText{EA}{DI}{3}, \pLinkDisplayText{EA}{DI}{4}, \pLinkLocalDisplayText{2} and \pLinkLocalDisplayText{3} show that
                                                                          \DI{8}{0}{\OpGr{PS}\Of{\Vol{v,v;s,s}}\Equals\Vol{s,v;t,t}\Comma}{4}
                                 \DI{8}{0}{\OpGr{PS}\Circ\OpGr{PS}\Of{\Vol{v,v;s,s}}\Equals\Vol{v,s;s,\OpGr{ps}\Of{s}}\Circ\OpGr{PS}\Of{\Vol{v,v;s,s}}\Circ\Vol{v,s;s,\OpGr{ps}\Of{s}}}{5}
                                             \DII{8}{-8}{\setRelation{\OpGr{ps}}\Of{\Mabf{s}{t}{v}}\Equals\Mabf{\Op{\gR{ps}}\oF{s}}{\Op{\gR{ps}\cIrc\gR{ps}}\oF{s}}{v}}{and}{6}
 \Par That \pLinkLocalDisplayText{1} holds now follows from \pLinkLocalDisplayText{4}, \pLinkLocalDisplayText{5} and \pLinkLocalDisplayText{6}. Furthermore \OpGr{PS}\ is the progenitor of \OpGr{ps}\Period\ \QED
  \Item{T2}{Theorem} Let M be an exponential meridian. Then the underlying field has characteristic \Zero.
 \Proof Let \Function\OpGr{w}sendsin\Mab{\Infty}{\Infty}toin\Mabf{\zEro}{\oNe}{\iNfty}\end\ be a \Two-exponential function originating at \Zero\ with issue \One\Period\  Assume that the underlying field were of characteristic p different from 0 (and 2 of course). Then p were a prime number different from 2. We introduce field operations for the basis \Set{\Zero,\One,\Infty}\ of M. Choose \vruleh{14}x\In\MAB{\zEro}{\iNfty}\ such that
                                                                  \Dc{$\hbox{\OpGr{w}\Of{x}\Equals}\overbrace{\hbox{\One\Plus\One\Plus \dots\Plus\One}}^{\RMBx{7}{p\thinspace-\oNe\ times}}$\Period}
  \Par Then
   \Dc{\OpGr{w}\Of{x\Plus x}\Equals\OpGr{w}\Of{x}\Cdot\OpGr{w}\Of{x}\Equals$\overbrace{\hbox{\One\Plus\One\Plus\dots\Plus\One}}^{\RMBx{7}{(p\cDot p-\tWo\cDot p-\oNe) times}}$\Equals\One\Equals\OpGr{w}\Of{o}\quad\Implies\quad x\Plus x\Equals o}
 \Par which would imply that the field had characteristic 2, which is absurd. \QED
 \Item{L1}{Lemma} Suppose that x is an element of a meridian M and that \Set{\OpGr{a},\OpGr{b}}\Sin\MAB{x}{x}\Period\ Then
                                                             \DIc{\ThereIs{\OpGr{c}\In\MAB{x}{x}}\quad\OpGr{a}\Equals\OpGr{c}\Circ\OpGr{b}\Circ\OpGr{c}\Period}{1}
 \Proof Since \OpGr{a}\ is an involution with a fixed point x, there is another fixed point a of \OpGr{a}\Period\ Similarly, there is a point b distinct from x which \OpGr{b}\QUad fixes. Letting \OpGr{d}\Equiv\Vol{x,x;a,b}\Circ\OpGr{b}\Circ\Vol{x,x;a,b}\Comma\ we see that
                                                               \Dc{\OpGr{d}\Of{x}\Equals x\Equals\OpGr{a}\Of{x}\Andd\OpGr{d}\Of{a}\Equals a\Equals\OpGr{a}\Of{a}}
 \Par Since these two involutions fix two common distinct points, it follows from \pLinkItemText{E}{L3} that they must be identical. Letting \OpGr{c}\Equiv\Vol{x,x;a,b}\Comma\ we see that \OpGr{c}\ is in \MAB{x}{x}\ and that \pLinkLocalDisplayText{1} holds. \QED
 \Item{T4}{Theorem} Let \Pr{a,t,b}\ be an ordered basis for the exponential meridian M and let
                                                                                        \DIc{u\Equiv\Vol{a,a;b,b}\Of{t}\Period}{000}
 \Par Then
                                                                    \DI{8}{0}{\Mabf{a}{t}{b}\Cap\Mabf{a}{u}{b}\Equals\Void}{0}
                                                     \DII{8}{8}{\kn{-12}M\Equals\Mabf{a}{t}{b}\Cup\Mabf{a}{u}{b}\Cup\Set{a,b}\Period}{and}{1}
 \Proof Assume that u were in \Mabf{a}{t}{b}. Then
                                                                           \DIc{\ThereIs{\OpGr{s}\In\MABF{a}{b}}\quad\OpGr{s}\Of{t}\Equals u.}{3}
 \Par We would have
                      \Dc{\OpGr{s}\Circ\Vol{a,b;t,t}\Circ\OpGr{s}\Of{u}\Eq{\hLinkLocalDisplayText{000}}\OpGr{s}\Circ\Vol{a,b;t,t}\Of{t}\Equals\OpGr{s}\Of{t}\Eq{\hLinkLocalDisplayText{000}}u\Period}
 \Par Since \OpGr{s}\Circ\Vol{a,b;t,t}\Circ\OpGr{s}\And\Vol{a,b;t,t}\ would agree at u, a and b, it would follow that
                                                                        \DIc{\OpGr{s}\Circ\Vol{a,b;t,t}\Circ\OpGr{s}\Equals\Vol{a,b;t,t}\Period}{4}
 \Par By Theorem \pLinkItemText{EA}{T1} there would exist a \Two-exponential \OpGr{w}\ from \Mab{b}{b}\ with base a and issue t onto \Mabf{a}{t}{b}\Period\ Let \OpGr{W}\ \thinspace be the progenitor for \OpGr{w}\ and let \OpGr{th}\In\MAB{b}{b}\ be such that \OpGr{W}\Of{\OpGr{th}}\Equals\OpGr{s}\Period\ We would have
           \D{8}{0}{\OpGr{W}\Of{\OpGr{th}}\Circ\OpGr{W}\Of{\Vol{b,b;a,a}}\Circ\OpGr{W}\Of{\OpGr{th}}\Eq{\hLinkDisplayText{EA}{DI}{2}}\OpGr{W}\Of{\OpGr{th}}\Circ\Vol{b,b;t,t}\Circ\OpGr{W}\Of{\OpGr{th}}\Equals}
                                  \D{8}{8}{\OpGr{s}\Circ\Vol{a,b;t,t}\Circ\OpGr{s}\Eq{\hLinkLocalDisplayText{4}}\Vol{a,b;t,t}\Eq{\hLinkDisplayText{EA}{DI}{2}}\OpGr{W}\Of{\Vol{b,b;a,a}}}
 \Par which would imply that
                                                                       \DIc{\OpGr{th}\Circ\Vol{b,b;a,a}\Circ\OpGr{th}\Equals\Vol{b,b;a,a}\Period}{5}
 \Par The only way equation \pLinkLocalDisplayText{5} could hold in the context of \MAB{b}{b}\ would be for \OpGr{th}\ to be \Vol{b,b;a,a}\Comma\ which, in view of \pLinkDisplayText{EA}{DI}{2} and the fact that t is the issue of \OpGr{w}\Comma\ would imply that \OpGr{s}\Equals\Vol{a,b;t,t}\Period\ Consequently
                                                 \Dc{\Vol{b,b;a,a }\Of{t}\Equals u\Equals\OpGr{s}\Of{t}\Equals\Vol{a,b;t,t}\Of{t}\Equals t\quad\Implies\quad t\In\Set{b,a}}
 \Par which would be absurd: it follows that
                                                                                            \DIc{u\NIn\Mabf{a}{t}{b}\Period}{6}
 \PaR Now assume that there existed an element x\In\Mabf{a}{t}{b}\Cap\Mabf{a}{u}{b}\Period\ Then
                                                  \Dc{\ThereIs{\Set{\OpGr{f},\OpGr{h}}\Sin\MABF{a}{b}}\quad\OpGr{f}\Of{t}\Equals x\Equals\OpGr{h}\Of{u}\Period}
 \Par The meridian M being an exponential meridian, \MABF{a}{b}\ is a libra and so we have \OpGr{h}\Circ\OpGr{f}\Circ\Vol{a,b;t,t}\In\MABF{a}{b}\Period\ Consequently
                                 \Dc{\Mabf{a}{t}{b}\Ni\OpGr{h}\Circ\OpGr{f}\Circ\Vol{a,b;t,t}\Of{t}\Equals\OpGr{h}\Circ\OpGr{f}\Of{t}\Equals\OpGr{h}\Of{x}\Equals u}
 \Par which violates \pLinkLocalDisplayText{6}. This establishes \pLinkLocalDisplayText{0}.
 \PaR Let p be any element of \Mab{a}{b}. Let q be the fixed point of \Vol{a,b;p,p}\ which is distinct from p. Evidently \Vol{a,b;t,t}\Circ\Vol{a,a;b,b}\Circ\Vol{a,b;t,t}\ is in \Mii{M}\ and leaves both a and b fixed: hence
                                                                       \Dc{\Vol{a,b;t,t}\Circ\Vol{a,a;b,b}\Circ\Vol{a,b;t,t}\Equals\Vol{a,a;b,b}\Period}
  \Par It follows that
          \DIc{\Vol{a,b;t,t}\Circ\Vol{a,a;b,b}\Equals\Vol{a,a;b,b}\Circ\Vol{a,b;t,t}\quad\Implies\quad\Vol{a,b;t,t}\Circ\Vol{a,a;b,b}\Of{t}\Equals\Vol{a,a;b,b}\Circ\Vol{a,b;t,t}\Of{t}\Equals\Vol{a,a;b,b}\Of{t}}{7}
 \Par Since \MABF{a}{b}\ is isomorphic as a libra to \MAB{b}{b}\Comma\ it follows from Lemma \pLinkItemText{EA}{L1} that there exists \OpGr{u}\In\MABF{a}{b}\ such that
                                                                                \Dc{\Vol{a,b;p,p}\Equals\OpGr{u}\Circ\Vol{a,b;t,t}\Circ\OpGr{u}}
 which implies both
                 \DIc{\Vol{a,b;p,p}\Of{\OpGr{u}\Of{t}}\Equals\OpGr{u}\Circ\Vol{a,b;t,t}\Circ\OpGr{u}\Of{\OpGr{u}\Of{t}}\Equals\OpGr{u}\Of{t}}{8}
 and
       \D{8}{0}{\Vol{a,b;p,p}\Of{\OpGr{u}\Of{\Vol{b,b;a,a}\Of{t}}}\Equals\OpGr{u}\Circ\Vol{a,b;t,t}\Circ\OpGr{u}\Of{\OpGr{u}\Of{\Vol{b,b;a,a}\Of{t}}}\Equals}
                                                                           \DI{0}{0}{}{9}
                \D{0}{8}{\OpGr{u}\Circ\Vol{a,b;t,t}\Circ\Vol{b,b;a,a}\Of{t}\Eq{\hLinkLocalDisplayText{7}}\OpGr{u}\Of{\Vol{a,a;b,b}\Of{t}}\Period}
 \Par By definition of p and q it now follows from \pLinkLocalDisplayText{8} and \pLinkLocalDisplayText{9} that
                                                                     \Dc{\Set{\OpGr{u}\Of{t}\Comma\OpGr{u}\Of{\Vol{a,a;b,b}\Of{t}}}\Equals\Set{p,q}\Period}
 \Par If \OpGr{u}\Of{t}\Equals p, then p\In\Mabf{a}{t}{b}\Period\ If \OpGr{u}\Of{\Vol{a,a;b,b}\Of{t}}\Equals p, then p\In\Mabf{a}{u}{b}\Period\ This establishes \pLinkLocalDisplayText{1}. \QED
 \Item{N}{Definitions and Notation} Let \Op{b}\ be any ordered basis for M: that is, an ordered triple of pairwise distinct elements. We shall adopt the following notation:
                                                                               \DIc{\Pr{\BZero{b},\BOne{b},\BInfty{b}}\ \Equiv\ \Op{b}\Period}{2}
 \Par The field operations relative to \Op{b}\ will be denoted by \BPlus{b}, \BMinus{b}, \BTimes{b}, and \BMinus{}\kn{1}\Vbox{\hbox{\Op{\RMBx{5}{b}}}\kn{1.5}}: thus
                                    \DIc{x\BPlus{b}y\Equals\Quin{x}{y}{\BInfty{b}}{\BInfty{b}}{\BZero{b}}\Comma\quad x\BMinus{b}y\Equals\Quin{x}{\BZero{b}}{\BInfty{b}}{\BInfty{b}}{y}\quad%
                                       x\BTimes{b}y\Equals\Quin{x}{y}{\BInfty{b}}{\BZero{b}}{\BOne{b}}\Andd \BOver{x}{y}{b}\Equals\Quin{x}{\BOne{b}}{\BInfty{b}}{\BZero{b}}{y}\Period}{3}
 \Par We shall denote by
                                                                                             \DIc{F$_{\hbox{\Op{\RMBx{7}{b}}}}$}{4}
 \Par the corresponding field \Mab{\iNfty}{\iNfty}\Period
 \PaR When it is clear from the context which ordered basis \Op{b}\ is under consideration, we shall occasionally suppress the symbol \Quotes{\Op{\RMBx{7}{b}}} from the field operations, shall occasionally suppress the symbol \Quotes{\BTimes{b}} completely, and shall occasionally introduce the exponential notation
\DIc{x$^n$\ \Equals\ $\overbrace{\hbox{x\Cdot x\Cdot\dots\Cdot x}}^{\RMBx{7}{n\ times}}$\Period}{5}
  \Item{F2}{Lemma} Let \Mabf{\Zero}{\One}{\Infty}\ be an ordered basis for a meridian M\Period\ Let \Two \kn{2}\Equiv\One \Plus \One\Period\ Then
                                                                \DIc{\Mabf{\One}{\Two}{\Infty}\Equals\SetSuch{\One \Plus t$^2$}{t\In\Mab{\iNfty}{\iNfty}}\Period}{1}
  \Proof Let x be an element of \Mabf{\One}{\Two}{\Infty}\Period\ Then there exists \OpGr{th}\In\MABF{\iNfty }{\EUBx{7}{1}}\ such that \OpGr{th}\Of{\Two }\Equals x. Since \OpGr{th}\Of{\One }\Equals\Infty \Comma\ there exists b\In\Mab{\iNfty}{\iNfty} such that
                                                                             \Dc{\ForAll{s\In\Mab{\iNfty }{\iNfty }}\quad\OpGr{th}\Of{s}\Equals\Over{s\Plus b}{s\Minus\One }\Period}
  \Par Since \OpGr{th}\ has a fixed point p, it follows from the quadratic formula that b\Plus \One\  has a square root t. Thus
                    \Dc{x\Equals\OpGr{th}\Of{\Two }\Equals\Over{\Two \Plus b}{\Two\Minus\One }\Equals\One \Plus b\Plus \One \Equals\One \Plus t$^2$\Period}
  \PaR Now suppose that x\Equals\One \Plus t$^2$ for some t. Setting b\Equiv\One\Minus t$^2$, direct calculation shows that the function
                                            \Dc{\Function sendssinFto\Over{s\Minus b}{s\Minus\One }inM\end\Quad is in \MABF{\iNfty }{\EUBx{7}{1}}\ and sends \Two \ to x\Comma}
  \Par whence follows that x is in \Mabf{\One}{\Two}{\Infty}\Period\ \QED
  \Item{T6}{Theorem} Let \Pr{a,t,b}\ be an ordered basis for an exponential meridian M and let \OpGr{f}\ be an element of \MABF{a}{b}\Period\ Then one fixed point of \OpGr{f}\ is in \Mabf{a}{t}{b}\ and one is not.
 \Proof Since \OpGr{f}\Circ\Vol{a,a;b,b}\Circ\OpGr{f}\ and \Vol{a,a;b,b}\ both fix the distinct points a and b, it follows that they are the same. Let p and q be the fixed points of \OpGr{f}\Period\ Then
                                            \Dc{\OpGr{f}\Of{\Vol{a,a;b,b}\Of{p}}\Equals\OpGr{f}\Of{\OpGr{f}\Circ\Vol{a,a;b,b}\Circ\OpGr{f}\Of{p}}\Equals\Vol{a,a;b,b}\Of{p}\Period}
 \Par Since p is neither a nor b, it follows that \Vol{a,a;b,b}\Of{p}\NEq p and so
                                                                                       \DIc{\Vol{a,a;b,b}\Of{p}\Equals q\Period}{1}
 \Par From \pLinkDisplayText{EA}{T4}{3} follows that p is in either \Mabf{a}{t}{b}\ or in \Mabf{a}{\Vol{a,a;b,b}\Of{t}}{b}\Period\ By Theorem \pLinkItemText{EA}{III} \Mabf{a}{p}{b}\ is either \Mabf{a}{t}{b}\ or \Mabf{a}{\Vol{a,a;b,b}\Of{t}}{b}\Period\ By \pLinkLocalDisplayText{1} and \pLinkDisplayText{EA}{T4}{1}, whichever one \Mabf{a}{p}{b}\ is, \Mabf{a}{q}{b}\ is the other. \QED
  \Item{FF2}{Lemma} Let \Mabf{\Zero}{\One}{\Infty}\ be an ordered basis for a meridian M\Period\ Then
                                        \DIc{\Mabf{\zEro}{\One}{\Infty}\Equals\SetSuch{f$^2$}{t\In\Mab{\iNfty}{\iNfty}}\Equals\SetSuch{f$^2$}{t\In\Mabf{\zEro}{\oNe}{\iNfty}}\Period}{1}
  \Proof Let t be an element of \Mab{\iNfty}{\iNfty}\Period\ Then \Vol{\Zero,t$^2$,\One,\Zero}\ is in \MAB{\Zero}{\Infty}\Comma\ leaves t fixed, and
                                                      \DIc{\Vol{\Zero,t$^2$,\One,\Zero}\Of{\One}\Equals t$^2$\quad\Implies\quad t$^2$\In\Mabf{\Zero}{\One}{\Infty}\Period}{2}
  \PaR Now let x be an element of \Mabf{\Zero}{\One}{\Infty}\Period\ There exists \OpGr{f}\In\MABF{\zEro}{\iNfty}\ such that \OpGr{f}\Of{\One}\Equals x\Period\ It follows from \pLinkItemText{EA}{T6} that one of the fixed points of \OpGr{f}\ is in \Mabf{\Zero}{\One}{\Infty}\Colon\ we shall denote that point by t\Period\ We have
                                   \DIc{\OpGr{f}\Equals\Vol{\Zero,t$^2$,\One,\Zero}\quad\Implies\quad x\Equals\OpGr{f}\Of{\One}\Equals\Over{t$^2$}{\One}\quad\Implies\quad x\Equals t$^2$}{3}
 \PaR That \pLinkLocalDisplayText{1} holds is now a consequence of \pLinkLocalDisplayText{2} and \pLinkLocalDisplayText{3}.
 \Item{PlusTimes}{Lemma} Let M be an exponential meridian and \Mii{M}\ the meridian family of involutions on M. Let \Pr{\Zero,\One,\Infty}\ be an ordered basis for M and let \Quotes{+} and \Quotes{\Cdot} be the associated  binary operations on the field \Mab{\iNfty}{\iNfty}\Period\ Then,
                          \DIc{\ForAll{\Set{x,y}\Sin\Mabf{\zEro}{\oNe}{\iNfty}}\quad \Set{x\Cdot y\Comma\QUad\Over{y}{x}\Comma\Quad x\Plus y}\Sin\Mabf{\zEro}{\oNe}{\iNfty}\Period}{1}
 \Proof From \pLinkItemText{EA}{FF2} follows that there exists \Set{s,t}\Sin\Mabf{\zEro}{\oNe}{\iNfty}\ such that
                                                                                       \DIc{x\Equals s$^2$\Andd y\Equals t$^2$\Period}{2}
 \Par In particular
                                                               \DIc{xy\Equals(st)$^2$\quad\Imp{\hLinkItemText{EA}{FF2}}\quad xy\In\Mabf{\zEro}{\oNe}{\iNfty}}{3}
 \Par and
                                              \DIc{\Over{y}{x}\Equals(\Over{t}{s})$^2$\quad\Imp{\hLinkItemText{EA}{FF2}}\quad \Over{y}{x}\In\Mabf{\zEro}{\oNe}{\iNfty}\Period}{4}
 \PaR By \pLinkItemText{EA}{T1} there exists a \Two-exponential \OpGr{w}\ of \Mab{\iNfty}{\iNfty}\ onto \Mabf{\zEro}{\oNe}{\iNfty}\ with base \Zero\ and issue \One\Period\ We have by \pLinkDisplayText{EA}{DI}{4}
                   \Dc{\Mabf{\oNe}{\tWo}{\iNfty}\Equals\setRelation{\OpGr{w}}{(\Mabf{\zEro}{\oNe}{\iNfty})}\Sin\Range{\OpGr{w}}\Equals\Mabf{\zEro}{\oNe}{\iNfty}\quad\Imp{\hLinkDisplayText{EA}{F2}{1}}\quad%
                                        \One+\Over{y}{x}\In\Mabf{\zEro}{\oNe}{\iNfty}\quad\Imp{\hLinkLocalDisplayText{3}}\quad x(\One+\Over{y}{x})\In\Mabf{\zEro}{\oNe}{\iNfty}\Period}
 \Par This, with \pLinkLocalDisplayText{3} and \pLinkLocalDisplayText{4}, implies \pLinkLocalDisplayText{1}. \QED
  \Item{Defff}{Notation} Let \Op{b}\Equals\Pr{\Zero,\One,\Infty}\ be a basis for a meridian M. For \Set{a,b,c,d}\Sin F$_{\hbox{\Op{\RMBx{7}{b}}}}$\Equals\Mab{\iNfty}{\iNfty}\Comma\ we adopt the notation
                                                                                       \DIc{\BDet{a}{b}{c}{d}{b}\ \Equiv\ ad$-$bc\Period}{1}
  \PaR Furthermore we denote
                                                                  \DIc{\MNF\ \Equiv\ \SetSuch{\OpGr{f}\In\Mii{M}}{\OpGr{f}\ has no fixed point}\Period}{2}
  \Item{T7}{Theorem} Let \Pr{\Zero,\One,\Infty}\ be a basis for an exponential meridian M. Let \OpGr{f}\ be an element of \Mii{M}\ and suppose that \OpGr{f}\Equals\Vol{a,b,c,-a}\ as in \pLinkItemText{ED}{EI}. Then
                                                          \DIc{\OpGr{f}\In\MNF\Quad\Iff\Quad\BDet{a}{b}{c}{-a}{b}\In\Mabf{\zEro}{\oNe}{\iNfty}\Period}{1}
 \Proof It follows from the quadratic equation that \vruled{5}\Vol{a,b,c,-a}\ has a fixed point if, and only if, a$^2$\Plus b\Cdot c has a square root in \Mab{\iNfty}{\iNfty}\Comma\ which, in view of Lemma \pLinkItemText{EA}{FF2}, is true if, and only if, a$^2$\Plus b\Cdot c is in \Mabf{\zEro}{\oNe}{\iNfty}\Period\ It follows from Theorem \pLinkItemText{EA}{T4} that a$^2$\Plus b\Cdot c is in \Mabf{\zEro}{\oNe}{\iNfty}\ if, and only if, \Minus a$^2$\Minus b\Cdot c is not in \Mabf{\zEro}{\oNe}{\iNfty}\Period\ But
                                                       \Dc{$-\hbox{a}\Cdot\hbox{a}-\hbox{b}\Cdot\hbox{c}$\ \Equals\ \BDet{a}{b}{c}{-a}{b}}
\Par from which \pLinkLocalDisplayText{1} holds. \QED
 \Item{L8}{Lemma} Let M be an exponential meridian and let \OpGr{t}\ be a translation in \LLL{\Mii{M}}\Period\Foot{That is, \OpGr{t}\ has exactly one fixed point in M.} Then \OpGr{t}\ is the composition of nine elements of \MNF\Period\
 \Proof Let \Op{b}\ be the ordered basis \Pr{\Zero,\One,\Infty}\ where \Infty\ is the fixed point of \OpGr{t}\Comma\ \Zero\ is any other point of M and \One\Equiv\Inv{\OpGr{t}}\Of{\Zero}\Period\ Then
                                                                                      \Dc{\OpGr{t}\Equals\Vol{\One,\Minus\One,\Zero,\One}}
 \Par and direct calculation\Foot{One constructs the corresponding \Two\Cross\Two\ matrices and uses matrix multiplication.} yields
                     \DIc{\OpGr{t}\Equals\Vol{\Zero,\Minus\Six,\One,\Zero}\Circ\Vol{\Zero,\Minus\Five,\One,\Zero}\Circ\Vol{\One,\Minus\One,\Two,\Minus\One}\Circ\Vol{\Four,\Minus\Five,\Five,\Minus\Four}%
        \Circ\Vol{\Zero,\Minus\Four,\One,\Zero}\Circ\Vol{\Zero,\Minus\One,\One,\Zero}\Circ\Vol{\One,\Minus\One,\Two,\Minus\One}\Circ\Vol{\Zero,\Minus\Five,\One,\Zero}\Circ\Vol{\Zero,\Minus\Six,\One,\Zero}\Period}{1}
 \Par That each of the factors in \pLinkLocalDisplayText{1} is in \MNF\ follows from \pLinkItemText{EA}{T7}. \QED
  \Item{L9}{Lemma} Let \Pr{\Zero,\One,\Infty}\ be an ordered basis for an exponential meridian. Then
                                           \DIc{\Mabf{\zEro}{\oNe}{\iNfty}\Cap\Mabf{\iNfty}{\zEro}{\oNe}\Cap\Mabf{\oNe}{\iNfty}{\zEro}\Equals\Void}{1}
  \Proof Assume that there existed some p in \Mabf{\zEro}{\oNe}{\iNfty}\Cap\Mabf{\iNfty}{\zEro}{\oNe}\Cap\Mabf{\oNe}{\iNfty}{\zEro}\Period\ Then
        \Dc{\ThereIs{\Pr{\OpGr{a},\OpGr{b},\OpGr{c}}\In\MABF{\zEro}{\iNfty}\Cross\MABF{\zEro}{\oNe}\Cross\MABF{\oNe}{\iNfty}}\quad\OpGr{a}\Of{\One}\Equals\OpGr{b}\Of{\Infty}\Equals\OpGr{c}\Of{\Zero}\Equals p\Period}
 \Par It follows that
                               \Dc{\OpGr{a}\Equals\Vol{\Zero,p,\One,\Zero}\Comma\quad\OpGr{b}\Equals\Vol{p,-p,\One,-p}\Andd\OpGr{c}\Equals\Vol{-\One,p,-\One,\One}\Period}
 \Par Applying the quadratic equation to each of the homographies \OpGr{a}\Comma\ \OpGr{b}\ and \OpGr{c}, respectively, we obtain the formulae
                                                            \Dc{$\pm\sqrt{\RmBx{p}},\quad\RmBx{p}\pm\sqrt{\RmBx{p}\Cdot\RmBx{p}-\RmBx{p}}\Andd\One\pm\sqrt{\One-p}$}
 \Par respectively, for their fixed points. Since these three homographies would have fixed points, there would exist \Set{a,b,c}\In\Mab{\iNfty}{\iNfty}\ such that
                                                             \Dc{a$^2$\Equals p,\quad b\Cdot b\Equals p\Cdot p$-$p\Andd c\Cdot c\Equals\One$-$p\Period}
 \Par It would follow that
                                                                                      \Dc{(a\Cdot c)\Cdot(a\Cdot c)\Equals\Minus b\Cdot b}
 \Par which by \pLinkDisplayText{EA}{T4}{000} is impossible. \QED
  \Item{T9}{Theorem} Let \Pr{\Zero,\One,\Infty}\ be a basis for an exponential meridian M. Then
                                                  \DIc{\Set{\Set{\Zero},\Set{\One},\Set{\Infty},\Mabf{\zEro}{\vOl{\zEro,\zEro;\iNfty,\iNfty}\Of{\oNe}}{\iNfty},%
                      \Mabf{\oNe}{\vOl{\oNe,\oNe;\iNfty,\iNfty}\Of{\zEro}}{\iNfty},\Mabf{\oNe}{\vOl{\oNe,\oNe;\zEro,\zEro}\Of{\iNfty}}{\zEro}}\Quad is a partition of M.}{1}
  \Proof We first show that
                                                    \DIc{M\Equals(\Set{\Set{\Zero},\Set{\One},\Set{\Infty},\Mabf{\zEro}{\vOl{\zEro,\zEro;\iNfty,\iNfty}\Of{\oNe}}{\iNfty},%
                                              \Mabf{\oNe}{\vOl{\oNe,\oNe;\iNfty,\iNfty}\Of{\zEro}}{\iNfty},\Mabf{\oNe}{\vOl{\oNe,\oNe;\zEro,\zEro}\Of{\iNfty}}{\zEro}})}{2}
   \Par Let t be any element of M\Cop\Set{\Zero,\One,\Infty}\Period\ If t\NIn\vruled{12}\Mabf{\zEro}{\vOl{\zEro,\zEro;\iNfty,\iNfty}\Of{\oNe}}{\iNfty}, then \pLinkItemText{EA}{T4} implies that t\In\Mabf{\Zero}{\One}{\Infty}\Period\ \vruleh{13}Similarly, if t is not in \Mabf{\oNe}{\vOl{\oNe,\oNe;\iNfty,\iNfty}\Of{\zEro}}{\iNfty}, then t is in  \vruled{13}\Mabf{\iNfty}{\zEro}{\oNe}, and if t is not in \Mabf{\oNe}{\vOl{\oNe,\oNe;\zEro,\zEro}\Of{\iNfty}}{\zEro}, then t is in \Mabf{\oNe}{\iNfty}{\zEro}\Period\ It follows from \pLinkItemText{EA}{L9} that these three things cannot occur concurrently. We have established
 \pLinkLocalDisplayText{2}.
 \PaR Assume that the constituents of \pLinkLocalDisplayText{1} were not pairwise disjoint. Without loss of generality, we may assume that \Mabf{\zEro}{\vOl{\zEro,\zEro;\iNfty,\iNfty}\Of{\oNe}}{\iNfty}\ and \Mabf{\oNe}{\vOl{\oNe,\oNe;\iNfty,\iNfty}\Of{\zEro}}{\iNfty}\ had a point x in common: which is to say that
                                                                          \DIc{x\In\Mabf{\zEro}{-\oNe}{\iNfty}\Cap\Mabf{\oNe}{\tWo}{\iNfty}\Period}{3}
\Par Let \OpGr{w}\ be a \Two-exponential from \Mab{\iNfty}{\iNfty}\ to \Mabf{\zEro}{\oNe}{\iNfty}\Period\ From \pLinkDisplayText{EA}{DI}{4} follows that there would exist u\In\Mabf{\zEro}{\oNe}{\iNfty}\ such that \OpGr{w}\Of{u}\Equals x\Period\ But \OpGr{w}\Of{u}\ is in \Mabf{\zEro}{\oNe}{\iNfty}, which by \pLinkDisplayText{EA}{T4}{0} would contradict \pLinkLocalDisplayText{3}. \QED
 \Item{CO9}{Corollary} Let \Op{b}\Equals\Pr{\Zero,\One,\Infty}\ be an ordered basis for an exponential meridian M. Then
                                      \DI{8}{0}{\Set{\Mabf{\zEro}{\Over{\oNe}{\tWo}}{\oNe},\Set{\One},\Mabf{\oNe}{\tWo}{\iNfty}}\ is a partition of \Mabf{\zEro}{\oNe}{\iNfty}\Period}{1}
 \Proof From \pLinkItemText{EA}{T9} follows that
                                                                          \DIc{\Set{\Set{\Zero},\Set{\One},\Set{\Infty},\Mabf{\zEro}{-\oNe}{\iNfty},%
                                                                 \Mabf{\oNe}{\tWo}{\iNfty},\Mabf{\zEro}{\Over{\oNe}{\tWo}}{\oNe}}\Quad is a partition of M}{2}
 \Par and from \pLinkItemText{EA}{T4} follows that
                                              \DIc{\Set{\Set{\Zero},\Set{\Infty},\Mabf{\zEro}{\oNe}{\iNfty},\Mabf{\zEro}{-\oNe}{\iNfty}}\Quad is also a partition of M\Period}{3}
 \Par That \pLinkLocalDisplayText{1} holds follows from \pLinkLocalDisplayText{2} and \pLinkLocalDisplayText{3}. \QED
 \Item{L10}{Lemma} Let M be an exponential meridian and let \OpGr{r}\In\LLL{\Mii{M}}\Cop\Mii{M}\ be a pure rotation.\Foot{That is, \OpGr{r}\ has no fixed point.} Then there exists \OpGr{p}\In\MNF\ and a translation \OpGr{t}\ such that
                                                                                      \DIc{\OpGr{r}\Equals\OpGr{p}\Circ\OpGr{t}\Period}{1}
 \Proof Let \Infty\ be any element of M. Let a\Equiv\OpGr{r}\Of{\Infty}. Since \OpGr{r}\ is not in \Mii{M}\Comma\ it follows from Corollary \pLinkItemText{ED}{TTP} that \OpGr{r}\Of{a}\NEq\Infty\Period\ Let \One\Equiv\Inv{\OpGr{r}}\Of\Infty\Period\ It follows from \pLinkItemText{ED}{HI} and \pLinkItemText{ED}{Har} that there exists \Zero\In M such that \Set{\Set{\Zero,\Infty},\Set{a,\One}}\ is a harmonic pair. Introducing the canonical field operations on F\Equiv\Mab{\Infty}{\Infty}\Comma\ we see that a\Equals\Minus\One\ and so, if b\Equiv\OpGr{r}\Of{\Zero}, we have
                                                                                  \DIc{\OpGr{r}\Equals\Vol{\One,b,$-$\One,\One}\Period}{2}
 \Par Thus the equation for \OpGr{r}\ to have a fixed point x is
                                                                                            \Dc{x$^2$\Plus b\Equals\Zero\Period}
 \Par Since \OpGr{r}\Comma\ being a rotation, has no fixed point, we have b\In\Mabf{\zEro}{\oNe}{\iNfty}\Period\ From \pLinkItemText{EA}{F2} follows that \One+b\ is in \Mabf{\oNe}{\tWo}{\iNfty}\Period\ Since M is an exponential meridian, it follows from \pLinkDisplayText{EA}{PlusTimes}{1} that \One+b\ is in \Mabf{\zEro}{\oNe}{\iNfty}, whence follows by Lemma \pLinkItemText{EA}{FF2} that there exists m\In\Mabf{\zEro}{\oNe}{\iNfty}\ such that
                                                                                            \DIc{m$^2$\Equals\One+b\Period}{2.5}
  \Par Let
                                                                                       \DIc{p\ \Equiv\ \Over{b+\Two+\Two m}{b$^2$}\Period}{3}
 \Par We have
                            \D{8}{0}{p$^2$b$^2$$-$\Two pb+\One$-$\Four p\ \Equals\ \LPAREN{14}\hbox{\lower4pt\hbox{\Over{b+\Two+\Two m}{b$^2$}}}\RPAREN{14}\hbox{\raise5pt\hbox{$^2$}}b$^2$$-$%
                        \Two\LPAREN{14}\hbox{\lower4pt\hbox{\Over{b+\Two+\Two m}{b$^2$}}}\RPAREN{14}b+\One$-$\Four\LPAREN{14}\hbox{\lower4pt\hbox{\Over{b+\Two+\Two m}{b$^2$}}}\RPAREN{14}\ \Equals}
             \DI{8}{0}{\Over{b$^2$+\Four+\Four m$^2$+\Four b+\Four m+\Two bm$-$\Two b$^2$$-$\Four b$-$2bm+b$^2$$-$\Four b$-$\Eight$-$\Four m}{b$^2$}\ \Equals\ \Over{$-$\Four+\Four m$^2$$-$\Four b}{b$^2$}}{4}
                                                         \D{8}{8}{\Eq{\hLinkLocalDisplayText{2.5}}\ \Over{$-$\Four+\Four(\One+b)$-$\Four b}{b$^2$}\ \Equals\ \Zero}
\Par Since the left hand side of \pLinkLocalDisplayText{4} is the discriminant of the quadratic polynomial
                                                                                      \DIc{px$^2+$(pb$-$1)x+\One\Comma}{5}
\Par it follows from \pLinkLocalDisplayText{4} that that polynomial has a single root. This in turn implies that the homography \Vol{-\One,\One,-p,-pb}\ has exactly one fixed point, and thus is a translation. Letting
                                                         \Dc{\OpGr{t}\ \Equiv\ \Vol{-\One,\One,-p,-pb}\Andd\OpGr{p}\ \Equiv\ \Vol{\Zero,\One,-p,\Zero}\Comma}
\Par we see that from \pLinkLocalDisplayText{2} that \pLinkLocalDisplayText{1} holds.
\PaR It remains to show that \OpGr{p}\ has no fixed point. But the fixed point equation for \OpGr{p}\ is
                                                                                                 \Dc{\One\Equals-px$^2$}
 \Par and, since \One\ and px$^2$ are both in \Mabf{\zEro}{\oNe}{\iNfty}\Comma\ it follows from \pLinkDisplayText{EA}{T4}{0} that this equation has no solution. \QED
 \Item{T11}{Theorem} Let M be an exponential meridian and \OpGr{f}\ an element of \LLL{\Mii{M}}\Period\ Let \Op{b}\Equals\Pr{\Zero,\One,\Infty}\ be any ordered basis for M and suppose that \OpGr{f}\Equals\Vol{a,b,c,d}\ relative to that basis. Then
                                                 \DIc{\BDet{a}{b}{c}{d}{b}\In\Mabf{\zEro}{\oNe}{\iNfty}\quad\Iff\quad f is a composition of involutions without fixed points.}{1}
 \Proof [\Implied]\Quad If \OpGr{f}\ is a composition of involutions without fixed points, then it follows from \pLinkItemText{EA}{T7} that the matrix of each component has a determinant in \Mabf{\zEro}{\oNe}{\iNfty}\Period\ It follows from Lemma \pLinkItemText{EA}{PlusTimes} that the matrix of \OpGr{f}\ would have determinant in \Mabf{\zEro}{\oNe}{\iNfty} \Period
\Par\Quad[\Implies]\Quad Let then \Op{b}\ be an ordered basis for M and suppose that the determinant of the matrix of \OpGr{f}\Equals\Vol{a,b,c,d}\ is in \Mabf{\zEro}{\oNe}{\iNfty}\Period\ If \OpGr{f}\ is an involution, then it follows from Theorem \pLinkItemText{EA}{T7} that the determinant of the matrix is in \Mabf{\zEro}{\oNe}{\iNfty}. If \OpGr{f}\ is a translation, then it follows from \pLinkItemText{EA}{L8} that it is the composition of an involution without fixed elements of \MNF\Comma\ and so its matrix has determinant in \Mabf{\zEro}{\oNe}{\iNfty} \Period\ If \OpGr{f}\ is a pure rotation, then Lemma \pLinkItemText{EA}{L10} implies that it is the composition of an involution without fixed points and a translation. Since the determinants of the matrices of both these constituents are in \Mabf{\zEro}{\oNe}{\iNfty}\Comma\ it follows that the matrix of the rotation is as well.
\PaR Thus we may presume  that \OpGr{f}\ is a non-involutive dilation and that the determinate of the matrix of \Vol{a,b,c,d}\ is in \Mabf{\zEro}{\oNe}{\iNfty}\Period\ Let \OpGr{th}\ be any element of \LLL{\Mii{M}}\ which sends \Zero to one of the fixed points of \OpGr{f}\ and \Infty\ to the other. The \OpGr{th}\Circ\OpGr{f}\Circ\Inv{\OpGr{th}}\ is a non-involutive dilation and so has a matrix relative to \Op{b}\ of the form \Vol{r,\Zero,\Zero,s}. Since \BDet{a}{b}{c}{d}\ is in \Mabf{\zEro}{\oNe}{\iNfty}\Comma\ it follows that rs\Equals\BDet{r}{\Zero}{\Zero}{s}\ is as well. Thus we can choose r to be \One\ and m\In\Mabf{\zEro}{\oNe}{\iNfty}\ such m$^2$\Equals s\Period\ Thus
                                                              \DIc{\Vol{r,\Zero,\Zero,s}\Equals\Vol{\Zero,\One,-m,\Zero}\Circ\Vol{\Zero,-m,\One,\Zero}\Period}{2}
\Par Evidently \BDet{\Zero}{-m}{\One}{\Zero}\ and \BDet{\Zero}{\One}{-m}{\Zero}\ are in \Mabf{\zEro}{\oNe}{\iNfty}\Comma\ whence follows that the determinants of matrices of \Inv{\OpGr{th}}\Circ\Vol{\Zero,-m,\One,\Zero}\Circ\OpGr{th}\ and \Inv{\OpGr{th}}\Circ\Vol{\Zero,\One,-m,\Zero}\Circ\OpGr{th}\ are as well. Since
          \Dc{\OpGr{f}\Equals\Inv{\OpGr{th}}\Circ\OpGr{th}\Circ\OpGr{f}\Circ\Inv{\OpGr{th}}\Circ\OpGr{th}\Equals\Inv{\OpGr{th}}\Circ\Vol{\Zero,\One,-m,\Zero}\Circ\Vol{\Zero,-m,\One,\Zero}\Circ\OpGr{th}\Equals%
                                              (\Inv{\OpGr{th}}\Circ\Vol{\Zero,\One,-m,\Zero}\Circ\OpGr{th})\Circ(\Inv{\OpGr{th}}\Circ\Vol{\Zero,-m,\One,\Zero}\Circ\OpGr{th})}
\Par and so \OpGr{f}\ is the composition of two involutions without fixed points. \QED
 \Item{D5}{Definition} Recall that \Mii{G}\ is the notation for \LLL{\Mii{M}}\Colon\ the group of homographies of M. We shall write
                                                                                                      \DIc{\Mii{G}$^+$}{0}
 \Par for the set of elements of \Mii{G}\ which are compositions of involutions without fixed points. This terminology is motivated by the characterization \pLinkItemText{EA}{T11}. In view of \pLinkItemText{EA}{T11}, \Mii{G}$^+$ is a normal subgroup of \Mii{G}\Period\Foot{A subgroup N of a group G be \Def{normal} if, for each x\In G and n\In N, x\Cdot n]\Cdot\Inv{x}\In N\Period}
 \PaR We define an equivalence relation \OSim\ on the family of all ordered bases of M by, for any two such \Pr{a,b,c}\ and \Pr{r,s,t},
                                \DIc{\Pr{a,b,c}\OSim\Pr{r,s,t}\quad\Iff\quad\ThereIs{\OpGr{f}\In\Mii{G}$^+$ }\Quad\Pr{r,s,t}\Equals\Pr{\OpGr{f}\Of{a},\OpGr{f}\Of{b},\OpGr{f}\Of{c}}\Period}{1}
 \PaR There are two equivalence classes and we shall denote this set of two elements by
                                                                                              \DIc{\Set{\COreo,\CCOreo\Period}}{2}
 \Par These equivalence classes will be called \Def{orientations}.
 \PaR Let \COreo\ be an orientation of an exponential meridian M. Let a and b be distinct elements of M. We write
                                                                      \DI{8}{0}{\Arc{a,b}\quad\Equiv\quad\SetSuch{t\In\Mab{a}{b}}{\Pr{a,t,b}\In\COreo}}{3}
                                                                      \DII{8}{8}{\kn{-12}\CArc{a,b}\quad\Equiv\quad\Set{a,b}\Cup\Arc{a,b}\Period}{and}{4}
 \Par These two sets \Arc{a,b}\And\CArc{a,b}, respectively, will be called the \COreo-\Def{open} \Def{arc} \Def{associated} \Def{with} \Pr{a,b}\ and the \COreo-\Def{closed} \Def{arc} \Def{associated} \Def{with} \Pr{a,b}\Comma\ respectively. Thus, for \vruleh{12}t\In\Arc{a,b}, we have
                                                                \DIc{\Mabf{a}{t}{b}\Equals\Arc{a,b}\Andd\Set{a,b}\Cup\Mabf{a}{t}{b}\Equals\CArc{a,b}\Period}{5}
 \PaR When the orientation \COreo\ is evident from the context, we shall sometimes omit the symbol \COreo\ from the terms \COreo-open and \COreo-closed.
 \Item{P11}{Proposition} In an exponential meridian M
                                                                                     \DI{8}{0}{translations are in \Mii{G}$^+$,}{1}
                                                                                       \DI{8}{0}{pure rotations are in \Mii{G}$^+$,}{2}
                                                                             \DI{8}{0}{involutions with fixed points are not in \Mii{G}$^+$,}{3}
                                             \DI{8}{0}{\ForAll{\Pr{a,b,t} an ordered basis}\ForAll{u\In\Mabf{a}{b}{t}}\quad\Pr{a,b,t}\And\Pr{a,b,u}\ have the same orientation,}{4}
                                            \DI{8}{0}{\ForAll{\Pr{a,b,t} an ordered basis}\ForAll{u\NIn\Mabf{a}{t}{b}}\quad\Pr{a,t,b}\And\Pr{a,u,b}\ have opposite orientations,}{5}
                                                          \DI{8}{0}{\ForAll{\Pr{a,t,b} an ordered basis}\quad\Pr{a,b,t}\And\Pr{b,a,t}\ have opposite orientations}{6}
                                                    \DII{8}{8}{\kn{12}\ForAll{\Pr{a,t,b} an ordered basis}\quad\Pr{a,b,t}\And\Pr{b,t,a}\ have the same orientation.}{and}{7}
 \Proof [\Implies\pLinkLocalDisplayText{1}]\Quad Follows from Lemma \pLinkItemText{EA}{L8}.
 \Par\Quad[\Implies\pLinkLocalDisplayText{2}]\Quad Follows from Lemma \pLinkItemText{EA}{L10}.
 \Par\Quad[\Implies\pLinkLocalDisplayText{3}]\Quad Follows from Theorem \pLinkItemText{EA}{T7}.
 \Par\Quad[\Implies\pLinkLocalDisplayText{4}]\Quad The involutions \Vol{a,b;t,u}\ and \Vol{a,b;u,u}\ both have fixed points and so the determinants of their matrices are not positive. It follows that the matrix of their composition is positive and so \Vol{a,b;u,u}\Circ\Vol{a,b;t,u}\ is in \In\Mii{G}$^+$. Furthermore \Vol{a,b;u,u}\Circ\Vol{a,b;t,u}\Circ\Pr{a,b,t}\Equals\Pr{a,b,u}\Period
 \Par\Quad[\Implies\pLinkLocalDisplayText{5}] The involution \Vol{a,b;t,u}\ has no fixed point by Theorem \pLinkItemText{EA}{T4} and \Pr{b,a,t}\Equals\Vol{a,b;t,u}\Circ\Pr{a,b,t}\Period
  \Par\Quad[\Implies\pLinkLocalDisplayText{6}] The involution \Vol{a,b;t,t}\ has a fixed point and \Pr{b,a,t}\Equals\Vol{a,b;t,t}\Circ\Pr{a,b,t}\Period
 \Par\Quad[\Implies\pLinkLocalDisplayText{7}] The involutions \Vol{a,b;t,t}\ and \Vol{a,t;b,b} have fixed points and so the determinants of their matrices are not positive. It follows that the matrix of their composition is positive and so \Vol{a,t;b,b}\Circ\Vol{a,b;t,t}\ is in \In\Mii{G}$^+$. Furthermore \Vol{a,t;b,b}\Circ\Vol{a,b;t,t}\Circ\Pr{a,b,t}\Equals\Pr{b,t,a}\Period\ \QED
 \Item{C11}{Corollary} Let M be an exponential meridian with an orientation \COreo\ and let a and b be distinct elements of M. Then
                                                                       \DI{8}{8}{\Set{\Set{a},\Set{b},\Arc{a,b},\Arc{b,a}}\ is a partition of M\Comma}{1}
                                                                       \DI{0}{8}{\ForAll{t\In\Arc{a,b}}\quad\Set{\Arc{a,t},\Set{t},\Arc{t,b}}\ is a partition of \Arc{a,b}}{2}
                                                                       \DII{0}{8}{\kn{-12}\ForAll{t\In\Arc{b,a}}\quad\Set{\Arc{b,t},\Set{t},\Arc{t,a}}\ is a partition of \Arc{b,a}\Period}{and}{3}
\Proof This follows from \pLinkItemText{EA}{CO9} and \pLinkItemText{EA}{P11}. \QED
  \Item{T12}{Theorem} Let \Set{a,b}\And\Set{p,q}\ be subsets of a circular meridian M, both of cardinality 2. Then
   \DI{8}{0}{if p\Equals a or q\Equals b, either\quad\Arc{a,b}\Sin\Arc{p,q}\Orr\Arc{p,q}\Sin\Arc{a,b}\Semicolon}{1}
   \DI{8}{0}{if p\Equals b, then\quad(\Arc{a,b}\Cap\Arc{p,q})\In\Set{\Void,\Arc{a,q}}\Semicolon}{2}
   \DI{8}{0}{if q\Equals a, then\quad(\Arc{a,b}\Cap\Arc{p,q})\In\Set{\Void,\Arc{p,b}}\Semicolon}{3}
   \DI{8}{0}{if \Set{p,q}\Sin\Arc{a,b}\Comma\ then\quad\Arc{p,q}\Sin\Arc{a,b}\Orr\Arc{a,b}\Cap\Arc{p,q}\Equals\Arc{p,b}\Cup\Arc{a,q}\Semicolon}{4}
   \DI{8}{0}{if \Set{a,b}\Sin\Arc{p,q}\Comma\ then\quad\Arc{a,b}\Sin\Arc{p,q}\Orr\Arc{p,q}\Cap\Arc{a,b}\Equals\Arc{q,a}\Cup\Arc{b,p}\Semicolon}{5}
   \DI{8}{0}{if p\In\Arc{a,b}\ and q\NIn\Arc{a,b}\Comma\ then\quad\Arc{a,b}\Cap\Arc{p,q}\Equals\Arc{p,b}\Semicolon}{6}
   \DII{8}{0}{\kn{-12}if p\NIn\Arc{a,b}\ and q\In\Arc{a,b}\Comma\ then\quad\Arc{a,b}\Cap\Arc{p,q}\Equals\Arc{a,q}\Period}{and}{7}
 \Proof\Quad[\pLinkLocalDisplayText{1}]\Quad We presume that p\Equals a\Comma\ the other case being susceptible to analogous treatment. If q\In\Arc{a,b}\Comma\ it follows from \pLinkDisplayText{EA}{C11}{2} that
                                                 \DIc{\Arc{a,b}\Equals\Arc{a,q}\Cup\Set{q}\Cup\Arc{q,b}\quad\Implies\quad\Arc{p,q}\Equals\Arc{a,q}\Sin\Arc{a,b}\Period}{8}
  \Par Suppose, on the other hand that q\NIn\Arc{a,b}\Period\ If q\Equals b\Comma\ then \Arc{p,q}\Equals\Arc{a,b}\ so by \pLinkDisplayText{EA}{C11}{1} we may presume that q is in \Arc{b,a}\Period\ From \pLinkDisplayText{EA}{C11}{2} follows that
                                                           \DIc{\Arc{q,p}\Sin\Arc{b,a}\quad\Imp{\hLinkDisplayText{EA}{C11}{1}}\quad\Arc{a,b}\Sin\Arc{p,q}\Period}{9}
  \Par From \pLinkLocalDisplayText{8} and \pLinkLocalDisplayText{9} follows \pLinkLocalDisplayText{1}.
  \Par\Quad[\pLinkLocalDisplayText{2}]\Quad Here we presume that p\Equals b\Period\ Suppose that q\In\Arc{a,b}\Period\ From \pLinkDisplayText{EA}{C11}{2} follows that
                       \DIc{a\NIn\Arc{q,b}\Equals\Arc{q,p}\quad\Imp{\hLinkDisplayText{EA}{C11}{1}}\quad a\In\Arc{p,q}\quad\Imp{\hLinkDisplayText{EA}{C11}{2}}\quad\Arc{p,a}\Sin\Arc{p,q}\Period}{10}
  \Par Now suppose that q\NIn\Arc{a,b}\Period\ If q\Equals a\Comma\ then \Arc{a,b}\Cap\Arc{p,q}\Equals\Void\ by \pLinkDisplayText{EA}{C11}{1} so we shall presume that q\In\Arc{b,a}\Period\ From \pLinkDisplayText{EA}{C11}{2} follows that
                                           \DIc{\Arc{p,q}\Equals\Arc{b,q}\Sin\Arc{b,a}\quad\Imp{\hLinkDisplayText{EA}{C11}{1}}\quad\Arc{a,b}\Cap\Arc{p,q}\Equals\Void\Period}{11}
  \Par From \pLinkLocalDisplayText{10} and \pLinkLocalDisplayText{11} follows \pLinkLocalDisplayText{2}.
  \Par\Quad[\pLinkLocalDisplayText{3}]\Quad The proof of \pLinkLocalDisplayText{3} is analogous to that of \pLinkLocalDisplayText{2}.
  \Par\Quad[\pLinkLocalDisplayText{4}]\Quad We presume that \Set{p,q}\Sin\Arc{a,b}\Period\ From \pLinkDisplayText{EA}{C11}{2} follows that \Arc{a,p}\Cup\Set{p}\Cup\Arc{p,b}\Equals\Arc{a,b}\Period\ If q is in \Arc{p,b}\Comma\ then \Arc{p,q}\Cup\Set{q}\Cup\Arc{q,b}\Equals\Arc{p,b}\ and so \Arc{p,q}\Sin\Arc{a,b}\Period\ Thus we may presume that q\In\Arc{a,p}\Period\ From \pLinkDisplayText{EA}{C11}{2} follows that \Arc{a,q}\Cup\Set{q}\Cup\Arc{q,p}\Equals\Arc{a,p}\ and so \Set{\Arc{a,q},\Set{q},\Arc{q,p},\Set{p},\Arc{p,b}}\ is a partition of \Arc{a,b}\Period\ Since \pLinkDisplayText{EA}{C11}{1} implies that
  \Arc{p,q}\Cap\Arc{q,p}\Equals\Void\Comma\ we have \Arc{a,b}\Cap\Arc{p,q}\Equals\Arc{p,b}\Cap\Arc{a,q}\Comma\ which establishes \pLinkLocalDisplayText{4}.
  \Par\Quad[\pLinkLocalDisplayText{5}]\Quad The proof of \pLinkLocalDisplayText{5} is analogous to that of \pLinkLocalDisplayText{4}.
  \Par\Quad[\pLinkLocalDisplayText{6}]\Quad We presume that p\In\Arc{a,b}\ and q\NIn\Arc{a,b}\Period\ From \pLinkDisplayText{EA}{C11}{1} follows that \Set{\Set{a},\Arc{a,b},\Set{b},\Arc{b,a}}\ is a partition of M. From \pLinkDisplayText{EA}{C11}{2} follows that \Set{\Arc{a,p},\Set{p},\Arc{p,b}}\ is a partition of \Arc{a,b}\ and that \Set{\Arc{b,q},\Set{q},\Arc{q,a}}\ is a partition of \Arc{b,a}\Period\ From these three facts follows that
                                                    \DIc{\Set{\Set{a},\Arc{a,p},\Set{p},\Arc{p,b},\Set{b},\Arc{b,q},\Set{q},\Arc{q,a}}\ is a partition of M\Period}{12}
 \Par Evidently \pLinkLocalDisplayText{12} implies \pLinkLocalDisplayText{6}\Period
 \Par\Quad[\pLinkLocalDisplayText{7}]\Quad The proof of \pLinkLocalDisplayText{7} is analogous to that of \pLinkLocalDisplayText{6}. \QED
 \Item{D13}{Definition} Let M be an exponential meridian with an orientation \COreo\Period\ It follows from \pLinkItemText{EA}{T12} that the family of all open arcs is a base for a topology, which we shall call the \Def{arc} \Def{topology}.\Foot{We here adopt the convention that the empty set \Void\ is an arc.} A necessary and sufficient condition for a subset of an exponential meridian to be open relative to the arc topology is for it to be a union of arcs. We recall that a topological space is \Def{compact} if every open covering has a finite sub-covering.
 \Item{T14}{Theorem} Let M be an exponential meridian which is compact relative to the arc topology. Then M is isomorphic to a circle meridian. In particular, relative to the arc topology, M is homeomorphic to a circle. Furthermore, if \Op{b}\ is any basis for M, then F$_{\hbox{\Op{\RMBx{7}{b}}}}$ is isomorphic to the field of real numbers.
 \Proof Let \COreo\ be the orientation of M of which \Op{b}\Equals\Pr{\Zero,\One,\Infty}\ is an element. On the field F$_{\hbox{\Op{\RMBx{7}{b}}}}$ we define the relation \Order\Sin(F$_{\hbox{\Op{\RMBx{7}{b}}}}$\Cross F$_{\hbox{\Op{\RMBx{7}{b}}}}$) by
                                                  \DIc{\ForAll{\Set{x,y}\Sin F$_{\hbox{\Op{\RMBx{7}{b}}}}$ }\quad x\Order y\Quad\Iff\Quad y$-$x\In\CArc{\Zero,\Infty}\Period}{1}
 \Par For x\In\BField\ we have
                                                                     \DIc{x$-$x\Equals\Zero\quad\Implies\quad x$-$x\In\CArc{\Zero,\Infty}\Period}{2}
 \Par For \Set{x,y}\Sin\BField\Comma\ it follows from \pLinkItemText{EA}{T4} that
                                                                      \DIc{either y$-$x or x$-$y is in \Mabf{\zEro}{\oNe}{\iNfty}\Comma\ but not both.}{3}
 \Par If \Set{x,y,z}\Sin\BField\ and both x\Order y and y\Order z\Comma\ then
                         \DIc{\Set{y$-$x,z$-$y}\Sin\Mabf{\Zero}{\One}{\Infty}\Imp{\hLinkItemText{EA}{PlusTimes}}((y$-$x$)+$(z$-$y)\In\Mabf{\Zero}{\One}{\Infty}\quad\Implies\quad x\Order z\Period}{4}
 \Par It follows from \pLinkLocalDisplayText{2}, \pLinkLocalDisplayText{3} and \pLinkLocalDisplayText{4} that \Order\ is a total order in the sense of \pLinkItemText{Graphs}{Order}. For \Set{x,y,z}\In\BField\ such that x\Order y\Comma\ we have
                                                       \DIc{y$-$x\In\Mabf{\zEro}{\oNe}{\iNfty}\Andd(z$+$y)-(z$+$x)\Equals y$-$x\quad\Implies\quad(z$+$x)\Order(z$+$y)}{5}
 \Par and for \Set{x,y}\Sin\Mabf{\zEro}{\oNe}{\iNfty}\ we have by \pLinkItemText{EA}{PlusTimes}
                                                                   \DIc{xy\In\Mabf{\zEro}{\oNe}{\iNfty}\quad\Implies\quad\Zero\Order\thinspace xy\Period}{6}
 \PaR It is a classical result of real analysis that if a field with a total order satisfies \pLinkLocalDisplayText{2} through \pLinkLocalDisplayText{6}, then it is isomorphic with the field of real numbers if, and only if, to each subset of \BField\ with an upper bound corresponds a least upper bound. Let then S\Sin\BField\ and b\In\BField\ be such that
                                                                                           \Dc{\ForAll{x\In S}\quad x\Order b\Period}
 \Par For \Set{x,y}\Sin S\Comma\ we have that x\Order b and y\Order b\Comma\ whence follows that S is a directed set. Thus the identity function \Id{S}\ is a net and, since M is compact, there is a subnet of \Id{S}\ which converges to some element m of M. Specifically, this means that there exists a directed set D with direction \ORder\ and a co-final function \Function\OpGr{w}sendsinDtoinS\end\ such that \Id{S}\Circ\OpGr{w}\ converges to m. \PaR Assume that m were not an upper bound for S. Then there would exist s\In S such that m\Order s and m\NEq s. In that case there would exist d\In D such that, for all e\In D with d\ORder e,
                                                          \Dc{s\Order\OpGr{w}\Of{e}\quad\Implies\quad\Id{S}\Circ\OpGr{w}\ is eventually outside of \Arc{\Two m$-$s,s}}
 \Par which, since m is in \Arc{\Two m$-$s,s}\Comma\ would contradict the fact that \Id{S}\Circ\OpGr{w}\ converges to m. It follows that m is an upper bound of S.
 \PaR Let u be any upper bound of S such that u\Order m. If m\NEq u, then m would be in \Arc{u,\Two m}\ and so \Id{S}\Circ\OpGr{w}\ would be eventually in \Arc{u,\Two m}\Period\ But, u being an upper bound of S, this would be absurd. It follows that m is a least upper bound for S. \QED
 \vfill\eject
                                                                                  \Section{SetTheory}{Appendix I: Mathematical Notation and Terminology}\xrdef{pageSetTheory}
\def\Sing#1{\Bf{#1}}
\def\OrderShriek{\hbox{\Stack{\BF{7}{!}}{\Order}{-1pt}}}%
\PAR The language of modern mathematics is set theory. We set down here the many of the basic notions used in this paper and its sequel.
 \Item{Logic}{Logical Notation} We shall write
                                                                                                      \DIc{\MiiChBx{71}}{i}
 for the \Def{existential} \Def{quantifier} which is read as \Quotes{there exists} or \Quotes{for some}: the notation
                                                                                                 \DIc{\MiiChBx{71}\kn{1}!}{ii}
 is read as \Quotes{there exists a unique} or \Quotes{for exactly one}. We write
                                                                                                     \DIc{\MiiChBx{70}}{iii}
 for the \Def{universal quantifier} which is read as \Quotes{for all} or \Quotes{for each} or \Quotes{for every}.
 \Par The notation
                                                                                                       \DIc{\Implies}{iv}
 is read \Quotes{implies}, the notation
                                                                                                        \DIc{\Implied}{v}
 is read \Quotes{follows from} and the notation
                                                                                                         \DIc{\Iff}{vi}
 is read \Quotes{is a statement equivalent to the statement}.
 \Item{Sets}{Sets} A \Def{set} is a collection of objects. Such objects are called \Def{members} or \Def{elements} of the set. We denote that an object x is a member of a set X by
                                                                                                \DIc{x\In X\Orr X\Ni x\Period}{1}
We denote than an object x is not a member of X by
                                                                                                     \DIc{x\Nin X\Period}{2}
If X and Y are sets, then X is said to be a \Def{subset} \Def{of} Y unless
                                                                                          \DIc{\ThereIs{x\In X}\quad x\Nin Y\Period}{3}
If X is a subset of Y, we say that Y is a \Def{superset} \Def{of} X: we express this in symbols by
                                                                                               \DIc{X\Sin Y\Orr Y\Sout X\Period}{4}

 \PaR We convene that there is exactly one set with no elements, which we call \Def{void}, or the \Def{empty} \Def{set} or the \Def{null} \Def{set}:
                                                                                                      \DIc{\Void\Period}{5}
For any set X, it follows from \pLinkLocalDisplayText{3} that
                                                                                                   \DIc{\Null\Sin X\Period}{6}
If x,\dots,y is a list of the elements of a set X, we express this fact by
                                                                                           \DIc{X\ \Equals\ \Set{x,\dots,y}\Period}{7}
\Par A set with a single element is called a \Def{singleton}, one with exactly two elements a \Def{doubleton} \It{et}\ \It{cetera}.
\PaR If P\Of{x} is a condition on an object x, we denote the subset of a set X consisting of all those elements of X which satisfy the condition P\Of{x} by
                                                                                            \DIc{\SetSuch{x\In X}{P\Of{x}}\Period}{8}
 \Par The symbol
                                                                                                        \DIc{\Equiv}{vii}
 is used to define a symbol placed to the left of \Quotes{\Equiv}, as the set which is placed to the left. For instance the \Def{union} \Def{of} \Def{two} \Def{subsets} X \Def{and} Y of a third set Z is denoted and defined by
                                                                              \DIc{X\Cup Y\ \Equiv\ \SetSuch{t\In Z}{t\In X\Or t\In Y}\Period}{9}
The \Def{intersection} \Def{of} \Def{two} \Def{subsets} X and Y of a third set Z is denoted and defined by
                                                                              \DIc{X\Cap Y\ \Equiv\ \SetSuch{t\In Z}{t\In X\And t\In Y}\Period}{10}
\Par The terms \Def{family} and \Def{collection} are  synonyms for the term set, although family is often reserved for sets, the elements of which are sets themselves. The family of all subsets of a given set X is denoted
                                                                                                 \DIc{\Subsets{X}\Period}{11}
Thus
                                                                                         \DIc{S\In\Subsets{X}\Iff S\Sin X\Period}{12}
For subsets X and Y of a set Z, the \Def{symmetric} \Def{difference} \Def{of} X \Def{and} Y is defined and denoted by
                                                                             \DIc{X\Cop Y\ \Equiv\ \SetSuch{t\In X\Cup Y}{t\NIn X\Cap Y}\Period}{13}
When Y\Sin X, the symmetric difference X\Cop Y is often called the \Def{complement} \Def{of} Y \Def{in} X.
 \PaR The two \Def{distributive} \Def{laws} hold for any three subsets A, B and C of a set X:
 \D{10}{0}{A\Cap(B\Cup C)\Equals(A\Cap B)\Cup(A\Cap C)\Semicolon}
 \DI{0}{0}{}{14}
 \D{0}{7}{A\Cup(B\Cap C)\Equals(A\Cup B)\Cap(A\Cup C)\Period}
 A third law, which is symmetric in \Cup\ and \Cap\ is
                                                                 \DIc{(A\Cup B)\Cap(B\Cup C)\Cap(C\Cup A)\Equals(A\Cap B)\Cup(B\Cap C)\Cup(C\Cap A)\Period}{15}
 The \Def{De} \Def{Morgan} \Def{Laws} for subsets A and B of a set X are
 \D{10}{0}{X\Cop(A\Cup B)\Equals(X\Cop A)\Cap(X\Cop B)\Semicolon}
 \DI{0}{0}{}{16}
 \D{0}{7}{X\Cop(A\Cap B)\Equals(X\Cop A)\Cup(X\Cop B)\Period}
 \PaR For any subfamily \Mii{S}\Sin\Mii{P}\Of{X} of subsets of a set X, we denote
        \DIc{\Union{\RM{7}{S}\iN\MII{7}{S}}{S}\ \Equiv\ \SetSuch{x\In X}{\ThereIs{S\In\Mii{S}}\ x\In S}\Andd\Intersection{\RM{7}{S}\iN\MII{7}{S}}{S}\ \Equiv\ \SetSuch{x\In X}{\ForAll{S\In\Mii{S}}\ x\In S}\Period}{17}
  \Item{CarProducts}{Cartesian Products} A \Def{pair} is a set of the form
                                                                                                       \DIc{\Set{x,y}}{1}
where x\NEq y\Period\ An \Def{ordered} \Def{pair} \Def{of} \Def{elements} x \Def{and} y is defined and denoted by
                                                                                       \DIc{\Pr{x,y}\ \Equiv\ \Set{x,\Set{x,y}}\Period}{2}
We say that x is the \Def{first} \Def{coordinate} of the ordered pair \Pr{x,}\And y the \Def{second} \Def{coordinate}. For subsets X and Y of a set Z, the \Def{Cartesian} \Def{product} \Def{of} X \Def{with} Y is defined and denoted by
                                                                            \DIc{X\Cross Y\ \Equiv\ \SetSuch{\Pr{x,y}}{x\In X\And y\In Y}\Period}{3}
\PaR More generally, a \Def{triple} is a set \Set{x,y,z}\ with exactly 3 elements, and an \Def{ordered} \Def{triple} is defined and denoted by
                                                                               \DIc{\Pr{x,y,z}\ \Equiv\ \Set{x,\Set{x,y},\Set{x,y,z}}\Period}{4}
\Def{Quadruples}, \Def{ordered} \Def{quadruples}, \Def{quintuples}, \Def{ordered} \Def{quintuples} \It{et} \It{cetera} are defined analogously.
 \Item{Finite}{Finite Sets} A subset S of a set X is a \Def{proper} \Def{subset} if S\NIn\Set{\Void,X}. A set X is said to be \Def{infinite} if there exists a proper subset S of X and a subset \Order\ of X\Cross S such that
                                                                         \DIc{\ForAll{x\In X}\ThereIsShriek{s\In S}\quad\Pr{x,s}\In\Order\Period}{1}
 A set which is not infinite is called \Def{finite}. Every subset of a finite set is finite.
  \Item{PropertiesAndClasses}{Conditions And Classes} A \Def{condition} \Def{on} \Def{a} \Def{set} is a statement about a set which can be determined to either be \Def{true} or \Def{false}, but not both. A \Def{class} is the aggregate of those sets for which some property is true. A set of a class is called a \Def{member} of that class.
  \PaR A family is always a class, but a class may not be a family. We demonstrate this to be true by what is known as \Quotes{Russell's paradox}, published by Bertrand Russell in 1901, but first enunciated by Ernst Zermelo in 1900:
  \PaR Assume that all classes are families. Let \Mii{F}\ be the family of all infinite sets. Since \Mii{F}\ is infinite itself, it follows that \Mii{F}\In\Mii{F}. Let \Mii{C}\ be the family of all sets which are not elements of themselves. It is trivial that \Void\ is an element of \Mii{C}\ and we have seen that \Mii{F}\NIn\Mii{C}. If \Mii{C}\In\Mii{C}, then \Mii{C}\ is not in \Mii{C} by definition of \Mii{C}. But if \Mii{C}\NIn\Mii{C}, then \Mii{C}\ is in \Mii{C}, again by the definition of \Mii{C}. Thus our assumption cannot have been correct.
  \PaR So how can we know if a class is a family (set)? If a class has been constructed from sets using unions, intersections, subsets or Cartesian products, we can be sure it is a set. Ernst Zermelo suggested another method, called the \Def{axiom} \Def{of} \Def{choice}: to wit
                                \DIc{\ForAll{\Mii{S} a family of sets}\ThereIs{a set X}\ForAll{S\In\Mii{S}}\ThereIs{x$_{\RM{7}{S}}$\In S}\quad X\Cup S\Equals\Set{x$_{\RM{7}{S}}$}\Period}{1}
  Such a function
                                                                                    \Dc{\Function sendsSin\Mii{S}tox$_{\RM{7}{S}}$inS\end}
  is called a \Def{choice} \Def{function}. It is known that if one constructs new sets using the axiom of choice, one is not thereby led into logical contradictions. However there is no way to deduce that the axiom of choice holds from the other axioms of set theory. One can use it or not, depending on ones taste.
 \Item{Relations}{Equivalence Relations} The terms \Quotes{relation} and \Quotes{graph} are in some places treated as synonyms. Here we shall treat the term graph as being more restrictive than a relation, and shall reserve the next section for its presentation.
 \PaR Let \Mii{A}\And\Mii{B}\ be two classes or sets. By a \Def{relation} \Def{from} \Mii{A}\ to \Mii{B}, we shall mean a class of ordered pairs \Pr{A,B}\ where each A is a member of \Mii{A}\ and each B a member of \Mii{B}. \PaR if \Mii{A}\Equals\Mii{B}\And the following three conditions hold, we say that a relation \Sim\ from \Mii{A}\ to \Mii{B}\ is an \Def{equivalence} \Def{relation}:
                                                                             \DI{8}{0}{\ForAll{A\In\Mii{A}}\quad\Pr{A,A}\In\Sim\Semicolon}{1}
                                                           \DI{8}{0}{\ForAll{\Set{A,B}\Sin\Mii{A}}\quad\Pr{A,B}\In\Sim\ \Iff\  \Pr{B,A}\In\Sim\Semicolon}{2}
                                                     \DIc{\ForAllSuch{\Set{A,B,C}\Sin\Mii{A}}{\Pr{A,B}\In\Sim\And\Pr{B,C}\In\Sim}\quad\Pr{A,C}\In\Sim\Period}{3}
 \PaR An \Def{equivalence} \Def{class} \Def{associated} \Def{with} \Def{the} \Def{equivalence} \Def{relation} \Sim\ \Def{and} \Def{a} \Def{member} A \Def{of} \Mii{A}\ is
                                                              \DIc{the class of all members B of \Mii{A}\ such that \Pr{A,B}\ is a member of \Mii{A}\Period}{4}
 The \Def{partition} \Def{of} \Mii{A}\ \Def{associated} \Def{with} \Def{the} \Def{equivalence} \Def{relation} \Sim\ is the class  of all equivalence classes associated with \Sim.
 \PaR Conversely, for a class or set \Mii{A}\ we define a \Def{partition} \Def{of} \Mii{A}\ to be a class \Mii{P}\ of sub-classes of \Mii{A}\ such the such that each member of \Mii{A}\ is a member of a unique member of \Mii{P}. The \Def{equivalence} \Def{relation} \Sim \Def{associated} \Def{with} \Def{the} \Def{partition} \Mii{P}\ is the class of all ordered pairs \Pr{X,Y}\ where X and Y are members of a common member of \Mii{P}.
 \vfill\eject
                                                                                                   \Section{Graphs}{Appendix II: Graphs}\xrdef{pageGraphs}
\PAR Many things can be described in the context of a subset of the Cartesian product of two sets. We present some of the more important ones in the present section.
\Item{Definitions}{Definitions} A subset \Order\  of a Cartesian product X\Cross Y of sets is said to be a \Def{graph} \Def{of} X \Def{with} Y. For x\In X and y\In Y, we shall sometimes write
                                                                                            \DIc{x\thinspace\Order\thinspace y}{0}
 to indicate that \Pr{x,y}\ is an element of \Order\ and
                                                                                        \DIc{x\thinspace\NOrder\thinspace y\Comma}{1}
 to indicate that \Pr{x,y}\ is not an element of \Order. To indicate that \Pr{x,y}\Sin\Order\ but x\NEq y, we shall write
                                                                                                \DIc{x\OrderShriek y\Period}{2}
 \PaR The \Def{domain} \Def{of} \Def{a} \Def{graph} \Order\Sin X\Cross Y is defined and denoted by
                                                                \DIc{\Domain{\Order} \Equiv\ \SetSuch{x\In X}{\ThereIs{y\In Y}\ \Pr{x,y}\In\Order}}{3}
and its \Def{range} by
                                                             \DIc{\Range{\Order} \Equiv\ \SetSuch{y\In Y}{\ThereIs{x\In X}\ \Pr{x,y}\In\Order}\Period}{4}
 The \Def{world} \Def{of} \Def{a} \Def{graph} \Order\ is defined and denoted by
                                                                          \DIc{\World{\Order}\ \Equiv\ \Domain{\Order}\Cup\Range{\Order}\Period}{5}
 \PaR The \Def{inverse} \Def{of} \Def{a} \Def{graph} \Order\ is defined and denoted by
                                                                         \DIc{\Inv{\Order}\ \Equiv\ \SetSuch{\Pr{y,x}}{\Pr{x,y}\In\Order}\Period}{6}
 In particular
                                                                                      \DIc{\Inv{(X\Cross Y)}\ \Equals\ Y\Cross X\Period}{7}
 \PaR The \Def{complement} \Def{of} \Def{a} \Def{graph} \Order\ is the set
                                                \DIc{\Cop\Order\ \Equiv\ \SetSuch{\Pr{x,y}\In\Domain{\Order}\Cross\Range{\Order}}{\Pr{x,y}\NIn\Order}\Period}{8}
 \PaR If \Order\And\OpGr{s}\ are graphs, then the \Def{composition} \Def{of} \Order\ \Def{with} \OpGr{s}\ is defined and denoted by
                                      \DIc{\Order\Circ\OpGr{s}\ \Equiv\ \SetSuch{\Pr{x,z}}{\ThereIs{y\In\Domain{\OpGr{s}}}\ \Pr{x,y}\In\Order\And\Pr{y,z}\In\OpGr{s}}\Period}{9}
 \Item{FDefs}{Functions} A graph \OpGr{f}\ is a \Def{function} if
                                                      \DIc{\ForAll{x\In\Domain{\OpGr{f}}}\ThereIsShriek{y\In\Range{\OpGr{f}}}\ \Pr{x,y}\In\OpGr{f}\Period}{10}
For a function \OpGr{f}\And\Pr{x,y}\In\OpGr{f}, we define the notation
                                                                                            \DIc{\OpGr{f}\Of{x} \Equiv\ y\Period}{11}
A function \OpGr{f}\Sin X\Cross Y is said to be \Def{surjective} if \Range{\OpGr{f}}\Equals Y and \Domain{\OpGr{f}}\Equals X\thinspace. A function is said to be \Def{injective} if its inverse is also a function.
A function is \Def{bijective} if it is both injective and surjective.
 \PaR If \Range{\OpGr{f}}\Equals\Domain{\OpGr{f}}\ for a bijective function, we say that \OpGr{f}\ is a \Def{permutation} \Def{of} \Def{its} \Def{domain}. A subset S of the domain of a permutation \OpGr{f}\ such that \OpGr{f}\Of{s}\In S for each s\In S is said to be \Def{invariant}. An invariant set containing no proper invariant subsets is said to be an \Def{orbit} \Def{of} \OpGr{f}. If \OpGr{f}\ and \OpGr{th}\ are permutations with a common domain, and if each orbit of \OpGr{th}\ is a subset of some orbit of \OpGr{f}, we shall say that \OpGr{th}\ is a \Def{refinement} \Def{of} \OpGr{f}.
 \PaR The notation
                                                                                           \DIc{\Function\OpGr{f}sendsinXtoinY\end}{12}
means that \OpGr{f}\ is a function such that \Domain{\OpGr{f}}\Equals X and \Range{\OpGr{f}}\Sin Y. This is sometimes described by saying that \OpGr{f}\ \Def{is} \Def{a} \Def{function} \Def{from} X \Def{into} Y (when \Range{\OpGr{f}}\Equals Y the term \Quotes{into} can be replaced by \Quotes{onto}). If E\Of{x}\ denotes an expression in terms of x, then
                                                                                      \DIc{\Function\OpGr{f}sendsxinXtoE\Of{x}inY\end}{13}
indicates in addition to \pLinkLocalDisplayText{12} that
                                                                               \DIc{\ForAll{x\In X}\quad\OpGr{f}\Of{x}\Equals E\Of{x}\Period}{14}
A statement of the form
                                                                                          \DIc{\Function sendsxinXtoE\Of{x}inY\end}{15}
defines a function without giving it a name.
  \PaR The \Def{identity} \Def{function} \Def{on} \Def{a} \Def{set} X is the function
                                                                           \DIc{\Id{X}\ \Equiv\ \SetSuch{\Pr{x,x}}{x\In X}\Period}{16}
 For any bijective function \OpGr{f}\Sin X\Cross Y, we have
                                                         \DIc{\Inv{\OpGr{f}}\Circ\OpGr{f}\ \Equals\ \Id{X}\Andd\OpGr{f}\Circ\Inv{\OpGr{f}}\ \Equals\ \Id{Y}\Period}{17}
 \PaR A function of which the range is a singleton is called a \Def{constant} \Def{function}: if X\Equals\Domain{\OpGr{f}}\ and \Set{c}\Equals\Range{\OpGr{f}}, we write
                                                                   \DIc{\Constant{X}{c}\quad for the function\quad \Function sendsxinXtocin\end\Period}{18}
 \PaR Associated with any graph \Order\Sin X\Cross Y\ are four \Quotes{set functions}:
                                  \DI{10}{5}{\Function\setRelation{\Order}sendsSin\Mii{P}\Of{X}to\SetSuch{y\In Y}{\ThereIs{x\In S}\ \Pr{x,y}\In\Order}in\Mii{P}\Of{Y}\end\Comma}{19}
                                   \DI{5}{5}{\Function\setRelationBack{\Order}sendsSin\Mii{P}\Of{Y}to\SetSuch{x\In X}{\ThereIs{y\In S}\ \Pr{x,y}\In\Order}in\Mii{P}\Of{X}\end\Comma}{20}
                                        \DI{5}{10}{\Function\RightPolar{\Order}sendsSin\Mii{P}\Of{X}to\SetSuch{y\In Y}{\ForAll{x\In S}\ \Pr{x,y}\In\Order}in\Mii{P}\Of{Y}\end}{21}
 and
                                        \DIc{\Function\LeftPolar{\Order}sendsSin\Mii{P}\Of{Y}to\SetSuch{x\In X}{\ForAll{y\In S}\ \Pr{x,y}\In\Order}in\Mii{P}\Of{X}\end\Period}{22}
We have in general
                           \DIc{\setRelation{\Order}\ \Equals\ \setRelationBack{\Inv{\Order}}\Comma\quad\setRelationBack{\Order}\ \Equals\ \setRelation{\Inv{\Order}}\Comma\quad
                                               \RightPolar{\Order}\ \Equals\ \RightPolar{\Inv{\Order}}\Andd\LeftPolar{\Order}\ \Equals\ \RightPolar{\Inv{\Order}}\Period}{23}
These two set functions are actually determined by their values at the singletons: to exploit this fact, we define
              \DIc{\Function\setRel{\Order}sendsxinXto\setRelation{\Order}\Of{\Sing{x}}\Equals\RightPolar{\Order}\Of{\Sing{x}}\Equals\SetSuch{y\In Y}{\Pr{x,y}\In\Order}in\Mii{P}\Of{Y}\end}{24}
and
              \DIc{\Function\setRelBack{\Order}sendsyinYto\setRelationBack{\Order}\Of{\Sing{y}}\Equals\LeftPolar{\Order}\Of{\Sing{y}}\Equals\SetSuch{x\In X}{\Pr{x,y}\In\Order}in\Mii{P}\Of{X}\end\Period}{25}
 \PaR We have
                                             \DIc{\ForAll{S\Sin X}\quad\setRelation{\Order}\Of{S}\Equals\Union{\RM{7}{x}\iN\RM{7}{S}}{\setRel{\Order}\Of{x}}\Comma\quad
                                                 \ForAll{S\Sin Y}\quad\setRelationBack{\Order}\Of{S}\Equals\Union{\RM{7}{y}\iN\RM{7}{S}}{\setRelBack{\Order}\Of{y}}\Comma}{26}
                                               \DIc{\ForAll{S\Sin X}\quad\RightPolar{\Order}\Of{S}\Equals\Intersection{\RM{7}{x}\iN\RM{7}{S}}{\setRel{\Order}\Of{x}}\Andd
                                               \ForAll{S\Sin Y}\quad\LeftPolar{\Order}\Of{S}\Equals\Intersection{\RM{7}{y}\iN\RM{7}{S}}{\setRelBack{\Order}\Of{y}}\Period}{27}
 \PaR If S is any subset of the domain of an order \Order, then the \Def{restriction} \Def{of} \Order\ to S is defined and denoted by
                                                                                            \DIc{\Order\Restriction{S}\Period}{28}
  \PaR For sets X and Y, the set of all functions from X into Y is sometimes denoted by
                                                                                      \DIc{\Power{X}{Y}}{30}
 \Par We shall extend this notation to
          \DIc{\Bijections{X}{Y}\ \Equiv\ \SetSuch{\OpGr{f}\In\Functions{X}{Y}}{\OpGr{f}\ is bijective}\Period}{30.5}
 \PaR Two functions \Function\OpGr{f}sendsinXtoinY\end\And\Function\OpGr{th}sendsinXtoinY\end\ will be said to be \Def{compatible} if
                                                   \DIc{\ForAll{\Set{x,y}\Sin X}\quad\OpGr{f}\Of{x}\Equals\OpGr{f}\Of{y}\Iff\OpGr{th}\Of{x}\Equals\OpGr{th}\Of{y}\Period}{31}
 Compatibility of functions \OpGr{f}\And\OpGr{th}\ can be expressed more succinctly as
                                                                                     \DIc{\OpGr{f}\Circ\Inv{\OpGr{th}}\ is a function.}{32}
 \PaR A graph \Order\ will be said to be \Def{symmetric} provided that there exists a bijection \Function\OpGr{f}sendsin\Domain{\Order}toin\Range{\Order}\end\ such that
              \DIc{\ForAll{\Pr{a,b}\In \Domain{\Order}\Cross\Range{\Order}}\quad\Pr{a,b}\In\Order\Quad\Iff\Quad\Pr{\Inv{\OpGr{f}}\Of{b},\OpGr{f}\Of{a}}\In\Order\Period}{33}
 \Item{Equivalence}{Equivalence} The concept of an equivalence relation was introduced in \pLinkItemText{SetTheory}{Relations}. When the symbol \Mii{A}\ there represents a set, the equivalence relation \Sim\ is a graph. Such an equivalence relation is distinguished by the following:
                                                           \DIc{\Id{\MII{5}{A}}\Sin\Sim ,\quad\Sim \Equals\Inv{\Sim }\Andd(\Sim \Circ\Sim)\Sin\Sim \Period}{1}
 The three parts of \pLinkLocalDisplayText{1} may also be expressed by
                         \Dc{\ForAll{x\In \Mii{A}}\ x\Sim x\ ,\quad\ForAll{\Set{x,y}\Sin\Mii{A}}\ x\Sim y\Iff y\Sim x\Andd\ForAllSuch{\Set{x,y,z}\Sin\Mii{A}}{x\Sim  y\And y\Sim z}\ \ x\Sim z\Period}
  \Item{Order}{Order} An \Def{order} \Def{on} \Def{a} \Def{set} X is a relation \Order\Sin X\Cross X such that
                                                                       \DIc{\Id{X}\Equals\Order\Cap\Inv{\Order}\Andd(\Order\Circ\Order)\Sin\Order\Period}{1}
 Another way of expression the two conditions of \pLinkLocalDisplayText{1} is
                    \Dc{\ForAll{x\In X}\ x\Order x\Comma\quad\ForAllSuch{\Set{x,y}\Sin X}{x\Order y\And y\Order x}\ x\Equals y\Andd\ForAllSuch{\Set{x,y,z}\Sin X}{x\Order y\And y\Order z}\ x\Order z\Period}
An order is sometimes called a \Def{partial} \Def{ordering} (or simply an \Def{ordering}) and, relative to \Order, X is called an \Def{ordered} \Def{set}, or \Def{partially} \Def{ordered} \Def{set} or \Def{poset}.
 \PaR An order \Order\ is said to be \Def{linear} or \Def{total} if
                                                                            \DIc{X\Cross X\Equals(\Order\Cup\Inv{\Order})}{2}
 or, what amounts to the same thing, if
                                                                      \DIc{\ForAll{\Pr{x,y}\In X\Cross X}\quad either \ x\Order y\ or y\Order x\Period}{3}
 A subset \OpGr{s}\ of an order \Order\ such that \Domain{\OpGr{s}}\Cross\Domain{\OpGr{s}}\Equals\OpGr{s}\Cup\Inv{\OpGr{s}}\ which is called a \Def{chain} in \Order. The \Def{Hausdorff} \Def{Maximality} \Def{Principle}, which is equivalent to the Axiom of Choice, asserts that each chain is a subset of some maximal chain.
 \PaR The \Def{supremum} of a subset S of \Domain{\Order}, if it exists, is an element \Ft{sup}\Of{S} such that
                                              \DIc{\ForAll{s\In S}\ s\Order\Ft{sup}\Of{S}\Andd\ForAllSuch{x\In S}{\ForAll{s\In S}\ s\Order x}\quad\Ft{sup}\Of{S}\Order x\Period}{4}
 The \Def{infimum} of a subset S of \Domain{\Order}, if it exists, is an element \Ft{inf}\Of{S} such that
                                              \DIc{\ForAll{s\In S}\ \Ft{inf}\Of{S}\Order s\Andd\ForAllSuch{x\In S}{\ForAll{s\In S}\ x\Order s}\quad x\Order\Ft{inf}\Of{S}\Period}{5}
 For x and y in an ordered set, we define (if they exist)
                                                             \DIc{x\Wedge y\ \Equiv\ \Ft{inf}\Of{\Set{x,y}}\Andd x\Vee y\ \Equiv\ \Ft{sup}\Of{\Set{x,y}}\Period}{6}
 \Item{Cardinality}{Cardinality} If \OpGr{f}\ is any bijective function, then \Domain{\OpGr{f}}\ and \Range{\OpGr{f}}\ are said \Def{to} \Def{have} \Def{the} \Def{same} \Def{cardinality}. Because each identity function is bijective,
 each set has the same cardinality as itself. Furthermore, if sets X and Y have the same cardinality and Y and Z have the same cardinality, then X and Z have the same cardinality. Thus it is possible, for any set X, to attach a
 symbol to it, and by extension, attach this same symbol to every other set with the same cardinality as X. There are many equivalent ways to do this for common sets. Here are some typical examples:
 \PaR We attach the symbol \Def{\Zero} to the empty set \Void\ and the symbol \Def{\One} to the singleton \Set{\Void}. Attach the symbol \Def{\Two} to the ordered pair \Pr{\One,\Void}, \Def{\Three} to the ordered pair \Pr{\Two,\Void}, \Def{\Four} to the
 the ordered pair \Pr{\Three,\Void} and so on \It{ad}\ \It{infinitum}. The set
                                                                                                        \DIc{\Aii{N}}{1}
 of all such symbols distinct from 0 is called the set of \Def{natural} \Def{numbers}. We attach the symbol $\aleph_0$ to the set \Aii{N}. Any set which (by extension) has one of these symbols attached to it is said to be
 \Def{countable}. If such symbol is attached to some set, this symbol is said to be its \Def{cardinality}. These symbols are called \Def{cardinal} \Def{numbers}. If the cardinality of a non-void set is an element of \Aii{N}, it is
 said to be \Def{finite}. This definition of \Quotes{finite} is consistent with the one given in \pLinkItemText{SetTheory}{Finite}.
 \PaR A set of which the cardinality is $\aleph_0$ is said to be \Def{countably} \Def{infinite}. We denote the cardinality of any set X by
                                                                                                       \DIc{\Card{X}}{2}
 and, for any natural number n, we denote by
                                                                                                     \DIc{\Underbar{n}}{3}
 the set of natural numbers less than or equal to n. Any function with domain of the form \Underbar{n}\ is said to be a \Def{finite} \Def{sequence}. An \Def{infinite} \Def{sequence} is a function with domain \Aii{N}. For a sequence
 \OpGr{s}\ and m\In\Domain{\OpGr{s}}, the m-th \Def{element} \Def{of} \OpGr{s}\ is defined and denoted by
                                                                                      \DIc{\OpGr{s}$_m$\ \Equiv\ \OpGr{s}\Of{m}\Period}{4}
 \PaR For any n\In\Aii{N}\ and any set X, the elements of
                                                                            \DIc{\SEQUENCE{X}{n}}{5}
 are called \Def{n-tuples}: the sequences from \Underbar{n}\ into X. Such an n-tuple x is frequently written
                                                                                 \DIc{\Nr{x$_{\RMBx{5}{1}}$,\dots,x$_{\RMBx{7}{n}}$}\Period}{6}
 \PaR Beginning with 1, every other finite cardinal number is said to be \Def{odd}. The other finite cardinal numbers are \Def{even}. A set is said to be \Def{even} if its cardinality is even -- \Def{odd} if its cardinality is odd. Two finite sets \Def{have} \Def{equal} \Def{parity} if they are both even, or both odd.\vfill\eject
                                                                                 \Section{Topology}{Appendix III: Topology}\xrdef{pageTopology}
 \Item{D}{Definitions} A \Def{topology} is a family \MMiiBx{T}\ of subsets of a set X containing X, the empty set \Void, the union of the members of any subfamily of \MMiiBx{T}\ and the intersection of the members of any finite subfamily of \MMiiBx{T}\Period The members of a topology are said to be \Def{open} subsets of X and complements of open sets are said to be \Def{closed} subsets of X. A subset K of X is said to be \Def{compact} if, whenever K is a subset of the union of a subfamily of \MMiiBx{T}\Comma\ then it is also subset of some finite subfamily of that \Quotes{\Def{covering}} subfamily. A subset of X is \Def{disconnected} if it is the union of two disjoint elements of \MMiiBx{T}\Comma\ and \Def{connected} if it is not disconnected.
 \Par A subset N of X is said to be a \Def{neighborhood} of one of its elements x if there exists O\In\MMiiBx{T}\ such that x\In O and O\Sin N. A topology is said to be a \Def{hausdorff} topology if each, for all distinct x and y in X, there exist disjoint open sets U and V such that x\In U and y\In V.
 \Item{RT}{Relative Topology} If \MMiiBx{T}\ is a topology on a set X and if S is a subset of X, then the family \SetSuch{O\Cap S}{O\In\MMiiBx{T}}\ is a topology on S called the \Def{relativized} \Def{topology}. The set S is compact for its relativized topology if, and only if, it is compact as a subset of X. The same is true for the concepts \Quotes{connected} and \Quotes{hausdorff}.
 \Item{BASE}{Base for a Topology} Let X be a set and \Mii{B}\ a family of subsets of X such that the intersection of any two elements of \Mii{B}\ is an element of \Mii{B}. Suppose that the null set is in \Mii{B}\ and that each element of x is in some element of \Mii{B}\Period\ Let \MMiiBx{T}\ be the family of all subsets S\Sin X such that S is a union of elements of \Mii{B}\Period\ Then \MMiiBx{T}\ is a topology for X. We say that \Mii{B}\ is a \Def{base} \Def{for} \Def{the} \Def{topology} \MMiiBx{T}\Period
 \Item{Nets}{Nets} Let \MMiiBx{T}\ be a topology on a set X. An order \Order\ on a set D is said to be a \Def{direction} \Def{on} D provided
                                                                     \DIc{\ForAll{\Set{x,y}\Sin D}\ThereIs{z\In D}\quad x\Order z\And y\Order z\Period}{1}
 \PaR Any function with domain a directed set is said to be a \Def{net}. A net \OpGr{w}\ from one directed set D with direction \Order\ into another directed set E with direction \ORder\ is said to be \Def{co-final} if
                                                            \DIc{\ForAll{e\In E}\ThereIs{d\In D}\ForAllSuch{x\In D}{d\Order x}\quad e\ORder\OpGr{w}\Of{x}\Period}{2}
 \Par If D and E are as above, \OpGr{w}\ is a co-final net in E and if \OpGr{n}\ is a net with domain E, then \OpGr{n}\Circ\OpGr{w}\ is said to be a \Def{subnet} of \OpGr{n}. A net \OpGr{v}\ defined on a directed set D is said to be \Def{eventually} \Def{in} \Def{a} \Def{set} S if
                                                                        \DIc{\ThereIs{d\In D}\ForAllSuch{x\In D}{d\Order x}\quad \OpGr{n}\Of{x}\In S}{3}
and to be \Def{frequently} \Def{in} \Def{a} \Def{set} S if
                                                                     \DIc{\ForAll{d\In D}\ThereIs{x\In D}\quad d\Order x\Andd\OpGr{n}\Of{x}\In S\Period}{4}
 \PaR A net \OpGr{n}\ is said to \Def{converge} to an element x of X if
                                                  \DIc{\ForAllSuch{O\In\MMiiBx{T}}{x\In O}\quad\OpGr{n}\ is eventually in O\Period}{5}
 \Par The point x is said to be the \Def{limit} of the net \OpGr{n}\Period
 \PaR The proofs of the following facts are pedestrian:
  \DI{8}{0}{\ForAll{S\Sin X}\quad S is open\quad\Iff\quad each net convergent to a point in S is eventually in S,}{6}
  \DI{8}{0}{\ForAll{S\Sin X}\quad S is closed\quad\Iff\quad the limit of any convergent net in S is in S,}{7}
  \DI{8}{8}{\ForAll{S\Sin X}\quad S is compact\quad\Iff\quad each net in S has a subnet convergent to an element of S.}{8}
 \Item{Con}{Continuity} Let \Mii{T}$_1$ be a topology on a set X$_1$ and \Mii{T}$_2$ be a topology on a set X$_2$. A function \Function\Op{f}sendsin X$_1$toin X$_2$\end\ is said to be \Def{continuous} provided
                                                                     \DIc{\ForAll{O\In\Mii{T}$_2$}\quad\setRelationBack{\Op{f}}\Of{O}\In\Mii{T}$_1$}{1}
 \Par or, equivalently
                                                  \DIc{\ForAll{\OpGr{w}\ convergent to w\In X$_1$}\quad\Op{f}\Circ\OpGr{w}\ converges to \Op{f}\Of{w}\ in X$_2$\Period}{2}
 \PaR The proofs of the following facts are pedestrian:
                                                                        \DI{8}{0}{the composition of any two continuous functions is continuous,}{3}
                          \DI{8}{0}{\ForAllSuch{S\Sin X$_1$}{S is compact}\ForAll{\Function\Op{f}sendsin X$_1$toin X$_2$\end\ continuous}\quad\setRelation{\Op{f}}\Of{S}\ is compact in X$_2$}{4}
                   \DII{8}{8}{\ForAllSuch{S\Sin X$_1$}{S is connected}\ForAll{\Function\Op{f}sendsin X$_1$toin X$_2$\end\ connected}\quad\setRelation{\Op{f}}\Of{S}\ is connected in X$_2$\Comma}{and}{5}
 \Item{Metric}{Metrics} A \Def{metric} \Op{d}\ is a non-negative function on the cartesian product of a set X with itself such that\Foot{The value d\Of{\Pr{x,y}}\ of d at an ordered pair \Pr{x,y}\ is usually abbreviated to d\Of{x,y}.}
                                                          \DI{8}{0}{\ForAll{\Set{x,y}\Sin X}\quad d\Of{\Pr{x,y}}\Equals\Zero\Quad\Iff\Quad x\Equals y\Comma}{1}
                                                                          \DI{9}{0}{\ForAll{\Set{x,y}\Sin X}\quad d\Of{\Pr{x,y}}\Equals d\Of{\Pr{y,x}}}{2}
                                                               \DII{8}{8}{\ForAll{\Set{x,y,z}\Sin X}\quad d\Of{\Pr{x,y}}+d\Of{\Pr{y,z}}$\le$d\Of{\Pr{x,z}}\Period}{and}{3}
 \PaR The \Def{ball} \Def{of} \Def{radius} r \Def{around} \Def{a} \Def{point} x in X, for r$\ge\Zero$, is the set
                                                                                   \DIc{\SetSuch{t\In X}{d\Of{\Pr{x,t}}$\le r$}\Period}{4}
 \Par The family of all such balls forms a base for a \Def{metric} \Def{topology}.
 \PaR A metric topology is Hausdorff. A subset of a metric space is compact if, and only if, each sequence has a convergent subsequence.
\vfill\eject
                                                                                                        \Section{Notation}{Notation}\xrdef{pageNotation}
 \PAr\baselineskip=12pt
\setbox\BoxI\hbox{\vbox{\halign{#\dotfill&\dots\dotfill\LinkDisplayText#\cr
         \Ft{Homograph}\Of{X}&{P}{FT1}{1}\cr
         V$^{(\EUBx{7}{n})}$&{P}{FT4}{3}\cr
         V$_{(\EUBx{7}{n})}$&{P}{FT4}{4}\cr
         \hbox{\OpGr{l}$_{\EUBx{7}{n}}$\Of{x}}&{P}{FT4}{5}\cr
         \hbox{\Op{\Gr{n}$_{\EUBx{7}{n}}$}}&{P}{FT4}{6}\cr
         \hbox{\OpGr{a}$^{(\EUBx{7}{n})}$}&{P}{FT4}{8}\cr
         A$^{oo}$&{P}{FT6}{2}\cr
         \LLL{\OpGr{a},\OpGr{b},\OpGr{c}}&{ED}{L}{1}\cr
         \LLL{\Mii{F}}&{ED}{L}{2}\cr
         \vruleh{11}\Vol{a,b;c,d}&{ED}{MD}{5}\cr
         \Vol{a,b,c,d}&{ED}{EI}{3}\cr
         \Aii{P}&{ED}{Circle}{1}\cr
         \hbox{\Op{p\Lower{4pt}{\RMBx{7}{C}}}}\Of{x}&{ED}{Circle}{2}\cr
         \Line{x,y}&{ED}{Circle}{2}\cr
         L\Wedge M&{ED}{LineM}{1}\cr
         \Infinity{L}&{ED}{SphereM}{1}\cr
         \Vol{S;L}&{ED}{SphereM}{5}\cr
         \Hol{a}{b}{c}{u}{v}{w}&{E}{D}{2}\cr
         \Mii{G}$_{\RMBx{7}{n}}$&{E}{D}{4}\cr
         \LL{a,b,c}&{LL}{LD}{1}\cr
         \LL{a$_1$,a$_2$,\dots,a$_n$}&{LL}{LD}{2}\cr
         \LibraP{a}{b}&{LL}{D4}{2}\cr
         \LibraR{a}{b}&{LL}{D4}{2}\cr
         \LibraL{a}{b}&{LL}{D4}{2}\cr
         \Ru{K}&{CL}{BD}{1}\cr
         \Es{F}&{CL}{BD}{2}\cr
         \GreenStar, \Spade, \Club, \RedHeart, \RedDiamond&{CL}{BD}{4}\cr
         \hbox{\GIn, \BIn, \RIn}&{CL}{BD}{1}\cr
         \WideHat{\Eb{F}}&{CL}{BD}{1}\cr
         \hbox{\WWurf{a,b,c,d}}&{CL}{BD}{2}\cr
         \Quin{a}{b}{c}{d}{e}&{CL}{PMD}{3.5}\cr
         \hbox{\PB{Q}, \OB{Q}}&{W}{C}{1}\cr
         \hbox{\PB{a}, \PB{b}, \PB{c} ,\PB{d}}&{W}{C}{1}\cr
         \hbox{\OB{a}, \OB{b}, \OB{c}, \OB{d}}&{W}{C}{1}\cr
         \hbox{\GInv, \BInv, \RInv}&{W}{C}{5}\cr
         \hbox{\GadTurn}&{W}{C}{6}\cr
         \hbox{\ABInv}&{W}{C}{8}\cr
         \hbox{\AVInv}&{W}{C}{10}\cr
         \hbox{\AVbcInv}&{W}{C}{12}\cr
         \Ft{K}&{W}{C}{14}\cr
         \Ft{K}$_4$&{W}{C}{14}\cr
         \hbox{\GSquare, \BSquare, \RSquare}&{W}{W}{2}\cr
         \hbox{\wurf{\Rb{b}}{\Rb{d}}{\Rb{l}}{\Rb{q}}{a}{b}{c}{d}}&{W}{W}{8}\cr
         \hbox{\Wurf{a,b,c,d}{W}{W}}&{W}{W}{9}\cr
         \hbox{\Wurf{a,b,c,d}{W}{W}$^{\sIm}$}&{W}{W}{10}\cr
         \Wuerf{x}{y}{z}&{W}{DW}{2}\cr
         \Ft{G}&{W}{CI}{1}\cr
         \hbox{\GSquare,\BSquare,\RSquare}&{W}{HW}{2}\cr}}}
\setbox\BoxII\hbox{\vbox{\halign{#\dotfill&\dots\dotfill\LinkDisplayText#\cr
         \MAB{a}{d}&{EA}{II}{1}\cr
         \Mab{a}{d}&{EA}{II}{2}\cr
         \MABF{a}{d}&{EA}{II}{3}\cr
         \Mabf{a}{l}{d}&{EA}{II}{6}\cr
         \hbox{\BZero{b},\BOne{b},\BInfty{b}}&{EA}{N}{2}\cr
         \hbox{\BPlus{b},\BMinus{b},\BTimes{b},\BMinus{}}&{EA}{N}{3}\cr
         \BField&{EA}{N}{4}\cr
         \BDet{a}{b}{c}{d}{b}&{EA}{Defff}{1}\cr
         \MNF&{EA}{Defff}{2}\cr
         \Mii{G}$^+$&{EA}{D5}{0}\cr
         \OSim&{EA}{D5}{1}\cr
         \hbox{\COreo,\CCOreo}&{EA}{D5}{2}\cr
         \Arc{a,b}&{EA}{D5}{3}\cr
         \CArc{a,b}&{EA}{D5}{4}\cr
         \MiiChBx{71}&{SetTheory}{Logic}{i}\cr
         \MiiChBx{71}\kn{1}!&{SetTheory}{Logic}{ii}\cr
         \MiiChBx{70}&{SetTheory}{Logic}{iii}\cr
         \Implies&{SetTheory}{Logic}{iv}\cr
         \Implied&{SetTheory}{Logic}{v}\cr
         \Iff&{SetTheory}{Logic}{vi}\cr
         x\In X, X\Ni x&{SetTheory}{Sets}{1}\cr
         x\Nin X&{SetTheory}{Sets}{2}\cr
         \Void&{SetTheory}{Sets}{5}\cr
         \Set{x,\dots,y}&{SetTheory}{Sets}{7}\cr
         \SetSuch{x\In X}{P\Of{x}}&{SetTheory}{Sets}{8}\cr
         \Equiv&{SetTheory}{Sets}{vii}\cr
         X\Cup Y&{SetTheory}{Sets}{9}\cr
         X\Cap Y&{SetTheory}{Sets}{10}\cr
         \Subsets{X}&{SetTheory}{Sets}{11}\cr
         X\Cop Y&{SetTheory}{Sets}{13}\cr
         \Union{\RM{7}{S}\iN\MII{7}{S}}{S}&{SetTheory}{Sets}{17}\cr
         \Intersection{\RM{7}{S}\iN\MII{7}{S}}{S}&{SetTheory}{Sets}{17}\cr
         \Pr{x,y}&{SetTheory}{CarProducts}{2}\cr
         X\Cross Y&{SetTheory}{CarProducts}{3}\cr
         \Sim&{SetTheory}{Relations}{2}\cr
         \hbox{x\NOrder y}&{Graphs}{Definitions}{1}\cr
         \hbox{x\OrderShriek y}&{Graphs}{Definitions}{2}\cr
         \hbox{\Domain{\Order}}&{Graphs}{Definitions}{3}\cr
         \hbox{\Range{\Order}}&{Graphs}{Definitions}{4}\cr
         \hbox{\World{\Order}}&{Graphs}{Definitions}{5}\cr
         \hbox{\Inv{\Order}}&{Graphs}{Definitions}{6}\cr
         \hbox{\Cop\Order}&{Graphs}{Definitions}{8}\cr
         \hbox{\Order\Circ\OpGr{s}}&{Graphs}{Definitions}{9}\cr
         \hbox{\Function\OpGr{f}sendsinXtoinY\end}&{Graphs}{FDefs}{12}\cr
         \hbox{\Function\OpGr{f}sendsxinXtoE\Of{x}inY\end}&{Graphs}{FDefs}{13}\cr
         \hbox{\Function sendsxinXtoE\Of{x}inY\end}&{Graphs}{FDefs}{15}\cr
         \hbox{\Id{X}}&{Graphs}{FDefs}{16}\cr}}}
\setbox\BoxIII\hbox{\vbox{\halign{#\dotfill&\dots\dotfill\LinkDisplayText#\cr
         \hbox{\Constant{X}{c}}&{Graphs}{FDefs}{18}\cr
         \hbox{\setRelation{\Order}}&{Graphs}{FDefs}{19}\cr
         \hbox{\setRelationBack{\Order}}&{Graphs}{FDefs}{20}\cr
         \hbox{\RightPolar{\Order}}&{Graphs}{FDefs}{21}\cr
         \hbox{\LeftPolar{\Order}}&{Graphs}{FDefs}{22}\cr
         \hbox{\Order\Restriction{S}}&{Graphs}{FDefs}{28}\cr
         \hbox{X\kn{.5}{\Vbox{\RMBx{7}{Y}\kn{5}}}}&{Graphs}{FDefs}{30}\cr
         \Ft{sup}\Of{S}&{Graphs}{Order}{4}\cr
         \Ft{inf}\Of{S}&{Graphs}{Order}{5}\cr
         \Aii{N}&{Graphs}{Cardinality}{1}\cr
         \Card{X}&{Graphs}{Cardinality}{2}\cr
         \Underbar{n}&{Graphs}{Cardinality}{3}\cr
         \Nr{x$_{\RMBx{5}{1}}$,\dots,x$_{\RMBx{7}{n}}$}&{Graphs}{Cardinality}{6}\cr}}}
         \Dimen=\hsize\advance\Dimen by-10pt%
         \advance\Dimen by-\wd\BoxI\advance\Dimen by-\wd\BoxII\advance\Dimen by-\wd\BoxIII\divide\Dimen by2%
         \CountI=1\CountII=1\CountIII=1\DimenI=\ht\BoxI\DimenII=\ht\BoxII\DimenIII=\ht\BoxIII%
         \ifdim\DimenI>\DimenII\advance\CountI by1\else\advance\CountII by1\fi%
         \ifdim\DimenI>\DimenIII\advance\CountI by1\else\advance\CountIII by1\fi%
         \ifdim\DimenII>\DimenIII\advance\CountII by 1\else\advance\CountIII by1\fi%
         \ifnum\CountI=3\DimenO=\DimenI\else\ifnum\CountII=3\DimenO=\DimenII\else\ifnum\CountIII=3\DimenO=\DimenIII\fi\fi\fi%
         \advance\DimenI by-\DimenO\advance\DimenII by-\DimenO\advance\DimenIII by-\DimenO%
         \hbox{\Vbox{\box\BoxI\kern-\DimenI}\kern\Dimen\Vbox{\box\BoxII\kern-\DimenII}\kern\Dimen\Vbox{\box\BoxIII\kern-\DimenIII}}
\vfill\eject
                                                                                                     \Section{Index}{Index}\xrdef{pageIndex}
\PAr\vskip-0in
\setbox\BoxI\hbox{\vbox{\halign{#\dotfill&\dots\dotfill\LinkItemText#\cr
        \Two-exponential isomorphism&{EA}{DI}\cr
        \Three-transitive Group&{E}{I}\cr
        Abel, Henrik&{P}{H}\cr
        Abelian&{ED}{EI}\cr
        Adjacent&{ED}{LineM}, \LinkItemText{CL}{BD}, \LinkItemText{W}{C}\cr
        Affine Space&{LL}{E1}\cr
        Alberti, Leone Batista&{P}{H}\cr
        Arc&{EA}{D5}\cr
        Arc Topology&{EA}{D13}\cr
        Axiom of Choice&{SetTheory}{PropertiesAndClasses}\cr
        Axis&{W}{C}\cr
        Balanced Triple of Pairs&{ED}{LineM}\cr
        Base for a Topology&{Topology}{BASE}\cr
        Basis&{P}{FT6}, \LinkItemText{E}{D}\cr
        Bijective&{Graphs}{FDefs}\cr
        Bounded Interval&{Graphs}{Order}\cr
        Brunelleshi, Filipo&{P}{H}\cr
        Cardinality&{Graphs}{Cardinality}\cr
        Cartesian Product&{SetTheory}{CarProducts}\cr
        Chain&{Graphs}{Order}\cr
        Chevalier, Auguste&{P}{H}\cr
        Choice Function&{SetTheory}{PropertiesAndClasses}\cr
        Circle Meridian&{ED}{Circle}\cr
        Class&{SetTheory}{PropertiesAndClasses}\cr
        Closed&{Topology}{D}\cr
        Co-final&{Topology}{Nets}\cr
        Collection&{SetTheory}{Sets}\cr
        Compact&{Topology}{D}\cr
        Compatible Functions&{Graphs}{FDefs}\cr
        Complement&{SetTheory}{Sets}, \LinkItemText{Graphs}{Definitions}\cr
        Complete Quadrilateral&{ED}{LineM}\cr
        Composition&{Graphs}{Definitions}\cr
        Condition on a Set&{SetTheory}{PropertiesAndClasses}\cr
        Connected&{Topology}{D}\cr
        Continuous&{Topology}{Con}\cr
        Convergent&{Topology}{Nets}\cr
        Constant Function&{Graphs}{FDefs}\cr
        Countable&{Graphs}{Cardinality}\cr
        Countably Infinite&{Graphs}{Cardinality}\cr
        Cross Ration&{W}{I}\cr
        Cube&{W}{C}\cr
        Da Vinci, Leonardo&{P}{H}\cr
        De Morgan Laws&{SetTheory}{Sets}\cr
        Dedekind, Richard&{P}{H}\cr
        Desargue, G\'erard&{P}{H}\cr
        Dilation&{E}{Trans}, \LinkItemText{E}{Trans}\cr
        Dimension&{P}{FT6}\cr
         Direction&{Topology}{Nets}\cr
        Distributive Laws&{SetTheory}{Sets}\cr
        Domain&{Graphs}{Definitions}\cr}}}
\setbox\BoxII\hbox{\vbox{\halign{#\dotfill&\dots\dotfill\LinkItemText#\cr
        Doubleton&{SetTheory}{Sets}\cr
        D\"urer, Albrecht&{P}{H}\cr
        Empty Set&{SetTheory}{Sets}\cr
        Even&{Graphs}{Cardinality}\cr
        Erlanger Programm&{ED}{EP}\cr
        Equivalence&{Graphs}{Definitions}\cr
        Equivalence Class&{SetTheory}{Relations}\cr
        Equivalence Relation&{SetTheory}{Relations}\cr
        Euclidean Plane&{ED}{Circle}\cr
        Existential Quantifier&{SetTheory}{Logic}\cr
        Exponential Meridian&{EA}{DI}\cr
        Family&{SetTheory}{Sets}\cr
        Field Meridian&{ED}{EI}\cr
        Finite&{SetTheory}{Finite}, \LinkItemText{Graphs}{Cardinality}\cr
        First Form of Fundamental Theorem&{P}{FT5}\cr
        Function&{Graphs}{FDefs}\cr
        Function Libra&{ED}{L}\cr
        Fundamental Theorem&{P}{FT6}, \LinkItemText{ED}{FT}\cr
        Galois, \"Evariste&{P}{H}\cr
        Gaussian Plane&{ED}{SphereM}\cr
        Graph&{Graphs}{Definitions}\cr
        Group&{ED}{EI}\cr
        Group Libra Operator&{LL}{D2}\cr
        Harmonic Pair&{ED}{Har}\cr
        Harmonic Wurf&{W}{W}\cr
        Hausdorff&{Topology}{D}\cr
        Hausdorff Maximality Principle&{Graphs}{Order}\cr
        Homography&{P}{FT1}, \LinkItemText{ED}{MD}, \LinkItemText{ED}{EI}\cr
        Identity&{ED}{EI}\cr
        Identity Function&{Graphs}{FDefs}\cr
        Independent Set&{P}{FT6}\cr
        Infimum&{Graphs}{Order}\cr
        Injective&{Graphs}{FDefs}\cr
        Inner Involution Libra&{LL}{DV}\cr
        Interval&{Graphs}{Order}\cr
        Invariant&{Graphs}{FDefs}\cr
        Intersection&{SetTheory}{Sets}\cr
        Inverse&{ED}{EI}, \LinkItemText{Graphs}{Definitions}\cr
        Involution&{E}{Trans}\cr
        Isomorphism of Graphs&{Graphs}{Definitions}\cr
        Issue (of a \Two-exponential)&{EA}{DI}\cr
        Kepler, Johannes&{P}{H}\cr
        Klein, Felix&{ED}{EP}\cr
        Klein 4-group&{W}{C}\cr
        Libra&{LL}{LD}\cr
        Libra Generated by \Mii{F}&{ED}{L}\cr
        Libra Homomorphism&{LL}{D3}\cr
        Libra Isomorphism&{LL}{D3}\cr
        Libra Left Translation&{LL}{D4}\cr
        Libra Operator&{LL}{LD}\cr}}}
\setbox\BoxIII\hbox{\vbox{\halign{#\dotfill&\dots\dotfill\LinkItemText#\cr
        Libra Right Translation&{LL}{D4}\cr
        Limit&{Topology}{Nets}\cr
        Line at Infinity&{ED}{Circle}\cr
        Line Dual to a Line&{ED}{SphereM}\cr
        Line Meridian&{ED}{LineM}\cr
        Linear Ordering&{Graphs}{Order}\cr
        Load&{CL}{BD}\cr
        Member&{SetTheory}{Sets}, \LinkItemText{SetTheory}{PropertiesAndClasses}\cr
        Meridian Basis&{ED}{EI}\cr
        Meridian Cycle&{ED}{HI}\cr
        Meridian Family of Involutions&{ED}{MD}\cr
        Meridian family of loads&{CL}{PMD}\cr
        Meridian Field&{ED}{EI}\cr
        Meridian Group of Permutations&{E}{MD}\cr
        Meridian Libra Isomorphism&{EA}{DI}\cr
        Meridian Isomorphism&{ED}{MD}\cr
        Meridian Orbit&{ED}{HI}\cr
        Meridian Quinary Operator&{CL}{PMD}\cr
        Monge, Gaspard&{P}{H}\cr
        Moore, Eliakim Hastings&{P}{H}\cr
        n-tuple&{Graphs}{Cardinality}\cr
        Natural Numbers&{Graphs}{Cardinality}\cr
        Net&{Topology}{Nets}\cr
        Odd&{Graphs}{Cardinality}\cr
        Open&{Topology}{D}\cr
        Open Interval&{Graphs}{Order}\cr
        Opposite&{ED}{LineM}, \LinkItemText{CL}{BD}, \LinkItemText{W}{C}\cr
        Orbit&{Graphs}{FDefs}\cr
        Order&{Graphs}{Order}\cr
        Ordering&{Graphs}{Order}\cr
        Ordered Pairs&{SetTheory}{CarProducts}\cr
        Oriented Circular Meridian&{EA}{D5}\cr
        Orientation&{EA}{D5}\cr
        Originating at&{EA}{DI}\cr
        Pair&{SetTheory}{CarProducts}\cr
        Pappus of Alexandria&{P}{H}\cr
        Parity&{Graphs}{Cardinality}\cr
        Partial Ordering&{Graphs}{Order}\cr
        Partition&{SetTheory}{Relations}\cr
        Pascal, Blaise&{P}{H}\cr
        Peano, Giuseppe&{P}{FT1}\cr
        Permutation&{Graphs}{FDefs}\cr
        Perspectivity&{P}{FT7}\cr
        Piero della Francesca&{P}{H}\cr
        Plane at Infinity&{ED}{SphereM}\cr
        Poset&{Graphs}{Order}\cr
        Pre-meridian family of loads&{CL}{PMD}\cr
        Projective 3-Space&{ED}{SphereM}\cr
        Projective Automorphism&{P}{FT3}\cr
        Projective Isomorphism&{P}{FT3}\cr
        Projective Plane&{ED}{Circle}\cr
        Projective Space&{P}{FT1}\cr
        Proper Interval&{Graphs}{Order}\cr}}}
\setbox\BoxIV\hbox{\vbox{\halign{#\dotfill&\dots\dotfill\LinkItemText#\cr
        Proper Subset&{SetTheory}{Finite}\cr
        Poncelet,Jean-Victor&{P}{H}\cr
        Progenitor (of a \Two-exponential)&{EA}{DI}\cr
        Projectivity&{P}{FT7}\cr
        Quadriad&{W}{C}\cr
        Range&{Graphs}{Definitions}\cr
        Refinement&{Graphs}{FDefs}\cr
        Regular&{CL}{BD}, \LinkItemText{W}{W}\cr
        Relative Topology&{Topology}{RT}\cr
        Relation&{SetTheory}{Relations}\cr
        Representation Field&{P}{FT4}\cr
        Restriction&{Graphs}{FDefs}\cr
        Rieman Sphere&{ED}{SphereM}\cr
        Rotation&{E}{Trans}, \LinkItemText{E}{Trans}\cr
        Separable Order Interval&{Graphs}{Order}\cr
        Sequence&{Graphs}{Cardinality}\cr
        Set&{SetTheory}{Sets}\cr
        Simplex&{P}{FT6}\cr
        Singleton&{SetTheory}{Sets}\cr
        Singular&{CL}{BD}, \LinkItemText{W}{W}\cr
        Sphere Meridian&{ED}{SphereM}\cr
        Steiner, Jakob&{P}{H}\cr
        Stereographic Projection&{P}{H}\cr
        Subset&{SetTheory}{Sets}\cr
        Superset&{SetTheory}{Sets}\cr
        Supremum&{Graphs}{Order}\cr
        Surjective&{Graphs}{FDefs}\cr
        Symmetric Difference&{SetTheory}{Sets}\cr
        Symmetric Graph&{Graphs}{FDefs}\cr
        Tits, Jacques&{ED}{MC}\cr
        Topology&{Topology}{D}\cr
        Translation&{E}{Trans}, \LinkItemText{E}{Trans}, \LinkItemText{LL}{E1}\cr
        Triple&{SetTheory}{CarProducts}\cr
        Union&{SetTheory}{Sets}\cr
        Universal Quantifier&{SetTheory}{Logic}\cr
        Vector Representation&{P}{FT1}\cr
        Void&{SetTheory}{Sets}\cr
        Von Staudt, Karl Georg Christian&{P}{H}\cr
        Weyl, Hermann&{P}{H}\cr
        World&{Graphs}{Definitions}\cr
        Wurf&{P}{H}\cr}}}
 \MaxIV{\wd\BoxI}{\wd\BoxII}{\wd\BoxII}{\wd\BoxIV}%
 \DimenI=\DimenReturn\advance\DimenI by-\wd\BoxI%
 \DimenII=\DimenReturn\advance\DimenI by-\wd\BoxII%
 \DimenIII=\DimenReturn\advance\DimenI by-\wd\BoxIII%
 \DimenIV=\DimenReturn\advance\DimenI by-\wd\BoxIV%
 \Dimen=\hsize\advance\Dimen by-2\DimenReturn%
 \DimenArg=\ht\BoxIII\advance\DimenArg by-\ht\BoxIV%
 \setbox\Box\vbox{\box\BoxIV\kern\DimenArg}%
 \hbox{\box\BoxI\kern\Dimen\box\BoxII}%
 \vfill
 \eject%
 \hbox{\box\BoxIII\kern\Dimen\box\Box}%
 \vfill\eject%
 \baselineskip=14pt
                                                                                           \Section{Eponymy}{Eponymy}\xrdef{pageEponymy}
\PAR\It{Abelian} (group and libra): Henrik Abel was born in Nedstrand, Norway in 1802 and died in 1829.
\PAR\It{Euclidean} (plane and space): Euclid of Alexandria lived in Greece and Egypt \It{circa} 300 BC.
\PAR\It{Gaussian} (plane): Johann Carl Friedrich Gauss was born in Braunschweig, Germany in 1777 and died in 1855.
\PAR\It{Hausdorff} (topological space and maximality principle): Felix Hausdorff was born in Breslau, Germany in 1868 and died in 1942.
\PAR\It{Klein} (4-group): Christian Felix Klein was born in D\"usseldorf, Prussia in 1849 and died in 1925.
\PAR\It{Riemann} (sphere): Bernhard Riemann was born in Brelenz, Hanover in 1826 and died in 1866.

\vfill\eject%
                                                                                           \Section{Acknowledgement}{Acknowledgement}\xrdef{pageAcknowledgement}
\PAR This work, in part, was prosecuted at the Division of Information Sciences at Academia Sinica in Taipei, Taiwan. The author wishes to thank that institution for providing facilities and specifically Professors Hsu Tsan-Sheng and Kao Ming-Tat for their help, friendship and hospitality.
\vfill\eject
                                                                                              \Section{Bibliography}{Bibliography}\xrdef{pageBibliography}
\Par \Mark{BibHandS}{HEWITT, EDWIN AND KARL STROMBERG}: Real and Abstract Analysis. Springer Verlag 1969.
\PAr \Mark{BibSeidenberg}{SEIDENBERG,ABRAHAM}: Lectures in Projective Geometry. Princeton: D. Van Nostrand Company, Inc. 1962
\PAr \MarkWithText{BibTITS}{Jacques Tits}TITS, JACQUES: G\'en\'eralisation des groupes projectifs bas\'es sur la notion de transitivit\'e. Mem. Acad. Roy. Beg. \Bf{27} (1952), 115 p.
\PAr \Mark{BibVandY}VEBLEN, OSWALD AND JOHN WESLEY YOUNG: Projective Geometry, Vol. I. Boston: Ginn and Company, 1910.
\PAr \MarkWithText{BibSTAUDT}{\Underbar{Geometrie der Lage}}VON STAUDT, KARL: Geometrie der Lage. N\"urnberg: Verlag der Friedr. Korn'schen Buchhandlung, 1847.
 \vfill\eject\end